\documentclass[11pt,reqno,tbtags,a4paper]{amsart}
\usepackage{amssymb}
\usepackage{xpunctuate}
\usepackage{url}
\usepackage[square,numbers]{natbib}
\bibpunct[, ]{[}{]}{;}{n}{,}{,}
\usepackage{graphicx}
\usepackage{tikz}
\graphicspath{{./}}

\title[Frege's election method]
{Frege's election method: discrepancy, invariant polygons and Farey geometry}

\date{10 November 2026}

\newcommand{\MissingGraphicBox}[2][6.0cm]{%
  \begingroup
  \setlength{\fboxsep}{0pt}%
  \fbox{%
    \begin{minipage}[c][#1][c]{0.82\textwidth}
      \centering\small Missing figure: \texttt{#2}
    \end{minipage}}%
  \endgroup}
\makeatletter
\let\OrigIncludegraphics\includegraphics
\renewcommand{\includegraphics}[2][]{%
  \IfFileExists{#2}{\OrigIncludegraphics[#1]{#2}}{\MissingGraphicBox{#2}}}
\def\pgfdeclareimage{\@ifnextchar[{\pgfdeclareimage@opt}{\pgfdeclareimage@noopt}}
\def\pgfdeclareimage@opt[#1]#2#3{%
  \expandafter\gdef\csname pgfimageopts@#2\endcsname{#1}%
  \expandafter\gdef\csname pgffile@#2\endcsname{#3}}
\def\pgfdeclareimage@noopt#1#2{%
  \expandafter\gdef\csname pgfimageopts@#1\endcsname{}%
  \expandafter\gdef\csname pgffile@#1\endcsname{#2}}
\def\pgfuseimage#1{%
  \begingroup
  \edef\pgf@opts{\csname pgfimageopts@#1\endcsname}%
  \edef\pgf@file{\csname pgffile@#1\endcsname}%
  \IfFileExists{\pgf@file}{\expandafter\OrigIncludegraphics\expandafter[\pgf@opts]{\pgf@file}}{\MissingGraphicBox{\pgf@file}}%
  \endgroup}
\makeatother

\author{Jean-Pierre Conze}
\author{Svante Janson}
\author{Anders \"Oberg}
\thanks{Supported by the Knut and Alice Wallenberg Foundation
and the Swedish Research Council}
\address{Jean-Pierre Conze, IRMAR, Campus Beaulieu, 35042 Rennes Cedex, France}
\email{jean-pierre.conze@univ-rennes.fr}
\address{Svante Janson and Anders \"Oberg, Department of Mathematics, Uppsala University, PO Box 480,
SE-751~06 Uppsala, Sweden}
\email{svante.janson@math.uu.se}
\email{anders@math.uu.se}
\newcommand\urladdrx[1]{{\urladdr{\def~{{\tiny$\sim$}}#1}}}
\urladdrx{http://www2.math.uu.se/~svante/}

\numberwithin{equation}{section}

\renewcommand\le{\leq}
\renewcommand\ge{\geq}
\renewcommand\leq{\leqslant}
\renewcommand\geq{\geqslant}

\allowdisplaybreaks

\theoremstyle{plain}
\newtheorem{theorem}{Theorem}[section]
\newtheorem{fregeintrotheorem}{Theorem}

\newtheorem{lemma}[theorem]{Lemma}
\newtheorem{proposition}[theorem]{Proposition}
\newtheorem{corollary}[theorem]{Corollary}

\newtheorem*{method}{\MethodName}

\theoremstyle{definition}

\newcommand\xqed[1]{%
    \leavevmode\unskip\penalty9999 \hbox{}\nobreak\hfill
    \quad\hbox{#1}}

\newtheorem{exampleqqq}[theorem]{Example}
\newenvironment{example}{\begin{exampleqqq}}
  {\xqed{$\triangle$}\end{exampleqqq}}

\newtheorem{remarkqqq}[theorem]{Remark}
\newenvironment{remark}{\begin{remarkqqq}}
  {\xqed{$\triangle$}\end{remarkqqq}}

\newtheorem{definition}[theorem]{Definition}
\newtheorem{problem}[theorem]{Problem}

\theoremstyle{remark}

\newenvironment{romenumerate}[1][-10pt]{
\addtolength{\leftmargini}{#1}\begin{enumerate}
 \renewcommand{\labelenumi}{\textup{(\roman{enumi})}}%
 \renewcommand{\theenumi}{\labelenumi}%
 }{\end{enumerate}}

\newcounter{oldenumi}
{\setcounter{oldenumi}{\value{enumi}}
\begin{romenumerate} \setcounter{enumi}{\value{oldenumi}}}
{\end{romenumerate}}

\newcounter{thmenumerate}

\newcounter{xenumerate}   

\newcounter{cases}

\newcommand\pfitemx[1]{\par#1:}
\newcommand\pfitemref[1]{\pfitemx{\ref{#1}}}
\newcommand\pfitemnr{\smallskip\noindent\refstepcounter{cases}%
(\roman{cases}). } 

\newcommand\pfcase[2]{\smallskip\noindent\emph{Case #1: #2} \noindent}
\newcounter{steps}
\newcommand\stepx[1]{\smallskip\noindent\refstepcounter{steps}%
\emph{Step \arabic{steps}: \ignorespaces#1}\noindent}

\newcommand{\refT}[1]{Theorem~\ref{#1}}
\newcommand{\refTs}[1]{Theorems~\ref{#1}}
\newcommand{\refC}[1]{Corollary~\ref{#1}}
\newcommand{\refCs}[1]{Corollaries~\ref{#1}}
\newcommand{\refL}[1]{Lemma~\ref{#1}}

\newcommand{\refR}[1]{Remark~\ref{#1}}

\newcommand{\refS}[1]{Section~\ref{#1}}
\newcommand{\refSs}[1]{Sections~\ref{#1}}
\newcommand{\refSS}[1]{Section~\ref{#1}}

\newcommand{\refE}[1]{Example~\ref{#1}}

\newcommand{\refFig}[1]{Figure~\ref{#1}}
\newcommand{\refFigs}[1]{Figures~\ref{#1}}

\newcommand{\refStep}[1]{Step~\ref{#1}}
\newcommand{\refSteps}[1]{Steps~\ref{#1}}

\begingroup
  \divide\count255 by 60
  \multiply\count255 by -60
  \advance\count255 by \time
  \ifnum \count255 < 10 \xdef\klockan{\the\count1.0\the\count255}
  \else\xdef\klockan{\the\count1.\the\count255}\fi
\endgroup

\newcommand{\sumim}{\sum_{i=1}^m}
\newcommand{\sumjm}{\sum_{j=1}^m}

\newcommand\set[1]{\ensuremath{\{#1\}}}
\newcommand\bigset[1]{\ensuremath{\bigl\{#1\bigr\}}}

\newcommand\biggset[1]{\ensuremath{\biggl\{#1\biggr\}}}
\newcommand\lrset[1]{\ensuremath{\left\{#1\right\}}}
\newcommand\xpar[1]{(#1)}
\newcommand\bigpar[1]{\bigl(#1\bigr)}
\newcommand\Bigpar[1]{\Bigl(#1\Bigr)}

\newcommand\lrpar[1]{\left(#1\right)}

\newcommand\bigsqpar[1]{\bigl[#1\bigr]}

\newcommand\xcpar[1]{\{#1\}}

\newcommand\bigabs[1]{\bigl\lvert#1\bigr\rvert}

\def\rompar(#1){\textup(#1\textup)}    
\newcommand\xfrac[2]{#1/#2}

\newcommand\innprod[1]{\langle#1\rangle}

\def\xexp(#1){e^{#1}}
\newcommand\ceil[1]{\lceil#1\rceil}
\newcommand\lrceil[1]{\left\lceil#1\right\rceil}
\newcommand\floor[1]{\lfloor#1\rfloor}

\newcommand\round[1]{\langle#1\rangle}
\newcommand\frax[1]{\{#1\}}

\newcommand\ntoo{\ensuremath{{n\to\infty}}}
\newcommand\Ntoo{\ensuremath{{N\to\infty}}}

\newcommand\mtoo{\ensuremath{{m\to\infty}}}

\newcommand\punkt{\xperiod}    
    
\newcommand\ie{i.e\punkt}
\newcommand\eg{e.g\punkt}

\newcommand\cf{cf\punkt}

\newcommand{\aex}{a.e\punkt}

\newcommand{\tend}{\longrightarrow}

\newcommand\wto{\overset{\mathrm{w}}{\tend}}

\newcommand\bbR{\mathbb R}

\newcommand\bbN{\mathbb N}
\newcommand\bbT{\mathbb T}
\newcommand\bbQ{\mathbb Q}
\newcommand\bbZ{\mathbb Z}

\newcounter{CC}
\newcounter{cc}

\newcommand\E{\operatorname{\mathbb E{}}}
\renewcommand\P{\operatorname{\mathbb P{}}}

\newcommand\argmax{\operatorname{argmax}}
\newcommand\argmin {\operatorname{argmin}}

\newcommand\supp{\operatorname{supp}}

\newcommand\inter{\ensuremath{^\circ}}

\newcommand\ga{\alpha}
\newcommand\gb{\beta}
\newcommand\gd{\delta}
\newcommand\gD{\Delta}

\newcommand\gG{\Gamma}

\newcommand\kk{\kappa}
\newcommand\gl{\lambda}

\newcommand\go{\omega}
\newcommand\gO{\hat{\cO}}

\newcommand\gS{\Sigma}

\newcommand\gU{\Upsilon}
\newcommand\eps{\varepsilon}

\renewcommand\phi{\xxx}  

\newcommand\cA{\mathcal A}
\newcommand\cB{\mathcal B}
\newcommand\cC{\mathcal C}
\newcommand\cD{\mathcal D}

\newcommand\cH{\mathcal H}

\newcommand\cL{{\omega}}
\newcommand\hcL{\hat{\cL}}

\newcommand\LL{{\mathcal L}}
\newcommand\LLL{{\overline\LL}} 
\newcommand\LG{{\mathcal G}} 
\newcommand\HG{\mathcal L} 

\newcommand\cN{\mathcal N}
\newcommand\cO{\mathcal O}
\newcommand\cP{\mathcal P}

\newcommand\cR{{\mathcal R}}
\newcommand\cS{{\mathcal S}}

\newcommand\cU{{\mathcal U}}
\newcommand\cV{\mathcal V}

\newcommand\cX{{\mathcal X}}
\newcommand\cY{{\mathcal Y}}

\newcommand\indic[1]{\boldsymbol1\xcpar{#1}}

\newcommand\etta{\boldsymbol1}

\newcommand\qw{^{-1}}

\newcommand\oi{\ensuremath{[0,1]}}

\newcommand\dtv{d_{\mathrm{TV}}}

\newcommand\dd{\,\mathrm{d}}
\newcommand\ddx{\mathrm{d}}

\newcommand\rhs{right-hand side}

\newcommand\dH{D'Hondt}
\newcommand\SL{Saint-Lagu\"e}
\newcommand\Ho{H_{\0}}
\newcommand\hHo{\hH_{\0}}

\newcommand\ve{\mathbf{e}}

\newcommand\vi{{\mathbf{i}}}

\newcommand\vn{\mathbf{n}}
\newcommand\vp{\mathbf{p}}
\newcommand\hvp{\hat{\mathbf{p}}}
\newcommand\vq{\mathbf{q}}
\newcommand\vs{\mathbf{s}}
\newcommand\vu{\mathbf{u}}
\newcommand\hvu{{\hat{\mathbf{u}}}}
\newcommand\vv{\mathbf{v}}
\newcommand\hvv{{\hat{\mathbf{v}}}}
\newcommand\vw{\mathbf{w}}
\newcommand\vx{\mathbf{x}}
\newcommand\vy{\mathbf{y}}
\newcommand\vz{\mathbf{z}}

\newcommand\hvx{\hat{\mathbf{x}}}

\newcommand\hvz{\hat{\mathbf{z}}}

\newcommand\hnu{\hat{\nu}}
\newcommand\cp{^{\mathsf c}}
\newcommand\bs{\overline{s}}
\newcommand\bsk{\overline{s}'_k}
\newcommand\hx{\hat x}

\newcommand\hz{\hat z}
\newcommand\etal{{et al\xperiod}}
\newcommand\Rp{R_\vp}
\newcommand\bcO[1]{\overline{\cO(#1)}}
\newcommand\bgO[1]{\overline{\gO(#1)}}
\newcommand\hH{\widehat{H}}
\newcommand\bbTxm{\widetilde{\bbT}^m}
\newcommand\bbNo{\bbN_0}
\newcommand\hcN{\widehat{\cN}}
\newcommand\bgD{\overline{\gD}}
\newcommand\bgDx{\overline{\gD}^*}
\newcommand\sfs{{\mathsf{s}}}
\newcommand\cLx{\cL^*}
\newcommand\Kx{K^*}
\newcommand\glKx{\gl_{\Kx}}

\newcommand\stj{_*}
\newcommand\glhHo{\gl_{\hHo}}
\newcommand\dW{\partial W}
\newcommand\mmx{[m]^{(2)}}
\newcommand\hB{\widehat B}
\newcommand\hBo{\hB\inter}
\newcommand\Bo{B\inter}
\newcommand\SA{^{(1)}}
\newcommand\SB{^{(2)}}
\newcommand\SC{^{(3)}}
\newcommand\tgG{\widetilde{\gG}}
\newcommand\chgG{\check{\gG}}

\newcommand\uv{_{u/v}}
\newcommand\uvi{_{u'/v'}}

\newcommand\uvx[1]{_{u_{#1}/v_{#1}}}
\newcommand\ttt{\widetilde t}
\newcommand\jx{{j^*}}
\newcommand\jox{j_0^*}
\newcommand\kx{{k^+}}
\newcommand\kxx[1]{{#1^+}}
\newcommand\oo{^0}
\newcommand\linje[1]{\overline{#1}}
\newcommand\Tn{T^{(n)}}
\newcommand\sx{s^*}

\newcommand{\splitatcommas}[1]{
  \begingroup
  \begingroup\lccode`~=`, \lowercase{\endgroup
    \edef~{\mathchar\the\mathcode`, \penalty0 \noexpand\hspace{0pt plus 1em}}%
  }\mathcode`,="8000 #1%
  \endgroup
}

\newcommand{\R}{{\mathbb{R}}}

\newcommand{\Z}{{\mathbb{Z}}}

\newcommand{\T}{{\mathbb{T}}}

\newcommand\0{{\boldsymbol 0}}
\newcommand\1{{\boldsymbol 1}}

\begin{document}

\begin{abstract}
We give a comprehensive study of Frege's modified election method as a
proportional allocation rule and as a dynamical system. 
The method was introduced by Harrenstein, Lackner and  Lackner  (2022) 
as a modification of Frege's original proposal (1918), and it
coincides with the greedy algorithm for the chairman assignment problem
studied by Coppersmith et al. (2011). We consider several aspects of this
method.

First, we prove
that the harmonic bound on the discrepancy earlier shown to be sharp for 
time-dependent vote proportions remains sharp  for 
fixed vote proportions; we further quantify the restrictions on votes needed
for the discrepancy to come close to this bound. 

Secondly, the method defines a dynamical system with a
piecewise translation  
which, for almost every vote vector, 
has a unique compact attractor made up
of finitely many polytopes, and every orbit is equidistributed on it,
whatever choices are made at ties. 

Thirdly, for three parties, the dynamics is
two-dimensional. We construct the invariant polygons explicitly and give
a complete description of the geometry governing their bifurcations
as the votes vary in the two-dimensional parameter space;
this space 
is partitioned into a countably infinite number of ``Farey cells'' indexed by
reduced fractions in $(0,\frac12]$, 
and each of them is further partitioned into a finite number of
stability domains, separated by stability walls.
In an open Farey cell indexed by a reduced fraction
$u/v$, the invariant polygon has $4v+3$ vertices inside each stability
domain and $4v+2$ on an internal stability wall. 
We construct the Farey cells explicitly and determine all stability walls and
their endpoints; we also show that stability walls in different Farey cells
combine to form a pattern of triangles.
The Farey construction is defined
directly by neighbouring fractions in the Stern--Brocot tree
and their mediants, independently of the Frege map.
\end{abstract}

\maketitle
\tableofcontents

\section{Introduction and summary of main results}\label{S:intro}
Frege's election method was proposed in 1918 by Gottlob Frege \cite{Frege};
it is based on an idea which seems to be unique to Frege's method
among the multitude of election methods:
the method aims to obtain proportionality
on the average over a sequence of elections.
To this end, ``unused'' votes are saved from one election
to the next; this makes the election method a dynamical system.
The method has not been used in practice, and
it was for a long time forgotten, but it was studied in 2022
by \citet{HLL}, who also pointed out that a minor modification of Frege's
original method is mathematically nicer, and also seems to better satisfy
Frege's objectives.
We will in the present paper only consider this modification, which we
simply call \emph{Frege's method};
see \refR{RFrege} for a brief description of the original version.
For background on election methods, see for example
\citet{BY} and \citet{Pukelsheim}.

A detailed mathematical definition of the (modified) Frege method is given
in \refS{SFrege}. Briefly, we have a finite number $m$ of candidates (or
parties) that at each integer time take part in an election where one seat is
assigned. (See further \refR{Rprop}.)
Each candidate has a score, initially zero.
In each round, the vote proportions are added to the scores;
the candidate with the largest score wins and has one unit deducted.
The scores after the round thus record the discrepancy between the sum of
the vote proportions received so far and the number of times each candidate has
won.

In general, the votes may be time-dependent, but an important special case
is the time-homogeneous case, when the vote proportions are the same in
every round. Much of our analysis concerns this case.
For fixed positive vote proportions $\vp=(p_1,\dots,p_m)$ with $\sum_i p_i=1$,
Frege's method is described by the map acting on the hyperplane
$\Ho:=\{\vx\in\bbR^m:\sum_i x_i=0\}$ by
\begin{equation}\label{tvx1}
 T(\vx)=
 \bigl\{\vx+\vp-\ve_i:x_i+p_i=\max_j(x_j+p_j)\bigr\},
\end{equation}
where $\ve_i$ is the $i$th standard basis vector.
We define $T$ in \eqref{tvx1} as a multi-valued map:
in case of a tie, $T(\vx)$ contains all possible resulting states.
In many of our results we do not assume any particular tie-breaking rule,
but retain all admissible images.
We may also fix a tie-breaking rule and obtain a single-valued version of $T$.
Although one expects similar results for the multi-valued and single-valued
versions, this does not always hold (\refR{Rnoback}).
We will therefore often consider both versions.

Several other authors have studied the same or closely related dynamical
systems in different contexts, at least in single-valued versions.
We describe some of them here, together with summaries of our main results;
see \refSs{SSChair} and \ref{SSAdler} for further details and references.

\medskip
\emph{First,} 
the \emph{chairman assignment problem} 
\cite{Meijer73,Tijdeman80,Coppersmith}
has mathematically the same setting: a representative (chairman)
is chosen repeatedly, with the objective of making the choices 
approximately proportional to some given weights
(which in general may be time-dependent).
Much of the literature concerns the optimization problem;
in particular,
\citet{Coppersmith} study a greedy
algorithm, which coincides with  Frege's method;
see further \refSS{SSChair}.
The main results in these papers concern the discrepancy from perfect
proportionality. In particular, consider for simplicity the time-homogeneous
case (see \refS{SFrege} for the general, time-dependent case), denote the
vote proportions as above by $p_i$, and
let $r_j(n)$ be the number of seats assigned to candidate $j$ in the first
$n$ rounds; then the discrepancy is
\begin{equation}
 z_j(n):=np_j-r_j(n).
\end{equation}
\citet{Coppersmith} prove, even for time-dependent votes, the
upper bound $z_j(n)\le h_m-1$, where 
$h_m:=\sum_{i=1}^m\xfrac1i$,
and Tijdeman \cite{Tijdeman80} had earlier shown
that this bound can be attained when the votes are allowed to be
time-dependent. We prove that this bound is sharp also
in the time-homogeneous case.

\begin{fregeintrotheorem}[Fixed-vote discrepancy]\label{TIntroA}
For every $m\ge2$ and $\varepsilon>0$, there is a fixed positive vote
vector $(p_i)_1^m$ for which $z_j(n)>h_m-1-\varepsilon$ for some $j$ and $n$.
For $m\ge3$, however, every fixed vote vector satisfies
\begin{equation}\label{tia1}
 \sup_{\substack{1\le j\le m\\ n\ge1}}z_j(n)<h_m-1.
\end{equation}
Moreover, if $z_j(n)\ge h_m-1-\varepsilon$, the candidates can be
relabelled so that
\begin{equation}\label{tia2}
 p_{i+1}\le C_m\varepsilon p_i,\qquad 2\le i\le m-1,
\end{equation}
where $C_m$ depends only on $m$.
\end{fregeintrotheorem}

These assertions are proved in Theorems~\ref{T>}, \ref{TR<<}, and \ref{Tchar}.
The estimates \eqref{tia2} show that near equality in the bound \eqref{tia1}
requires each of $p_3,\dots,p_m$ to be small compared with the preceding
vote proportion. In particular, if $m\ge3$ and
the vote vector is chosen uniformly in the simplex, then
the probability that the discrepancy ever comes within $\varepsilon$ of the
upper bound is 
$O(\varepsilon^{\binom{m-1}{2}})$
(\refT{Tchar2}).

\medskip
\emph{Secondly,}
there is also a closely related dynamical literature.
The dynamical system \eqref{tvx1} describing Frege's method for fixed votes
is a piecewise translation.
Adler \etal{} \cite{adler0,adler10,adler17,adler18} study
a class of piecewise translations which includes the Frege map as a special
case, in the context of convex dynamics and error diffusion.
They consider a single-valued map with a fixed rule for ties,
and show the existence of absorbing invariant
tiles that are finite unions of polytopes.
G\"unt\"urk and Thao \cite{gun}
give more general tiling and regularity results for piecewise-affine maps.
We prove corresponding results for the multi-valued Frege map \eqref{tvx1},
where compact invariance includes every admissible image at a tie.
We also determine the limiting behaviour of every
\emph{branch orbit}, i.e., an orbit of $T$ obtained by choosing one of the
alternatives at each step when there is a tie.
See \refSS{SSAdler} for a more detailed comparison.

The Frege map $T$ in \eqref{tvx1} projects modulo $\bbZ^m$
to translation by $\vp$ on the subtorus
\begin{equation}
\hHo:=\bigl\{\widehat{\vx}\in\bbT^m:
\textstyle\sum_i\widehat x_i=0\bigr\}
\cong\bbT^{m-1}.
\end{equation}
We call $\vp$ \emph{ergodic} if
$p_1,\dots,p_m$ are rationally independent over $\bbQ$; this is equivalent
to ergodicity of the translation on $\hHo$.
(This condition holds for almost every vote vector.)

\begin{fregeintrotheorem}[The compact attractor]\label{TIntroB}
For every positive ergodic vote vector $\vp$, the multi-valued Frege map
has a unique non-empty compact invariant set $\Kx\subset\Ho$.
It absorbs every compact subset of $\Ho$ in finite time.
The set $\Kx$ is a finite union of full-dimensional polytopes with facets
parallel to the tie hyperplanes $x_i=x_j$.
Normalized Lebesgue measure on $\Kx$ is the unique invariant probability
measure on $\Ho$. Every branch orbit, from every initial point in $\Ho$,
has limit set $\Kx$ and is equidistributed with respect to this measure.
\end{fregeintrotheorem}

These conclusions are proved in Theorems~\ref{Tinv}, \ref{Tgl},
\ref{TB}, and \ref{TX}.

Berth\'e, Carton, Chevallier, Steiner and Yassawi \cite{BCCSY24}
study Tijdeman's algorithm for the chairman assignment problem.
They show that the resulting sequences are natural codings of toral
translations with polytopal partition atoms.
This algorithm uses a different selection rule from Frege's method and
achieves Tijdeman's uniform bound on absolute discrepancy.
We compare this bound with the sharp bound for Frege's greedy rule
in \refSS{SSAdler}.

Chevallier \cite{Chev} relates geometric codings of ergodic translations
of $\bbT^2$ to fundamental
domains and bounded remainder sets, while Berth\'e, Steiner and
Thuswaldner \cite{BST23} develop constructions based on multidimensional
continued fractions. For Frege's method, the translation vectors and
winner regions are prescribed by the allocation rule.

\medskip
\emph{Thirdly,}
we study the time-homogeneous dynamical system in the case $m=3$ in detail.
(The simpler case $m=2$ is straightforward and treated in \refS{Sm=2}.)
For $m=3$, both the phase space and the parameter space are two-dimensional.
For given vote proportions $\vp$,
the invariant domain is a polygon with sides in three fixed directions.
We construct it explicitly and describe how it changes as $\vp$ varies;
in particular, its number of vertices is unbounded.

Write $\vp=(p,q,r)$, where $r=1-p-q$. By symmetry, for distinct vote
proportions it suffices to consider the parameter triangle
\begin{equation}
 \Delta:=\{(p,q):0<q<p<1-p-q\}.
\end{equation}
Up to boundaries, we partition this triangle into open \emph{Farey cells}
$\cR_{u/v}$, indexed by reduced fractions $u/v\in(0,1/2]$.
Inside each cell, a \emph{stability domain} is an open parameter region
on which the symbolic description of the polygon, including its division
into winner regions, is constant.

\begin{fregeintrotheorem}[Planar classification]\label{TIntroC}
\begin{romenumerate}
\item
For $(p,q)$ in any open Farey cell $\cR_{u/v}$, we construct explicitly
a compact polygonal domain $\cD$ with $T(\cD)=\cD$.
If\/ $p,q,1$ are rationally independent (the ergodic case), then $\cD=\Kx$,
the invariant set in \refT{TIntroB}.
\item
Each Farey cell $\cR\uv$ is divided by finitely many disjoint
\emph{stability walls}, each lying on a straight line, into
open convex \emph{stability domains}. On every such domain, the vertices of
$\cD$ depend affinely on $(p,q)$ and their number is exactly $4v+3$.
On a wall the number of vertices is $4v+2$.
\item
The walls are the intersections of the open Farey cells with
explicit line segments $\LG_{u/v}$, indexed by reduced fractions
$u/v\in(0,1/2)$ and starting at $(u/v,0)$.
Furthermore, these segments form an infinite family of triangles with base on
$q=0$:
each fraction $u/v$ with odd denominator is
paired with its unique child with even denominator in the Stern--Brocot
tree, and the two corresponding segments have a common endpoint.
\end{romenumerate}
\end{fregeintrotheorem}

The Farey partition is given in \refSS{SSm3farey};
the polygon construction in \refSS{SSm3construction};
the stability domains and vertex counts in
\refSS{SSSD} (Theorems~\ref{TSD}, \ref{TW2}, and \ref{Tm3-4v+3});
and the wall segments and their pairing in
Sections~\ref{SSCW} and~\ref{SSm3-wall-pattern}
(Theorems~\ref{Twalls} and~\ref{ConjFarey+lines}).

The construction also applies to non-ergodic parameters in the open
cells. Every ergodic parameter in $\Delta$ lies in a stability domain,
so the classification gives its invariant attractor explicitly.

The Farey cells and their further subdivision into stability domains
can be studied independently of Frege's map:
the cells and wall segments are defined directly
by simple arithmetic relations using neighbouring
fractions in the Stern--Brocot tree
and their mediants.
Here we use this arithmetic and geometric
construction to classify the invariant polygons of the Frege map and
their changes with the parameters.

The paper is organized as follows. Section~\ref{SFrege} gives the precise
definitions and the comparison with earlier work.
Sections~\ref{Sbounds} and \ref{Sopt} prove the discrepancy bounds and
their optimality for fixed votes. Sections~\ref{Sorb}--\ref{Sinter}
develop the dynamical theory, treating rational, ergodic, and intermediate
parameters separately; the compact attractor results in the ergodic case are in
Section~\ref{Sergodic}. Section~\ref{Sbias} discusses
the asymptotic bias of the election method and
relates it to the invariant measure.
Section~\ref{Sm=2}
treats the one-dimensional case $m=2$, and
Section~\ref{Sm=3} contains the planar classification for the case $m=3$,
from the Farey partition and polygon construction to the
wall pattern. 
Finally, Section~\ref{Sdiv} briefly describes related dynamical
systems that are defined by standard
divisor election methods.

\section{Frege's election method and related mathematics}\label{SFrege}

\subsection{Some notation}\label{SSnot}
We generally use boldface for vectors in $\bbR^m$ (for some $m\ge2$), 
and write for example $\vp=(p_i)_{i=1}^m$. 
We let $\0=(0,\dots,0)$ and $\1=(1,\dots,1)$.
$\gS_m$ denotes the  $(m-1)$-dimensional simplex
\begin{align}
\gS_m:=\biggset{\vp\in\bbR^m:p_i\ge0\,\forall i   \text{ and } \sumim p_i=1},
\end{align}
i.e., the set of probability vectors in $\bbR^m$.

We let $[m]:=\set{1,\dots,m}$.

An important role is played by the \emph{harmonic numbers}
(used already in \refT{TIntroA})
\begin{align}
  h_m:=\sum_{i=1}^m \frac{1}{i}, \qquad m\ge1. 
\end{align}

The indicator of a set $A$ is denoted $\etta_A$.

Lebesgue measure (in any dimension, and sometimes normalized by a constant) 
is denoted by $\gl$.
Occasionally we use $\gl_K$ for Lebesgue measure on a set $K$.

The inner product (scalar product) of two vectors $\vx$ and $\vy$ is denoted
$\innprod{\vx,\vy}$.

For a real number $x$, $\floor x$, $\ceil x$, and $\round x$ denote $x$
rounded to the largest integer $\le x$, the smallest integer $\ge x$, 
and the nearest integer, respectively. We regard $\round x$ as taking both
values $x\pm\frac12$ when $x$ is halfway between two integers.

We let $\bbNo:=\set{0,1,\dots}$. 

\subsection{The Frege election method}\label{SSFrege}

The (modified) Frege election method introduced  in \cite{HLL}
can be defined as follows.
(See \refR{RFrege} for Frege's original version, which will not be
considered in the present paper.)
We will throughout the paper use the notations below.

We are given $m\ge2$ candidates (or parties) $1,\dots,m$.
At each integer time $t\ge1$, we are given a vector 
$\vp(t)=(p_j(t))_1^m\in\gS_m$ of vote proportions;
thus $p_j(t)\ge0$ and
\begin{align}\label{a1}
  \sumjm p_j(t)=1.
\end{align}
(In other words, $\vp(t)$ is a probability distribution on \set{1,\dots,m}.)

Let $E(t)\in[m]:=\set{1,\dots,m}$ be the candidate elected at time $t\ge1$.
$E(t)$ is determined as follows.

\def\MethodName{Frege's method}
\begin{method}
Define recursively aggregate scores $s_j(t)$ and $z_j(t)$ for the candidate $j$
before and after electing representative $t\ge1$ by
\begin{align}\label{z0}
  z_j(0)&:=0;\\
s_j(t)&:=z_j(t-1)+p_j(t);\label{sj}
\end{align}
let $E(t)$ be the candidate $j$ that maximizes $s_j(t)$
(treating ties by any rule);
continue the recursion by
\begin{align}
z_j(t)&:=s_j(t)-\indic{E(t)=j}.\label{zj}
\end{align}
\end{method}

From \eqref{a1} and \eqref{z0}--\eqref{zj} we have
\begin{align}\label{sumsum}
\sumjm s_j(t)=1
\text{\quad and\quad}
  \sumjm z_j(t)=0,
\qquad t\ge1.
\end{align}

Let $r_j(t)$ be the number of times up to time $t$ that $j$ is elected, and
let
\begin{align}
  \label{Sj}
P_j(t):=\sum_{s=1}^t p_j(s).
\end{align}
Note that it follows from \eqref{z0}--\eqref{zj} that
\begin{align}\label{rj}
  z_j(t) = P_j(t)-r_j(t).
\end{align}
Hence, $z_j(t)$ measures how much the outcome of Frege's method deviates
from perfect proportionality up to time $t$.
Many of our results are stated in terms of the deviations $z_j(t)$.

Note that in the definition above, we let the vote proportions $p_i(t)$
depend on $t$.
This is obviously the case in the situation considered by Frege,
with a sequence of elections at regular intervals.
However, the special case (which we call \emph{time-homogeneous})
when the probabilities $p_i(t)=p_i$ are independent of $t$ is of obvious
theoretical interest, 
and in this paper we are mainly concerned with this case.
Note also that in the case of real elections, the proportions $p_i(t)$ are
rational numbers, but we make no such assumption and allow arbitrary 
real $p_i(t)\in\oi$ satisfying \eqref{a1}.

\begin{remark}\label{RFrege}
  Frege's original proposal \cite{Frege} does not use vote proportions $p_j$
  but the actual numbers of votes, $\pi_j(t)$ say, for the different candidates.
The method can be described as in \eqref{z0}--\eqref{zj} with the difference
that $\pi_j(t)$ is used instead of $p_j(t)$, and that in \eqref{zj}, if
$E(t)=j$, we subtract (instead of 1) the integer part of the average score,
\begin{align}
 \left\lfloor\frac 1m \sumjm s_j(t)\right\rfloor.
\end{align}
If the total number of voters is constant, this cost eventually equals
that number; see \cite[Lemma~1]{HLL}. After rescaling, the recurrence then
has the form of the modified method, with initial scores determined by
the preceding elections.
\end{remark}

\begin{remark}\label{Rprop}
  \citet{Frege} stated his method for a sequence of different elections in a
  single-member constituency.
\citet{HLL} points out that the (modified) method can be used also to elect
several representatives in a multi-member constituency (perhaps the entire
country) in a single election where each voter votes for a party
(each party having a list of candidates):
we let the vector $\vp$ be the proportions of votes for the different parties,
and use the algorithm repeatedly with the same vector $\vp$ until the
desired number are elected. 
Regarded in this way (which was not Frege's own), the  Frege's  method
becomes a new proportional election method (also called apportionment method).
This aspect is investigated further in \cite{HLL}, with comparisons to 
some other well-known proportional election methods.
Note that in this application of Frege's method,
only the time-homogeneous case is relevant; 
this is another motivation for studying that case.
\end{remark}

Although \citet{Frege} formulated his method for individual candidates and
not for a party system, it is obvious that the method can be applied to 
a system with parties as well.
(In the situation in \refR{Rprop}, it is necessary to consider parties.)
We will therefore use ``candidates'' and ``parties'' interchangeably below.

\subsection{Frege's method as a dynamical system}\label{SSdynamical}
Consider now the time-homogeneous case $\vp(t)=\vp$ for some fixed
probability vector $\vp$.
The recursion \eqref{sj}--\eqref{zj} then defines both
$\vs(t):=(s_j(t))_{j=1}^m$  and $\vz(t):=(z_j(t))_{j=1}^m$ as dynamical
systems.
These two dynamical systems are essentially the same; by \eqref{sj} they
differ only by a translation by the fixed vector $\vp$, and also by
different initial conditions. 
We focus on $\vz(t)$ and note that  all our results can be 
translated to $\vs(t)$.

Note that these two dynamical systems act on vectors $\vs$ and $\vz$ in
$\bbR^m$; however, with our initial condition $\vz(0)=\0$ in \eqref{z0},
$\vs(t)$ and $\vz(t)$ are confined to the hyperplanes \eqref{sumsum};
hence the dynamics is really in a $(m-1)$-dimensional space.

In dynamical system notation, we can reformulate the definition of Frege's
method as follows.

For a vector $\vx=(x_i)_1^m$, let
\begin{align}\label{theta}
  \theta(\vx) := \argmax_{1\le i\le m} x_i,
\end{align}
i.e., the index $i$ such that $x_i = \max_j x_j$;
we say that $\theta(\vx)$ is the winner among $\{1, \dots, m\}$.
If the maximum is tied, we may assume that $\theta$ is chosen by some rule
(for example the smallest $i$ such that $x_i$ attains the maximum),
but it is sometimes better to regard $\theta(\vx)$ as multi-valued in
such cases. 
Also, let $\ve_1,\dots,\ve_m$ denote the standard basis vectors in $\bbR^m$,
and let
$\1:=(1, 1, \dots, 1)^t\in\bbR^m$, the vector with all coordinates 1. 
Then \eqref{sj}--\eqref{zj} can be written as
\begin{align}\label{dyn}
  \vz(t) = T \bigpar{\vz(t-1)},
\end{align}
where $T$ is the (piecewise translation)  operator
\begin{align}\label{tx}
  T(\vx) = \vx+\vp-\ve_{\theta(\vx+\vp)}.
\end{align}

Since $\vp$ is a probability vector, we have $\langle \vp, \1 \rangle = 1$,
and it follows from \eqref{tx} that
$\langle T\vx, \1 \rangle = \langle \vx, \1 \rangle$ for all $\vx\in\bbR^m$.
Therefore
\begin{align}\label{tn1}
\langle T^t \vx, \1 \rangle = \langle \vx, \1 \rangle, 
\qquad t \ge 1,  
\end{align}
showing again that the dynamics is confined to a hyperplane.
Note that the map is not injective; for instance, 
we have $T(\ve_i-\vp)=\0$ for every basis vector $\ve_i$.

In Frege's election method, the recursion starts with $\vz(0)=\0$, but when
viewed as a dynamical system, it is of interest to study arbitrary initial
values, and we will often do so.
In general, we 
denote the initial value of the process by
$\vu = (u_1, \dots, u_m)$; i.e., we
let $\vz(0) = \vu$.
We may also write $\vz(t;\vu)$ and $z_j(t;\vu)$ to show the dependence on 
the initial state $\vu$; thus $\vz(0;\vu)=\vu$. By \eqref{tn1},
\begin{align}\label{tn2}
\langle \vz(t;\vu),\1 \rangle = \langle \vu, \1 \rangle, 
\qquad t \ge 0. 
\end{align}

Note also that for any $a\in\bbR$, 
we can translate the orbits by 
the vector $a\1=(a,\dots,a)^t$, 
because the translation by this vector commutes with the map $T$:
\begin{align}\label{tna}
  T(\vx+a\1) = T(\vx)+a\1,
\end{align}
since the set of possible winners remains the same under this translation.
For a single-valued version, \eqref{tna} requires that the tie rule be
unchanged when the same constant is added to all scores; any fixed priority
rule has this property.
Therefore, we can when convenient 
without loss of generality assume that the initial vector
$\vu\in\Ho$, where 
$\Ho$ is the hyperplane
\begin{align}\label{Ho}
\Ho:=\set{\vx:\innprod{\vx,\1}=0}  
=\biggset{\vx:\sumim x_i=0}.
\end{align}
(We will often assume this, but not always.)

\begin{remark}\label{Rmulti}
Unless we say otherwise, the results below are valid both if we regard $T$
as a single-valued map, with some arbitrary but fixed rule for ties, and 
if we regard $T$ as a multi-valued map.

If we want a single-valued map, the simplest choice
of a tie resolution rule
is to give priority to 
the candidates in a fixed order;
without loss of generality we may (as said above)
let the winner be the candidate (among the tying ones) with smallest index.
We call this the \emph{lexicographic tie resolution rule}.
More complicated rules are possible, but seem less interesting.

To regard $\theta$ and thus $T$ as multi-valued means formally that
we regard $T$ as mapping sets to sets; more precisely, $T$ maps the power set
$2^{\bbR^m}$ of all subsets of $\bbR^m$ to itself by, \cf{} \eqref{tx},
\begin{align}\label{txa}
  T(A) := \bigcup_{i=1}^m \bigpar{((A+\vp)\cap \cV_i) - \ve_i},
\end{align}
where, for $i\in[m]$,
\begin{align}
\cV_i:=\bigset{\vx\in\bbR^m: x_i=\max_{1\le j\le m}x_j}
=\bigset{\vx\in\bbR^m: x_i\ge x_j\text{ for every }j}.
\end{align}
This is the closed polyhedral cone where $i$ is the winner, or at
least one of the winners.
We will use this notation as a set-valued map sometimes for clarity, but we
often write $T$ as an ordinary function in $\bbR^m$, leaving 
the  interpretation in the multi-valued case to the reader when it is obvious.
(In any case, $T$ and $T^n$ are single-valued for \aex{} initial value
$\vu$, see \refR{RD} below.)

When regarding $T$ as a multi-valued operator, $\vz(n):=T^n\vu$ is formally
a set (finite and non-empty), but for simplicity we again often write
$\vz(n)$ as a point when this will not cause any confusion.
We may also consider a single \emph{branch}
$(\vz'(n))_{n\ge0}$, meaning a sequence of points 
$\vz'(n)=\vz'(n,\vu)\in\bbR^m$
such that $\vz'(0)=\vu$ and 
\begin{align}\label{branch}
\vz'(n+1)\in T\bigl(\vz'(n)\bigr),
\qquad n\ge0.
\end{align}

In real elections, ties are commonly not broken by a fixed rule but by lot,
i.e., by making a random choice. This is one reason for considering the
multi-valued version, which includes all possible outcomes of the random
choice(s).

There are also some mathematical advantages of the multi-valued version. 
Since every $\cV_i$ is closed, it follows from \eqref{txa} that
\begin{align}\label{compact}
\text{If $K\subset\bbR^m$
is compact, then $T(K)$ is compact.}
\end{align}
(If $T$ is regarded as a single-valued function with some rule for resolving
ties, then $T(K)$ is not necessarily compact; note that $T$ is discontinuous
at the boundaries of $\cV_i$.)

On the other hand, there are also obvious advantages in the analysis to
consider a single-valued $T$. We conclude that there are good reasons to
study both the single-valued and multi-valued versions.
(For many properties, there is no difference.)
\end{remark}

\begin{remark}\label{RD}
Let $\cD_\sfs$ be the set of all $\vx\in\bbR^m$ such that
there is a tie in \eqref{theta}, 
let $\cD:=\cD_\sfs-\vp$,
and let $\cU:=\bbR^m\setminus\cD$,
the set where $T\vx$ is uniquely defined by  \eqref{tx}.
The problem with ties and possible multiple values thus occurs only
when $\vx$ belongs to the exceptional set $\cD$,
which is contained in a union of a finite number of hyperplanes. 
(In fact, using the notation in \refR{Rmulti}, 
$\cD+\vp=\cD_\sfs=\bigcup_i\partial \cV_i$.)
In particular, $\gl(\cD)=0$, and thus
$T$ is uniquely determined, 
and thus single-valued,
for \aex{} $\vx\in\bbR^m$.
Furthermore, $T$ is continuous in $\cU$, but not at $\cD$.

We define also $\cD_0:=\cD\cap\Ho$ and $\cU_0:=\cU\cap\Ho$ and note that
$\gl(\cD_0)=0$ for Lebesgue measure $\gl$ in $\Ho$;
hence, $T$ is also uniquely determined (and continuous) for \aex{} $\vx\in\Ho$.
\end{remark}

\begin{remark}\label{Rothers}
Frege's  method is not unique among election methods
in being described by a dynamical system.
We show in \refS{Sdiv} that also the commonly used 
\dH{} and \SL{} election methods may be formulated as dynamical
systems, very similar to the one above.
\end{remark}

\subsection{The chairman assignment problem}\label{SSChair}
As mentioned in the introduction,
a related problem, often called the 
\emph{chairman assignment problem}, has been studied by a number of authors
including \citet{Meijer73}, \citet{Tijdeman80} and \citet{Coppersmith} 
(see these for further references). 
Mathematically,
the setting is the same as in \refSS{SSFrege} above. 
(The formulation in \cite{Tijdeman80} is that a number of states
form a union; each year they elect a chairman, with states having
different weights corresponding to our $p_i$ or $p_i(t)$.)
However, \cite{Meijer73} and \cite{Tijdeman80}
do not consider Frege's algorithm 
but the optimization problem of choosing a sequence $E=(E(t))_{t=1}^\infty$
minimizing (in our notation) 
\begin{align}\label{b1}
D(E):=  \sup_{j,t}|P_j(t)-r_j(t)|
=
  \sup_{j,t}|z_j(t)|.
\end{align}
The main result of \cite{Meijer73} and \cite{Tijdeman80} is that, 
assuming that the entire sequence $\vp(t)$ is known, 
it is always possible to achieve
$D(E) \le 1-(2m-2)\qw$, and that this bound is best possible,
i.e., 
\begin{align}\label{b2}
D_m:=\sup_{(\vp(t))} \inf_{(E(t))}D=1-\frac{1}{2m-2},
\end{align}
taking the supremum
over all sequences of probability vectors $(\vp(t))_{t\ge1}$.
\citet{Tijdeman80} notes also, see \refE{E1} below,
that if we only consider online algorithms,
i.e., algorithms that choose $E(t)$ knowing only $\vp(s)$ for $s\le t$,
then the bound \eqref{b2} cannot always be achieved; 
for any online algorithm we have
\begin{align}\label{b3}
\sup_{(\vp(t))}D\bigl((E(t)),(\vp(t))\bigr)
\ge\sup_{(\vp(t))} 
\sup_{j,t} z_j(t)
\ge \sum_{2}^m \frac{1}{i}=h_m-1.
\end{align}
Note that this bound is strictly larger than \eqref{b2} for every $m\ge3$, 
and unbounded as \mtoo.

\citet{Coppersmith} study only online algorithms, and in particular the
natural greedy algorithm, which  is the same as  Frege's method
introduced later (independently and in a different context) 
by \cite{HLL} and studied in the present paper.
It is shown in \cite{Coppersmith}  that this algorithm achieves the bound in
\eqref{b3}, i.e., that for 
Frege's method, for any sequence $\vp(t)$,
\begin{align}\label{b4<}
z_j(t)\le h_m-1,
\qquad j\in[m],\,t\ge1.
\end{align}
The corresponding sharp lower bound 
\begin{align}\label{b4>}
z_j(t)\ge -1+\frac1m,
\qquad  j\in[m],\,t\ge1,
\end{align}
(which is easier) 
is also shown in \cite{Coppersmith};  see \refE{E2} below for sharpness.
Frege's method is thus optimal among all online algorithms for the
chairman assignment problem.

\begin{remark}
\citet{HLL} studied, as said above, both
Frege's original method and the modified method.
In particular they prove  \cite[Theorem 2]{HLL} that
in the time-homogeneous case,
for the modified Frege method 
(and also for Frege's original method, provided the total number of voters
remains constant),
we have the following proportionality result:
\begin{equation}\label{prop}
\lim_{t\to \infty} \frac{r_j(t)}{t}= p_j.
\end{equation}
They prove this as a consequence of \eqref{rj} and weaker versions of
\eqref{b4<}--\eqref{b4>}. 
\end{remark}

We give here for easy reference the well-known examples mentioned above.
\begin{example}\label{E1}
We repeat the following basic example that appears also in
\cite{Tijdeman80}, \cite{Coppersmith}, \cite{HLL}:
Let $m\ge2$.
At time $t=1,\dots,m$, split the votes evenly among the candidates that
have not yet been elected.
We may assume that the candidates are numbered in the order that they are
elected; we then thus let, for $1\le t\le m$,
\begin{align}
  p_j(t):=
  \begin{cases}
0, & 1\le j \le t-1,    
\\
\frac{1}{m-t+1}, & t\le j \le m
.  \end{cases}
\end{align}
We then have, by \eqref{rj},
\begin{align}
z_m(m-1)=
  P_m(m-1)=\sum_{i=2}^m\frac{1}i
=h_m-1.
\end{align}
As noted by \citet{Tijdeman80}, this example applies to any online algorithm.
\end{example}

\begin{example}\label{E2}
 The following almost trivial example \cite{Coppersmith} 
shows that \eqref{b4>} is sharp for every $m\ge2$, also in
the time-homogeneous case:
Let the $m$ candidates get equally many votes, so $p_i=1/m$.
If
$j$ wins the first election, then 
$z_j(1)=s_j(1)-1=\frac{1}{m}-1$.
\end{example}

\subsection{Related dynamical results}\label{SSAdler}
Frege's method is the greedy algorithm for the Chairman assignment problem
studied by \citet{Coppersmith}. The corresponding map also appears in
convex dynamics and error diffusion, where an input point in a polytope is added
to the current error vector and a nearest vertex is then subtracted.
For the standard simplex, the nearest-vertex rule is precisely the
largest-coordinate rule of Frege's method. With constant input, the map
is a piecewise translation which projects to a torus translation.

Adler \etal{} \cite{adler0} prove boundedness results 
for this class of maps
in the more general
setting of time-dependent inputs. For constant input, the subsequent
papers \cite{adler10,adler17,adler18} establish the existence of absorbing
invariant tiles and describe the polyhedral structure of their closures,
first in low dimensions and then for higher-dimensional acute simplices
under an ergodicity assumption. In these papers a rule is chosen to
resolve ties, so that the maps are single-valued. An invariant tile for
such a map need not be closed, and its closure need not be invariant.
Their results therefore do not directly give the compact invariant
attractor for the multi-valued map considered here.

For our purposes there is no need to specify how ties are resolved.
We retain all tied winners, and then $T$ maps compact sets to compact
sets, see \refR{Rmulti}. In the ergodic case we prove that there is a
unique non-empty compact invariant set $\Kx$, that it absorbs every
compact subset of $\Ho$ in finite time, and that it is a finite union of
full-dimensional polytopes. Here invariance means $T(\Kx)=\Kx$, 
for the multi-valued map
with all
admissible images of boundary points included. Moreover, every branch
orbit starting in $\Ho$ has limit set $\Kx$ and is equidistributed with
respect to normalized Lebesgue measure on $\Kx$.

Since the earlier theory was developed in several papers over a number
of years, and our treatment of the boundaries is different, we give
full proofs of the results needed here in Section~\ref{Sergodic}.

G\"unt\"urk and Thao \cite{gun} prove tiling and Jordan regularity
results for a more general class of piecewise-affine maps. Their results
are formulated in terms of essential invariance, that is, invariance up
to sets of measure zero. This is natural in their setting
but differs from our perspective;  
we consider pointwise invariance, including
every boundary point.

For related work on geometric codings of torus translations, see
Chevallier \cite{Chev} and Berth\'e, Steiner and Thuswaldner
\cite{BST23}. Particularly close to our setting is the work 
by Berth\'e, Carton, Chevallier, Steiner and Yassawi \cite{BCCSY24}, who
consider the chairman assignment problem 
but with a different algorithm than we do; 
they study
Tijdeman's algorithm \cite{Tijdeman80} 
and realize it
as a natural coding of a torus translation
with polytopal partition atoms. 

For general results on coding of piecewise isometries,
see \citet{Goetz}.

\section{The bounds for $z_j(t)$}\label{Sbounds}
As a background, and because we will reuse  the arguments later,
we begin by giving proofs of the bounds \eqref{b4<}--\eqref{b4>}.

\begin{theorem}[\citet{Coppersmith}]\label{T<}
  For every $j\in[m]$ and $t\ge0$,
  \begin{align}\label{t<}
-1+\frac{1}{m}\le    z_j(t) \le h_m-1=\sum_{i=2}^m\frac{1}{i}.
  \end{align}
\end{theorem}

\begin{proof}
Our proofs are similar to \cite{Coppersmith}.
We begin with the lower bound.

Observe that since $E(t)$ maximizes $s_j(t)$,
\begin{align}\label{a2}
  s_{E(t)}(t) \ge \frac{\sumjm s_j(t)}{m}=\frac{1}{m}.
\end{align}

  We use induction in $t$. The claim holds for $t=0$ by \eqref{z0}.
Suppose that it holds for $t-1$.
If $E(t)=j$, then, by \eqref{a2},
\begin{align}\label{a21}
z_j(t)=s_j(t)-1=s_{E(t)}-1\ge\frac{1}{m}-1. 
\end{align}
On the other hand, if
$E(t)\neq j$, then 
\begin{align}\label{a22}
  z_j(t)=s_j(t)=z_j(t-1)+p_j(t)\ge z_j(t-1)\ge-1+\frac{1}{m}
\end{align}
by the induction hypothesis.

We turn to the upper bound.
Fix $j\in[m]$ and a time $T\ge0$.
For $0\le t\le T$,
let 
\begin{align}\label{at}
  A_t:=\set{j}\cup\set{E(s):t<s\le T},
\end{align}
i.e., the set of candidates that will be elected at least once
until $T$, together with our chosen candidate $j$.
Obviously, $A_{t}\subseteq A_{t-1}$ for all $1\le t \le T$, with
either $A_{t}=A_{t-1}$ or $A_{t}=A_{t-1}\setminus\set{E(t)}$.

For any set $A$ of candidates, let
\begin{align}\label{qu1}
  z_A(t):=\sum_{j\in A} z_j(t),
\qquad
  s_A(t):=\sum_{j\in A} s_j(t),
\qquad
  p_A(t):=\sum_{j\in A} p_j(t)
.\end{align}

We will prove, by induction on $t$, that
\begin{align}\label{qu2}
  z_{A_t}(t) \le |A_t|\sum_{i=|A_t|+1}^m\frac{1}{i},
\qquad 0\le t\le T.
\end{align}
The base case $t=0$ is trivial, since $z_A(0)=0$ for every $A$.

Suppose that $t\ge1$ and that \eqref{qu2} holds with $t$ replaced by $t-1$.
By the definition \eqref{at}, we have $E(t)\in A_{t-1}$, and thus, 
by \eqref{zj}, 
\begin{align}\label{qu3}
  z_{A_{t-1}}(t)&
=s_{A_{t-1}}(t)-\indic{E(t)\in A_{t-1}}
= s_{A_{t-1}}(t)-1 
.\end{align}
Furthermore, by \eqref{sj} and \eqref{a1},
\begin{align} \label{qu6}
 s_{A_{t-1}}(t)
= z_{A_{t-1}}(t-1)+ p_{A_{t-1}}(t)
\le z_{A_{t-1}}(t-1)+1
.\end{align}
There are two cases:

\pfcase1{$A_t=A_{t-1}$.}
In this case, \eqref{qu3} and \eqref{qu6} yield
$z_{A_{t}}(t)\le z_{A_{t-1}}(t-1)$, and thus \eqref{qu2} follows from the
induction hypothesis.

\pfcase2{$A_t\subsetneq A_{t-1}$.}
In this case $A_t=A_{t-1}\setminus\set{k}$, where
$k:=E(t)$.
Let 
\begin{align}\label{qu4}
\bs:=\frac{s_{A_{t-1}}(t)}{|A_{t-1}|},
\end{align}
the average of $s_j(t)$ for $j\in A_{t-1}$.
Since $k$ was elected in this round, we must have $s_k(t)\ge s_j(t)$ for
every $j$, and thus also $s_k(t)\ge \bs$.
Furthermore, $E(t)=k\notin A_t$, and thus $z_j(t)=s_j(t)$ for every $j\in
A_t$ by \eqref{zj}. 
Hence,
\begin{align}\label{qu5}
z_{A_t}(t)=s_{A_t}(t)=s_{A_{t-1}}(t)-s_k(t)
\le s_{A_{t-1}}(t)-\bs
=\frac{|A_{t-1}|-1}{|A_{t-1}|}s_{A_{t-1}}(t)
.
\end{align}
Consequently, by \eqref{qu5}, \eqref{qu6},
and the induction hypothesis,
\begin{align}\label{qu7}
\frac{z_{A_t}(t)}{|A_t|}
&
=
\frac{z_{A_t}(t)}{|A_{t-1}|-1}
\le
\frac{1}{|A_{t-1}|}s_{A_{t-1}}(t)
\le \frac{1}{|A_{t-1}|}z_{A_{t-1}}(t-1)
+\frac{1}{|A_{t-1}|}
\notag\\&
\le \sum_{i=|A_{t-1}|+1}^m\frac{1}i
+\frac{1}{|A_{t-1}|}
= \sum_{i=|A_{t}|+1}^m\frac{1}i,
\end{align}
which shows \eqref{qu2} in this case. 
This completes the induction step.

Hence, \eqref{qu2} holds, for all $t\le T$. Take $t=T$.
We have, by \eqref{at},
$A_T=\set j$, and thus $z_{A_T}(T)=z_j(T)$ and $|A_T|=1$. 
Hence, \eqref{qu2} yields \eqref{t<}.
\end{proof}

We focus on the right inequality in \eqref{t<} 
and consider next the case of equality there.
The next theorem shows that \refE{E1} gives essentially the only way to
achieve
equality.
\begin{theorem}\label{TL=}
  Suppose that $z_j(T)=h_m-1$ for some $T\ge1$ and $j\in[m]$.
Then the candidates can be renumbered such that for some sequence of times
$1\le T_1<\dotsm< T_{m-1}= T$,
the vector $\vz(T_j)$ equals the vector
$\vz(j)$ in \refE{E1} for $j=1,\dots,m-1$.
\end{theorem}
\begin{proof}
Note first, by inspecting the proof above, that
if we have strict inequality in \eqref{qu2} for $t-1$, 
then so do we for $t$.
(In the case $A_t\subsetneq A_{t-1}$, this is seen in \eqref{qu7}.)
Hence, if equality holds for $t=T$, then it holds for all $t\in[0,T]$.
In particular, for $t=0$ the sum must vanish, and thus $|A_0|=m$.

Furthermore, the proof shows also that if we have equality for some $t\ge1$,
then we must have equality in \eqref{qu6}, and thus $p_{A_{t-1}}(t)=1$.
In other words, at time $t$ all votes are given to the members of $A_{t-1}$,
i.e., to the chosen $j$ and to the candidates that will be elected in $[t,T]$.

Let $T_\ell=\min\set{t:|A_t|=m-\ell}$ for $\ell=0,\dots,m-1$.
Note first that if $t\in[T_{m-1},T]$, then $|A_t|=1$ and thus $A_t=\set j$,
and equality in \eqref{qu2} shows that $z_j(t)=h_m-1$.
Moreover, for $t>T_{m-1}$, we thus have $p_j(t)=p_{A_{t-1}}(t)=1$.
It follows that for $t>T_{m-1}$ all votes go to $j$ who wins, and thus
\eqref{sj}--\eqref{zj} yield
$\vz(t)=\vz(t-1)$. Hence, $\vz(T)=\vz(T_{m-1})$, and we may replace $T_{m-1}$ by $T$.

If $t=T_\ell$ for some $\ell\ge1$, then $|A_t|<|A_{t-1}|$, so Case 2 in the proof
above applies. Furthermore, to have equality in \eqref{qu7}
we have to have equality in \eqref{qu5}. 
This means that if we again let $k=E(t)$, then $s_k(t)=\bs$, and thus by
\eqref{qu4} and the fact that $k$ is elected, we must have $s_i(t)=\bs$ for
all $i\in A_{t-1}$. 
Hence, by \eqref{zj}, $z_k(t)=\bs-1$ and $z_i(t)=\bs$ for all $i\in A_t$.
Furthermore, since we have equality in \eqref{qu2}, we have
$\bs=\sum_{i=m-\ell+1}^m\frac{1}{i}$. Note also that for $i\notin A_{t-1}$, 
$i$ has neither received any votes nor been elected in the interval
$(T_{\ell-1},T_\ell]$, and hence $z_i(T_\ell)=z_i(T_{\ell-1})$. 
It follows, by induction on $\ell$, that $\vz(T_\ell)$ is as $\vz(\ell)$ in
\refE{E1}. 
\end{proof}

In the time-homogeneous case, we have a sharper result.
(We will improve this further to an even sharper result in \refT{TR<<}.)
\begin{theorem}\label{TR<}
In the time-homogeneous case, if $m\ge3$,
then \eqref{t<} can be sharpened to
$z_j(t)<h_m-1$ for all $j$ and $t$.
\end{theorem}
\refE{E2} shows that we may have equality in the left inequality in
\eqref{t<} for any $m$. The same example shows also that we may have
equality in the right inequality for $m=2$.

\begin{proof}
In the time-homogeneous case, we see from the proof of \refT{TL=} that
if $z_j(T)=h_m-1$, then
$p_{A_{T-1}}=1$.
Since $|A_{T-1}|\le2$, this means that the probability vector $\vp$ has at
most two non-zero coefficients. If $m\ge3$, this means that $p_i=0$ for some
$i$.
But this is impossible, since then we can ignore $i$ and consider only the
remaining $m-1$ candidates, in which case \refT{T<} yields the stronger bound
$z_j(t)\le h_{m-1}-1$.
Consequently,  $z_j(T)=h_m-1$ is impossible in the time-homogeneous case when
$m\ge3$. 
\end{proof}

\begin{remark}
\citet{HLL} give weaker bounds for $z_j(t)$. They study also
whether Frege's method satisfies the  properties
\emph{variable upper quota}, meaning that for any time $t\ge 1$ we have
\begin{equation}\label{quo1}
r_j(t)\le \left \lceil\sum_{s=1}^t p_j(s) \right \rceil
\end{equation}
and 
\emph{variable lower quota}, which means
\begin{equation}\label{quo2}
r_j(t)\ge \left \lfloor \sum_{s=1}^t p_j(s) \right \rfloor.
\end{equation}
(In the time-homogeneous case these are called just
\emph{upper quota} and
\emph{lower quota}, see \cite{BY};
\cite{HLL} adds  \emph{variable} to the names to signify that
the probabilities $p_j(t)$ may be time-dependent.)
It is easily seen from \eqref{rj} that 
\eqref{quo1} is equivalent to $z_j(t)>-1$, and
\eqref{quo2} is equivalent to $z_j(t)<1$.
Hence, it follows from \eqref{b4<}--\eqref{b4>} and \eqref{b3} that, as shown by
\cite{HLL},
Frege's method satisfies variable upper quota for all 
$m$, but variable lower quota only for $m=2,3$.
It follows as special cases that in the time-homogeneous case, the method
satisfies upper quota for all $m$ and lower quota for $m=2,3$;
\cite[Example 5]{HLL} gives an example (with $m=6$)
that shows that lower quota does not hold in general.
See \refT{T>} below for a more precise result.
\end{remark}

\begin{remark}\label{Rleft}
  For the left inequality in \eqref{t<}, the argument in
  \eqref{a2}--\eqref{a22} shows that equality is possible only if 
all candidates tie, and thus $z_i(t)=1/m$ for all $i\neq j$,
or if this happened at some earlier time $t_1$ and $p_j(t')=0$ for
$t'\in(t_1,t]$. 
\end{remark}

\section{Optimality of the upper bound}\label{Sopt}
\subsection{Optimality for time-independent probabilities}\label{Snonconst}
The upper bound in \refT{T<} is  best possible by the simple \refE{E1}
if we allow
time-dependent votes.
It is less obvious that this upper bound in \refT{T<} is best possible also
in the time-homogeneous case; note that in this case, if $m\ge 3$, then
strict inequality holds  by \refT{TR<}.
In the present section we consider only the time-homogeneous case.
The following result shows that we in this case
may come arbitrarily close to the upper bound in \eqref{t<}.

\begin{theorem}\label{T>}
  Let $m\ge2$ and $\eps>0$. Then there exist (time-independent)
  probabilities $p_j>0$, $j\in[m]$, with $\sum_1^m p_j=1$
such that there exists a time $t$ with
\begin{align}\label{t>}
  z_j(t) > h_m-1-\eps
\end{align}
for some $j\in[m]$.

If $m=2$, we may moreover attain $z_j(t)=h_m-1=1/2$.
\end{theorem}

\begin{proof}
We assume
for convenience that ties are resolved by a deterministic rule;
more precisely, we assume that the candidate with highest number wins.
We may perturb slightly any such example with ties to give an example
without ties where \eqref{t>} still holds. Indeed, replace $p_j$ by
$p_j+\eta(j-(m+1)/2)$ for sufficiently small $\eta>0$.
Along the finite itinerary, strict score comparisons are preserved and
former ties become strict comparisons favouring the highest index. The terminal discrepancy
changes continuously, so \eqref{t>} is preserved.

We use induction on $m$, with the stronger induction hypothesis:
%
\begin{equation*}
\vbox{\em There exists $p_j$ as above, a time $t_1\ge1$ when $s_j(t_1)=1/m$ 
for all $j$, and a time $t=t_2\ge1$ such that \eqref{t>} holds for some $j$.}
\end{equation*}
The base case $m=2$ is trivial: we take $p_1=p_2=1/2$ and $t_1=t_2=1$
(as in the first round in \refE{E1} or \ref{E2}).
Note that in this case we achieve equality $z_1(1)=1/2=h_2-1$, as claimed.

Now assume that the induction hypothesis holds for some $m\ge2$.
It follows from \eqref{rj} and \eqref{Sj} that
\begin{align}\label{dix4}
  z_j(t)
\equiv  P_j(t)=tp_j \pmod \bbZ
\end{align}
for every $t$; thus, taking $t=mt_1$ and $t=t_1$, and recalling $s_j(t_1)=1/m$, 
\begin{align}\label{dix5}
  z_j(mt_1)\equiv mt_1p_j \equiv mz_j(t_1)
\equiv m s_j(t_1)
=1\pmod\bbZ,
\end{align}
\ie, $z_j(mt_1)$ is an integer for every $j$.
By the lower bound in \eqref{t<}, 
this implies $z_j(mt_1)\ge0$, and by \eqref{sumsum} we
conclude that $z_j(mt_1)=0$ for every $j$. This means that the election
becomes periodic, with elected candidates and the functions $z_j(t)$ and
$s_j(t)$ periodic with period $mt_1$.

Let $N$ be a large integer, 
and let $t'_1:=(Nm+1)t_1$.
The periodicity just shown and the induction hypothesis imply that
\begin{align}\label{dix6}
  s_j(t'_1)=s_j(Nmt_1+t_1)=s_j(t_1)=\frac{1}{m},
\qquad j\in[m].
\end{align}
Let further
\begin{align}\label{dix7}
  \eta:=\frac{1}{m(m+1)t'_1}
=\frac{1}{m(m+1)(Nm+1)t_1},
\end{align}
and define new probabilities
$p'_1,\dots,p'_{m+1}$ by
\begin{align}\label{dix8}
  p'_j:=
  \begin{cases}
    p_j-\eta, & 1\le j\le m,
\\ m\eta, & j=m+1.
  \end{cases}
\end{align}
We assume that $N$ is so large that all $p'_j>0$.
Let $s'_j(t)$ and $z'_j(t)$ denote the quantities in \eqref{sj}--\eqref{zj}
for $p'_1,\dots,p'_{m+1}$.

For $t<t'_1$, we have
\begin{align}\label{dix9}
s'_{m+1}(t)\le P'_{m+1}(t) < t'_1m\eta =\frac{1}{m+1},             
\end{align}
which by \eqref{a2} shows that $m+1$ cannot be elected before time $t'_1$.
Since furthermore $p_j-p'_j$ is constant for $j\in[m]$,
it follows that for $t=1,\dots,t'_1-1$, we have the same sequence of elected
as for $p_1,\dots,p_m$, with
\begin{align}\label{dix10}
  s'_j(t)=s_j(t)-t\eta,
\qquad \text{for } j\le m,\,  t\le t'_1.
\end{align}
In particular, by \eqref{dix6} and \eqref{dix7},
\begin{align}\label{dix11}
  s'_j(t'_1)=s_j(t'_1)-t'_1\eta=\frac{1}{m}-\frac{1}{m(m+1)}=\frac{1}{m+1},
\qquad j\le m.
\end{align}
By \eqref{sumsum}, we then also have $s'_{m+1}(t'_1)=1/(m+1)$,
so the first part of the induction hypothesis is verified for $(p'_j)_1^{m+1}$.

At time $t'_1$, we thus have a tie between all candidates.
By our assumption,  candidate $m+1$ wins. Thus
\begin{align}
  z'_{m+1}(t'_1)=s'_{m+1}(t'_1)-1=-\frac{m}{m+1}.
\end{align}
Let $t_2':=t_1'+t_2$, and assume that $N$ is so large that
$t_2m\eta<m/(m+1)$. Then, for every $t\in(t'_1,t'_2]$, we have 
\begin{align}
s'_{m+1}(t) \le z'_{m+1}(t'_1)+(t-t'_1)p'_{m+1} 
\le  z'_{m+1}(t'_1)+t_2m\eta 
<0,
\end{align}
and thus candidate $m+1$ will not win again until after $t'_2$.
At time $t'_1$, candidates $1,\dots,m$ start with the same value
$z'_j(t'_1)=1/(m+1)$, and it follows that for $t\in(t'_1,t'_2]$,
the same sequence is elected as starting from time 0 with $(p_j)_1^m$,
with
\begin{align}\label{ki1}
  z'_j(t'_1+s)=\frac{1}{m+1}+z_j(s)-s\eta,
\qquad j\in[m],\, 0\le s\le t_2.
\end{align}
In particular, for $j\le m$,
\begin{align}\label{ki2}
  z'_j(t'_2)=\frac{1}{m+1}+z_j(t_2)-t_2\eta.
\end{align}
By the induction hypothesis, there exists $j$ such that
$z_j(t_2)>h_m-1-\eps$. Hence, by choosing $N$
large and thus $\eta$ small enough,
also
$z_j(t_2)-t_2\eta>h_m-1-\eps$ and then, by \eqref{ki2},
\begin{align}\label{ki3}
  z'_j(t'_2)>\frac{1}{m+1}+h_m-1-\eps=h_{m+1}-1-\eps.
\end{align}
This completes the induction step, and thus the proof.
\end{proof}

\subsection{Conditions for being close to the supremum}\label{Sclose}
We saw in \refT{TR<} that for $m\ge3$, there is always strict inequality
in the bound \eqref{t<} for $z_j(t)$. 
We show in \refT{T>} that we can come arbitrarily
close. The construction in the proof of \refT{T>} uses probabilities $p_i$
with
$p_m \ll p_{m-1} \ll\dots\ll p_3 \ll p_2$, 
where $\ll$ informally means ``much smaller than''.
We show that this is not just convenient in our construction, but actually
necessary to get close to the upper bound.
\begin{theorem}\label{Tchar}
  For every $m\ge3$ there exists a constant $C_m$ such that:
If $\eps>0$ and 
\begin{align}\label{fb1}
  z_j(T) \ge h_m-1-\eps
\end{align}
for some $j\in[m]$ and $T\ge1$, 
then we may number the parties such that
\begin{align}\label{fb2}
  p_{i+1} \le C_m \eps p_i, 
\qquad 2\le i\le m-1.
\end{align}
\end{theorem}
The proof yields $C_m=2m(m-1)$; we do not know the best possible value.

Note that the conclusion \eqref{fb2} shows that $p_3,\dots,p_m$ are small
if $\eps$ is small, and in particular 
$p_i=O(\eps^{i-2})$ for $i\ge3$, but note also that
the theorem does not say anything about the relation between $p_1$ and $p_2$; in
fact, in the construction in the proof of \refT{T>}, we have
$p_1=p_2\approx 1/2$.  

\begin{proof}
We fix $m\ge3$.
Fix also a time $T\ge0$, and suppose that \eqref{fb1} holds.
By symmetry, we may assume $j=1$.

We argue as in the proof of \refT{T<}, with some further precision.
Note that we may assume that $\eps$ is small when needed, 
since \eqref{fb2} is trivial for $\eps\ge C_m\qw$
if  we order $p_2,\dots,p_m$ in decreasing order.

We use the notation in the proof of \refT{T<}, now with $j=1$; 
in particular, recall $A_t$ defined in \eqref{at}, and
that we have $A_{t}\subseteq A_{t-1}$ for all $1\le t \le T$, with
either $A_{t}=A_{t-1}$ or $A_{t}=A_{t-1}\setminus\set{E(t)}$.

For $i=2,\dots,m$, let $T_i$ be the last time in $[1,T]$ that $i$ was
elected; if $i$ never was elected, let $T_i:=0$.
Number the parties except 1 such that $T_2\ge T_3\ge\dots\ge T_m$;
in other words, the parties are numbered according to the last time $\le T$
that they were elected (in inverse order).
Then, for every $t\in[1,T]$, we have $A_t=\set{1,\dots,|A_t|}$. 

Define again $z_A(t)$, $s_A(t)$, and $p_A(t)$ by \eqref{qu1},
and define, for $0\le t\le T$,
\begin{align}\label{<qu2}
L(t):=
  \sum_{i=|A_t|+1}^m\frac{1}{i}-  \frac{z_{A_t}(t)}{|A_t|}
= h_m - h_{|A_t|}-  \frac{z_{A_t}(t)}{|A_t|}
,
\end{align}
and note that \eqref{qu2} shows that $L(t)\ge0$ for all $t$.
Our aim is to find a more precise bound.
Note that the assumption \eqref{fb1} is equivalent to, since $A_T=\set1$,
\begin{align}\label{fb6}
  L(T) \le \eps.
\end{align}

 For every $t\in[1,T]$, we have \eqref{qu3}, and we may sharpen
\eqref{qu6} to
\begin{align} \label{<qu6}
 s_{A_{t-1}}(t)
= z_{A_{t-1}}(t-1)+ p_{A_{t-1}}
= z_{A_{t-1}}(t-1)+1-p_{A\cp_{t-1}}
.\end{align}
We consider the same two cases as before:

\pfcase1{$A_t=A_{t-1}$.}
In this case, \eqref{qu3} and \eqref{<qu6} yield
\begin{align} \label{fb3}
 z_{A_{t}}(t)
= z_{A_{t-1}}(t-1)-p_{A\cp_{t-1}}
\end{align}
and thus
\begin{align}\label{fb4}
  L(t)=L(t-1)+\frac{p_{A\cp_{t-1}}}{|A_{t-1}|}
.\end{align}

\pfcase2{$A_t\subsetneq A_{t-1}$.}
In this case, 
$t=T_k$ for some $k=E(t)\ge2$, and, 
with the numbering chosen above,
$A_t=\set{1,\dots,k-1}$ and  $A_{t-1}=\set{1,\dots,k}$.
Let, cf.\ \eqref{qu4}, 
\begin{align}\label{<qu4}
\bsk:=\frac{s_{A_{t}}(t)}{|A_{t}|}
=\frac{s_{A_{t}}(t)}{k-1}
,\end{align}
the average of $s_j(t)$ for $j\in A_{t}$.
Since $k$ was elected in this round, we must have $s_k(t)\ge s_j(t)$ for
every $j$, and thus also 
\begin{align}
  \label{fb9}
s_k(t)\ge \bsk.
\end{align}
Furthermore, $E(t)=k\notin A_t$, and thus $z_j(t)=s_j(t)$ for every $j\in
A_t$ by \eqref{zj}. 
Hence, $z_{A_t}(t)=s_{A_t}(t)$ and thus, using \eqref{<qu2},
\begin{align}\label{fb7}
 \bsk= \frac{z_{A_{t}}(t)}{|A_{t}|}
=h_m - h_{k-1} - L(T_k) 
.\end{align}
Also, by \eqref{<qu4},
\begin{align}\label{<qu5}
\frac{k}{k-1}z_{A_t}(t)&
=\frac{k}{k-1}s_{A_t}(t)
=s_{A_t}(t)+\bsk
=s_{A_{t-1}}(t)-s_k(t)+\bsk
.
\end{align}
Consequently, by \eqref{<qu5}, \eqref{<qu6}  and \eqref{<qu2},
recalling $|A_{t}|=k-1$ and $|A_{t-1}|=k$,
\begin{align}\label{<qu7}
\frac{z_{A_t}(t)}{|A_t|}
&
=
\frac{s_{A_t}(t)}{k-1}
=
\frac{s_{A_{t-1}}(t)}k-\frac{s_k(t)-\bsk}k
\notag\\&
= \frac{1}{k}\bigpar{z_{A_{t-1}}(t-1) 
+1
-p_{A\cp_{t-1}}}
-\frac{s_k(t)-\bsk}{k}
\notag\\&
= \sum_{i=k+1}^m\frac{1}i -L(t-1)
+\frac{1}{k}
-\frac{p_{A\cp_{t-1}}}k
-\frac{s_k(t)-\bsk}{k}
\notag\\&
= \sum_{i=k}^m\frac{1}i
-L(t-1)
-\frac{p_{A\cp_{t-1}}}k
-\frac{s_k(t)-\bsk}{k}
,\end{align}
and thus
\begin{align}\label{<qu8}
L(t)
&
=
L(t-1)
+\frac{p_{A\cp_{t-1}}}k
+\frac{s_k(t)-\bsk}{k}
.\end{align}

Combining \eqref{fb4} and \eqref{<qu8}, 
and recalling \eqref{fb9},
we see that in both cases
\begin{align}\label{<qu88}
L(t)
-
L(t-1)
\ge
\frac{p_{A\cp_{t-1}}}{|A_{t-1}|}
\ge0.
\end{align}
In particular, $L(t)$ is a weakly increasing sequence.
Furthermore, $z_{A}(0)=0$ for every set $A$, and thus \eqref{<qu2} yields
\begin{align}\label{fa1}
  L(0)=\sum_{|A_0|+1}^m\frac{1}{i}
\ge0
.\end{align}
Consequently, \eqref{fb6} yields 
\begin{align}\label{fb5}
0\le L(t)\le L(T)\le \eps,
\qquad 0\le t \le T
.\end{align}
Moreover,
if $|A_0|<m$, then \eqref{fa1} yields $L(0)\ge 1/m$. 
If we assume, as we may, 
$\eps<1/m$, this is impossible by \eqref{fb5}, and thus
we must have $|A_0|=m$, meaning that every party, except possibly 1, is
elected at some time. Thus $1\le T_m < \dots <T_2 \le T$.

Next, \eqref{fb9}, \eqref{fb7} and \eqref{fb5} yield
\begin{align}\label{fb8}
s_k(T_k)
\ge \bsk
\ge h_m-h_{k-1}-\eps
.\end{align}
On the other hand, \eqref{<qu8} and \eqref{fb5} imply that
\begin{align}\label{fa2}
s_k(T_k)- \bsk
\le k L(T_k) \le k\eps,
\end{align}
and thus for every $j$, using \eqref{fb7},
\begin{align}\label{fa3}
  s_j(T_k)\le s_k(T_k) 
\le \bsk + k\eps
\le h_m-h_{k-1} + k\eps
.\end{align}
In particular, if $2\le k<m$, then \eqref{fa3} with $k$ replaced by $k+1$ yields
\begin{align}\label{fa4}
  s_k(T_{k+1})\le h_m-h_k+(k+1)\eps.
\end{align}
Consequently, 
\eqref{fb8} and \eqref{fa4} yield
\begin{align}\label{fa5}
  s_k(T_k)-s_k(T_{k+1}) &\ge h_m-h_{k-1}-\eps - \bigpar{h_m-h_k+(k+1)\eps}
\notag\\&
= \frac{1}k-(k+2)\eps
\ge \frac{1}{m-1}-(m+1)\eps
\ge \frac{1}{2m},
\end{align}
provided $\eps\le 1/(2m(m-1))$, which we may assume.
It follows from \eqref{sj}--\eqref{zj} that, for every $t\ge2$, 
\begin{align}\label{fc1}
s_k(t)-s_{k}(t-1)\le s_k(t)-z_k(t-1)=p_k.  
\end{align}
Hence, \eqref{fa5} implies
\begin{align}\label{fc2}
  (T_k-T_{k+1})p_k \ge s_k(T_k)-s_k(T_{k+1})\ge \frac{1}{2m}.
\end{align}
Furthermore, 
for $t\ge T_{k+1}$
we have $k+1\in A\cp_{t}$ and $|A_t|\le k$,
and thus \eqref{<qu88}
implies
\begin{align}\label{fc3}
  L(T_k)-L(T_{k+1})
\ge \sum_{t=T_{k+1}+1}^{T_k}  \frac{p_{A\cp_{t-1}}}{|A_{t-1}|}
\ge(T_k-T_{k+1})\frac{p_{k+1}}k
.\end{align}
Finally, we find from \eqref{fc2}--\eqref{fc3} and \eqref{fb5}, for every
$2\le k\le m-1$,
\begin{align}
  \frac{p_{k+1}}{p_k}&
=  \frac{(T_k-T_{k+1})p_{k+1}}{(T_k-T_{k+1})p_k}
\le \frac{k(L(T_k)-L(T_{k+1}))}{1/(2m)}\le 2mk\eps
\notag\\&
\le 2m(m-1)\eps,
\end{align}
showing \eqref{fb2}. 
\end{proof}

We now reach another main result: 
we consider random probabilities $p_j$ and
obtain an estimate of the probability of ever having a
value $z_j(t)$ close to the upper bound. Together with \refT{Tchar} this
next result also gives an indication  
of the geometry of an attractor for large $m$. We begin our study of the existence and the properties of the
attractor in the next section. 
\begin{theorem}\label{Tchar2}
  Fix $m\ge3$ and let $\vp$ be uniformly random in $\gS_m$.
Then, as $\eps\to0$, for every $j\in[m]$,
\begin{align}\label{cchar}
  \P\Bigpar{\sup_{t\ge1}z_j(t) \ge h_m-1-\eps}
= O\bigpar{\eps^{\kk_m}}
,\end{align}
where $\kk_m=\binom{m-1}2$.
\end{theorem}

\begin{proof}
  Suppose that $\sup_t z_j(t)\ge h_m-1-\eps$.
For every $\eta>0$, \refT{Tchar} applies with $\eps+\eta$ in place of $\eps$
and some $T\ge1$.
Letting $\eta\downarrow0$ along a subsequence with a fixed numbering of the
parties yields \eqref{fb2}. Since there is only a fixed
  number $m!$ of numberings, we may consider each numbering separately, 
and estimate the probability that \eqref{fb2} holds for a given numbering.
By induction, \eqref{fb2} implies $p_i \le C_{m}^{i-2} \eps^{i-2}$ for $i\ge3$.
Since $(p_2,\dots,p_m)$ is uniformly distributed in the simplex 
$\set{(p_i)_2^m:p_i\ge0\;\forall i\ge2\text{ and }\sum_2^m p_i\le1}$,
which has volume $1/(m-1)!$, the density function is $(m-1)!$ in this simplex,
this yields
\begin{align}
  \P\bigpar{\text{\eqref{fb2} holds}}
\le (m-1)! \prod_{i=3}^m C_m^{i-2}\eps^{i-2}
= C_m' \eps^{\binom {m-1}2}
\end{align}
for some constant $C_m'$.
\end{proof}


We assume that the value for $\kk_m$ given in \refT{Tchar2} is not sharp,
but we do not know any better value.

\begin{problem}\label{Pchar}
  Find the best possible exponent $\kk_m$ for $m\ge3$.
\end{problem}

As another corollary, we may sharpen \refT{TR<},
which says that if $m\ge3$, then we cannot have equality in the upper bound
in \eqref{t<} for any $t$. 
The construction in the proof of \refT{T>} coming
close to the upper bound $h_m-1$
uses rational $p_j$; then the sequence $z_j(t)$ is periodic and therefore
also
$\sup_t z_j(t)$ is strictly less than $h_m-1$. However, this leaves the
possibility that for some irrational $p_j$ it might be possible to find a
sequence of times $t_n$ such that $z_j(t_n)\to h_m-1$.
We now show that this is impossible.

\begin{theorem}\label{TR<<}
In the time-homogeneous case, if $m\ge3$,
then \eqref{t<} can be sharpened to
$\sup_{j\in[m],\,t\ge 1}z_j(t)<h_m-1$.
\end{theorem}

\begin{proof}
  Suppose that for some $j$, $\sup_t z_j(t)=h_m-1$. 
Then, for every $\eps>0$ there is a time $T$ such that \eqref{fb1} holds,
and thus \eqref{fb2} holds for some numbering of the parties. 
In particular,
\begin{align}
  \min_{1\le i\le m} p_i \le p_3 \le C_m \eps p_2 \le C_m \eps.
\end{align}
Since $\eps$ is arbitrary, this yields $\min_i p_i=0$, so $p_i=0$ for some $i$.
But in this case we only have to consider the remaining $m-1$ parties,
and thus \refT{T<} yields $\sup_t z_j(t)\le h_{m-1}-1<h_m-1$.
\end{proof}

\section{Orbits, attractor and invariant measure for the dynamical system}\label{Sorb}
We consider in the remainder of the paper the dynamical system in
\refSS{SSdynamical} with  a given (time-independent) 
probability vector $\vp$ and an arbitrary starting point
$\vz(0)=\vu$.
Unless we say otherwise, we may consider either the multi-valued version of $T$
or a single-valued version defined by some fixed rule for resolving ties,
see \refR{Rmulti}.
To avoid some trivialities, we assume from now on that $p_j>0$
for every $j\in[m]$. 
Denote the \emph{orbit} of $\vu\in\bbR^m$ by
\begin{align}\label{Ou}
  \cO(\vu):=\set{T^n\vu:n\ge0}=\set{\vz(n;\vu):n\ge0}.
\end{align}
When we consider the multi-valued version of $T$,
this is interpreted as 
\begin{align}\label{Ou2}
  \cO(\vu):=\bigcup_{n\ge0}T^n(\vu).
\end{align}
By \eqref{tn1}, the orbit $\cO(\vu)$ is contained in the
hyperplane
\begin{align}\label{Hu}
H_\vu = \lrset{x\in\bbR^m: \sum_{i=1}^m x_i = \sum_{i=1}^m u_i}.
\end{align}

The
\emph{limit set} 
(also called \emph{$\go$-limit set}) 
of $\vu\in\bbR^m$ is
\begin{align}\label{Lu}
  \cL(\vu) := \bigcap_{n=0}^\infty\overline{\set{\vz(t;\vu):t\ge n}}
=\bigcap_{n=0}^\infty\overline{\cO(T^n \vu)}.
\end{align}
The limit set is thus a closed set.
We have $\cL(\vu)\subseteq H_\vu$, again by \eqref{tn2}.

Recall from \refR{Rmulti} that
in the multi-valued case,
we may also consider a single \emph{branch}
$\vz'=(\vz'(n,\vu))_{n\ge0}$ satisfying \eqref{branch}.
We then define the branch orbit $\cO(\vz')$
and branch limit set $\omega(\vz')$ 
in the same way as \eqref{Ou} and \eqref{Lu},
i.e.,
\begin{align}
\cO(\vz') :=\set{\vz'(n,\vu):n\ge0},
&&&
\omega(\vz'):=\bigcap_{n=0}^\infty\overline{\set{\vz'(t;\vu):t\ge n}}.
\end{align}
Note that 
$\cO(\vz')\subseteq\cO(\vu)$
and
$\go(\vz')\subseteq\cL(\vu)$.

 A set $A\subseteq \bbR^m$ is 
\emph{forward invariant}  if $T(A)\subseteq A$,
\emph{backward invariant}  if $A\subseteq T(A)$,
and
\emph{invariant} if it is both forward and backward invariant, i.e.,
if $T(A)=A$. 

Note that the definition \eqref{Ou} shows that 
\begin{align}\label{cO2}
T(\cO(\vu))=\cO(T\vu)\subseteq\cO(\vu),  
\end{align}
\ie, $\cO(\vu)$ is forward invariant.
For a continuous dynamical system, it follows immediately that also the
limit set is forward invariant, but $T$ is not continuous everywhere (see
\refR{RD}).
We will show later (\refC{CTX}) that, nevertheless,
the limit set $\cL(\vu)$ is forward invariant
in the ergodic case (defined in \refS{Sergodic}).
Outside this case, forward invariance can fail even for the multi-valued
map, as \refE{Einter1} shows.

We will prove backward invariance for the multi-valued version
in Lemma \ref{Lback} below.

\begin{remark}\label{Rnoback}
  The simple \refE{Em=2I}
shows that 
for a single-valued version of $T$,
$\cL(\vu)$ is in general (and presumably typically)
\emph{not} backward invariant.
This is one instance where the multi-valued version has nicer properties. 
\end{remark}

A set $A\subseteq\Ho$ is \emph{absorbing} 
(in $\Ho$) 
for $T$ if for every initial value
$\vu\in\Ho$ there exists $N\ge0$ such that
\begin{align}
  T^n\vu \subseteq A \qquad \text{for all }n\ge N,
\end{align}
\ie, equivalently, $\cO(T^N\vu)\subseteq A$ for some $N\ge0$.
Thus every orbit eventually stays in $A$.
If $A$ is a closed absorbing set, then $\cL(\vu)\subseteq A$ for every
$\vu\in\Ho$.

A compact invariant set $K\subseteq\Ho$ is called a \emph{global attractor}
if it is absorbing. 
For a general dynamic system, it is not obvious that such a set exists, nor
that it is  unique, but
we will show later (\refT{Tinv})
that for the Frege map,
at least in the ergodic case, 
a unique global attractor exists; we then may call it
\emph{the global attractor}.

\subsection{Orbits are bounded}\label{Sorbit}

We  show that every orbit is bounded by extending the estimates
\eqref{t<} to arbitrary initial points $\vu$. 
\begin{theorem}\label{Tu}
Suppose $\vz(0)=\vu=(u_j)_1^m\in\bbR^m$.
\begin{romenumerate}
  
\item \label{Tu1}
Then,
  for every $j\in[m]$ and $t\ge0$,
  \begin{align}\label{tu1}
-1+\frac{1}{m} + \min_j u_j
\le    z_j(t) 
\le h_m-1 + \max_j u_j.
  \end{align}
\item \label{Tu2}
If $\vu\in\Ho$, then for every $j\in[m]$ and all sufficiently large $t$,
\begin{align}\label{tu2}
  z_j(t)\ge -1+\frac{1}{m}.
\end{align}
Furthermore, if $K\subset\Ho$ is compact, then there exists $t_0$ such that
\eqref{tu2} holds for all $\vu\in K$ and all $t\ge t_0$.
\item \label{Tu3}
Moreover, if $\vu\in\Ho$ and \eqref{tu2} holds for $t=0$ and all
$j\in[m]$,
then \eqref{tu2} holds for all $t\ge0$ and all $j\in[m]$.
\end{romenumerate}
\end{theorem}
\begin{proof}
\pfitemref{Tu1}
  This follows by simple modifications of the proof of \refT{T<}.
It suffices to prove the bounds for every branch of the multi-valued map.
Recall that we may translate such a branch by any vector of the form $a\1$, see
\eqref{tna}.

For the left inequality, we may by translation assume that $\min_j u_j=0$.
Then \eqref{tn1} implies that $\innprod{\vz(t),\1}=\innprod{\vz(0),\1}\ge0$
and thus $\innprod{\vs(t),\1}\ge1$ for every $t$. Consequently, we still have
$s_{E(t)}(t)\ge 1/m$ as in \eqref{a2}, and the proof is concluded by
induction as before, noting that the case $t=0$ is trivial.

For the right inequality, we may by translation assume that $\max_ju_j=0$.
Thus $z_j(0)=u_j\le0$ for every $j$, and it follows that $z_A(0)\le0$ for
every set $A\subset[m]$. Consequently \eqref{qu2} holds for $t=0$.
The induction step in the proof of \refT{T<} then holds without any changes,
and consequently \eqref{qu2} holds and thus $z_j(T)\le h_m-1$.

\pfitemref{Tu2}
Note first that every party $j$ has to win sometimes, because otherwise
$z_j(t)=z_j(0)+tp_j$ would increase to infinity, contradicting \eqref{tu1}.
Moreover, if $K$ is compact, this holds uniformly for $\vz(0)=\vu\in K$:
there exists $t_0$ such that every party wins at least 
once before $t_0$.

The assumption $\vu\in\Ho$ implies that $\innprod{\vz(t),\1}=0$
and $\innprod{\vs(t),\1}=1$ for every $t$, and thus \eqref{a2} holds.
It follows by induction that $z_j(t)\ge \frac{1}{m}-1$ for all times $t$
after the first time party $j$ wins, and hence for all $t\ge t_0$.

\pfitemref{Tu3}
As in \ref{Tu2}, now by induction on $t\ge0$.
\end{proof}

\begin{corollary}\label{Cu}
The closure $\overline{\cO(\vu)}$ of the orbit is a compact subset of $H_\vu$
for every $\vu\in\bbR^m$.
\end{corollary}
\begin{proof}
  Obvious from \eqref{tu1} and \eqref{tn1}.
\end{proof}

\begin{corollary}\label{CLu}
The limit set $\cL(\vu)$  is a non-empty compact subset of $H_\vu$
for every $\vu\in\bbR^m$.
\end{corollary}
\begin{proof}
By \refC{Cu} and \eqref{cO2},
$\overline{\cO(T^n \vu)}$ is a decreasing sequence of non-empty
compact subsets of $H_\vu$;
thus the result follows from \eqref{Lu}.
\end{proof}

Of course, the compact set in \refC{Cu} depends on the starting point $\vu$,
but we show next that there is a fixed compact set $K_0$ that contains
$\cL(\vu)$ for every $\vu \in \Ho$. 
Moreover, $K_0$ is forward invariant and absorbing.

\begin{lemma}\label{LK0}
  There exists a compact forward invariant absorbing set $K_0\subseteq \Ho$.
In fact, we may take
\begin{align}\label{jb3}
  K_0 := \Ho\cap [-1+\tfrac1m,m-1]^m.
\end{align}
In particular, $\cL(\vu)\subseteq K_0$ for every $\vu\in\Ho$.
\end{lemma}
\begin{proof}
Let $\vu\in\Ho$.
By \refT{Tu}\ref{Tu2},
there exists $t_0$ such that
\begin{align}\label{jb1}
  z_j(t) \ge -1+\frac{1}{m}>-1,
\qquad j\in[m],\, t\ge t_0.
\end{align}
Since $\vz(t)\in \Ho$, 
it follows from \eqref{jb1} that also 
\begin{align}\label{jb2}
  z_j(t) < m-1,
\qquad j\in[m],\, t\ge t_0.
\end{align}
Consequently, 
$\vz(t)\in K_0$ 
for $t\ge t_0$.
Hence, $K_0$ is absorbing.
As noted above, this implies $\cL(\vu)\subseteq K_0$ for $\vu\in\Ho$.

Moreover, if $\vu\in K_0$, then \refT{Tu}\ref{Tu3} shows that we may take
$t_0=0$ above. In particular, $T\vu\in K_0$.
Thus $K_0$ is forward invariant.
\end{proof}

We next prove backward invariance of the limit sets for the multi-valued
version; recall that this may fail for single-valued versions
(\refR{Rnoback} and \refE{Em=2I}).

\begin{lemma}\label{Lback}
  Consider the  multi-valued version of $T$. Then the limit set $\cL(\vu)$
  is backward invariant for every $\vu\in\bbR^m$.
\end{lemma}

\begin{proof}
Suppose that $\vx\in\cL(\vu)$. Let $n\ge0$. Then $\vx\in \bcO{T^{n+1}\vu}$.
By \refC{Cu} and \eqref{compact}, 
$T(\bcO{T^n\vu})$ is compact, and since it contains
$T(\cO(T^n\vu))=\cO(T^{n+1}\vu)$, it follows that 
\begin{align}
  \vx\in \bcO{T^{n+1}\vu}
= \overline{T(\cO(T^n\vu))}
\subseteq
T(\bcO{T^n\vu}).
\end{align}
Thus, there exists $\vy_n\in \bcO{T^n\vu}$ such that $\vx\in T\vy_n$.
By \eqref{tx}, this means that $\vy_n=\vx-\vp+\ve_i$ for some $i\in[m]$;
thus there is only a finite number of possible $\vy_n$, and we may choose
$\vy$ as one of them such that $\vy_n=\vy$ for infinitely many $n$.
Then $\vx\in T\vy$. Moreover, $\vy\in\bcO{T^n\vu}$ for arbitrarily large $n$,
and since the sequence of sets $\bcO{T^n\vu}$ is decreasing, 
it follows that $\vy\in\bcO{T^n\vu}$ for all $n$, and thus $\vy\in\cL(\vu)$.
This proves that $\vx\in T\cL(\vu)$, and thus $\cL(\vu)\subseteq T\cL(\vu)$.
\end{proof}

\subsection{Projecting onto $\bbT^m$}\label{SSproj}
Let $\pi$ be the projection $\vx \mapsto \hvx := \vx \mod \Z^m$ from $\R^m$ onto the
torus $\T^m = \R^m / \Z^m$. The image $\pi H_\vu$ is the submanifold  
\begin{align}\label{Su}
\hH_\vu=
\biggset{(\hx_i)_1^m\in\bbT^m: \sum_{i=1}^m \hx_i =  \sumim u_i \pmod{ \Z}}
\subset\bbT^m. 
\end{align}
In particular, 
\begin{align}\label{Su0}
\hHo=
\biggset{(\hx_i)_1^m\in\bbT^m: \sum_{i=1}^m \hx_i = 0}.
\end{align}

Let $\Rp:\bbT^m\to\bbT^m$ be the translation (rotation) $\hvx\mapsto\hvx+\hvp$
in $\bbT^m$, where $\hvp:=\pi(\vp)$.
Since $\pi(\ve_j)=\0$ for every $j$, it is obvious from \eqref{tx} that the
operator $T$ projects to $\Rp$, in other words
\begin{align}\label{piT}
  \pi \bigpar{T(\vx)} = \Rp\bigpar{\pi(\vx)},
\qquad \vx\in\bbR^m.
\end{align}
Hence, the projection of the dynamical system to the torus 
$\bbT^m$ has a very simple,
and well-known, structure. 
(Note also that the problem of ties disappears in the projection.)

In particular, the orbit $\cO(\vu)\subset\bbR^m$ of $\vu$ under $T$ 
projects to the orbit $\gO(\hvu)\subseteq\bbT^m$ 
of $\hvu:=\pi(\vu)$ under $\Rp$. Note that
$\gO(\hvu)\subseteq \hH_\hvu\subset \bbT^m$.  
The closure $\bgO{\hvu}$ has a well-known simple structure.
We collect some simple facts  in a lemma; 
these results are partly well-known,
but we give proofs for completeness;
see also \cite[Section 2]{gun}.
We let $\0$ denote also the zero in $\bbT^m$,
and we denote the limit set of $\hvu\in\bbT^m$ under $\Rp$ by
\begin{align}
  \label{hcL}
\hcL(\hvu):=\bigcap_{n\ge0}\bgO{\Rp^n\hvu}.
\end{align}

\begin{lemma}\label{LP}
Let $\vu\in\bbR^m$ and let $\hvu:=\pi(\vu)\in\bbT^m$.
Then
\begin{align}\label{lp}
  \pi\bigpar{\bcO{\vu}}
=\pi\bigpar{\cL(\vu)}
=\bgO{\hvu}
=\hcL(\hvu).
\end{align}
The projected limit set in \eqref{lp} is non-empty and compact.
Furthermore, let
\begin{align}\label{G}
G:=\bgO{\0}=\hcL(\0). 
\end{align}
Then $G$ is a closed subgroup of\/ $\hHo\subset\bbT^m$
which is invariant for $\Rp$: 
$\Rp G = G$; furthermore, $G$ is a minimal closed invariant subset of\/ 
$\bbT^m$.
For every $\vu\in\bbR^m$, $\hcL(\hvu)$ is the coset
\begin{align}\label{lp2}
  \pi\bigpar{\cL(\vu)}=\hcL(\hvu)=G+\hvu.
\end{align}
\end{lemma}

\begin{proof}
First, consider the action of $\Rp$ on $\bbT^m$.
It is easy to see that every neighbourhood of $\0$ in $\bbT^m$ contains
$\Rp^n(\0)=n\hvp$ for arbitrarily large $n$, and thus 
\eqref{hcL} implies
\begin{align}
  \label{hcl0}
\0\in\hcL(\0).
\end{align}
Since $\Rp$ is continuous, it is easily seen that
$\hcL(\hvu)$ is forward invariant, \ie,
$\Rp(\hcL(\hvu))\subseteq\hcL(\hvu)$ for every $\hvu\in\bbT^m$;
thus it follows from \eqref{hcl0} that $\Rp^n(\0)\in\hcL(\0)$ for every $n\ge0$.
Hence, $\gO(\0)=\set{\Rp^n(\0):n\ge0}\subseteq\hcL(\0)$,
which together with \eqref{Lu} yields $\hcL(\0)=\bgO{\0}$
as claimed in \eqref{G}.
Consequently, we can define $G$ by \eqref{G}.

Furthermore, it follows from \eqref{hcl0} and \eqref{hcL} 
that for every $n\ge0$, we have
$\0\in\hcL(\0)=\hcL(\Rp^n\0)$. 
Thus by translation invariance,
$-n\hvp=\Rp^{-n}\0\in\hcL(\0)$.
Consequently, $\set{n\hvp:n\in\bbZ}\subseteq\hcL(\0)$, and thus
\begin{align}\label{lp3}
G=
\bgO{\0}=
  \hcL(\0)
=\overline{\set{n\hvp:n\in\bbZ}}.
\end{align}
It follows that $G=\hcL(\0)$ is a closed subgroup of $\bbT^m$, 
and $\Rp G = G$. Since every non-empty closed invariant subset has to contain
$\hvx+\set{n\hvp:n\ge0}$ for some $\hvx\in G$, 
it follows from \eqref{G} that $G$ is minimal.
Translation invariance and \eqref{G} yield
\begin{align}\label{lp4}
\bgO{\hvu}=
\hcL(\hvu)=\hcL(\0)+\hvu=G+\hvu,
\qquad \hvu\in\bbT^m.   
\end{align}

We now go back to $T$, and an arbitrary $\vu\in\bbR^m$.
It follows easily from 
$\pi\bigpar{\cO\xpar{\vu}}=\gO\xpar{\hvu}$ and
the compactness of $\bcO{\vu}$ (\refC{Cu}) that
\begin{align}\label{lp4b}
\pi\bigpar{\bcO{\vu}}=\bgO{\hvu}.  
\end{align}
Hence, using also \eqref{Lu} and \eqref{lp4},
\begin{align}\label{lp5}
\pi\bigpar{\cL(\vu)}\subseteq \pi\bigpar{\bcO{T^n\vu}}=\bgO{\Rp^n\hvu}
=G+n\hvp+\hvu=G+\hvu=\hcL(\hvu).  
\end{align}
Conversely, 
if $\hvx\in\hcL(\hvu)$, then for every $n\ge0$, we have by \eqref{hcL}
and \eqref{lp4b}
\begin{align}\label{lp6}
  \hvx\in \bgO{\Rp^n\hvu}=\pi\bigpar{\bcO{T^n\vu}}.
\end{align}
Hence, using \eqref{Lu} again,
\begin{align}
  \pi\qw(\hvx)\cap\cL(\vu)
=  \bigcap_{n\ge0}\Bigpar{\pi\qw(\hvx)\cap\bcO{T^n\vu}}
\neq\emptyset,
\end{align}
since it
is the intersection of a decreasing sequence of non-empty
compact sets. Consequently, $\hvx\in\pi\bigpar{\cL(\vu)}$, and thus
$\hcL(\hvu)\subseteq\pi\bigpar{\cL(\vu)}$.

Thus $\pi\bigpar{\cL(\vu)}=\hcL(\hvu)$, 
which together with \eqref{lp4} and \eqref{lp4b} shows both \eqref{lp}
and \eqref{lp2}.
Finally,
that $\hcL(\hvu)$ is compact and non-empty 
follows from \eqref{lp2}
(or by \eqref{hcL} and the argument for  \refC{CLu}).
\end{proof}

It follows that 
every coset of $G$ also
is a minimal closed invariant set for $\Rp$;
furthermore, $\bbT^m$ is the disjoint union of
these compact invariant sets.
Moreover, $G$ is, as every closed subgroup of $\bbT^m$, 
isomorphic to a direct product $F\times\bbT^\mu$ for some finite Abelian
group $F$ and some $0\le\mu\le m$.
In our case, this can be seen more explicitly as follows.
We say that a finite indexed family $(a_1,\dots,a_k)$ of real numbers is
\emph{rationally independent} if it is linearly independent over $\bbQ$,
or equivalently, if
\begin{align}\label{pi3}
  \sum_{i=1}^kn_ia_i=0
\end{align}
for integers $n_i\in\bbZ$ implies $n_i=0$ for all $i$.
We let in the sequel $\mu$ denote the number
\begin{align}\label{mu}
  \mu:=
\max\bigset{|I|:I\subset[m] \text{ with }
\bigl((p_i)_{i\in I},1\bigr) \text{ rationally independent}}.
\end{align}
Since $1=\sum_i p_i$, it is easy to see that 
\begin{align}\label{mu2}
\mu+1 = \dim_\bbQ \bigpar{\bbQ p_1+\dots+\bbQ p_m}
,\end{align}
the dimension of the $\bbQ$-vector space spanned by $p_1,\dots,p_m$.   
In particular, $0\le \mu\le m-1$.
 
\begin{lemma}\label{LG}
Let $I\subset[m]$ be a set of indices achieving the maximum in \eqref{mu}.
Then the projection $\Pi_I:\bbT^m\to\bbT^\mu$ 
onto the $\mu$ coordinates in $I$,
\ie, $\Pi_I(\hx_i)_1^m:=(\hx_i)_{i\in I}$,
maps $G$ onto $\bbT^\mu$. Moreover, if $G_0$ is the connected component of $G$
containing $\0$, then $G_0\cong\bbT^\mu$, the projection $\Pi_I$ maps $G_0$
onto $\bbT^\mu$ (but is not necessarily a bijection), and $G\cong F\times
G_0\cong F\times\bbT^\mu$ for some finite cyclic group $F$.
\end{lemma}

\begin{proof}
This is well known,
but for completeness we sketch a proof.
We assume for notational convenience that $I=\set{1,\dots,\mu}$,
so $\Pi_I$ is the projection onto the first $\mu$ coordinates. 
The assumption means that Weyl's criterion applies to
$\Pi_I(\hvp)$, which implies that the sequence $(n\Pi_I(\hvp))$ is
equidistributed
in $\bbT^\mu$. In particular, this sequence is dense in $\bbT^\mu$. 
Moreover, every $p_j$ is a rational linear combination of $p_1,\dots,p_\mu,1$.
This implies that there exists an integer $N\ge1$ such that
for every $j>\mu$ and some integers $n_{ji}\in\bbZ$,
\begin{align}
  Np_j = \sum_{i=1}^\mu n_{ji}p_i + n_{j0}.
\end{align}
It follows that
\begin{align}
\hvp\in\bbTxm:=\lrset{(x_i)_1^m\in\bbT^m: Nx_j= \sum_{i=1}^\mu n_{ji}x_i,
\qquad j=\mu+1,\dots,m}.
\end{align}
Since $\bbTxm$ is a closed subgroup of $\bbT^m$, it follows that
$G\subseteq\bbTxm$. Hence, the kernel $K_1$ of the projection
$\Pi_I:G\to\bbT^\mu$ is finite.
Let $G_0$ be the connected component of $\0$ in $G$.
Since a closed subgroup of a torus has finitely many connected components,
$G/G_0$ is finite. The subgroup $\Pi_I(G_0)$ is compact and hence closed.
Since $\Pi_I$ maps $G$ onto $\bbT^\mu$, the quotient
$\bbT^\mu/\Pi_I(G_0)$ is a homomorphic image of $G/G_0$, and is therefore
finite. It is also connected, being a quotient of the connected torus
$\bbT^\mu$, so it is trivial. Thus $\Pi_I(G_0)=\bbT^\mu$.
The projection $\Pi_I:G_0\to\bbT^\mu$ has finite kernel
$K_2\subseteq K_1$ and is therefore a finite covering and a local
isomorphism. It follows that $G_0\cong\bbT^\mu$.
Finally, the quotient $G/G_0$ is isomorphic to $K_1/K_2$ and is thus also
a finite Abelian group, and thus a (finite) product of finite cyclic groups.
Since $G_0\cong \bbT^\mu$ is a divisible group, it follows easily
that there is a finite subgroup $F\subseteq G$ such that $G\cong F\times G_0$.
Finally, by projecting $G$ to $F$, it follows from \eqref{lp3} that $F$ is
a cyclic group generated by the projection of $\hvp$ to $F$.
\end{proof}

The most interesting cases are the two extreme cases $\mu=0$ and $\mu=m-1$,
which are  discussed in the following two sections.
The general case combines features from these extreme cases;
we give some comments and simple examples in \refS{Sinter}.

\section{Orbits in the rational case}
\label{Srational}

One extreme case is  $\mu=0$. By \eqref{mu}, this holds if and only if
\emph{all $p_j$ are rational}, and 
we call this the \emph{rational case}. 
By \refL{LG}, in this (and only this) case, $G$ is a finite (cyclic) group;
equivalently, the orbit $\gO(\hvu)$ is finite for every $\hvu\in\bbT^m$.
Furthermore, for every $\vu\in\bbR^m$,
\refC{Cu} implies that 
$\cO(\vu)\cap\pi\qw(\vx)$ is finite
for every $\hvx\in \gO(\pi(\vu))$,
and thus every orbit $\cO(\vu)$ is finite. 
It follows that if we consider $T$ as a single-valued operator with any
fixed rule for ties, then for every initial point $\vu$, the dynamical system
$T^n\vu$ is ultimately periodic, with eventual period (cycle length) $|\cL(\vu)|$.
More explicitly, we have the following.

\begin{theorem}\label{TQ}
  Suppose that $p_1,\dots,p_m$ all are positive rationals, and write
  them on a least common denominator as
  \begin{align}\label{tq}
    \vp=\lrpar{\frac{q_1}{Q},\dots,\frac{q_m}Q}
  \end{align}
for some positive integers $q_1,\dots,q_m$ and $Q$, with 
$\gcd(q_1,\dots,q_m)=1$
and $Q=\sum_i q_i$.
Consider a single-valued version of $T$.
Then, if we start with $\vz(0)=\0$, the sequence $\vz_n=T^n\0$ is periodic
with period exactly $Q$. More generally, the same holds for any initial
value $\vu\in[0,\frac1m)^m$.

Furthermore, 
for every initial value,
the projected 
sequence $\hvz_n:=\pi(\vz_n)$ in $\bbT^m$ has period $Q$.
\end{theorem}

\begin{proof}
  Consider first $\Rp$. We have, for any $\hvu\in\bbT^m$,
  \begin{align}\label{tq1}
 \Rp^n\hvu=\hvu
&\iff n\hvp =\0\text{ in }\bbT^m
\iff n\vp\in\bbZ^m
\iff \frac{nq_i}{Q}\in\bbZ\;\forall i\in[m]
\notag\\&
\iff
Q\mid\gcd(nq_1,\dots,nq_m)
.  \end{align}
Furthermore,
\begin{align}\label{tq2}
  \gcd(nq_1,\dots,nq_m)=n\gcd(q_1,\dots,q_m)=n,
\end{align}
and thus $\Rp^n\hvu=\hvu$ if and only if $n$ is a multiple of $Q$.

It follows that, for any $\vu\in\bbR^m$, 
$Q|n$ is a necessary condition for $T^n\vu=\vu$. 
Conversely, for any $\vu$, \eqref{tx} shows that
\begin{align}\label{tq3}
  z_i(Q;\vu)\in u_i + Qp_i + \bbZ = u_i+\bbZ.
\end{align}
If $u_i\in[0,\frac{1}m)$ for every $i$, then
\refT{Tu} shows that $z_i(Q;\vu)\ge -1+1/m>-1+u_i$, and since \eqref{tq3} shows
that
$z_i(Q;\vu)-u_i$ is an integer, we have $z_i(Q;\vu)- u_i\ge0$.
On the other hand, $\sum_i z_i(Q;\vu)=\sum_i u_i$ by \eqref{tn2}, and thus
we must have $z_i(Q;\vu)- u_i=0$ for all $i$, \ie, $T^Q\vu=\vz(Q;\vu)=\vu$.
\end{proof}

\begin{remark}\label{RQ}
  It seems natural to consider a single-valued version of $T$ here, but we may
also consider the multi-valued version. 
For it, it follows by the same argument that 
$T^Q\0=\set\0$, and thus $T$ has period exactly $Q$ also as a set-valued
mapping, if we start with $\vu=\0$ (or $\vu\in[0,\frac{1}m)^m$).
\end{remark}

\begin{problem}\label{PQ}
Under the rationality assumptions of \refT{TQ}, fix a single-valued version
of $T$. We conjecture that for every initial value $\vu\in\bbR^m$, the
sequence $T^n\vu$ is ultimately periodic with least eventual period $Q$.
We have been unable to prove this in general and leave it as an open problem.
\end{problem}

As a substitute,
we state  the following weaker result.

\begin{theorem}
  Suppose that $p_1,\dots,p_m$ all are positive rationals, and write
  them on a least common denominator $Q$ as in \eqref{tq}.
Fix a single-valued version of $T$.
Then, for every initial point $\vu\in\bbR^m$, 
the orbit
$\vz(n;\vu)=T^n\vu$ is ultimately periodic with some period $P$ such that
\begin{align}\label{fq1}
Q\mid P.
\end{align}
Moreover, the projected orbit on the torus $\bbT^m$ has exact period $Q$.
\end{theorem}

\begin{proof}
%
The first part of the proof of \refT{TQ} is valid for every $\vu$, and shows
that the projected orbit in $\bbT^m$ has period $Q$.

Moreover, every orbit of $T$ is bounded. Since the projected orbit in
$\bbT^m$ is finite, and each fibre of the projection intersects a bounded subset
of $\bbR^m$ in only finitely many points, it follows that the orbit
\begin{align}
\{T^n\vu:n\ge0\}
\end{align}
is finite. Therefore the sequence $\vz(n;\vu)$ is ultimately periodic with
some period $P$.

If $T^{n+P}\vu=T^n\vu$ for all sufficiently large $n$, then the projected
orbit also repeats with period $P$. Since the projected orbit has exact period
$Q$, it follows that $Q$ divides $P$.
\end{proof}

\begin{remark}
The divisibility \eqref{fq1} comes only from the torus factor. To deduce the
stronger conclusion $P=Q$, one would need to control
the $Q$-step displacement
\begin{align}
T^Qx-x\in\bbZ^m
\end{align}
along the eventual cycle, with the aim of showing that it vanishes. 
Note that this would hold if $T$
commuted with integer translations; however,  
the Frege map does not have this property:
shifting $x$ by an integer vector may change the relative order of the
coordinates, and hence the sequence of winners during the next $Q$ steps.
\end{remark}

\section{Orbits in the ergodic case}\label{Sergodic}

The other extreme case is $\mu=m-1$.
By \eqref{mu2}, this occurs if and only if $p_1,\dots,p_m$ span an
$m$-dimensional vector space over $\bbQ$, \ie, 
if and only if
\emph{the indexed family $(p_1,\dots,p_m)$ is rationally independent}.

In this case, we may in \refL{LG} take $I$ as any subset of $[m]$ with
$|I|=m-1$.
It then follows from \eqref{Su0} that
the projection $\Pi_I$ is 
an isomorphism of $\hH_\0$ onto $\bbT^\mu$. 
Hence, $G\subseteq\hHo$ and \refL{LG} imply that  $G=G_0=\hHo\cong\bbT^{m-1}$.
Furthermore, it is well known that the translation by $\Pi_I(\hvp)$ in
$\bbT^{m-1}$ then is ergodic, and thus $\Rp$ is an ergodic translation in
$\hH_\0$. We therefore call this the \emph{ergodic} case.
(It may also be called the \emph{rationally independent case}.)

\begin{remark}
 The ergodic case is generic, in the sense of both measure and
topology.
I.e., the
complement in $\gS_m$ of the set of ergodic $\vp$ has measure 0
and is a meagre set. 
\end{remark}

In the ergodic case, 
the sequence $n\hvp$ is
equidistributed in $\hHo$, and thus dense in $\hHo\cong \bbT^{m-1}$, 
\cf{} the proof of \refL{LG}; hence the same holds for every orbit
$\Omega(\hvu)=\set{\Rp^n\hvu:n\ge0}$. 
The main result of the present section shows how this lifts to a result on
invariant sets in $\Ho$.

Note that if $A$ is invariant for $T$, then, by \eqref{piT},
\begin{align}\label{Rpinv}
  \Rp(\pi(A))=\pi(T(A))=\pi(A)
\end{align}
and thus $\pi(A)\subseteq\bbT^m$ is invariant for $\Rp$,
and similarly for backward and forward invariant.

In the following lemma, we allow both the single-valued and the multi-valued
interpretation of $T$.
We let $\gl$ denote the Lebesgue measure on both $\Ho$ and
$\hHo$, in both cases normalized  by defining
it to be  $\ddx x_1 \dotsm \ddx x_{m-1}$.
(This differs by an unimportant constant from the $(m-1)$-dimensional
Hausdorff measure.) 
Note that on the compact group $\hHo$, $\gl$ is thus the
normalized Haar measure.

\begin{lemma}\label{LC}
  Suppose that $\vp$ is ergodic.
  \begin{romenumerate}
    
  \item \label{LC1}
Let $A\subset\Ho$ be a measurable bounded backward invariant set,
and define the function $\psi:\hHo\to\bbNo$ by 
\begin{align}\label{psi} 
  \psi(\hvx):=\bigabs{A\cap\pi\qw(\hvx)},
\end{align}
i.e., the number of points in $A$ that cover $\hvx$ for the projection $\pi$.
Then there is a constant $\chi\in\bbNo$ such that $\psi(\hvx)=\chi$ for \aex{}
$\hvx\in\hHo$.
Furthermore, the measure $\gl(A)=\chi$.

\item \label{LC2}
In particular, if $A$ is a non-empty compact backward invariant set, then
$\pi(A)=\hHo$, and $\psi(\hvx)=\chi\ge1$ for \aex{} $\hvx\in\hHo$.
  \end{romenumerate}
\end{lemma}

\begin{proof}
\pfitemref{LC1}
  We have
  \begin{align}\label{tc1}
    \psi(\pi(\vx)):=\sum_{\vz\in\bbZ^m}\etta_A(\vx+\vz),
\qquad \vx\in\Ho.
  \end{align}
Since $A$ is bounded and measurable, this shows that $\psi$ maps 
$\pi\Ho=\hHo$ into
$\bbNo$ and that $\psi$ is measurable.

Let $\hvx\in\hHo$ and let $\pi\qw(\hvx)\cap A=\set{\vx_1,\dots,\vx_\ell}$, where
$\ell=\psi(\hvx)$.
Thus every $\vx_i\in A$, and since $A$ is backward invariant, there exist 
$\vy_i\in A$ such that $T(\vy_i)=\vx_i$. 
By \eqref{piT}, we have
\begin{align}\label{tc2}
\pi(\vy_i)=\pi(\vx_i)-\hvp=\hvx-\hvp.   
\end{align}
Since the exceptional subset $\cD_0$ of $\Ho$ where $T$ may take multiple
values (in the multi-valued interpretation) has Lebesgue measure 0, see
\refR{RD}, 
it follows that so has its projection $\pi(\cD_0)\subset\hHo$.
Consequently, for \aex{} $\hvx\in\hHo$, $T$ is single-valued on
$\pi\qw(\hvx-\hvp)$, and thus at every $\vy_i$. Since the images $T(\vy_i)$
are distinct, it follows that $\vy_1,\dots,\vy_\ell$ are distinct, and thus,
using again \eqref{tc2}, 
\begin{align}\label{tc3}
\psi(\hvx-\hvp)\ge 
|\set{\vy_i:1\le i\le \ell}|
=\ell = \psi(\hvx).  
\end{align}
Thus $\psi(\hvx-\hvp)-\psi(\hvx)\ge0$ \aex. 
On the other hand, $\psi$ is bounded and measurable, and thus integrable.
Since Lebesgue measure (= Haar measure) on $\hHo\cong\bbT^{m-1}$ is
invariant, it follows that
\begin{align}\label{tc4}
  \int_{\hHo} \bigpar{\psi(\hvx-\hvp)-\psi(\hvx)}\dd \hvx 
=
  \int_{\hHo} \psi(\hvx-\hvp)\dd \hvx  -   \int_{\hHo} \psi(\hvx)\dd \hvx 
=0.
\end{align}
Consequently, 
\begin{align}\label{tc5}
\psi(\hvx-\hvp)=\psi(\hvx) 
\qquad \text{a.e.}
\end{align}
 Since $\vp$ is ergodic, it follows from \eqref{tc5}, for example by
considering the Fourier transform, that $\psi(\hvx)$ is \aex{} 
equal to some constant $\chi$.

Finally, the projection $\pi$ is a local bijection, and 
locally measure preserving
(by our normalization of the measures). 
Hence it follows from \eqref{tc1} that
\begin{align}
  \gl(A)=\int_{\hHo}\psi(\hvx)\dd\hvx = \chi.
\end{align}

\pfitemref{LC2}
In this case, $\pi(A)$ is non-empty and compact, and backward
invariant for $\Rp$ by
\eqref{piT}. 
Let $\hvx_0\in\pi(A)$. Since $\pi(A)$ is backward invariant, we have
$\Rp^{-n}(\hvx_0)=\hvx_0-n\hvp\in\pi(A)$ for every $n\ge0$.
On the other hand, as said above, 
when $\vp$ is ergodic, 
then the sequence $n\hvp$ is
dense in $\hHo$, 
Hence $\pi(A)$ contains a subset that is dense in $\hHo$, 
and since $\pi(A)$ is compact, it follows that $\pi(A)=\hHo$.

Consequently, $\psi(\hvx)\ge1$ for every $\hvx\in\hHo$, and 
the result follows from \ref{LC1}.
\end{proof}

We state the next lemma for the multi-valued version of $T$, so we regard
$T^n\vu$ as a set; this implies trivially the same result for any
single-valued version.
Recall that $\cU$ is the set where $T$ is uniquely determined without ties,
see \refR{RD}.
\begin{lemma}\label{Luni}
Suppose that $\vp$ is ergodic. Then,
  for every $\vu\in\bbR^m$, there exists $N$ such that $T^n(\vu)\subset \cU$
  for every $n\ge N$.
\end{lemma}
\begin{proof}
  Suppose that $\vz(n)=T^n\vu\notin\cU$. Then $z_i(n)+p_i=z_j(n)+p_j$ for some
  $1\le i<j\le m$. 
Hence we have, since \eqref{tx} implies $z_i(n)-u_i-np_i\in\bbZ$, 
  \begin{align}\label{luni1}
u_i-u_j + (n+1)(p_i-p_j) 
\in z_i(n)+p_i-z_j(n)-p_j+\bbZ
= \bbZ.
  \end{align}
This cannot happen for two different values of $n$ unless $p_i-p_j$ is
rational, which is impossible when $\vp$ is ergodic.
Consequently, there is at most one such $n$ for each pair $(i,j)$, and thus
$T^n\vu\notin\cU$ happens at most $\binom m2$ times.
\end{proof}

In the next theorem we consider only the multi-valued version of $T$ which 
leads to a nice formulation with a compact global attractor. 
(Recall \refR{Rmulti}, and
in particular that $T$ then maps compact sets to compact.)
The limit set $\cL(\vu)$ is by 
\refC{CLu} and \refL{Lback} 
compact and
backward invariant. Forward invariance in the ergodic case will be
proved below in \refC{CTX}.
Therefore, we define also
\begin{align}\label{cLx}
  \cLx(\vu):=\bigcup_{n=0}^\infty T^n(\cL(\vu))
.\end{align}
Then $\cLx(\vu)\supseteq\cL(\vu)$, and
$\cLx(\vu)$ is forward invariant by definition. 
\begin{remark}\label{RPinv}
The auxiliary set $\cLx(\vu)$ is introduced because, at this stage,
forward invariance of $\cL(\vu)$ has not yet been established.
We will later show in \refC{CTX} that,
in the ergodic case at least, 
$\cL(\vu)$ actually is forward invariant and thus
$\cLx(\vu)=\cL(\vu)$.
Thus the definition \eqref{cLx} 
is only temporary.
(For non-ergodic $\vp$, forward invariance of $\cL(\vu)$ may fail;
see \refE{Einter1}.)
\end{remark}

At least if $\vp$ is ergodic, $\cLx(\vu)$
is also backward invariant 
as shown in the proof below.
(Recall that $\cL(\vu)$ is backwards invariant for every $\vp$.)
It follows from \refL{LK0} that $\cLx(\vu)\subseteq K_0$ for $\vu\in\Ho$.
It is less obvious that $\cLx(\vu)$ is compact, but  that is part of the
theorem below.
We will show in \refC{CTX} that also $\cL(\vu)=\Kx$ for $\vu\in\Ho$.
Further results on $\Kx$ will be given in \refTs{Tgl}, \ref{TB}, and \ref{TX}.

\begin{theorem}\label{Tinv}
Suppose that $\vp$ is ergodic and
consider the multi-valued version of $T$. 
Then $T$ has a unique non-empty compact invariant set $\Kx\subset\Ho$.
Moreover, $\Kx$ is absorbing, and is thus the (unique)
global attractor for the dynamics on $\Ho$.
In particular, for every $\vu\in\Ho$ there exists $N(\vu)$ such that
\begin{align}
T^n(\vu)\subseteq \Kx \qquad\text{for every }n\ge N(\vu),  
\end{align}
More generally, 
for every compact set $K\subseteq\Ho$,
there exists $N(K)<\infty$ such that
\begin{align}\label{tinvK}
  T^n(K)\subseteq \Kx
\qquad\text{for every }n\ge N(K).
\end{align}
Furthermore, 
there exists $N<\infty$ such that
for every $\vu\in\Ho$, 
\begin{align}\label{tinv}
  \Kx=\cLx(\vu)=\bigcup_{n=0}^NT^n(\cL(\vu))
.\end{align}
The Lebesgue measure $\gl(\Kx)=1$.
The projection $\pi$ onto $\bbT^m$ maps $\Kx$ onto $\hHo$, and it is
\aex{} injective in the sense that $|\pi\qw(\hvx)\cap \Kx|=1$
for \aex{} $\hvx\in\hHo$.
\end{theorem}

\begin{proof}
Since some arguments are a bit subtle because of the possibility of multiple
values of $T$, we will in this proof be careful to regard $T\vx$ as a set;
however, we will identify a singleton set $\set{\vx}$ and the point $\vx$.

Let $\vu\in\Ho$.
Recall the definition of the orbit $\cO(\vu)$ and the limit set $\cL(\vu)$
in \eqref{Ou2} and \eqref{Lu}.
We consider only initial points $\vu\in\Ho$.
Then \refCs{Cu} and \ref{CLu} show that 
$\bcO{\vu}$ and $\cL(\vu)$ are non-empty compact subsets of
$\Ho\subset\bbR^m$.

\stepx{The projection $\pi(\cL(\vu))=\hHo$.}\label{Tinv2}
$\cL(\vu)$ is compact and non-empty, and backward invariant by
\refL{Lback}. 
Hence, \refL{LC}\ref{LC2} applies.

\smallskip
For the remainder of the proof, 
fix $0<\gd\le 1/m$ and
let $\cN_\gd$ and
$\hcN_\gd$ be the neighbourhoods of $\0$ in $\Ho$ and $\hHo$, respectively,
defined by $\max_j|x_j|<\gd$ and $\max_j|\hx_j|<\gd$.
Recall $K_0$ defined in \eqref{jb3}.

\stepx{
$\cN_\gd\subseteq\cL(\vu)$.
Furthermore, 
if $\vx\in\cL(\vu)$ with $\pi(\vx)\in\hcN_\gd$, 
  then $\vx\in\cN_\gd$.
In other words,}\label{Tinv3} 
\begin{align}\label{tinv3}
  \cL(\vu)\cap\pi\qw(\hcN_\gd)=\cN_\gd.
\end{align}
To see this, suppose that $\vx=(x_j)_1^m\in\cL(\vu)$ with 
$\hvx=(\hx_j)_1^m:=\pi(\vx)\in\hcN_\gd$.
 Identify $\hx_j\in\bbT$ with the corresponding real number in
$(-\frac12,\frac12]$. Hence the assumption means $|\hx_j|<\gd$ for every
$j\in[m]$.  

By \refL{LK0}, 
$x\in K_0$ and thus
$x_j\ge -1+\frac1m$ for every $j\in[m]$;
hence
\begin{align}\label{tinv3b}
x_j-\hx_j \ge -1+\frac{1}{m}-|\hx_j| > -1+\frac{1}{m}-\gd   \ge-1,
\end{align}
Furthermore, $x_j-\hx_j\in\bbZ$,
and consequently $x_j-\hx_j\ge0$.
Moreover, since $\vx\in\cL(\vu)\subset\Ho$ 
and thus $\sum_j x_j=0$,
we have 
$\sum_j\hx_j\in\bbZ$;
since also  $\sum_j|\hx_j|<m\gd\le1$, it follows that
$\sum_j\hx_j=0$, and thus
\begin{align}\label{tinv3c}
 \sumjm\bigpar{ x_j-\hx_j}
=  \sumjm x_j - \sumjm \hx_j
=0.
\end{align}
Consequently, $x_j-\hx_j=0$ for every $j$, and thus $\vx\in\cN_\gd$.

This proves the inclusion $\subseteq$ in \eqref{tinv3}.
Conversely, if $\vx\in\cN_\gd$, then $\pi(\vx)\in\hcN_\gd$, and
\refStep{Tinv2} shows that $\pi(\vx)=\pi(\vy)$ for some $\vy\in\cL(\vu)$;
it follows by the already proved inclusion that  $\vy\in\cN_\gd$, and
then necessarily $\vy=\vx$.

\stepx{
If\/ $\vu\in K_0$ defined in \eqref{jb3} and, for some $n\ge0$, 
$\Rp^n\hvu\in\hcN_\gd$, then 
$T^n\vu$ is single-valued and $T^n\vu\in\cN_\gd$.}\label{Tinv4}
This is very similar to \refStep{Tinv3}.
Let $\vz(n)=(z_j(n))_1^m$
be any element of $T^n\vu$, and 
let $\hvz(n)=(\hz_j(n))_1^m=\pi(\vz(n))=\Rp^n\hvu$.
Hence the assumption means $|\hz_j(n)|<\gd$ for every $j\in[m]$.  
By \refL{LK0}, $T^n\vu\in K_0$ and thus
$z_j(n)\ge -1+\frac1m$. Furthermore,
$T^n\vu-\Rp^n\hvu\in\bbZ^m$ by \eqref{piT}, and thus
$z_j(n)-\hz_j(n)\in\bbZ$ for every $j\in[m]$.
It now follows as in \eqref{tinv3b}--\eqref{tinv3c} that
$z_j(n)-\hz_j(n)=0$ for every $j$, and thus
\begin{align}\label{tinv3d}
  z_j(n)=\hz_j(n).
\end{align}
Hence each $z_j(n)$ is uniquely determined, and thus $T^n\vu$ is
single-valued. Furthermore, \eqref{tinv3d} shows also that $\vz(n)\in\cN_\gd$.

\stepx{There exists $N$ such that if\/ $\vu,\vv\in \Ho$, then}\label{Tinv5} 
\begin{align}\label{tinv5}
 \cLx(\vv)\subseteq \bigcup_{n=0}^NT^n(\cL(\vu)). 
\end{align}
Since $\vp$ is ergodic, the open sets $\Rp^{k}\hcN_\gd$, $k\ge0$, 
cover $\hHo$, so
by compactness of $\hHo$, there exists $N$ such that
\begin{align}\label{tinv5a}
\bigcup_{k=0}^N\Rp^{k}\hcN_\gd=\hHo.  
\end{align}
Let $\vx\in\cL(\vv)$ and let $n\ge N$. Write as usual $\hvx:=\pi(\vx)$.
Then \eqref{tinv5a} implies that 
$\Rp^n\hvx=\pi(T^n\vx)\in \Rp^{k}\hcN_\gd$
for some $0\le k\le N$, and thus 
\begin{align}\label{tinv5b}
\Rp^{n-k}\hvx \in \hcN_\gd.
\end{align}
Furthermore,  $\vx\in K_0$ by \refL{LK0}.
Consequently, \refStep{Tinv4} 
applies to $\vx$, and shows together with \refStep{Tinv3} that
\begin{align}\label{tinv5c}
  T^{n-k}\vx\in \cN_\gd \subseteq \cL(\vu).
\end{align}
Hence,
\begin{align}\label{tinv5d}
  T^{n}\vx\in T^k \bigpar{\cL(\vu)}
\subseteq \bigcup_{k=0}^NT^k(\cL(\vu)). 
\end{align}
Consequently, 
\begin{align}\label{tinve}
T^n(\cL(\vv))\subseteq \bigcup_{k=0}^NT^k(\cL(\vu)).   
\end{align}
We have assumed $n\ge N$, but 
\refL{Lback}
shows that $\cL(\vv)\subseteq T(\cL(\vv))$,
and it follows that $T^n(\cL(\vv)) \subseteq T^N(\cL(\vv))$
when $0\le n<N$; hence \eqref{tinve} holds for $0\le n<N$ too.
The definition \eqref{cLx} yields  \eqref{tinv5},

\stepx{For $N$ as in \refStep{Tinv5}, and every\/ $\vu\in\Ho$,} \label{Tinv6}
\begin{align}\label{tinv6}
 \cLx(\vu) =\bigcup_{n=0}^NT^n(\cL(\vu)). 
\end{align}
The inclusion $\subseteq$ follows by taking $\vv=\vu$ in \refStep{Tinv5}. 
The opposite inclusion is immediate from the definition \eqref{cLx}.

\stepx{The set $\cLx(\vu)$ is the same for all $\vu\in\Ho$.}\label{Tinv7}
For any $\vu,\vv\in\Ho$,
\refSteps{Tinv5} and \ref{Tinv6} yield $\cLx(\vv)\subseteq\cLx(\vu)$, and
the result follows by symmetry.

\smallskip
By \refSteps{Tinv6} and \ref{Tinv7}, we can define $\Kx$ by \eqref{tinv}
for any $\vu\in\Ho$, and then \eqref{tinv} holds for every $\vu$.

\stepx{$\Kx$ is compact and invariant.}\label{Tinv8}
We have  \eqref{tinv}, and the definition \eqref{cLx} thus yields
\begin{align}\label{tinv8a}
  T(\Kx)=T\Bigpar{\bigcup_{n\ge0}T^n\bigpar{\cL(\vu)}}
=\bigcup_{n\ge1}T^n\bigpar{\cL(\vu)}.
\end{align}
Forward invariance of $\Kx=\cLx(\vu)$, \ie, $T(\Kx)\subseteq \Kx$, follows 
immediately, and backward invariance follows by \eqref{tinv8a} and
\refL{Lback}.
Finally, every $T^n(\cL(\vu))$ is compact by \eqref{compact}, and thus 
\eqref{tinv6} shows that $\Kx=\cLx(\vu)$ is compact.

\stepx{
$\Kx$ is absorbing.
Moreover,
for every compact set $K\subseteq\Ho$, there exists
$N<\infty$ such that for every $\vx\in K$ and $n\ge N$,
we have $T^n\vx\subseteq \Kx$.}\label{Tinv9}
Assume first $K=K_0$.
Similarly as in \refStep{Tinv5},
since $\vp$ is ergodic, the open sets $\Rp^{-k}\hcN_\gd$, $k\ge0$, cover
$\hHo$, so by compactness of $\hHo$ there exists $N$ such that
\begin{align}\label{tinv9a}
\bigcup_{k=0}^N\Rp^{-k}\hcN_\gd=\hHo.  
\end{align}
Let $\vx\in K_0$. 
Then \eqref{tinv9a} implies that 
$\hvx:=\pi(\vx)\in \Rp^{-k}\hcN_\gd$
for some $0\le k\le N$, and thus $\Rp^k\hvx\in\hcN_\gd$.
By \refSteps{Tinv4} and \ref{Tinv3}, it follows that
$T^k\vx\in\cN_\gd\subseteq\cL(\vu)\subseteq \cLx(\vu)=\Kx$. 
Since $\Kx$ is invariant by \refStep{Tinv8}, it follows that
$T^n\vx\subseteq \Kx$ for every $n\ge k$, and in particular for every
$n\ge N$. 

For a general compact set $K$,  \refT{Tu}\ref{Tu2} shows that for some
$N_1=N_1(K)$, we have $T^{N_1}\vx\subseteq K_0$ for every $\vx\in K$.
It follows that $T^n\vx=T^{n-N_1}T^{N_1}\vx\subseteq\Kx$ for every $n\ge N+N_1$.
This proves \eqref{tinvK} with $N(K):=N+N_1(K)$.

\stepx{
$\Kx$ is the unique non-empty compact invariant set in $\Ho$.}\label{Tinv10}
Suppose  that $K$ is any non-empty compact invariant set.
By \refStep{Tinv9}, for some $N$ we have
$T^N(K)\subseteq\Kx$. Since $K$ is invariant, this means $K\subseteq\Kx$.
Conversely, let $\vx\in K$. 
Since $K$ is invariant, the orbit $\cO(\vx)\subseteq K$, and since $K$ is
compact, this implies $\bcO{\vx}\subseteq K$. 
Hence, the definition \eqref{Lu} implies
$\cL(\vx)\subseteq\bcO{\vx}\subseteq K$, and thus, again using the
invariance of $K$, we have
$T^n(\cL(\vx))\subseteq T^n K = K$ for every $n\ge0$.
Hence, $\Kx\subseteq K$ by \eqref{tinv} (with $\vu=\vx$).
Thus $K=\Kx$.


\stepx{The projection $\pi$ maps $\Kx$ onto $\hHo$ and is \aex{} injective.} 
\label{Tinv11}
\refStep{Tinv8} and 
\refL{LC}\ref{LC2} (with $A=\Kx$) show that
$\pi(\Kx)=\hHo$, and also that $\psi(\hvx)$ (defined in \eqref{psi}) 
is equal to some constant $\chi$ \aex{}  in $\hHo$.

Suppose now that $\vu\in\Kx$ and $\hvu:=\pi(\vu)\in\hcN_\gd$.
Since $\Kx\subseteq K_0$, it follows from \refStep{Tinv4} (with $n=0$)
that $\vu\in\cN_\gd$, and thus \eqref{tinv3d} holds;
hence $\vu$ is uniquely determined by $\hvu$.
This shows that $\psi(\hvu)=1$ for $\hvu\in\hcN_\gd$, and thus $\chi=1$.
Consequently, $|\pi\qw(\hvx)\cap \Kx|=\psi(\hvx)=1$ \aex{} in $\hHo$.

\stepx{$\gl(\Kx)=1$}.\label{Tinv12}
This follows by \refL{LC} and $\chi=1$ in \refStep{Tinv11}.

This completes the proof of \refT{Tinv}.
\end{proof}

\begin{remark}
We have shown in \refT{Tinv} that
$\psi(\hvx):=|\pi\qw(\hvx)\cap \Kx|=1$ for \aex{} $\hvx\in\hHo$.
However, this cannot hold for every $\hvx$, i.e., there must be exceptional
points.
To see this, note that otherwise $\pi$ would be a continuous bijection of
the compact set
$\Kx$ onto $\hHo$, and thus a homeomorphism $\Kx\cong\hHo$;
this is impossible since there is no subset of $\bbR^{m-1}$ homeomorphic to
$\bbT^{m-1}$.
(By the invariance of domain theorem applied to $\pi\qw$, every point in
$\Kx$ would be an interior point, and thus $\Kx$ would be open, non-empty,
and compact, a contradiction.)
\end{remark}

\subsection{The invariant measure in the ergodic case}\label{Smeasure}

The invariant set $\Kx$
comes with an invariant measure, which
is just Lebesgue measure.
Recall that if $f$ is any measurable function from one measurable space
$\cX$ to another $\cY$, and $\nu$ is a measure on $\cX$, then $f\stj\nu$ is
the measure $f\stj\nu(A):=\nu(f\qw(A))$ on $\cY$.
For a multi-valued $f$, we  use the same definition but consider only measures
$\nu$ such that $f(\vx)$ is single-valued for $\nu$-\aex{} $\vx\in\cX$;
then we may ignore the set where $f(\vx)$ takes more than one value and
regard $f$ as single-valued.
With this interpretation, the following theorem is valid for
both the single-valued and multi-valued versions  of $T$.

Let $\gd_\vx$ denote the Dirac measure (point mass) at $\vx$.

\begin{theorem}\label{Tgl}
Suppose that $\vp$ is ergodic.
Let $\glKx:=\etta_{\Kx}\gl$ be Lebesgue measure on $\Ho$ restricted to $\Kx$.
\begin{romenumerate}
  
\item \label{Tgla}
Then $\glKx$ is a probability measure that is invariant for $T$.
Moreover, $\glKx$ is the unique invariant probability measure on $\Ho$.

\item \label{Tglb}
If\/ $\vu\in\Ho$, and $\vz'(n)$ is any element in $T^n\vu$, 
for each $n\ge0$, then the empirical distribution
\begin{align}\label{tgl1b}
  \nu_N:=\frac{1}{N}\sum_{n=1}^N \gd_{\vz'(n)}
\wto\glKx,
\qquad \text{as \Ntoo}, 
\end{align}
with convergence in the usual weak topology.
\end{romenumerate}  
\end{theorem}

We may express \eqref{tgl1b} by saying that the sequence $\vz'(n)$ is
\emph{asymptotically uniformly distributed in $\Kx$}.

\begin{proof}
\pfitemref{Tgla}
Suppose that $\nu$ is an invariant  probability measure on $\Ho$. 
Then \eqref{piT} implies
\begin{align}\label{69a}
(\Rp)\stj(\pi\stj(\nu))
=\pi\stj(T\stj(\nu))=\pi\stj(\nu),
\end{align}
and thus $\pi\stj\nu$ is invariant for $\Rp$ on $\pi(\Ho)=\hHo\subset\bbT^m$.
Since $\vp$ is ergodic, this implies that $\pi\stj\nu$ is the Lebesgue
measure $\gl$ on $\hHo$; we denote this measure by $\glhHo$ for clarity.

Furthermore, let $\eps>0$; then there exists a compact set $K\subseteq\Ho$
such that $\nu(K)>1-\eps$. 
It follows from the invariance that
\begin{align}\label{69b}
\nu(T(K))=  (T\stj\nu)(T(K)) = \nu\bigpar{T\qw(T(K))}
\ge \nu(K) > 1-\eps,
\end{align}
and by induction
\begin{align}\label{69c}
  \nu(T^n(K))\ge \nu(K)>1-\eps
\end{align}
for every $n\ge0$.
By \refStep{Tinv9} in the proof of \refT{Tinv}, 
for $n$ large enough we have $T^n(K)\subseteq\Kx$, and thus
\eqref{69c} implies $\nu(\Kx)>1-\eps$. Since $\eps$ is arbitrary, this
yields $\nu(\Kx)=1$, \ie, $\nu$ is concentrated on $\Kx$.

\refT{Tinv} implies that $\pi:\Kx\to\hHo$ is onto and has a Borel
measurable right inverse $\iota:\hHo\to\Kx$, so that
$\pi(\iota(\hvx))=\hvx$ for every $\hvx\in\hHo$;
for example, we may choose the smallest inverse image in lexicographic
order.
Let $E\subseteq\hHo$ be the set of points with more than one preimage
in $\Kx$. By \refT{Tinv}, $\glhHo(E)=0$.
For any probability measure $\rho$ on $\Kx$ with
$\pi\stj\rho=\glhHo$, we have
\begin{align}
 \rho\bigpar{\{\vx\in\Kx:\pi(\vx)\in E\}}
 =(\pi\stj\rho)(E)=\glhHo(E)=0.
\end{align}
Outside this set the lift is unique, so $\iota(\pi(\vx))=\vx$ for
$\rho$-a.e. $\vx\in\Kx$. Thus every such measure $\rho$ equals
$\iota\stj\glhHo$.
Applying this to $\nu$ gives
\begin{align}\label{69d}
  \nu=\iota\stj(\pi\stj\nu)
=\iota\stj\glhHo.
\end{align}
Hence $\nu$ is unique.

Conversely, define $\nu:=\iota\stj\glhHo$ on $\Kx$.
Then $\pi\stj\nu=(\pi\circ\iota)\stj\glhHo=\glhHo$.
The set $\pi(\cD_0)$ has Lebesgue measure zero, as noted in the proof of
\refL{LC}; hence
$\nu(\cD_0)\le(\pi\stj\nu)(\pi(\cD_0))=0$.
Thus $T\stj\nu$ is well defined also in the multi-valued case.
Using \eqref{piT}, we obtain
\begin{align}\label{69e}
  \pi\stj(T\stj(\nu))=(\Rp)\stj(\pi\stj(\nu))=(\Rp)\stj(\glhHo)
=\glhHo
=\pi\stj\nu.
\end{align}
Moreover, $T\stj\nu$ is supported on $\Kx$, since $T(\Kx)=\Kx$.
The preceding lift argument, applied to $\rho=T\stj\nu$, therefore gives
$T\stj\nu=\iota\stj\glhHo=\nu$.

It remains to identify this measure $\nu$ as $\glKx$.
For $\vn=(n_1,\dots,n_{m-1})\in\bbZ^{m-1}$, let 
\begin{align}\label{69f}
  Q_\vn:=\set{\vx\in\Ho: \floor{x_i}=n_i \text{ for } i=1,\dots,m-1}.
\end{align}
Then $\set{Q_\vn:\vn\in\bbZ^{m-1}}$ is a partition of $\Ho$, and
  $\pi:Q_\vn\to\hHo$ is a measure preserving bijection for every $\vn$.
Let 
$\psi(\hvx):= |\pi\qw(\hvx)\cap \Kx|$
be as in \refL{LC}, and 
recall that \refT{Tinv} shows that
$\psi(\hvx)=1$ \aex{} in $\hHo$.
It follows that for any measurable set $B\subseteq\hHo$,
\begin{align}\label{69g}
  \pi\stj(\glKx)(B)&
=\glKx\bigpar{\pi\qw(B)}
=\sum_{\vn\in\bbZ^{m-1}}\gl\bigpar{\pi\qw(B)\cap\Kx\cap Q_\vn}
\notag\\&
=\sum_{\vn\in\bbZ^{m-1}}\int_B\etta_{\Kx}\bigpar{\pi|_{Q_\vn}\qw(\hvx)}
\dd\glhHo(\hvx)
=\int_B\psi(\hvx)\dd\glhHo(\hvx)
\notag\\&
=\glhHo(B).
\end{align}
Hence, $\pi\stj(\glKx)=\glhHo$, and thus, using \eqref{69d},
$\nu=\iota\stj\glhHo=\glKx$.

\pfitemref{Tglb}
Suppose first that $\vu\in\Kx$. Since $\Kx$ is invariant by \refT{Tinv},
then every $\vz'(n)\in\Kx$, and thus $\nu_N$ is  a probability measure on
$\Kx$.

Clearly, $\pi\stj\nu_N$ is the empirical distribution of the sequence
$\pi(\vz'(n))=\Rp^n\hvu$. Since $\vp$ is ergodic, 
Weyl's criterion applies and shows that the sequence $\Rp^n\hvu$ is
equidistributed in $\hHo$, \ie,
\begin{align}\label{69h}
\hnu_N:=  \pi\stj(\nu_N)=\frac{1}{N}\sum_{n=1}^N\gd_{\Rp^n\hvu}\wto \glhHo.
\end{align}

Let $K\subseteq\Kx$ be compact. Then $\pi(K)\subseteq\hHo$ is compact, and 
$\hnu_N(\pi(K)) =\nu_N\bigpar{\pi\qw(\pi(K))}\ge\nu_N(K)$.
Hence, \eqref{69h} and the Portmanteau theorem
\cite[Theorem 2.1]{Billingsley1968}
yield
\begin{align}\label{69i}
  \limsup_\Ntoo \nu_N(K) &
\le
  \limsup_\Ntoo \hnu_N(\pi(K))
\le \glhHo(\pi(K))
=(\pi\stj\glKx)(\pi(K) )
\notag\\&
=\glKx\bigpar{\pi\qw(\pi(K) )}
=\glKx(K),
\end{align}
using again that $\pi$ is \aex{} injective on $\Kx$ by \refT{Tinv}.
By the Portmanteau theorem again,
\eqref{69i} shows that $\nu_N\wto\glKx$ as probability measures on $\Kx$,
which implies $\nu_N\wto\glKx$ also as probability measures on $\Ho$ (or
$\bbR^m$). 

For a general $\vu\in \Ho$, \refT{Tinv} shows that 
$T^k\vu\subseteq\Kx$ for some $k$.
The proof above shows that
\begin{align}\label{69j}
\nu_N':=  \frac{1}{N}\sum_{n=k+1}^{N+k}\gd_{\vz'(n)}\wto\glKx
\end{align}
as \Ntoo. Furthermore, the total variation distance satisfies
\begin{align}
  \dtv(\nu_N,\nu'_N)\le \frac{2k}{N}\to0
\end{align}
as \Ntoo, and thus $\nu_N\wto\glKx$ follows from \eqref{69j}.
\end{proof}

\begin{corollary}\label{Cgl}
 Suppose that $\vp$ is ergodic, and let, recalling \eqref{theta} and \eqref{tx},
  \begin{align}\label{clg1}
    \Kx_i:=\set{\vx\in\Kx: \theta(\vx+\vp)=i},\qquad i=1,\dots,m,
  \end{align}
where any tie convention is used on the tie hyperplanes.  
(This choice affects $\Kx_i$ only up to a null set.)
Equivalently, $\Kx_i$ is the subset of $\Kx$ where
$T\vx=\vx+\vp-\ve_i$.  Then
\begin{align}\label{clg2}
  \gl(\Kx_i)=p_i.
\end{align}
\end{corollary}

\begin{proof}
  For simplicity, choose $\vu\in\Kx$ such that the orbit $(T^n\vu)_{n\ge0}$
is unique, without any branching. (Recall that this holds for a.e.\
$\vu\in\Ho$, see \refR{RD}.) Then, by \eqref{prop},
\begin{align}\label{clg3}
\lim_{N\to\infty}\frac{\#\set{n\le N:T^n\vu\in \Kx_i}}{N}
= p_i.
\end{align}

The set $\Kx_i$ is a continuity set for $\gl_{\Kx}$.  Indeed, by
Theorem~\ref{TB} below, $\Kx$ is a finite union of polytopes, and thus
$\gl(\partial\Kx)=0$; apart from this boundary, the boundary of $\Kx_i$ is
contained in the finite union of tie hyperplanes
\begin{align}\label{clg35}
\bigcup_{j\ne k}\set{\vx\in\Ho:x_j+p_j=x_k+p_k},
\end{align}
which also has $\gl$-measure zero.  Hence
$\gl_{\Kx}(\partial \Kx_i)=0$.  The weak convergence in \eqref{tgl1b}
therefore gives
\begin{align}
\lim_{N\to\infty}\frac{\#\set{n\le N:T^n\vu\in \Kx_i}}{N}
=\gl_{\Kx}(\Kx_i)=\gl(\Kx_i),
\end{align}
and the result follows from \eqref{clg3}.
(As an alternative to using \refT{TB}, we may 
note that \eqref{tgl1b} holds also if we 
regard the measures $\nu_N$ and $\gl_\Kx$ as measures on the compact metric
space $\Kx$, and that
regarded as a subset of $\Kx$,
the set $\Kx_i$ is a continuity set for $\gl_{\Kx}$
since its boundary is
contained in \eqref{clg35}.)
\end{proof}

\begin{theorem}\label{TsupportK}
Suppose that $\vp$ is ergodic and consider the multi-valued version of $T$.
Let $\vu\in\Ho$. Then there exists $M<\infty$ such that $T^n(\vu)\subset \cU$
for all $n\ge M$. Writing
\begin{equation}\label{ts1}
T^M(\vu)=\{\vu^{(1)},\dots,\vu^{(r)}\}\subset \cU,
\end{equation}
we have
\begin{equation}\label{ts2}
\cL(\vu)=\bigcup_{a=1}^r \omega\bigl(\vu^{(a)}\bigr),
\end{equation}
where $\omega\bigl(\vu^{(a)}\bigr)$ is the limit set
of the single-valued orbit starting from $\vu^{(a)}$.
Moreover, for every $a=1,\dots,r$,
\begin{equation}\label{ts3}
\supp(\gl_{\Kx})
\subseteq
\omega\bigl(\vu^{(a)}\bigr)
\subseteq
\cL(\vu)
\subseteq
\Kx.
\end{equation}
Consequently, for every branch orbit $\vz'=(\vz'_n)_{n\ge0}$ starting from
an arbitrary $\vu\in\Ho$,
\begin{equation}\label{ts4}
\supp(\gl_{\Kx})
\subseteq
\omega(\vz')
\subseteq
\cL(\vu)
\subseteq
\Kx.
\end{equation}
\end{theorem}

\begin{proof}
The existence of $M$ is Lemma~\ref{Luni}. Since $T$ is single-valued on
$\cU$, \eqref{ts1} implies that for every $n\ge0$,
\begin{equation}\label{ts6}
T^{M+n}(\vu)=\{T^n(\vu^{(1)}),\dots,T^n(\vu^{(r)})\},
\end{equation}
where on the right-hand side $T$ is single-valued. For $a=1,\dots,r$ and
$n\ge0$, let
\begin{equation}\label{ts7}
B_{a,n}:=\overline{\{T^k(\vu^{(a)}):k\ge n\}}.
\end{equation}
Then each sequence $(B_{a,n})_{n\ge0}$ is decreasing, and finite initial
segments do not affect limit sets. Hence
\begin{equation}\label{ts8}
\cL(\vu)
=
\bigcap_{n\ge0}\overline{\bigcup_{k\ge n}T^{M+k}(\vu)}
=
\bigcap_{n\ge0}\bigcup_{a=1}^r B_{a,n}.
\end{equation}
Since the number of branches is finite and each sequence $(B_{a,n})_n$ is
decreasing, the intersection of these finite unions equals the union of the
intersections.  Indeed, if a point belongs to the right-hand side
of \eqref{ts8}, then for
infinitely many $n$ it belongs to $B_{a,n}$ for some fixed $a$; by
monotonicity it then belongs to $B_{a,n}$ for every $n$.  Thus
\begin{equation}\label{ts9}
\cL(\vu)=\bigcup_{a=1}^r \bigcap_{n\ge0} B_{a,n}
=
\bigcup_{a=1}^r \omega\bigl(\vu^{(a)}\bigr).
\end{equation}

Fix $a\in\{1,\dots,r\}$. By Theorem~\ref{Tgl}\ref{Tglb}, the empirical
measures
\begin{equation}\label{ts10}
\nu_N^{(a)}:=\frac1N\sum_{n=0}^{N-1}\delta_{T^n(\vu^{(a)})}
\end{equation}
converge weakly to $\gl_{\Kx}$.
Let $\vx\in\supp(\gl_{\Kx})$, and let $U$ be an open neighbourhood of
$\vx$ in $\Ho$. Since $\vx$ belongs to the support,
\begin{equation}\label{ts11}
\gl_{\Kx}(U)>0.
\end{equation}
By the Portmanteau theorem,
\begin{equation}\label{ts12}
\liminf_{N\to\infty}\nu_N^{(a)}(U)\ge \gl_{\Kx}(U)>0.
\end{equation}
Therefore $T^n(\vu^{(a)})\in U$ for infinitely many $n$, and since this
holds for every neighborhood $U$ of $\vx$, we conclude that
\begin{equation}\label{ts13}
\vx\in\omega\bigl(\vu^{(a)}\bigr).
\end{equation}
This proves
\begin{equation}\label{ts14}
\supp(\gl_{\Kx})\subseteq \omega\bigl(\vu^{(a)}\bigr).
\end{equation}
The inclusion $\omega\bigl(\vu^{(a)}\bigr)\subseteq \cL(\vu)$ follows from
\eqref{ts9}, and $\cL(\vu)\subseteq \Kx$ because $\Kx$ is 
closed and absorbing.  This proves \eqref{ts3}.

Finally, any branch orbit $\vz'=(\vz'_n)_{n\ge0}$ starting from $\vu$ chooses
one of the points $\vu^{(a)}$ at time $M$; hence
\begin{equation}\label{ts15}
\omega(\vz')=\omega\bigl(\vu^{(a)}\bigr)
\end{equation}
for some $a\in\{1,\dots,r\}$. Together with \eqref{ts3} this yields
\eqref{ts4}.
\end{proof}

\subsection{Polyhedral structure of the attractor in the ergodic case}\label{Spoly}

In the preceding subsections we developed the multi-valued dynamics of the
Frege map in the ergodic case and obtained the inclusions
\begin{align}
\supp(\gl_{\Kx})\subseteq \omega(\vz')\subseteq \cL(\vu)\subseteq \Kx.  
\end{align}
The purpose of the present section is to identify the geometry of $\Kx$ more
precisely. We shall prove that $\Kx$ is a finite union of full-dimensional
polytopes in $\Ho$. This implies that $\Kx$ is equal to the support of its
invariant Lebesgue measure, and therefore the inclusions above collapse to
\begin{align}
\omega(\vz')=\cL(\vu)=\Kx.  
\end{align}

Consider the $m(m-1)$ linear forms $\ell_{ij}(\vx):=x_i-x_j$ on $\bbR^m$
for $(i,j)\in\mmx:=\set{(i,j)\in[m]^2:i\neq j}$.
Thus $\ell_{ji}=-\ell_{ij}$, and
the hyperplane 
$\set{\vx:\ell_{ij}(\vx)=0}$
contains the set where $i$ and $j$ tie in
\eqref{theta}. The action of the multi-valued operator $T$ on an arbitrary 
set $A\subseteq\bbR^m$ can be described as follows, see \eqref{txa},
\begin{enumerate}
 \renewcommand{\labelenumi}{\textup{(A\arabic{enumi})}}%
 \renewcommand{\theenumi}{\labelenumi}%
\item\label{TA1} 
Translate $A$ to $A+\vp$.
\item\label{TA2} 
Cut $A+\vp$ into the $m$ pieces (possibly empty)
  \begin{align}\label{ta2}
A_i:=(A+\vp)\cap \cV_i =(A+\vp)\cap\bigcap_{j\neq i} \set{\vx:\ell_{ij}(\vx)\ge0}  
 . \end{align}
\item \label{TA3}
Translate each $A_i$ to $A_i-\ve_i$. 
Take the union: $T(A)=\bigcup_i(A_i-\ve_i)$.
\end{enumerate}

Since we here consider the action of $T$ on $\Ho$, and the translations in 
\ref{TA1} and \ref{TA3} bring us out of $\Ho$ and then back again, it will
be convenient to instead
use the following, obviously equivalent, description, for a
subset $B\subseteq \Ho$.

\begin{enumerate}
 \renewcommand{\labelenumi}{\textup{(B\arabic{enumi})}}%
 \renewcommand{\theenumi}{\labelenumi}%
\item\label{TB2} 
Cut $B$ into the $m$ pieces (possibly empty)
  \begin{align}\label{tb2}
B_i:=B\cap\bigcap_{j\neq i} \set{\vx:\ell_{ij}(\vx+\vp)\ge0}  
.  \end{align}
\item \label{TB3}
Translate each $B_i$ to $B_i+\vp-\ve_i$.
Take the union: $T(B)=\bigcup_i(B_i+\vp-\ve_i)$.
\end{enumerate}

We say that a  polytope  $B\subset\Ho$ is a
\emph{block} 
if it is a finite intersection
\begin{align}\label{tb3}
  B=\bigcap_{\nu=1}^N\set{ \vx\in\Ho:\ell_{i_\nu,j_\nu}(\vx)\ge a_\nu}
\end{align}
for some $N\in\bbN$,  $(i_\nu,j_\nu)\in\mmx$, and $a_\nu\in\bbR$.
$B$ is a \emph{good} block if it is a block with full dimension $m-1$,
i.e., $B$ has a non-empty interior  in $\Ho$.

\begin{remark}
 If $B$ is a good block given by \eqref{tb3}, then each facet
(i.e., maximal proper face) of $B$ is given by the intersection of $B$ 
with one 
of the hyperplanes $\set{\vx:\ell_{i_\nu j_\nu}(\vx)=a_\nu}$. 
(The converse does not hold:  
this intersection is not necessarily a facet for every $\nu$.)
In particular, every
facet is parallel to one of the hyperplanes 
$\set{\vx:\ell_{ij}(\vx)=0}$.  
\end{remark}

It is obvious that if $B$ is a block, then the pieces $B_i$ in \eqref{tb2}
are blocks (or empty), and thus $T(B)$ is a finite union of blocks.
However, if $B$ is a good block, then the pieces $B_i$ are not necessarily good
blocks, because it might happen that for some $i$ and $j$, we have
$\ell_{ij}(\vx+\vp)\le 0$ for all $\vx\in B$ with equality for some
$\vx$; then \eqref{tb2} implies that $B_i$ is included in the hyperplane
$\set{ \vx\in\Ho:\ell_{i,j}(\vx+\vp)=0}$ 
and thus has dimension $\le m-2$.
Nevertheless, we show that we can avoid such problems and preserve good
blocks in a special case. (It is easy to see that good blocks are preserved
in generic cases, but we will need a more specific result.)

\begin{lemma}\label{LTB}
  Suppose that $\vp$ is ergodic.
Let $\rho_{ij}$, for $1\le i<j\le m$, be positive real numbers such that
the indexed family consisting of $p_1,\dots,p_m$ and
$\rho_{ij}$, $1\le i<j\le m$, is rationally independent,
and define
\begin{align}\label{tb4}
B: =\bigset{\vx\in \Ho: -\rho_{ij}\le \ell_{ij}(\vx)\le \rho_{ij}
\text{ for all $1\le i<j\le m$}}.   
\end{align}
Then $T^n B$ is a finite union of good blocks for every $n\ge0$.
\end{lemma}

\begin{proof}
We obtain $T^n B$ from $B$ by performing \ref{TB2}--\ref{TB3} $n$ times.
It is obvious that $T^nB$ is the union of  the $m^n$ pieces $B_\vi$
obtained by taking  a sequence $\vi=(i_1,\dots,i_n)\in[m]^n$ and then
keeping only the $i_k$-th piece in \ref{TB3} the $k$-th time.
(Some of these pieces might be empty, and are ignored.)
By induction, every such piece $B_\vi$ is
an intersection of the translate $B+n\vp-\vq$ for some $\vq\in \bbZ^m$
and a number of halfspaces, which are translates of the halfspaces used for
cutting in \ref{TB2} in some step $k\le n$, and thus of the form,
for some $(i,j)\in\mmx$, $1\le k\le n$  and $\vq=(q_1,\dots,q_m)\in\bbZ^m$,
\begin{align}\label{tb5}
&  \bigset{\vx\in \Ho: \ell_{ij}(\vx+\vp)\ge0}+(n-k+1)\vp-\vq
\notag\\&\hskip4em
=
  \bigset{\vx\in \Ho: \ell_{ij}(\vx)\ge (n-k)(p_i-p_j)-q_i+q_j}.
\end{align}
Furthermore, by \eqref{tb4}, $B+n\vp-\vq$ is the intersection,
letting $\rho_{ji}:=\rho_{ij}$ for $i<j$,
\begin{align}\label{tb6}
&  \bigcap_{i\neq j} \bigset{\vx\in\Ho: \ell_{ij}(\vx)\ge -\rho_{ij}}
+n\vp-\vq
\notag\\&\hskip4em
= \bigcap_{i\neq j} \bigset{\vx\in\Ho: \ell_{ij}(\vx)
\ge -\rho_{ij}+n(p_i-p_j)-q_i+q_j}
\end{align}
Hence, the block $B_\vi$ is the intersection of the halfspaces in
\eqref{tb5} and \eqref{tb6}.

Suppose by induction that $B_\vi$ is a good block. To obtain $T(B_\vi)$, we
cut $B_\vi$ as in \eqref{tb2}. We select some index $i_{n+1}$ and make the
intersections with the halfspaces in \eqref{tb2} one at a time, giving us a
sequence of blocks $B'_0=B_\vi,B'_1, B'_2, \dots,B'_{m-1}$.
The final translated piece is
$B_{(\vi,i_{n+1})}=B'_{m-1}+\vp-\ve_{i_{n+1}}$;
translation preserves the property of being a good block.
At step $r=0,\dots,m-2$, we select some $j=j_r\neq i_{n+1}$, 
and intersect the current block $B'_r$
with, writing $i':=i_{n+1}$ and $j':=j_r$,
the halfspace 
$W_r:=\set{\vx\in \Ho:\ell_{i'j'}(\vx+\vp)\ge0}$, which is bounded by the
hyperplane 
\begin{align}\label{tb7}
\dW_r:=\set{\vx\in \Ho:\ell_{i'j'}(\vx+\vp)=0}.  
\end{align}
We assume by induction that $B'_r$ is a good block, and claim that
then $B'_{r+1}:=B'_r\cap W_r$ is either a good block or empty.
There are three cases:

\pfcase1{The hyperplane $\dW_r$ does not intersect the current block $B'_r$.}
Then either $B'_{r+1}=B'_r$ or $B'_{r+1}=\emptyset$;  both cases are
trivial. 

\pfcase2{The hyperplane $\dW_r$ intersects the interior of the block $B'_r$.}
Then the block $B'_{r+1}=W_r\cap B'_r$ has non-empty interior, 
so $B'_{r+1}$ is a good block.

\pfcase3{The hyperplane $\dW_r$ intersects the  block $B'_r$ but not its
interior.} 
(We will show that this cannot occur, so we assume that it does in order to
obtain a contradiction.)
In other words, $\dW_r$ is a supporting hyperplane of $B'_r$.
Then the intersection $\dW_r\cap B'_r$ is a face of $B'_r$ and thus contains
some vertex $\vy$ of $B'_r$.
The polytope $B'_r$ has interior points in the
$m-1$-dimensional space $\Ho$, and thus its vertex
$\vy$ is the intersection of some
set of $m-1$ hyperplanes
chosen from the hyperplanes bounding the halfspaces
in \eqref{tb5} and \eqref{tb6} defining $B'_0=B_\vi$ together with the
additional hyperplanes $\dW_0,\dots,\dW_{r-1}$ in \eqref{tb7}
used so far in the present step.
Say that these $m-1$ hyperplanes are \emph{marked}.
By \eqref{tb5}--\eqref{tb7},
each marked hyperplane is of the type 
\begin{align}\label{tb8}
\ell_{ij}(\vx)=a_{ij}  
\end{align}
for some 
$(i,j)\in\mmx$ and $a_{ij}\in\bbR$, 
and we then say that also the pair $(i,j)$ is marked;
furthermore, \eqref{tb5}--\eqref{tb7} show that
for some integer $n_{ij}\ge0$,  we have one of
\begin{align}\label{tb85}
a_{ij}&\equiv n_{ij}(p_i-p_j)\pmod\bbZ,
\\\label{tb86}
a_{ij}&\equiv n_{ij}(p_i-p_j)-\rho_{ij}\pmod\bbZ,
\\\label{tb87}
a_{ij}&= -(p_i-p_j),
\end{align}
where \eqref{tb87} only can occur for $i=i'$ and $j\neq j'$.

Let $G$ be the graph with vertices $[m]$ and an edge $ij$ for every marked pair
$(i,j)$.
The normal vector to the hyperplane in \eqref{tb8} 
is $\ve_i-\ve_j$, and since the
$m-1$ marked hyperplanes intersect in a single point, their normals have to
be linearly independent and span $\Ho$. Hence no marked pair $(i,j)$ can be
repeated; moreover, there cannot be any cycle in the graph $G$.
Since $G$ has $m$ vertices and $m-1$ edges, it follows that
$G$ is a tree. 
Hence, if $\dW_r$ is given by $\ell_{i'j'}$ as in
\eqref{tb7}, then there is a path $i'_0=i', i'_1, \dots, i'_s=j'$ in the graph
$G$ connecting the indices $i'$ and $j'$.

By summing \eqref{tb8} along this path
(using $\ell_{ij}=-\ell_{ji}$ when necessary), we obtain using
\eqref{tb85}--\eqref{tb87}, and writing 
$n'_\nu:=n_{i'_{\nu-1},i'_\nu}$,
\begin{align}\label{tb9}
  y_{i'}-y_{j'} = \sum_{\nu=1}^s (y_{i'_{\nu-1}}-y_{i'_\nu})
\equiv \sum_{\nu=1}^s n'_\nu(p_{i'_{\nu-1}}-p_{i'_\nu}) 
+\sum_{\nu=1}^s n''_\nu\rho_{i'_{\nu-1},i'_\nu}
\pmod\bbZ
\end{align}
where $n'_\nu\ge0$ and $n''_\nu\in\set{0,\pm1}$,
except that we have $n'_\nu=-1$ when \eqref{tb87} is used, 
which along this  path is possible only for $\nu=1$ (so $i'_{\nu-1}=i'$), 
and in this case we must have $i'_\nu=i'_1\neq j'$ and thus $s>1$.
On the other hand, $\vy\in\dW_r$, and thus by \eqref{tb7}
\begin{align}\label{tb10}
  y_{i'}-y_{j'}=-p_{i'}+p_{j'}.
\end{align}
Subtracting \eqref{tb10} from \eqref{tb9}, write the resulting congruence as
\begin{align}
 \sum_{i=1}^m c_i p_i+\sum_{i<j}d_{ij}\rho_{ij}=N,
 \qquad N\in\bbZ.  
\end{align}
Here $\sum_i c_i=0$, since the $p$-terms are sums of differences.
Replacing $N$ by $N\sum_i p_i$ and using rational independence gives
$d_{ij}=0$ and $c_i=N$ for every $i$. Summing yields $mN=\sum_i c_i=0$,
so $N=0$ and all coefficients vanish.
Thus the coefficients of each $p_i$ in \eqref{tb9} and \eqref{tb10}
are equal.
However, 
the coefficient of $p_{j'}=p_{i'_s}$ in \eqref{tb9} is $-n'_s\le0$,
while it is 1 in \eqref{tb10}.
This is the desired contradiction, which shows that this case 
cannot occur; this completes the proof of the lemma.
\end{proof}

\begin{theorem}\label{TB}
Suppose that $\vp$ is ergodic and consider the multi-valued version of $T$.
Then the compact invariant set $\Kx$
is a finite union of good blocks,
i.e., a finite union of polytopes in $\Ho$
that have non-empty interiors and furthermore have every facet parallel to
one of the hyperplanes  $\set{\vx\in\Ho:x_i=x_j}$. 
\end{theorem}
\begin{proof}

Fix some $0<\gd\le 1/m$ and let 
 $\cN_\gd$ and $\hcN_\gd$ be as in the proof of \refT{Tinv}.
Let
$B$ be the block in \eqref{tb4}, 
where we choose all $\rho_{ij}$ as in \refL{LTB} and so
small that $B$ is a subset of the neighbourhood
$\cN_\gd$  of $0$ in $\Ho$. 
Let $\Bo$ be the interior of $B$, and let
$\hB:=\pi(B)$ and 
$\hBo:=\pi(\Bo)$ be their projections in $\hHo$. 

Let $\hvu\in \hHo$.
Since $R_\vp$ is an ergodic
translation in $\hHo$, there exists $k\ge0$ such
that $R_\vp^{-k}\hvu\in \hBo$, and thus
$\hvu\in R_\vp^{k}(\hBo)$.
Hence, $\hHo=\bigcup_{k=0}^\infty R_\vp^k(\hBo)$.
Since $\hBo$ is open, so is every $R_\vp^k(\hBo)$,
and thus there exists a finite $N$ such that
$\hHo=\bigcup_{k=0}^N R_\vp^k(\hBo)$.

Let $\vu\in \Kx$. Then $\hvu:=\pi(\vu)\in\hHo$, 
and thus $\hvu=R_\vp^k(\hvv)$ for some $k\le N$ and $\hvv\in\hBo$.
We have shown in \refT{Tinv} that $\Kx$ is (backwards) invariant, and thus
there exists 
$\vv\in\Kx$ such that $T^k\vv=\vu$.
This means that $\hvu=\pi(T^k\vv)=R_\vp^k(\pi(\vv))$ and thus
$\pi(\vv)=R_\vp^{-k}\hvu=\hvv\in \hBo \subseteq\hcN_\gd$.
By \refStep{Tinv4} in the proof of \refT{Tinv}
(taking $n=0$), this implies $\vv\in\cN_\gd$, and since we
also have $\pi(\vv)\in\hBo$, and $\pi$ is injective on $\cN_\gd$, we have
$\vv\in\Bo$. Consequently, $\vu\in T^k(\Bo)\subseteq T^k(B)$.
This shows that
\begin{align}
  \Kx\subseteq\bigcup_{n=0}^N T^n B.
\end{align}
Conversely, by \refStep{Tinv3} in the proof of \refT{Tinv},
$B\subseteq\cN_\gd\subseteq\Kx$.  Since $\Kx$ is forward invariant, it follows
that $T^nB\subseteq\Kx$ for every $n\ge0$.  Hence
\begin{align}
  \Kx=\bigcup_{n=0}^N T^n B.
\end{align}
By \refL{LTB}, each $T^nB$ is a finite union of good blocks, and the result
follows.  
\end{proof}

\begin{remark}
  \refT{TB} shows that $\Kx$ is a finite union of (convex) polytopes, and
  thus $\Kx$ is a polyhedral set in the sense that its boundary is included in
  a finite set of hyperplanes. Note that $\Kx$ in general (and we believe
  essentially always) is \emph{not} convex; this is seen in the figures in
  \refS{Sm=3} (for $m=3$), see also \cite{adler10,adler17}. 
\end{remark}

\begin{remark}\label{Rapr27}
  Every facet of $\Kx$ lies in a translate of a cutting hyperplane
$\set{\vx\in\Ho:x_i+p_i=x_j+p_j}$.
The translation vector can be taken as
$\sum_{k=1}^m n_k(\vp-\ve_k)$ with non-negative integer coefficients.
Indeed, the proofs of \refT{TB} and \refL{LTB} use only such translations
of the cutting hyperplanes, apart from the hyperplanes in \eqref{tb6}
arising from the original facets of $B$.
To exclude the latter, repeat the construction with $\rho'_{ij}$ chosen
so that the indexed family consisting of $p_1,\dots,p_m$ and
$\rho_{ij},\rho'_{ij}$, $1\le i<j\le m$, is rationally independent.
A facet lying in a translated original facet for the first construction
would have to lie in a supporting hyperplane from the second construction.
Equality of these hyperplanes would give a rational relation with a
nonzero coefficient of $\rho_{ij}$, which is impossible; any integer term
can be absorbed using $1=\sum_i p_i$. Thus only translates of cutting
hyperplanes remain.
\end{remark}

\begin{lemma}\label{Lpoly-support}
Let $K\subseteq\Ho$ be a finite union
\begin{align}
K=\bigcup_{i=1}^N B_i,  
\end{align}
where each $B_i$ is a compact polytope with non-empty interior in $\Ho$.
Then
\begin{equation}\label{tx0}
\supp(\gl_K)=K.
\end{equation}
\end{lemma}

\begin{proof}
The inclusion $\supp(\gl_K)\subseteq K$ is immediate, since $\gl_K$ is
supported on $K$.
Conversely, let $\vx\in K$, and let $U$ be an open neighbourhood of $\vx$
in $\Ho$. Choose $i$ such that $\vx\in B_i$.
Because $B_i$ is a polytope with non-empty interior in $\Ho$, it is the
closure of its interior; hence
\begin{equation}
U\cap \operatorname{int}_{\Ho}(B_i)\neq\emptyset.
\end{equation}
Since $\operatorname{int}_{\Ho}(B_i)$ is open in $\Ho$, the set
$U\cap \operatorname{int}_{\Ho}(B_i)$ has positive Lebesgue measure in
$\Ho$, and therefore
\begin{equation}
\gl_K(U)>0.
\end{equation}
Thus every open neighbourhood of $\vx$ has positive $\gl_K$-measure, which
shows that $\vx\in\supp(\gl_K)$. This proves \eqref{tx0}.
\end{proof}

We can now improve \refT{Tinv} by adding some further conclusions.

\begin{theorem}\label{TX}
Suppose that $\vp$ is ergodic and consider the multi-valued version of $T$.
Then the compact invariant set $\Kx$
from Theorem~\ref{Tinv} is a finite
union of compact full-dimensional polytopes in $\Ho$.
Thus,
\begin{equation}\label{tx2}
\supp(\gl_{\Kx})=\Kx.
\end{equation}
Moreover, for every $\vu\in\Ho$ and every branch orbit
$\vz'=(\vz'_n)_{n\ge0}$ starting from $\vu$,
\begin{equation}\label{tx1}
\omega(\vz')=\cL(\vu)=\Kx.
\end{equation}
\end{theorem}

\begin{proof}
The polyhedral statement is Theorem~\ref{TB}. 
Hence, Lemma~\ref{Lpoly-support}
shows \eqref{tx2}.
Finally, \eqref{tx2}
and  Theorem~\ref{TsupportK} yield
\begin{equation}\label{ts24}
\Kx=\supp(\gl_{\Kx})
\subseteq
\omega(\vz')
\subseteq
\cL(\vu)
\subseteq
\Kx,
\end{equation}
which yields \eqref{tx1} and completes the proof.
\end{proof}

As a consequence, our temporary definition $\cLx(\vu)$ in \eqref{cLx} is no
longer needed.

\begin{corollary}\label{CTX}
Suppose that $\vp$ is ergodic and consider the multi-valued version of $T$.
Then for every $\vu\in\Ho$,
\begin{equation}\label{tx3}
\cLx(\vu)=\cL(\vu)=\Kx.
\end{equation}
In particular, the limit set is forward invariant in the ergodic case.
\end{corollary}

\begin{proof}
Theorem~\ref{TX} gives $\cL(\vu)=\Kx$, and $\Kx$ is invariant by
Theorem~\ref{Tinv}. Hence $\cL(\vu)$ is forward invariant and equal to
$\Kx$, so also $\cLx(\vu)=\cL(\vu)=\Kx$.
\end{proof}

\section{Orbits in the intermediate case}\label{Sinter}

The intermediate case $0<\mu<m-1$ combines features of the two extreme cases
discussed in \refSs{Srational} and \ref{Sergodic}.

As discussed in \refSS{SSproj},
the projection of the dynamical system to $\hHo\subset\bbT^m$ 
has a simple structure given by the translation $\Rp$.
When $\mu<m-1$, the closure of the projected orbit
$\bgO{\0}$ is the closed subgroup $G$ of $\hHo$ of dimension $\mu<\dim \hHo$.
In particular, $\bgO{\0}=G$ is a null set in $\hHo$.
It follows that no individual projected limit set $\hcL(\hvu)$ is a universal
limit set:
\refL{LP} shows that
every  initial point $\hvu\in\hHo$ has a limit set
$\hcL(\hvu)=\bgO{\hvu}\subset\hHo$ 
which is a coset of $G$, but different initial points typically yield different
cosets. Similarly,
the action of
$T$ on $\Ho$ splits up into cosets of the subgroup $\pi\qw(G)$,
and thus no single limit set is a universal attractor.

We believe that the restriction of the dynamical system to each of these
cosets has a structure that is similar to the ergodic case studied in
\refS{Sergodic}, but we give here only a partial result and 
two examples; we  leave a more 
detailed exploration for future research.

By \refL{LP}, the projection $\pi\bigpar{\bcO\vu}$ is a coset
of $G$.
It follows from  \refL{LG} that
the inverse image $\pi\qw(G)$ is a $\mu$-dimensional closed
subgroup of $\bbR^m$,
and is thus (\eg{} by Lie group theory)
the union of a countable number of $\mu$-dimensional (affine) planes.
Since $\bcO\vu$ is compact by \refC{Cu}, it follows that each orbit
$\cO(\vu)$ and each limit set
$\cL(\vu)$ is a subset of a finite union of $\mu$-dimensional planes in
$\bbR^m$. 

We have the following partial result.
Let $\gl_\mu$ denote $\mu$-dimensional Lebesgue (= Hausdorff) measure.
\begin{theorem}\label{Tmu}
Let $0<\mu<m-1$.
Consider the multi-valued version, and
  suppose that $\vp$ and $\vu$ are such that
$\gl_\mu\bigpar{\cL(\vu)\cap\cD}=0$.
Then there exists an integer $\chi\ge1$ such that if\/ $\psi(\hvx)$ is as in
\eqref{psi} with $A=\cL(\vu)$, then
$\psi(\hvx)=\chi$ for $\gl_\mu$-\aex{} $\hvx\in \pi(\cL(\vu))=G+\hvu$.
\end{theorem}
\begin{proof}
Let $m_{G+\hvu}$ denote normalized Haar measure on the compact coset
$G+\hvu$.  By \refL{LG}, $G\cong F\times\bbT^\mu$ with $F$ finite; on each
connected component this Haar measure is the corresponding $\mu$-dimensional
Lebesgue measure, up to the fixed normalization.  Since
$G=\overline{\{n\hvp:n\in\bbZ\}}$ by \eqref{lp3}, translation by $\hvp$ is
ergodic, indeed uniquely ergodic, on $G$, and therefore the sequence
$n\hvp+\hvu$ is equidistributed in the coset $G+\hvu$.

Put $A:=\cL(\vu)$.  The set $A$ is compact and backward invariant by
\refL{Lback}, and \refL{LP} gives $\pi(A)=G+\hvu$.  Define
$\psi(\hvx)=|A\cap\pi\qw(\hvx)|$ on $G+\hvu$, as in \eqref{psi}.  The proof of
\refL{LC} applies on the coset: backward invariance gives
$\psi(\hvx-\hvp)\ge\psi(\hvx)$ outside
$\pi(A\cap\cD)+\hvp$.
This exceptional set has $m_{G+\hvu}$-measure zero by the assumption
$\gl_\mu(A\cap\cD)=0$, since $\pi$ is locally Lipschitz.
Integrating over the Haar measure on
$G+\hvu$ gives equality a.e., and ergodicity of the translation by $\hvp$ on
this coset implies that $\psi$ is a.e.\ constant.  Since $\pi(A)=G+\hvu$, this
constant is at least one.
\end{proof}

\begin{example}\label{Einter1}
 Let $m=3$ and $\vp=(p,\frac13,\frac23-p)$ with $p\in(0,\frac23)$  irrational.
Then $\mu=1<m-1$. 
Consider for simplicity the orbit $\vz(n):=T^n(\0)$ with $\vu=\0$,
and recall that then $\vz(n)\in \Ho=\set{(x,y,-x-y)}$.
It follows immediately from the definition \eqref{tx}, or from \eqref{rj}, 
that $z_2(n)\in \frac13\bbZ$, and more precisely $z_2(n)\in \frac{n}3+\bbZ$.
Furthermore, \refT{T<} shows that $-\frac23\le z_j(n)\le \frac56$, 
and \refR{Rleft}
shows that
$z_j(n)=-\frac23$ is possible only with $z_1(n)=z_3(n)=\frac13$,
which is impossible since $p$ is irrational.
Consequently, for every $n\ge0$,
\begin{align}\label{ei1}
  z_2(n)\in\bigset{-\tfrac13,0,\tfrac13,\tfrac23},
\end{align}
and, more precisely,
\begin{align}\label{ei2}
&&  z_2(n) &= 0, && n\equiv 0 \pmod 3,&&\notag\\
&&  z_2(n) &= \tfrac13, && n\equiv 1 \pmod 3,&&\notag\\
&&  z_2(n) &\in\bigset{-\tfrac13,\tfrac23}, && n\equiv 2 \pmod 3.&&
\end{align}
For $y\in\bbR$, let
\begin{align}\label{eil}
  \cL_y:=\set{x\in\bbR:(x,y,-x-y)\in \cL(\0)}.
\end{align}
Then 
\eqref{ei1} shows that
$\cL_y=\emptyset$ if $y\notin\bigset{-\tfrac13,0,\tfrac13,\tfrac23}$,
and thus
\begin{align}\label{ei3}
  \cL(\0)=\bigcup_{y\in\set{-\tfrac13,0,\tfrac13,\tfrac23}}
\bigset{(x,y,-x-y):x\in\cL_y}.
\end{align}

If $n\equiv 1\pmod3$, then \eqref{ei2} shows $z_2(n)=\frac13$ and thus
$z_1(n)+z_3(n)=-\frac13$. Furthermore, $z_j(n)\ge-\frac23$ by \refT{T<}, and
it follows that $-\frac23\le z_j(n)\le-\frac13-(-\frac23)=\frac13$, $j=1,3$;
hence $\cL_{1/3}\subseteq[-\frac23,\frac13]$.
On the other hand, considering the projection to $\bbT=\bbR/\bbZ$,
$\set{z_1(n) \mod 1:n\ge0}$ is dense in $\bbT$,
and so is $\set{z_1(3n+1) \mod 1:n\ge0}$ since also $3p$ is irrational.
 Hence the projection of $\cL_{1/3}$ equals $\bbT$, and thus
$\gl(\cL_{1/3})\ge1$.
It follows that 
\begin{align}\label{ei4}
\cL_{1/3}=[-\tfrac23,\tfrac13].
\end{align}

Suppose now that $0<p<\frac13$. (The case $\frac13<p<\frac23$ follows by
interchanging $p_1$ and $p_3=\frac23-p_1$.)
If $z(3n+1)={(x,\frac13,-x-\frac13)}$, 
where thus $x\in(-\frac23,\frac13)$,
then
$s(3n+2)={(x+p,\frac23,-x-p+\frac13)}$, and  it follows that 
2 wins if $x+p>-\frac13$, and otherwise 3 wins. Thus
\begin{align}\label{ei5}
  z(3n+2) = 
  \begin{cases}
    (x+p,-\frac13,-x-p+\frac13), & x\in(-p-\frac13,\frac13),
\\
(x+p,\frac23,-x-p-\frac23), & x\in(-\frac23,-p-\frac13).
  \end{cases}
\end{align}
Since, as shown above, the possible values of $x$ form a dense subset of
$(-\frac23,\frac13)$, 
and each $\cL_y$ is closed,
it follows that
\begin{align}\label{ei6}
  \cL_{-1/3}&=\bigset{x+p: x\in[-p-\tfrac13,\tfrac13]}=[-\tfrac13,p+\tfrac13],
\\\label{ei7}
  \cL_{2/3}&=\bigset{x+p: x\in[-\tfrac23,-p-\tfrac13]}=[p-\tfrac23,-\tfrac13]
.\end{align}
Finally, it is seen in the same way that if 
$\vz(n)=(x,\frac23,-x-\frac23)$ with $x\in\cL_{2/3}$,
then 2 wins next, leading to $x+p\in\cL_0$ with 
$x+p\in(2p-\tfrac23,p-\tfrac13)$.
Similarly, if $\vz(n)=(x,-\frac13,-x+\frac13)$ with $x\in\cL_{-1/3}$,
then 1 wins (and $x+p-1\in\cL_0$)
if $x+p>\frac12$, and otherwise 3 wins (and $x+p\in\cL_0$);
note that the first case cannot occur if $2p+\frac13<\frac12$, \ie,
$p<\frac{1}{12}$. 
The density argument above then shows (we omit the details)
\begin{align}\label{ek1}
  \cL_0=
  \begin{cases}
    [2p-\frac23,p-\frac13]\cup[p-\frac13,2p+\frac13]
=   [2p-\frac23,2p+\frac13]
,&0<p< \frac{1}{12}
\\    [2p-\frac23,p-\frac13]\cup[-\tfrac12,2p-\frac23]\cup [p-\frac13,\frac12]
=[-\frac12,\frac12]
,& \frac{1}{12}<p<\frac13,
  \end{cases}
\end{align}
where the decompositions show the pieces corresponding to different winners
in the last round
(in the order 2,3 and 2,1,3).

The sets $\cL_y$ are thus compact intervals.
Away from finitely many endpoints and tie points, $T$ shuffles pieces of
these intervals. The  limit set $\cL(\0)$ need not be invariant.
For example, for any $p\in(0,\frac13)$,
\eqref{ek1} shows that $\frac13-p\in \cL_0$, and thus by \eqref{ei3}, 
\begin{align}\label{ekw}
 \vw:=(\tfrac13-p,0,p-\tfrac13)\in\cL(\0).
\end{align}
At this point in the limit set we have
$
 \vw+\vp=(\tfrac13,\tfrac13,\tfrac13),
$
and thus a three-way tie for $T$.
Choosing party 2 at this tie gives
$T_2\vw=(\tfrac13,-\tfrac23,\tfrac13)\notin\cL(\0)$.
Thus $\cL(\0)$ is not forward invariant for the multi-valued map,
although the orbit starting at $\0$ has no ties, as is easily seen.

Furthermore,
the point $\vv:=(\tfrac13,\tfrac13,-\tfrac23)\in\cL(\0)$ 
by \eqref{ei3}--\eqref{ei4},
and its only  preimage in $\cL(\0)$
is $\vw$ in \eqref{ekw}, which has $T_3\vw=\vv$ and thus
requires party 3 to win at $\vw+\vp$.
Hence,
for the single-valued map given by the lexicographic rule, 
the limit set $\cL(\0)$ is not backward invariant.

For the projection to $\bbT^3$, we see immediately  that 
$\Rp^{3n}\0\in \bbT':=\set{(\hx,0,-\hx):\hx\in\bbT}$.
Since also $3p$ is irrational, it follows that
$\set{\Rp^{3n}\0}=\set{3n\hvp}$ is dense in $\bbT'$, and thus 
$G=\hcL(\0)\supseteq \bbT'$. It follows that
\begin{align}
  \hcL(\0)=G=G_0\cup G_1\cup G_2,
\quad\text{with}\quad 
G_k:=\bigset{\bigpar{\hvx,\tfrac{k}{3},-\hvx-\tfrac{k}{3}}:\hvx\in\bbT}.
\end{align}
We see that $\pi$ maps $\cL_0\to G_0$, $\cL_{1/3}\to G_1$ and 
$\cL_{-1/3}\cup\cL_{2/3}\to G_2$; these mappings are bijective except at the
  endpoints of the intervals, and in particular $\chi=1$ in \refT{Tmu}.
Note how $\cL_{-1/3}$ and $\cL_{2/3}$ combine in the projection.
\end{example}

\begin{example}\label{Einter2}
   Let $m=4$ and $\vp=(p,\frac14,\frac14,\frac12-p)$ with $p\in(0,\frac14)$  
irrational.
(The case $p\in(\frac14,\frac12)$ follows by interchanging $p_1$ and $p_4$.)
Then $\mu=1<m-1$.
Consider $\vz(n):=T^n(\0)$.
This is similar to \refE{Einter1}, but ties will occur, and we will see how
this affects the limit set.

By \refT{T<} and \refR{Rleft}, $-\frac34< z_j(n)<\frac{13}{12}$, and it
follows that
\begin{align}
  z_2(n),z_3(n)\in\bigset{-\tfrac12,-\tfrac14,0,\tfrac14,\tfrac12,\tfrac34,1}
.\end{align}
Moreover, $z_2(n)-z_3(n)\in\bbZ$.

In particular, if $n\equiv 1\pmod4$, then $z_2(n)=z_3(n)=\frac14$.
Define, in analogy with \eqref{eil},
\begin{align}\label{ejl}
  \cL_{y,z}:=\set{x\in\bbR:(x,y,z,-x-y-z)\in \cL(\0)}.
\end{align}
Then $\cL_{\frac14,\frac14}\subseteq[-\frac34,\frac14]$, and by considering
the projection to $\bbT$, it follows similarly to \eqref{ei4} that
\begin{align}\label{ej4}
  \cL_{\frac14,\frac14}=[-\tfrac34,\tfrac14].
\end{align}

Consider first the single-valued version with
the lexicographic rule for ties, meaning that 2 wins over 3 (the only tie
possible along this orbit).
Arguing similarly to \refE{Einter1}, we then find
\begin{align}\label{ej5}
  \cL_{\frac12,\frac12}&=[p-\tfrac34,-\tfrac12],
\\\label{ej6}
  \cL_{-\frac12,\frac12}&=[-\tfrac12,p+\tfrac14],
\\\label{ej7}
  \cL_{-\frac14,\frac34}&=[2p-\tfrac34,p-\tfrac12]\cup[p-\tfrac12,-\tfrac14]
=[2p-\tfrac34,-\tfrac14],
\\\label{ej8}
  \cL_{-\frac14,-\frac14}&=[-\tfrac14,2p+\tfrac14],
\\\label{ej9}
\cL_{0,0}&=
           \begin{cases}
             [3p-\tfrac34,3p+\tfrac14], &0<p<\tfrac1{12},
\\
             [-\frac12,\frac12],
 &\tfrac1{12}<p<\tfrac14,
           \end{cases}
\end{align}
and all other $\cL_{y,z}$ empty.

Again, all $\cL_{y,z}$ are compact intervals, and
\eqref{ej4}--\eqref{ej9} show that $\pi:\cL(\0)\to G=\hcL(\0)$
is bijective except at the endpoints of the intervals.
The limit set $\cL(\0)$ is not invariant 
for the single-valued map defined by the lexicographic rule.
Indeed, $\cL(\0)$ contains $\vv:=(-\tfrac12,\tfrac12,\tfrac12,-\tfrac12)$,
whose only possible preimage in the displayed set is
\begin{align}
  \vw:=(-p-\tfrac12,\tfrac14,\tfrac14,p),\qquad
\text{with}\quad
 \vw+\vp=(-\tfrac12,\tfrac12,\tfrac12,\tfrac12).
\end{align}
The image $\vv$ requires $T_4$ and thus a win for party 4, 
whereas the lexicographic rule selects
party 2. Thus $\vv\notin T(\cL(\0))$.

However, if we instead consider the multi-valued version, then we also have
$\cL_{\frac12,-\frac12}=\cL_{-\frac12,\frac12}$
and 
$\cL_{\frac34,-\frac14}=\cL_{-\frac14,\frac34}$.
It follows that 
in this case one part of $\hcL(\0)$ is covered singly
by $\cL(\0)$,
but another part is covered doubly. 
More precisely, if $\hvx=(\hx_1,\hx_2,\hx_3,\hx_4)\in\hcL(\0)=G$, 
then $\hx_2=\hx_3\in\set{0,\tfrac14,\frac12,\frac34}$,
and, except for a finite number of points,
$\psi(\hvx)=2$ 
if $\hx_2=\frac12$ and $\hx_1\in\cL_{-\frac12,\frac12}$
or $\hx_2=\frac34$ and $\hx_1\in\cL_{-\frac14,\frac34}$,
while
$\psi(\hvx)=1$,
otherwise.
Hence, \refT{Tmu} does not hold for this example,
because of the (frequent) occurrence of ties.
\end{example}

\begin{problem}
  Explore the intermediate case further.
\end{problem}

\section{The asymptotic bias}\label{Sbias}

We may now study the asymptotic bias of Frege's election method,
still for the time-homogeneous case.
When $n$ seats have been elected, the \emph{seat excess} 
(over strict proportionality)
for party $j$ is, using \eqref{rj},  
\begin{align}
  r_j(n)-np_j = -z_j(n).
\end{align}
Of course, for any given vector $\vp$, these numbers will fluctuate as $n$
grows, so to see whether some parties are favoured systematically,
one takes some kind of average to obtain an (average) bias.
We will here use the approach in  \cite{SJ262},
which essentially was used already by \citet{StL2} to compare
his new election method with d'Hondt's:
we fix $\vp$ and let $n$ be uniformly random in \set{1,\dots,N};
we then let $N\to\infty$ 
and consider the limit
of the expectation, \ie, the limit of the average
\begin{align}\label{bias}
\bgD_j(N):= - \frac{1}N\sum_{n=1}^Nz_j(n).
\end{align}
(In the case of ties, we allow any choice.)
If this limit exists,
we call it the (asymptotic) \emph{bias} and
denote it by $\bgDx_j$.
In the ergodic case, we have a simple formula.

\begin{theorem}\label{Tbias}
Suppose that $\vp$ is ergodic. Let $\Kx$ and $\nu$ be the invariant attractor 
and the invariant probability measure in \refTs{Tinv} and \ref{Tgl}.
Then, as \Ntoo,
\begin{align}\label{tbias}
  \bgD_j(N) \to \bgDx_j=-\int_{\Kx} x_j\dd\vx.
\end{align}
\end{theorem}
\begin{proof}
By \refC{Cu}, the orbit is contained in a compact set $K\subset\Ho$.
Let $f_j$ be a bounded continuous function on $\Ho$ that agrees with
$-x_j$ on $K\cup\Kx$. It follows from \refT{Tgl}\ref{Tglb} that
\begin{align}
  \bgD_j(N)&
=\frac{1}{N}\sum_{n=1}^N f_j(\vz(n))
=\int f_j(\vx)\dd\nu_N(\vx)
\notag\\&
\to \int_{\Ho}f_j(\vx)\dd\glKx(\vx)
=\int_{\Kx}f_j(\vx)\dd\gl(\vx)
=-\int_{\Kx}x_j\dd\vx.
\end{align}
\end{proof}

\begin{example}\label{Ebias2}
Consider the simple case $m=2$, with $\vp=(p_1,1-p_1)$.
If $p_1$ is irrational, then $\vp$ is ergodic and $\Kx=[-\frac12,\frac12]$, see
\refE{Em=2I} below; thus \refT{Tbias} yields the asymptotic bias $\bgDx_j=0$
for $j=1,2$.

If $p_1$ is rational, with $p_1=q_1/Q$ in lowest terms, then it follows from
\refE{Em=2Q} below that if $Q$ is odd, then there are no ties and 
$z(n)$ is asymptotically uniformly distributed on the finite set \eqref{em=2q};
hence again $\bgDx_j=0$. 
If $Q$ is even and, say, ties are resolved in favour of party 1,
then \refE{Em=2Q} and \eqref{em=2q} similarly yield an asymptotic bias
\begin{align}\label{bias2}
\bgDx_1=  -\frac{1}Q\sum_{-Q/2}^{Q/2-1}\frac{k}Q = \frac{1}{2Q}.
\end{align}
If ties are resolved in favour of party 2, the bias \eqref{bias2} obviously
changes sign. Hence, for $m=2$, the bias is only a result of the
tie-breaking rule.
\end{example}

\begin{example}\label{Ebias3}
  Let $m=3$ and suppose that $\vp$ is ergodic.
We show in \refS{Sm=3} 
that the space of possible vectors $\vp$ is partitioned into an infinite
number of different stability domains with different shapes of the
attractor $\Kx$, and we give a construction
of $\Kx$  as an explicit polygon (\refT{TD}).
The bias then can be calculated by \refT{Tbias} 
for $\vp$ in any given stability domain.
For example, we obtain from \eqref{tbias} and 
the formulas in \refS{Sm=3} by straightforward
(computer-assisted) calculations the following 
formulas for the bias.

If $0<p_2<p_1<p_3<1$ and $2p_1+4p_2>1$ 
(the domain $S_{1/2}$ in \refSS{SSm3S12}, see \eqref{eq:S12}), 
then 
\begin{align}
\bgDx_1&= -\tfrac16 (p_1-p_2)(2-3p_1-3p_2)
\\
\bgDx_2&= -\tfrac13(p_1-p_2)(p_1+2p_2-1)
\\
\bgDx_3&= -\tfrac16{(p_1-p_2)^2}
.\end{align}

If $2p_1+4p_2<1$,
$3p_1+3p_2>1$, 
and $3p_1-3p_2<1$,
(the domain $S_{1/3}\SA$ in \refE{ER13}, 
see \eqref{DS13A}),  
then 
\begin{align}
\bgDx_1&= -\tfrac16 (p_1-p_2)(2-3p_1-3p_2)
\\
\bgDx_2&= -\tfrac13(p_1-p_2)(3p_1+6p_2-2)
\\
\bgDx_3&= -\tfrac16{(p_1-p_2)(2-3p_1-9p_2)}
.\end{align}

If $3p_1+3p_2<1$,
$4p_1+2p_2>1$, 
and $3p_1+6p_2>1$,
(the domain $S_{1/3}\SB$ in \refE{ER13},
see \eqref{DS13B}),  
then 
\begin{align}
\bgDx_1&= -\tfrac16\bigpar{6p_1^{2}+18p_1p_2+12p_2^{2}-4p_1-8p_2+1}
\\
\bgDx_2&= -\tfrac13(p_1-p_2)(3p_1+6p_2-2)
\\
\bgDx_3&= -\tfrac16(6p_1-1)(1-2p_1-4p_2)
.\end{align}

If $p_2<p_1$,
$4p_1+2p_2<1$, 
and $3p_1+6p_2>1$,
(the domain $S_{1/3}\SC$ in \refE{ER13},
see \eqref{DS13C}),  
then 
\begin{align}
\bgDx_1&= -\tfrac13 (p_1-p_2)(2-5p_1-4p_2)
\\
\bgDx_2&= -\tfrac13(p_1-p_2)(3p_1+6p_2-2)
\\
\bgDx_3&= -\tfrac23{(p_1-p_2)^2}
.\end{align}

It is easily checked that
on common boundaries of two of the domains above, the formulas agree; 
hence the asymptotic biases $\bgDx_j$ seem to be given by
some continuous functions of $\vp$.
However, the explicit formulas above are typically different in different
domains, 
although there are some coincidences,
and there is no obvious general pattern. 
The bias is thus much more complicated than for the election methods treated
in \cite{SJ262} where a single simple formula is found.
\end{example}

\begin{remark}
  In principle, the invariant attractor $K$ gives also the asymptotic
  distribution of the seat excess $-z_j(n)$ (with $n$ random in $[N]$ as
  above and letting \Ntoo). However, even for $m=3$, the
  rather complicated shape of $K$ suggests that the result is neither nice
  nor very illuminating. Cf.\ \cite{SJ262} where this asymptotic distribution is
  found for some other election methods.
\end{remark}

\begin{remark}
We have here kept $\vp$ fixed.
  Another, cruder, way to quantify the bias of an election method
is to instead regard $\vp$ as
random, and for example uniformly distributed over 
the subset $\set{p_1>p_2>\dots>p_m}$ of $\gS_m$,
and consider the expectation of the seat excess, letting \Ntoo.
This was  done already  by \citet{Polya1918a,Polya1918b,Polya1919}, and has more
recently
been studied by a number of authors, see \cite{SJ262} for references.
(Already \citet{StL2} did a closely related calculation.)
The bound \eqref{t<} is uniform in $\vp$ and $n$, and \refT{Tbias} gives
convergence of \eqref{bias} for almost every $\vp$.
Hence,
with the number of seats averaged uniformly over $\{1,\dots,N\}$ as
above, dominated convergence gives
\begin{align}
 \E_{\vp}\bigl[\bgD_j(N)\bigr]\longrightarrow
 \E_{\vp}\bigl[\bgDx_j\bigr].
\end{align}
Thus the asymptotic average bias exists at least as a Ces{\`a}ro limit
and is obtained by integrating the result of
\refT{Tbias} over $\vp$.
However, even for $m=3$, the formulas obtained above differ between domains,
and we do not have a general formula. We therefore leave the explicit
evaluation of this average bias as an open problem.
\end{remark}

\section{The attractor for the case $m=2$}\label{Sm=2}
The case $m=2$ is very simple, but we include it for completeness and to
illustrate some general features.
We have $\Ho=\set{(x,-x):x\in\bbR}$ and thus it suffices to consider the first
coordinate.
We define $T_1:\bbR\to\bbR$ as the action of $T$ on the first coordinate;
i.e.,  
\begin{align}
  \label{TT1}
T(x,-x)=(T_1x,-T_1x).
\end{align}
For any $\vx=(x,-x)\in\Ho$, we have 
$\vx+\vp=(x+p_1,-x-p_1+1)$, and it follows from \eqref{tx}
that $T_1$ is given by
\begin{align}\label{T1b}
  T_1 x =
  \begin{cases}
    x+p_1,& x+p_1\le 1/2,
\\
x+p_1-1, & x+p_1\ge 1/2,
  \end{cases}
\end{align}
with the usual interpretation at a tie, either as a single-valued map with
a fixed choice or as a multi-valued map.
The dynamical system then has the following description.
Recall that $\round{x}$ denotes $x$ rounded to the nearest integer,
and that we regard $\round x$ as double-valued when $x\in\frac12+\bbZ$.
\begin{theorem}\label{T2}
  Let $m=2$, and let $z(n):=T_1^nu$ for some $u\in\bbR$. 
Then
\begin{align}\label{t2a}
  z(n)=u+np_1-r_1(n),
\end{align}
where $r_1(n)$ is the number of wins by party 1 up to time $n$.
Furthermore:
  \begin{romenumerate}
  \item \label{T2a}
If\/ $u\in[-\frac12,\frac12]$,    
then $z(n)\in[-\frac12,\frac12]$ for all $n$,
and
we have 
\begin{align}\label{t2b}
  r_1(n)=\round{u+np_1}.
\end{align}
\item \label{T2b}
If\/ $u>\frac12$, then $r_1(n)=n$ and $z(n)=u+np_1-n=u-n(1-p_1)$ until $z(n)\le \frac12$;
the evolution for larger $n$ then is as in \ref{T2a}, 
including \eqref{t2b}.

\item \label{T2c}
If\/ $u<-\frac12$, then $r_1(n)=0$ and $z(n)=u+np_1$ until $z(n)\ge -\frac12$;
the evolution for larger $n$ then is as in \ref{T2a},
including \eqref{t2b}.
  \end{romenumerate}
\end{theorem}
\begin{proof}
First, \eqref{t2a} follows from \eqref{T1b}. 
(This is a version for $m=2$ of a general fact, cf.\ \eqref{rj} for $u=0$.)

\pfitemref{T2a}
  It is obvious from \eqref{T1b} that $T_1$ maps $[-\frac12,\frac12]$ into
  itself, and thus every $z(n)\in[-\frac12,\frac12]$.
Hence, \eqref{t2a} yields
\begin{align}
  \bigabs{(u+np_1)-r_1(n)} = \bigabs{z(n)} \le \tfrac12,
\end{align}
  which shows \eqref{t2b}.

\pfitemx{\ref{T2b},\ref{T2c}}
Follow directly from \eqref{T1b} and \ref{T2a}.
\end{proof}

\begin{remark}\label{RFrege2}
  Consider as in \refR{Rprop} Frege's election method as an election method
  for a multi-member constituency. Then \refT{T2} (with $u=0$) shows that
in the case of only 2 parties, with $n$ seats to distribute,   
the number of seats given to party $i$ is $\round{np_i}$.
This is exactly the result of the \emph{largest remainder method},
and it is well-known (and easy to verify) that in the case of two parties,
also \emph{\SL's method} gives the same result.
(Both methods are known under a number of different names, see \eg{}
\cite[pp.~76 and 70]{Pukelsheim}.)
For larger $m$, Frege's method differs from these two methods, and also from
other well-known methods; 
see \eg{} 
\cite[Example 7]{HLL}.
(In particular these two methods have asymptotic bias
$0$ for any $m$ \cite{SJ262}, 
which does not hold for Frege's method, see \refE{Ebias3}.)
\end{remark}

\begin{example}\label{Em=2I}
In the ergodic case (which for $m=2$ means that $p_1$ is irrational),
\refT{T2} illustrates \refTs{Tinv} and \ref{Tgl} in an explicit way.
We see that, projecting to the first coordinate, 
the invariant attractor $\Kx$ given by \refT{Tinv} equals
$[-\frac12,\frac12]$, and the convergence of the empirical measures
$\nu_N$ to Lebesgue measure on $[-\frac12,\frac12]$ in \refT{Tgl}
is equivalent to the well-known fact that the sequence of
fractional parts $\frax{np_1}:=np_1-\floor{np_1}$ is equidistributed.
We see also that in this case, the projection $\pi:\Kx\to\hHo\cong\bbT$
is injective except that $\pi(-\frac12)=\pi(\frac12)$, so the exceptional
null set of $\hvx$ in \refT{Tinv} consist of a single point.  

Note also that if we instead consider a single-valued version of $T$,
then 
for every $\vu\in\Ho$, 
the orbit $\cO(\vu)$ still eventually lies in $[-\frac12,\frac12]$ 
and is dense there,
and thus $\cL(\vu)=[-\frac12,\frac12]$.
However, if we use the lexicographic tie resolution rule giving priority to
party 1, then
$\frac12\notin T(\cL(\vu))$, and thus $\cL(\vu)$ is not backward invariant.
\end{example}

\begin{example}\label{Em=2Q}
  In the rational case, let $p_1=q_1/Q$ with $\gcd(q_1,Q)=1$.
For simplicity, start with $u=0$.
Consider first the single-valued version, where ties are resolved in favour
of party 1. Then we should round upwards in \eqref{t2b} when there is a tie,
and thus $z(n)\in[-\frac12,\frac12)$. The orbit $T_1^n(0)$ is periodic, with
period $Q$, and cycles through the $Q$ values
\begin{align}\label{em=2q}
  \lrset{\frac{k}{Q}: -Q/2 \le k< Q/2}
\end{align}
in some order. In particular, if $Q$ is odd there will never be a tie, while
if $Q$ is even, there will be a tie once in each period.
For even $Q$, 
any other deterministic tie resolution rule will give either the same
result, or replace $-1/2$ by $1/2$ once in each period.
We see also that if ties are resolved by lot, and $Q$ is even, 
then $\pm\frac12$ will appear with random signs; otherwise there is no
difference.
Similarly,
in the multi-valued version, we have exactly the same result when $Q$
is odd, while when $Q$ is even we get the set $\set{-\frac12,\frac12}$
instead of $-\frac12$. 
\end{example}

\section{The invariant polygon for the case $m=3$}\label{Sm=3}

We now specialize to the case $m=3$.
This is the first case after the simple one-dimensional case in \refS{Sm=2}.
The hyperplane $\Ho$ is two-dimensional, and thus
the orbits and limit sets are planar sets.  
In the
ergodic case, \refT{TB}
shows that the global attractor $\Kx$ of the multi-valued Frege map is a finite
union of polygons with non-empty interiors.  
We will see in the explicit construction
in \refSS{SSm3construction} that this union is connected and thus
a compact polygonal domain; furthermore it is non-convex.
Moreover, \refC{CTX} gives
\begin{align}
\cL(\vu)=\omega(\vz')=\Kx
\end{align}
for every initial point $\vu\in\Ho$ and every branch orbit $\vz'$.

In this section we describe this polygon explicitly and study its dependence
on the parameter $\vp$.  For any $\vp$, not necessarily ergodic, we construct
a compact polygonal domain $\cD$ satisfying $T(\cD)=\cD$.  If $\vp$ is
ergodic, Proposition~\ref{Pm3bridge} shows that $\cD$ is the global attractor
$\Kx$.

By symmetry, it suffices to consider
the lower part of the parameter triangle; we divide this into a 
countable number of
\emph{Farey cells} $\cR_{u/v}$, 
indexed by the reduced fractions $u/v\in(0,1/2]$
(see \refSS{SSm3farey}).
Each Farey cell is further subdivided 
by a finite number of straight lines, 
called \emph{(stability) walls},
into a finite number of 
\emph{stability domains};
for all parameter values in (the interior) of a given stability domain,
the invariant polygon $\cD$ 
can be described symbolically in the same way (see \refSS{SSsymbolic}),
and $\cD$ thus has the same combinatorial type.
More precisely, for all parameter values in one of the stability domains in a
Farey cell $\cR_{u/v}$, 
the polygon $\cD$  has exactly $4v+3$ vertices
(\refT{Tm3-4v+3}), 
and for each stability domain, 
the coordinates of the vertices are given by some affine functions of $\vp$.
Since the denominators $v$ are unbounded, infinitely many combinatorial types
occur. 
The boundaries of the Farey cells and 
the walls between stability domains in them
are straight lines with rational coefficients.
In particular, an ergodic $\vp$ belongs to the interior of exactly one
stability domain.
The stability walls are precisely the portions of the line segments
shown in \refFig{fig:m3-farey-triangles} that lie in open Farey cells;
see \refC{Cwalls}.

\begin{remark}\label{Rrand}
  Parameters on the boundary of a Farey cell or on a stability wall are
never ergodic.  
Parameters on a stability wall inside a Farey cell are included in \refT{TD}:
the construction in \refSS{SSm3construction} yields 
an invariant polygon which varies continuously, in the Hausdorff metric,
throughout each open Farey cell.
(In particular, for a parameter on a stability wall,
the same polygon is obtained in the limit when the parameter
is approached from either neighbouring stability domain, although the
combinatorial types differs.) 
Since the vertex coordinates are affine in each stability domain, the
formulas in \refSS{SSm3construction} 
define a limiting polygon also on the closure of the Farey cell;
however, in this limit different vertices or sides may merge,
so the number of vertices may decrease from the number in the open Farey cell.
The limiting polygon is invariant also on the boundary of the Farey cell;
this is verified after the proof of \refT{TD} below.
(It is not obvious because of the different branches that may occur at ties.)
Thus each closed Farey cell carries a continuously varying invariant polygon.
We conjecture that the limiting polygons from adjacent Farey cells agree on
their common boundary
and give a continuously varying invariant polygon $\cD$ for all parameters.
\end{remark}

\subsection{The planar form of the Frege map}\label{SSm3planar}

Let
\begin{align}\label{vp}
\vp=(p,q,r),\qquad r:=1-p-q.
\end{align}
Some results below are stated for the ergodic case when $p,q,r$ 
(or equivalently $p,q,1$) are rationally independent,
but we will not assume this unless explicitly stated.
However, we will, unless we explicitly consider boundary cases, 
throughout this section assume that $p,q,r$ are distinct and furthermore, 
without loss of generality by symmetry,  
\begin{equation}\label{m3-parameter-region}
0<q<p<r<1.
\end{equation}

We find it convenient to identify the hyperplane
\begin{align}
  \Ho=\{(x_1,x_2,x_3)\in\bbR^3:x_1+x_2+x_3=0\}
\end{align}
with $\bbR^2$ by
\begin{align}
  (x,y)\longleftrightarrow (x+\tfrac13-p,\;y+\tfrac13-q,\;-x-y+\tfrac13-r).
\end{align}
Under this identification, the three coordinates of $\vx+\vp$ are
\begin{align}\label{xy3}
   x+\tfrac13,\qquad y+\tfrac13,\qquad -x-y+\tfrac13.
\end{align}
Thus the winner regions are bounded by the three linear forms
\begin{align}
 f(x,y)&:=-x+y,\label{m3-f}\\
 g(x,y)&:=2x+y,\label{m3-g}\\
 h(x,y)&:=x+2y.\label{m3-h}
\end{align}
More precisely, if we denote by $\cA,\cB,\cC$ the regions where candidates
$1,2,3$ win, including ties, then
\begin{align}
\cA&=\{g\ge0,\ f\le0\},\label{m3-cA}\\
\cB&=\{f\ge0,\ h\ge0\},\label{m3-cB}\\
\cC&=\{h\le0,\ g\le0\}.\label{m3-cC}
\end{align}
The boundaries where there are ties and the regions overlap consist of
the half-lines
\begin{align}
\Delta_f&:=\{f=0,\ g\ge0\}
\qquad\text{(1 and 2 tie)},\\
\Delta_g&:=\{g=0,\ h\le0\}
\qquad\text{(1 and 3 tie)},\\
\Delta_h&:=\{h=0,\ f\ge0\}
\qquad\text{(2 and 3 tie)},
\end{align}
with a common intersection point (where all three tie)
\begin{align}\label{U00}
U :=(0,0).
\end{align}
The multi-valued Frege map is then the piecewise translation
\begin{equation}\label{m3-T}
T(x,y)=
\begin{cases}
T_1(x,y):=(x+p-1,y+q),& (x,y)\in\cA,\\
T_2(x,y):=(x+p,y+q-1),& (x,y)\in\cB,\\
T_3(x,y):=(x+p,y+q),& (x,y)\in\cC,
\end{cases}
\end{equation}
with two  values taken on the boundary half-lines, and three values at their
intersection $U$.

It is convenient to denote the translation vectors by
\begin{equation}\label{m3-XYZ}
X:=(p-1,q),\qquad Y:=(p,q-1),\qquad Z:=(p,q).
\end{equation}
We write translated points and polygons 
by a reduced notation, writing the translations as
\begin{align}\label{apT}
  T_1:V\mapsto VX, &&
  T_2:V\mapsto VY, &&
  T_3:V\mapsto VZ.
\end{align}
Thus, for a point $V\in\bbR^2$,
and non-negative integers $n_1,n_2,n_3$, we write
\begin{align}\label{fb11}
V X^{n_1}Y^{n_2}Z^{n_3}:=V+n_1X+n_2Y+n_3Z.
\end{align}
For translates of a polygon $Q$,
we use the same notation $QX$, $QY$,  and so on.

If $M_1,\dots,M_k$ is a finite sequence of points, then
the polygonal line through these points (in order) is denoted
\begin{align} \label{linje}
M_1\dotsm M_k
\quad\text{or sometimes (to avoid confusion)}\quad
\linje{M_1,\dots, M_k}
.\end{align}
We further let
\begin{align}\label{polyM}
[M_1,\dots,M_k]
\end{align}
denote the polygonal domain, including its boundary, bounded by the closed polygonal curve
$\linje{M_1,\dots,M_k,M_1}$
(assuming that this polygonal line has no self-intersection).

Finally, we write
\begin{align}\label{Delta}
\Delta:=\{(p,q):0<q<p<r\}
=\set{(p,q):0<q<p<1-p-q}
\end{align}
for the full parameter region that we consider.

\subsection{A consequence of the general theory}\label{SSm3bridge}

In the case $m=3$,
the general conclusions in Section~\ref{Sergodic}
take the following concrete planar form.
 
\begin{proposition}\label{Pm3bridge}
Assume that $p,q,1$ are rationally independent. Let $\cD\subset\Ho\cong\bbR^2$
be a compact invariant polygonal domain for the multi-valued Frege map, i.e.,
assume that $\cD$ is a polygon such that
\begin{align}
T(\cD)=\cD.
\end{align}
Then
\begin{align}
\cD=\supp(\gl_{\Kx})=\Kx=\cL(\vu)
\qquad\text{for every }\vu\in\Ho.
\end{align}
Moreover, for every initial point $\vu\in\Ho$ and every branch orbit $\vz'$
starting from $\vu$, we have
\begin{align}
\omega(\vz')=\cD,
\end{align}
and the empirical measures of the branch orbit converge weakly to Lebesgue
measure on $\cD$. In particular, every branch orbit eventually lies in $\cD$,
and every tail after entry is dense in $\cD$.
\end{proposition}

\begin{proof}
By \refT{Tinv}, $\Kx$ is the unique compact invariant set, and thus $\cD=\Kx$.
Furthermore, 
Theorem~\ref{TX} then gives 
$\cD=\Kx=\supp(\gl_{\Kx})=\omega(\vz')=\cL(\vu)$,
and Theorem~\ref{Tgl} yields the convergence of the empirical measures.
Eventual entry into $\cD$ follows from the absorption assertion in \refT{Tinv}.
\end{proof}

\begin{remark}\label{Rm3-nonergodic}
Outside the ergodic case, the identification $\cD=\cL(\vu)$ 
is not valid. In the rational and intermediate cases, each projected orbit on
the torus $\bbT^2$ is confined to a coset of a proper closed subgroup
(a finite coset in the rational case and otherwise a 1-dimensional one);
hence the orbit in $\Ho$ is either a finite set or lies on a finite set of
parallel lines, and thus the same holds for each
limit set $\cL(\vu)$.
In those
cases the correct conclusion from an explicit construction is only that $\cD$ is
a compact invariant polygon for the multi-valued map.
An individual limit set need not be forward invariant for the multi-valued map,
as Example~\ref{Einter1} shows.
See \cite{gun} for more on the dynamics in non-ergodic cases.
\end{remark}

\begin{remark}\label{Rparallel}
  \refT{TB} shows also that (in the ergodic case)
an invariant polygonal domain $\cD$ necessarily
  has all its sides parallel to one of $\gD_f, \gD_g,\gD_h$,
as can be seen in  the constructions below.
We call such sides $f$-sides, $g$-sides, and $h$-sides, respectively.
\end{remark}

\subsection{The first stability domain $S_{1/2}$}\label{SSm3S12}

To find the invariant polygon $\cD$,
we begin with $(p,q)$ in the triangular region $[J,K,M]$ with vertices
\begin{align}\label{JKM}
J:=\bigl(\tfrac13,\tfrac13\bigr),\qquad
K:=\bigl(\tfrac12,0\bigr),\qquad
M:=\bigl(\tfrac16,\tfrac16\bigr).
\end{align}
We will see that this is a stability domain in the sense defined in
\refSS{SSsymbolic} below, and
we denote it by
\begin{align}\label{eq:S12}
S_{1/2}:=[J,K,M]
=\bigsqpar{\bigl(\tfrac13,\tfrac13\bigr),
\bigl(\tfrac12,0\bigr),
\bigl(\tfrac16,\tfrac16\bigr)}.
\end{align}
More precisely, we 
define $S_{1/2}$ to be the open triangle $[J,K,M]$, and we
consider $(p,q)$ in this open triangle,
ignoring boundary cases.
(Note that this includes all cases with ergodic $(p,q)$ in the closed
triangle $[J,K,M]$.)

\begin{remark}\label{REX}
The vertex formulas below extend continuously to the closed triangle
$\overline S_{1/2}$ and define an invariant polygon, in which some of the
constructed points may coincide; see \refR{Rrand} and its justification
after the proof of \refT{TD}.
In particular,
for $(p,q)$ on the side $\linje{J,M}$ of $S_{1/2}$,
i.e., for $p=q\in[\frac16,\frac13]$,
the invariant polygon $\cD$ is the hexagon $\cH$ in \eqref{cH}.
\end{remark}

\begin{figure}[htbp]
\centering
\includegraphics[width=0.86\textwidth]{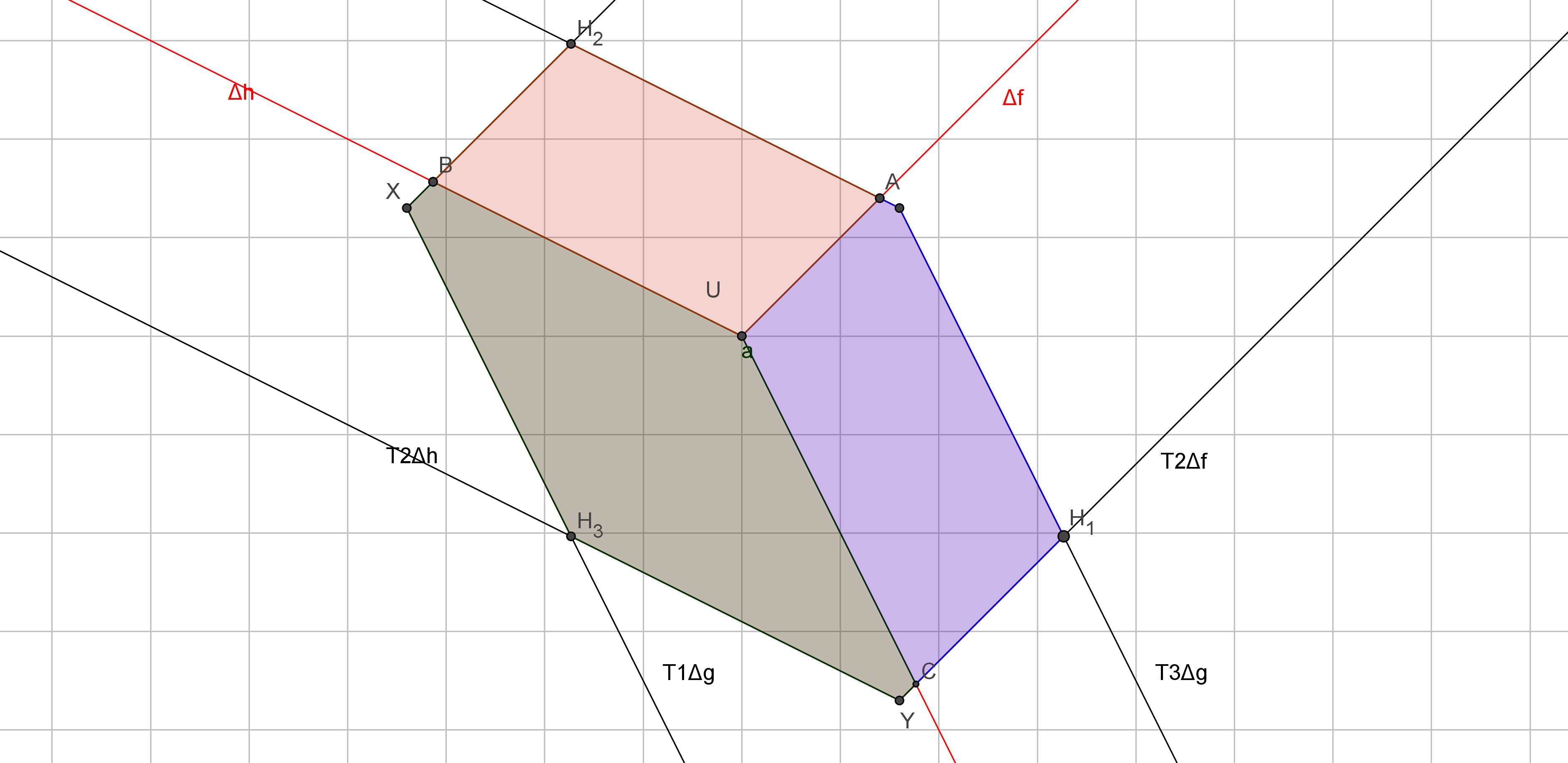}
\caption{The three boundary half-lines $\Delta_f$, $\Delta_g$, $\Delta_h$
  and the auxiliary hexagon $\cH$ 
  for $(p,q)$ in the first stability domain $S_{1/2}$.
($p=0.321, q=0.262$)
}
\label{fig:m3-hexagon}
\end{figure}

\begin{figure}[htbp]
\centering
\includegraphics[width=0.86\textwidth]{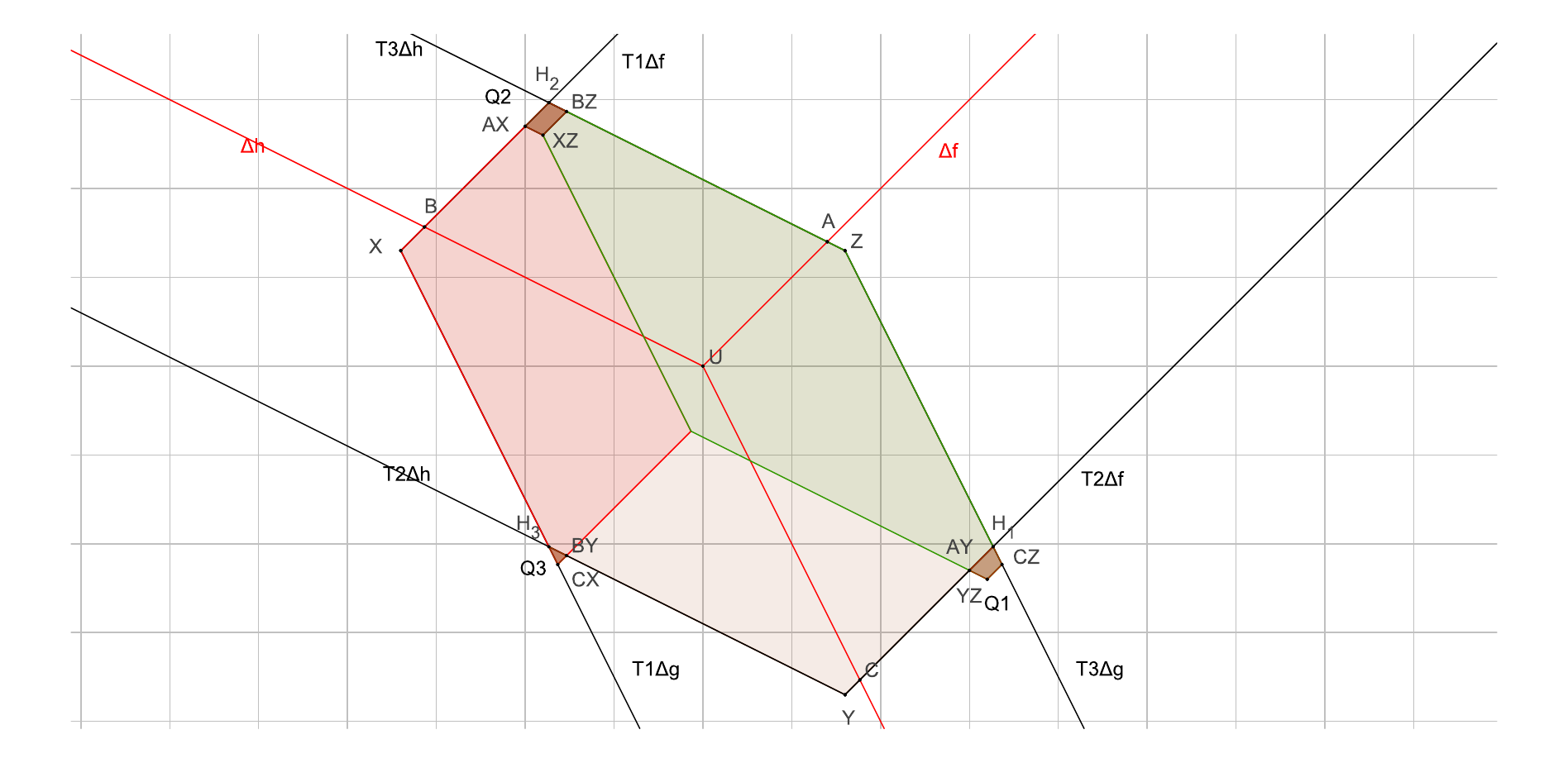}
\caption{Modification of the auxiliary hexagon $\cH$: the missing tetragon
  $Q_2$ is removed from one piece and the excess pieces $Q_1$ and $Q_3$ are
  added to the other pieces.
($p=0.321, q=0.262$)
}
\label{fig:m3-mod-hexagon}
\end{figure}

\begin{figure}[htbp]
\centering
\includegraphics[width=0.86\textwidth]{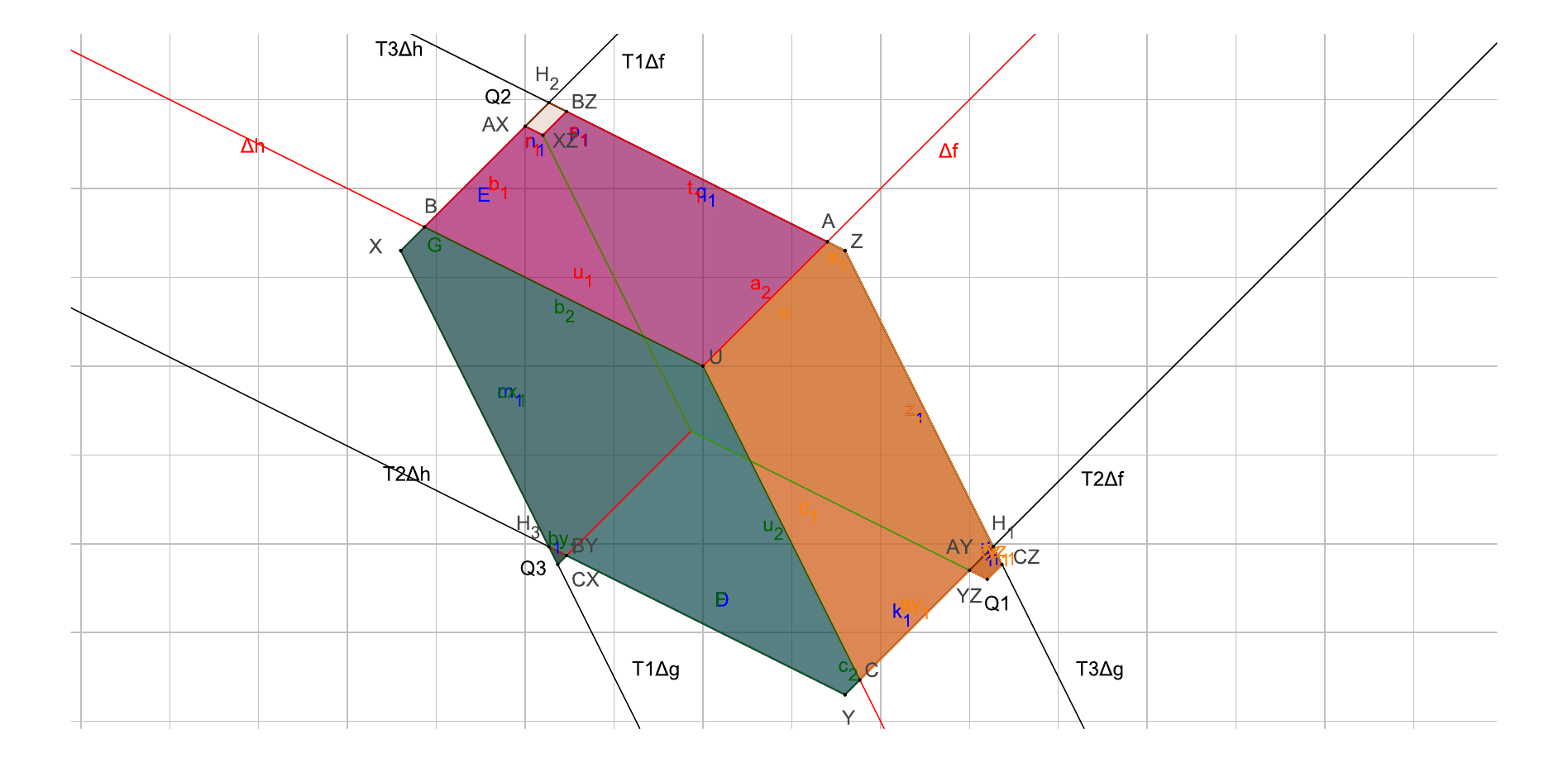}
\caption{The invariant polygon $\cD$ 
for $(p,q)$ in the domain $S_{1/2}$ after the modification of $\cH$.
($p=0.321, q=0.262$)
}
\label{fig:m3-S12-main}
\end{figure}

For $(p,q)\in S_{1/2}$ we define, besides 
$U=(0,0)$ in \eqref{U00} and
the translation vectors $X,Y,Z$ in
\eqref{m3-XYZ}, the points
\begin{align}
H_1&:=T_3\Delta_g\cap T_2\Delta_f=\bigl(p+\tfrac13,q-\tfrac23\bigr),\label{hi1}
\\
H_2&:=T_1\Delta_f\cap T_3\Delta_h=\bigl(p-\tfrac23,q+\tfrac13\bigr),\label{hi2}
\\
H_3&:=T_2\Delta_h\cap T_1\Delta_g=\bigl(p-\tfrac23,q-\tfrac23\bigr),\label{hi3}
\end{align}
and
\begin{align}
A&:=\Delta_f\cap T_3\Delta_h=\tfrac13(p+2q,p+2q),\label{pda}\\
B&:=\Delta_h\cap T_1\Delta_f=\tfrac13\bigl(2(p-q-1),-(p-q-1)\bigr),\label{pdb}\\
C&:=\Delta_g\cap T_2\Delta_f=\tfrac13\bigl(p-q+1,-2(p-q+1)\bigr).\label{pdc}
\end{align}
\refFig{fig:m3-hexagon} shows how
the hexagon 
\begin{align}\label{cH}
\cH:=[Z,H_1,Y,H_3,X,H_2]
\end{align}
splits into three pieces
\begin{align}
\cH_1&:=\cH\cap\cA=[U,A,Z,H_1,C],\\
\cH_2&:=\cH\cap\cB=[U,B,H_2,A],\\
\cH_3&:=\cH\cap\cC=[U,C,Y,H_3,X,B].
\end{align}
Translating these three polygons by $T_1,T_2,T_3$ yields, see \eqref{apT},
\begin{align}
T_1\cH_1&=[X,AX,XZ,H_1X,CX],\\
T_2\cH_2&=[Y,BY,H_2Y,AY],\\
T_3\cH_3&=[Z,CZ,YZ,H_3Z,XZ,BZ],
\end{align}
where we note the common vertex
\begin{align}
  H_1X=H_2Y=H_3Z=\bigpar{2p-\tfrac23,\;2q-\tfrac23}.
\end{align}
These translated polygons are disjoint up to their boundaries and almost cover
$\cH$, but three defects appear:
\begin{itemize}
\item a tetragon
\begin{align}
Q_2=[AX,H_2,BZ,XZ]
\end{align}
is left uncovered inside $\cH_2$;
\item a tetragon
\begin{align}
Q_1=[H_1,CZ,YZ,AY]
\end{align}
protrudes from $\cH_1$;
\item a triangle
\begin{align}
Q_3=[H_3,BY,CX]
\end{align}
protrudes from $\cH_3$.
\end{itemize}
To obtain an invariant set, we modify $\cH$ as follows, 
see \refFigs{fig:m3-mod-hexagon} and \ref{fig:m3-S12-main}.
We remove the uncovered tetragon $Q_2$ from $\cH_2$
and add the excess pieces $Q_1$ and $Q_3$ to $\cH_1$ and $\cH_3$,
respectively. 
This yields the three polygons
\begin{align}
\cP_1&=[U,A,Z,CZ,YZ,AY,C],\label{pd1}\\
\cP_2&=[U,B,AX,XZ,BZ,A],\label{pd2}\\
\cP_3&=[U,C,Y,BY,CX,X,B].\label{pd3}
\end{align}
Their union is the polygonal domain
\begin{equation}\label{m3-D12}
\cD=[CX,X,AX,XZ,BZ,Z,CZ,YZ,AY,Y,BY].
\end{equation}
The point $U$ lies in the interior of $\cD$; 
all other vertices of $\cP_1,\cP_2,\cP_3$ lie on the boundary of $\cD$.
We will call $\cP_2$ the \emph{middle piece}
and $\cP_1,\cP_3$ \emph{side pieces}.

\begin{figure}[htbp]
\centering
\includegraphics[width=0.6\textwidth]{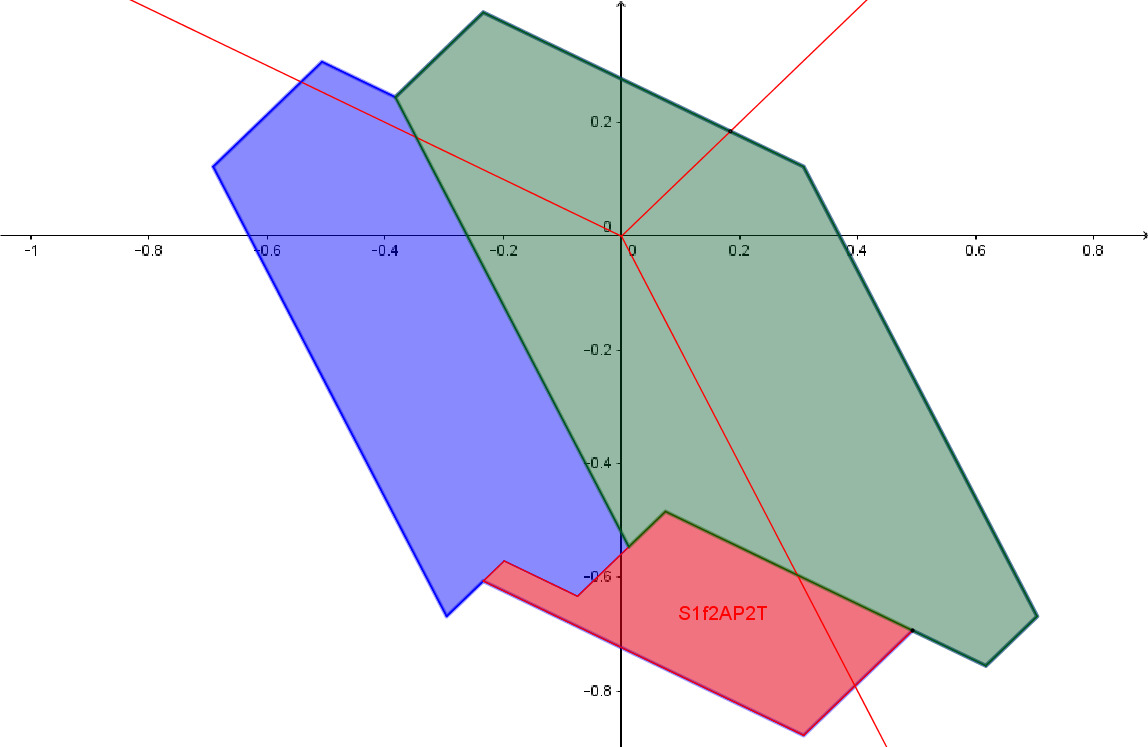} 
\caption{An invariant polygon in $S_{1/2}$ together with the translated pieces
  $T_i\cP_i$.}
\label{fig:m3-S12-translated}
\end{figure}

\begin{theorem}\label{Tm3-S12}
For every $(p,q)\in S_{1/2}$, the three polygons $\cP_1,\cP_2,\cP_3$ lie in the
winner regions of candidates $1,2,3$, respectively, and
\begin{align}
\cD=\cP_1\cup\cP_2\cup\cP_3.
\end{align}
Moreover,
\begin{align}
T_1\cP_1&=[X,AX,XZ,CXZ,XYZ,AXY,CX],\label{sv1}\\
T_2\cP_2&=[Y,BY,AXY,XYZ,BYZ,AY],\label{sv2}\\
T_3\cP_3&=[Z,CZ,YZ,BYZ,CXZ,XZ,BZ],\label{sv3}
\end{align}
these three images have pairwise disjoint interiors, and
\begin{align}
\cD=T_1\cP_1\cup T_2\cP_2\cup T_3\cP_3.
\end{align}
Consequently, 
$\cD$ is a compact invariant polygonal domain for the multi-valued
Frege map.
\end{theorem}

\begin{proof}
The formulas for the images $T_i\cP_i$ follow immediately from the translation
rules \eqref{apT}. The fact that the pieces $\cP_i$ lie in the correct
winner sectors is verified by checking the signs of the linear forms
\eqref{m3-f}--\eqref{m3-h} at the vertices; since the winner regions are
convex, this suffices. 
Note that this verification uses the assumption that
$(p,q)\in S_{1/2}=[J,K,M]$; for example, we have 
\begin{align}\label{hAX}
  h(AX)=h(XZ)=2p+4q-1>0,
\end{align}
while equality holds on the boundary segment $\linje{K,M}$.
Finally, it is easily verified that the three translated polygons
$T_1\cP_1, T_2\cP_2, T_3\cP_3$ fit together like jigsaw pieces in the
interior,
and that 
the outer boundary of
$T_1\cP_1\cup T_2\cP_2\cup T_3\cP_3$ equals the boundary
\eqref{m3-D12} of $\cD$.
(See \refFig{fig:m3-S12-translated}, and recall that
every side has one of 
the three directions in
\refR{Rparallel}.)
Hence the
union of the three translated polygons is exactly $\cD$.
\end{proof}

\begin{corollary}\label{Cm3-S12-ergodic}
Assume that $(p,q)\in S_{1/2}$ and that $p,q,1$ are rationally independent.
Then
\begin{align}
\cD=\Kx=\cL(\vu)\qquad\text{for every }\vu\in\Ho.
\end{align}
Moreover, for every branch orbit $\vz'$ of the multi-valued Frege map, one has
\begin{align}
\omega(\vz')=\cD.
\end{align}
\end{corollary}

\begin{proof}
This is an immediate consequence of Proposition~\ref{Pm3bridge} and 
\refT{Tm3-S12}.
\end{proof}

\subsection{Symbolic description, stability sets, and stability domains}
\label{SSsymbolic}

For $(p,q)$ outside the domain $S_{1/2}$, the invariant polygon $\cD$ is no
longer given by \eqref{m3-D12}, as is already seen from \eqref{hAX}.
The construction below nevertheless shows that 
all sides of $\cD$ lie on lines obtained from one of
the three tie-lines
$f=0$, $g=0$, $h=0$
by a translation of the type $X^{n_1}Y^{n_2}Z^{n_3}$ with $n_1,n_2,n_3\in\bbNo$.
(This is also shown by \refR{Rapr27} for the ergodic case, and this implies
the general case since the polygons $\cD$ that we construct depend continuously on $(p,q)$ within each closed Farey cell.)
Hence the vertices of $\cD$ can be defined as the intersections of certain
pairs of such lines.
It is thus possible to give a \emph{symbolic description} of the invariant
polygon as
a finite symbolic pattern  describing $\cD$ and its three pieces
$\cP_1$, $\cP_2$, $\cP_3$: a finite list of translated tie-lines that
support the
edges, the incidences that define the vertices, the cyclic order of the
vertices on each piece, and the winner sector containing each piece.
The pattern also records the order and side choices needed to specify
the pieces and their translated images; boundary equalities are recorded
separately. Thus alternative order or side choices are different patterns.
This can be formalized in different ways; we will 
use descriptions such as \eqref{m3-D12} and \eqref{pd1}--\eqref{pd3}
using supplementary definitions as in \eqref{pda}--\eqref{pdc}.
We will also refer to the data in the symbolic description as the
\emph{combinatorial type} of $\cD$.

\begin{definition}\label{Dm3-realization}
Fix a finite symbolic pattern $w$ describing one of the invariant polygons
constructed in this section and its
three pieces
$\cP_1$, $\cP_2$, $\cP_3$.
The \emph{stability set} $\Omega_w$ 
is the set of parameters $(p,q)\in\Delta$
for which the invariant polygon is described by $w$.
A \emph{stability domain} is 
a connected component with non-empty interior of a stability set;
it follows from the next proposition
that this is the same as an open component of a stability set.
\end{definition}

\begin{proposition}\label{Pm3-realization}
For every fixed symbolic pattern $w$,
the stability set $\Omega_w$ is cut
out in parameter space by finitely many affine relations and strict affine
inequalities. 
In particular, every stability domain is an open convex polygonal subset of
$\Delta$.  More generally, each stability set is relatively open in the affine
subspace determined by its equalities.  Moreover, on each stability set the
combinatorial type of the invariant polygon is constant and all designated
vertices depend affinely on $(p,q)$.
\end{proposition}

\begin{proof}
Each line supporting a side is obtained from one of $f=0$, $g=0$, $h=0$ by a word
in the translations $T_1,T_2,T_3$. Since the vectors $X,Y,Z$ depend affinely
on $(p,q)$, every such supporting line has an equation of the form
\begin{align}
F(x,y)=c(p,q),
\end{align}
where $F\in\{f,g,h\}$ and $c$ is affine in $(p,q)$
with rational coefficients.
Every prescribed vertex equals the intersection of two such lines, and is
(by simple linear algebra)
thus also affine in $(p,q)$. The prescribed order along a supporting line
and the prescribed side of another supporting line are determined by
differences of affine functions of $(p,q)$. Once these choices are fixed,
the remaining incidence, sector, and matching conditions are likewise
affine equalities or strict affine inequalities, with no alternative order
or side choices left to consider. Thus every connected
component is a convex polygonal set, relatively open in the affine subspace in
which its equalities hold. If the component has non-empty interior in $\Delta$,
then no non-trivial equality is active throughout it, and it is open in
$\Delta$. On such a component no incidence or order relation can change.
\end{proof}

\begin{remark}\label{Rboundaries}
As for $S_{1/2}$ above, we concentrate on stability domains, i.e.,
open stability sets, and ignore the modifications occurring at their
boundaries.
Formally, a stability domain is an open set, 
but we use the same polygonal notation \eqref{polyM} as for closed polygons;
the interpretation should be clear from the context. 
The same applies to the Farey cells defined below. 
Similarly,
we may say ``partition'' when we really mean a partition up to
the boundaries of the various parts; the boundaries may overlap or be absent.
\end{remark}

\subsection{Farey cells and first entry}
\label{SSm3farey}

The first stability domain $S_{1/2}=[J,K,M]$ occupies only the upper part of the
parameter triangle \eqref{Delta}.  We write
\begin{align}
\Delta=[J,K,L],
\end{align}
where we recall \eqref{JKM} and further define $L:=(0,0)$. 
In the complementary triangle $[K,L,M]$, new translated tie-lines give
polygons with arbitrarily many vertices. We describe them using Farey cells.
As we will see in \refSS{SSm3construction},
within a Farey cell the upper boundary and the
middle hexagon have a fixed symbolic form; the only change occurs when the
$\gD_g$ cut passes through a vertex of the lower boundary.  Since all
vertices depend affinely on $(p,q)$, these events lie on affine lines,
which we (as above) call (stability) walls.

We first partition $[K,L,M]$ into a countable family of parallelograms called
\emph{Farey cells}.  The name reflects their construction from the Farey, or
Stern--Brocot, tree.  We shall then construct the invariant polygon throughout
each open cell and obtain its internal stability walls from the construction.

\begin{figure}[htbp]
\centering
\includegraphics[width=0.86\textwidth]{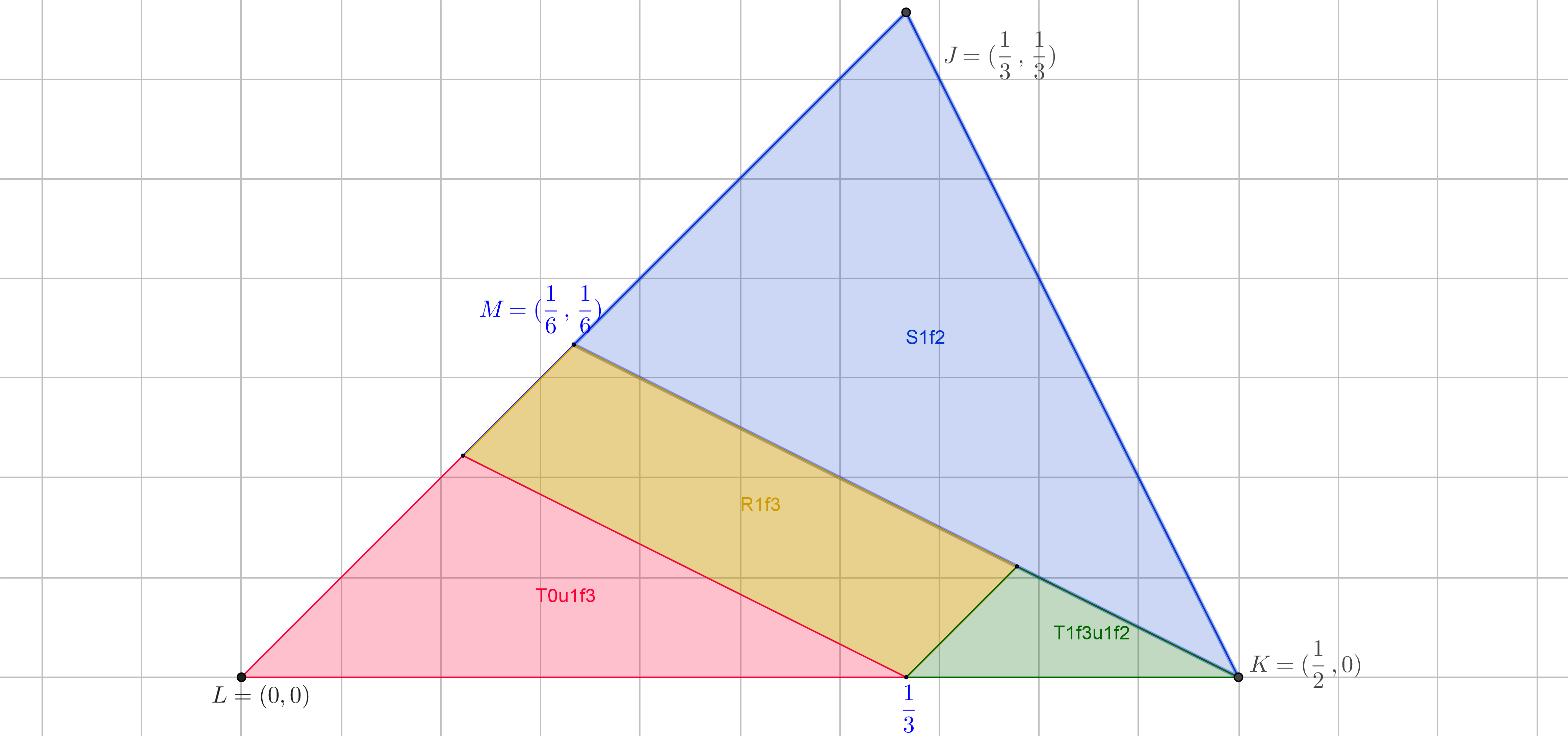}
\caption{The triangles $[J,K,M]$ and $[K,L,M]$ in the parameter space.}
\label{fig:m3-param-00}
\end{figure}

Recall the classical Farey tree, also called Stern--Brocot tree, 
of non-negative reduced fractions $a/b$ (with $a\ge0$ and $b>0$),
see e.g.\ \cite[p.~116--123]{CM}.
Two reduced fractions $0\le \frac ab<\frac cd$ are \emph{neighbouring}
if $bc-ad=1$;
this means that
\begin{align}\label{sb1}
  \frac cd - \frac ab = \frac 1{bd},
\end{align}
a relation that will be frequently used below (explicitly or implicitly). 
If $\frac ab$ and
$\frac cd$ are neighbouring reduced fractions, then
their \emph{mediant} (Farey sum) is defined as
\begin{align}\label{mediant}
\frac{a+c}{b+d}.
\end{align}
Note that $\frac ab < \frac{a+c}{b+d} < \frac cd$, and that
$(\frac ab,\frac{a+c}{b+d})$ and $(\frac{a+c}{b+d}, \frac cd)$
are two new neighbouring pairs.
We may thus continue taking mediants recursively, 
\emph{ad infinitum},
which gives the Farey tree; 
more precisely, we associate the interval $(\frac ab,\frac cd)$
to the mediant $\frac uv=\frac{a+c}{b+d}$; this interval is cut by 
$\frac uv$ into two, which are associated to their mediants
$\frac{a+u}{b+v}$ and $\frac{u+c}{v+d}$, and these two fractions are the
children of $\frac uv$.
We are interested in 
reduced fractions in $[0,\tfrac12]$ only,
so we start with the neighbouring pair $(\frac01,\frac12)$ and 
take successive mediants; this gives a part of the Farey tree
that contains every reduced fraction in $(0,\tfrac12)$
exactly once.
Hence, we have the following.

\begin{lemma}\label{Lm3-farey-unique}
For every reduced fraction $\frac uv\in(0,\tfrac12)$ there exists a unique
pair of neighbouring reduced fractions $\frac ab<\frac cd$ whose mediant is
$\frac uv$.
In other words, $\frac uv$ has a
unique representation
\begin{equation}\label{m3-farey}
\frac uv=\frac{a+c}{b+d},\qquad \text{with }bc-ad=1,
\end{equation}
where
$a\ge0$, $b,c,d>0$, and thus
\begin{align}\label{sb2}
0\le \frac ab<\frac uv<\frac cd\le \frac12.
\end{align}
\end{lemma}

We may then call $a/b$ and $c/d$ the \emph{(Farey) parents} of $u/v$,
but
note that the family analogy is not perfect: $u/v$ has two parents, but it
is a child of only one of them. (Since the Farey tree is a tree.)
It is easily checked that $u/v$ is the child of the parent 
with largest denominator ($b$ or $d$); the other ``parent'' is a more remote
ancestor in the Farey tree.
\begin{proof}
This is the standard uniqueness of the parents of $\frac uv$ in the
Farey tree, 
and may also be obtained from B\'ezout's identity;
see e.g.\ \cite[p.~116--118]{CM}.
\end{proof}

\begin{figure}[htbp]
\centering
\includegraphics[width=0.95\textwidth]{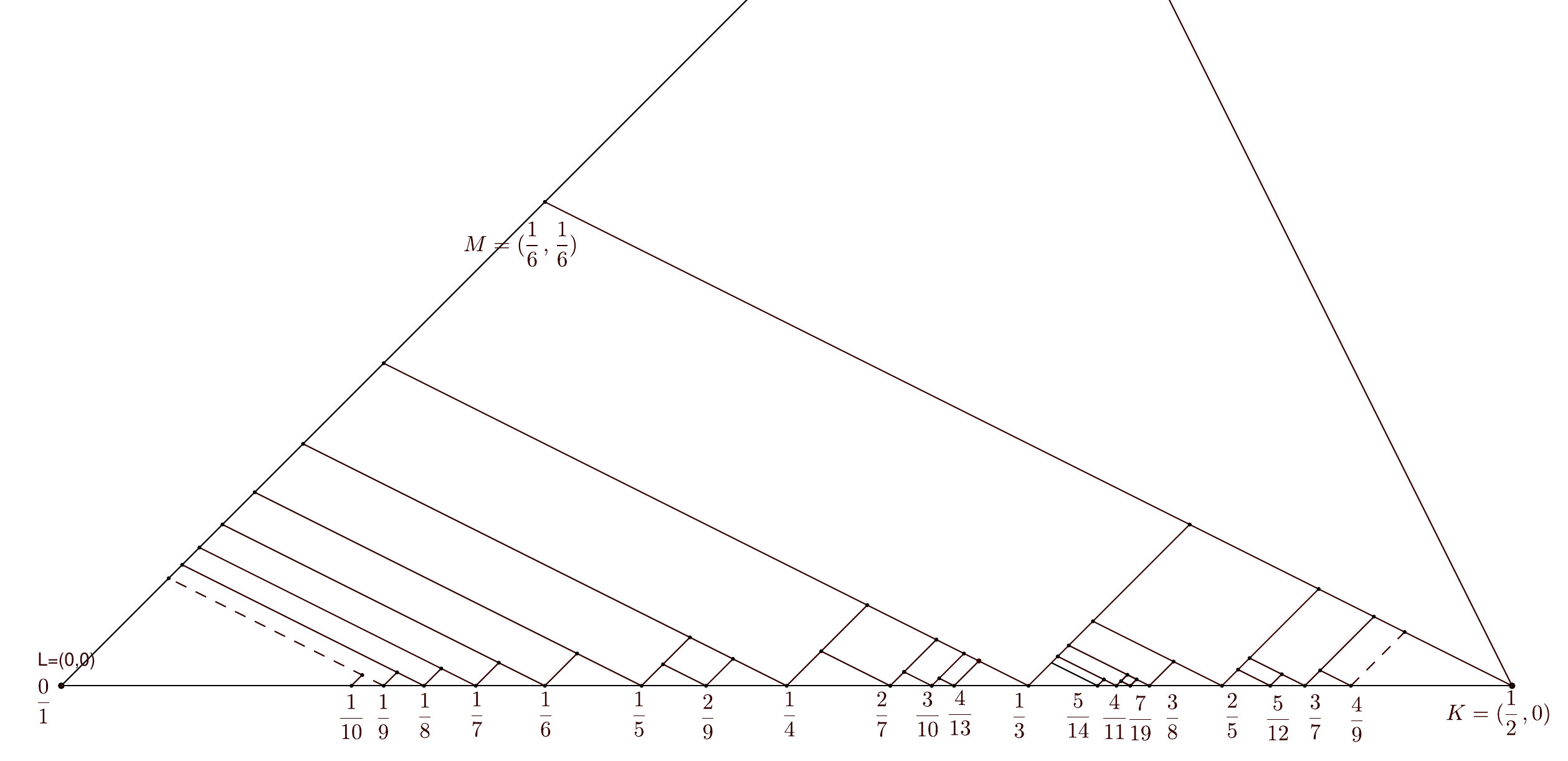}
\caption{Farey cells in the lower part of the parameter triangle.}
\label{fig:m3-farey-cells}
\end{figure}

To every reduced fraction $\frac uv\in(0,\tfrac12)$, or equivalently to
the corresponding neighbouring pair $\bigl(\frac ab,\frac cd\bigr)$ 
in \eqref{m3-farey},
we associate the 
points
\begin{align}\label{m3-fcu}
 U_{u/v}&:=U_{(\frac ab,\frac cd)}
:=\bigl(\tfrac ab+\tfrac1{3bd},\tfrac1{3bd}\bigr),
\\\label{m3-fcv}
 V_{u/v}&:=V_{(\frac ab,\frac cd)}
:=\bigl(\tfrac{a+c}{b+d},0\bigr)=\bigpar{\tfrac uv,0},
\\\label{m3-fcl}
 L_{u/v}&:=L_{(\frac ab,\frac cd)}
:=\bigl(\tfrac ab+\tfrac1{3b(b+d)},\tfrac1{3b(b+d)}\bigr)
=\bigl(\tfrac ab+\tfrac1{3bv},\tfrac1{3bv}\bigr)
,\\\label{m3-fcm}
 M_{u/v}&:=M_{(\frac ab,\frac cd)}
:=\bigl(\tfrac{a+c}{b+d}+\tfrac1{3d(b+d)},\tfrac1{3d(b+d)}\bigr)
=\bigl(\tfrac{u}{v}+\tfrac1{3dv},\tfrac1{3dv}\bigr)
.\end{align}
and the parallelogram 
which we call \emph{Farey cell}
(see \refFig{fig:m3-farey-cells}),
\begin{align}\label{m3farey2}
\cR_{u/v}:=\cR_{(\frac ab,\frac cd)}
:=[U_{u/v},L_{u/v},V_{u/v},M_{u/v}],
\end{align}
with sides  parallel to the lines $f=0$ and $h=0$.
(We mean by the Farey cell $\cR\uv$ the interior, i.e.,
the open parallelogram \eqref{m3farey2},
even when this is not said explicitly.)
We define also
the triangle 
\begin{align}\label{m3farey0}
\gD_{u/v}=\Delta_{(\frac ab,\frac cd)}
:=\bigsqpar{(\tfrac ab,0),U_{u/v},(\tfrac cd,0)}
=
\bigsqpar{(\tfrac ab,0),\bigl(\tfrac ab+\tfrac1{3bd},\tfrac1{3bd}\bigr),(\tfrac cd,0)}.
\end{align}
Note that the upper two sides of $\gD_{u/v}$ are 
parallel to the lines 
$f=0$ and $h=0$, and that $L_{u/v}$ and $M_{u/v}$ lie on these sides;
furthermore (identifying a point $(p,0)\in\bbR^2$ with $p\in\bbR$), 
$V_{u/v}$ is the mediant of the two lower vertices 
$V_{a/b}$ and $V_{c/d}$ of $\gD_{u/v}$.

\begin{proposition}\label{Pm3-farey-partition}
Up to boundaries, the triangle $[K,L,M]=\gD_{(\frac01,\frac12)}$ 
is partitioned into Farey cells
$\cR_{(\frac ab,\frac cd)}$ indexed by neighbouring reduced fractions
\begin{align}
0\le \frac ab<\frac cd\le \frac12.
\end{align}
Moreover,
\begin{align}\label{m3farey1}
\Delta_{(\frac ab,\frac cd)}
=
\Delta_{(\frac ab,\frac{a+c}{b+d})}
\cup \cR_{(\frac ab,\frac cd)}
\cup
\Delta_{(\frac{a+c}{b+d},\frac cd)},
\end{align}
again up to boundaries.
\end{proposition}

Equivalently, $[K,L,M]$ is partitioned into $\cR_{u/v}$ indexed by reduced
fractions $\frac uv\in(0,\frac12)$. We find it convenient to use both
notations.

\begin{proof}
The decomposition \eqref{m3farey1} follows from the definitions
\eqref{m3farey0} and \eqref{m3-fcu}--\eqref{m3farey2},
noting that $L_{(\frac ab,\frac cd)}=U_{(\frac ab,\frac{a+c}{b+d})}$
and
$M_{(\frac ab,\frac cd)}=U_{(\frac{a+c}{b+d},\frac cd)}$.
The decomposition into Farey cells $\cR_{(\frac ab,\frac cd)}$
then follows by applying \eqref{m3farey1} recursively.
Indeed, along any infinite branch of this recursion, $b+d\to\infty$ and
hence $bd\to\infty$. The corresponding triangles have height $1/(3bd)$,
which tends to zero. Thus no point with $q>0$ can remain in all of them.
\end{proof}

It is convenient to use the following coordinates:
\begin{equation}\label{m3-Phi}
\Phi(p,q):=(\alpha,\beta):=(p-q,\,p+2q).
\end{equation}
Then
\begin{align}\label{Phi1}
 p=\frac{2\alpha+\beta}{3},\qquad q=\frac{\beta-\alpha}{3},\qquad
 r=1-p-q=\frac{3-\alpha-2\beta}{3},
\end{align}
and therefore
\begin{align}\label{Phi2}
0<q<p<r
\iff 0<\alpha<\beta\quad\text{and}\quad \alpha+\beta<1.
\end{align}
Moreover,
\begin{align}\label{Phi3}
\Phi(K)=\bigl(\tfrac12,\tfrac12\bigr),\qquad
\Phi(L)=(0,0),\qquad
\Phi(M)=\bigl(0,\tfrac12\bigr),
\end{align}
so $[K,L,M]$ is carried to the
right-angled triangle
\begin{align}\label{Phi4}
\bigl\{(\alpha,\beta):0\le \alpha\le \beta\le \tfrac12\bigr\}.
\end{align}

\begin{proposition}\label{Pm3-farey-rect}
Let $\frac 01\le\frac ab<\frac cd\le\frac12$ 
be neighbouring reduced fractions with mediant
$\frac uv=\frac{a+c}{b+d}$. The closures of the corresponding regions satisfy
\begin{align}
\Phi\bigl(\overline{\Delta_{(\frac ab,\frac cd)}}\bigr)
&=
\Bigl\{(\alpha,\beta):\frac ab\le \alpha\le \beta\le \frac cd\Bigr\},
\label{m3-farey-rect1}
\\
\Phi\bigl(\overline{\cR_{(\frac ab,\frac cd)}}\bigr)
&=
\Bigl[\frac ab,\frac uv\Bigr]\times\Bigl[\frac uv,\frac cd\Bigr].
\label{m3-farey-rect2}
\end{align}
In particular, every Farey cell becomes a rectangle in the
$(\alpha,\beta)$-plane. 
Moreover, $(p,q)$ lies in the open Farey cell
$\cR_{u/v}$ if and only if $(\ga,\gb)=\Phi(p,q)$ satisfies
\begin{equation}\label{m3-open-farey}
        \frac ab<\alpha<\frac uv<\beta<\frac cd .
\end{equation}
\end{proposition}

\begin{proof}
Applying $\Phi$ to the three vertices in \eqref{m3farey0} gives,
using \eqref{sb1},
\begin{align}
\Bigl(\frac ab,\frac ab\Bigr),\qquad
\Bigl(\frac ab,\frac cd\Bigr),\qquad
\Bigl(\frac cd,\frac cd\Bigr).
\end{align}
This proves \eqref{m3-farey-rect1}.
For the Farey cell, the four vertices in 
\eqref{m3-fcu}--\eqref{m3-fcm} 
are mapped to
\begin{align}
\Bigl(\frac ab,\frac cd\Bigr),\qquad
\Bigl(\frac uv,\frac uv\Bigr),\qquad
\Bigl(\frac ab,\frac uv\Bigr),\qquad
\Bigl(\frac uv,\frac cd\Bigr),
\end{align}
and \eqref{m3-farey-rect2} follows, together with the version
\eqref{m3-open-farey} for the open Farey cell.
\end{proof}

For integers $k\ge0$ and $0\le m\le k$, put
\begin{equation}\label{m3-Pkm}
        P(k,m):=mX+(k-m)Z=(kp-m,kq).
\end{equation}
Then, in the Farey coordinates \eqref{m3-Phi},
\begin{equation}\label{m3-Pkm-sector}
        f(P(k,m))=m-k\alpha,
        \qquad
        h(P(k,m))=k\beta-m,
\end{equation}
and also, since $g=h-f$,
\begin{equation}\label{m3-Pkm-g}
        g(P(k,m))=k(\ga+\beta)-2m.
\end{equation}
Consequently,
recalling the definition \eqref{m3-cA}--\eqref{m3-cC} of $\cA,\cB,\cC$, 
and noting that $\ga<\gb$ by \eqref{m3-Phi}, we have
for $k\ge1$
\begin{align}
  P(k,m)\in\cA\inter &\iff  \frac{m}k < \ga, \label{qaA}
\\
  P(k,m)\in\cB\inter &\iff  \ga <\frac{m}k < \gb,\label{qaB}
\\
  P(k,m)\in\cC\inter &\iff  \gb <\frac{m}k;\label{qaC}
\end{align}
the equalities $m/k=\alpha$ and $m/k=\beta$ are the
boundary cases $P(k,m)\in\cA\cap\cB = \gD_f$
and $P(k,m)\in\cB\cap\cC = \gD_h$, respectively.

For $0\le k\le v$, define further
\begin{equation}\label{m3-nu-Pk}
        \nu(k):=\left\lceil (k-1)\frac uv\right\rceil,
        \qquad
        P_k:=P(k,\nu(k)).
\end{equation}
Note that $\nu(0)=0$ and thus $P_0=(0,0)=U$; this trivial case will be
convenient in later subsections, but here we will use $k\ge1$.

The step code $(\nu(k)-\nu(k-1))_{k=1}^v$, with $Z$ coded by $0$
and $X$ by $1$, is a cyclic shift of the upper Christoffel word with
$v-u$ zeros and $u$ ones \cite{BLRS09}.

We need another standard result for the Farey tree.
\begin{lemma}\label{Lm3-no-small}
If $r_1/s_1<r_2/s_2$ are neighbouring reduced fractions, then every reduced
fraction $m/k$ strictly between them has $k\ge s_1+s_2$.
\end{lemma}

\begin{proof}
We have
\begin{align}
        \frac{m}{k}-\frac{r_1}{s_1}
=\frac{ms_1-kr_1}{ks_1}
\ge\frac1{ks_1},
        \qquad
        \frac{r_2}{s_2}-\frac{m}{k}
=\frac{kr_2-ms_2}{ks_2}
\ge\frac1{ks_2}.
\end{align}
Adding and using
$r_2/s_2-r_1/s_1=1/(s_1s_2)$ gives $k\ge s_1+s_2$.
\end{proof}

\begin{theorem}
\label{Tm3-Farey-entry}
Let $u/v\in(0,1/2)$ be reduced, with Farey parents $a/b<c/d$, so that
$u=a+c$, $v=b+d$, and $bc-ad=1$.  
Assume that $(p,q)$ lies in the open Farey cell
$\cR_{u/v}$, or equivalently that \eqref{m3-open-farey} holds.
Then $U=P_0,P_1,P_2,\dots,P_v$ is (the beginning of) a branch orbit of $T$, 
defined by choosing $TU=P_1=Z$.
Furthermore, 
$P_v=P(v,u)=uX+(v-u)Z$ lies in $\cB\inter$, whereas $P_k\notin\cB$ for
$1\le k<v$.  More precisely, for $1\le k<v$,
if $\nu(k+1)=\nu(k)+1$, then
$P_k\in\cA\inter$ and $P_{k+1}=P_k+X$; if $\nu(k+1)=\nu(k)$, then
$P_k\in\cC\inter$ and $P_{k+1}=P_k+Z$; in both cases thus $TP_k=P_{k+1}$.
Thus the Farey cell 
$\cR_{u/v}$
is the parameter region in which the distinguished branch
from $U$ has no intermediate ties and first enters the sector $\cB\inter$
at the point $uX+(v-u)Z$.

Moreover the four sides of the Farey
cell are exactly the following parameter conditions, where one specified
point in this orbit
lies on a sector line, and thus there is a tie:
\begin{equation}
\begin{aligned}\label{ele1}
        \alpha=\frac ab &\Longleftrightarrow f(P_b)=0, &
        \alpha=\frac uv &\Longleftrightarrow f(P_v)=0,\\
        \beta =\frac uv &\Longleftrightarrow h(P_v)=0, &
        \beta =\frac cd &\Longleftrightarrow h(P_d)=0.
\end{aligned}
\end{equation}
\end{theorem}

\begin{proof}
First, \eqref{m3-nu-Pk} yields $\nu(v)=u$ and thus,
using also \eqref{m3-Pkm},
\begin{align}\label{jc2v}
  P_v=P(v,u)=uX+(v-u)Z=(vp-u,vq).
\end{align}
Hence, $P_v\in\cB\inter$ follows at once from
\eqref{qaB} and \eqref{m3-open-farey}.

If
$1\le k<v$, then $\nu(k)/k\ne u/v$ since $u/v$ is reduced. 
If $\nu(k)/k<u/v$, then $\nu(k+1)=\ceil{ku/v}>\nu(k)$
by \eqref{m3-nu-Pk},
and thus $\nu(k+1)=\nu(k)+1$ (since always $\nu(k+1)-\nu(k)\in\set{0,1}$);
hence $P_{k+1}=P_k+X$.
Furthermore, 
Lemma~\ref{Lm3-no-small}, applied to the neighbouring pair $a/b<u/v$,
implies $\nu(k)/k\le a/b<\alpha$, so 
by \eqref{qaA}, 
$P_k\in\cA\inter$ and thus $TP_k = P_k+X=P_{k+1}$.

If instead $\nu(k)/k>u/v$, then similarly $\nu(k+1)=\nu(k)$
by \eqref{m3-nu-Pk}
and thus $P_{k+1}=P_k+Z$; 
furthermore, 
Lemma~\ref{Lm3-no-small} applied to $u/v<c/d$
implies
$\nu(k)/k\ge c/d>\beta$, so 
$P_k\in\cC\inter$ by \eqref{qaC} and  thus $TP_k = P_k+Z=P_{k+1}$.

Hence the vertices $P_k$ form a
branch from $U$ that first enters $\cB$ at $P_v$.  Finally,
\begin{align}\label{jc1}
 \left\lceil (b-1)\frac uv\right\rceil
=
 \left\lceil \frac{av+1-u}v\right\rceil
=a,
 \qquad
 \left\lceil (d-1)\frac uv\right\rceil
= \left\lceil \frac {cv-1-u}v\right\rceil
=c,
\end{align}
because $bu-av=1$ and $du-cv=-1$, and $u<v/2<v-1$.  
Thus
\begin{align}
\label{jc2b}
  P_b&=P(b,a)= aX+(b-a)Z=(bp-a,bq),
\\\label{jc2d}
  P_d&=P(d,c)=  cX+(d-c)Z=(dp-c,dq).
\end{align}

Conversely, if the specified branch first enters $\cB\inter$ at $P_v$
without intermediate ties, then 
\eqref{jc2v} and \eqref{qaB} yield $\alpha<u/v<\beta$. 
Furthermore, since $b,d<v$,
the points $P_b$ and $P_d$ are outside $\cB$; 
thus it follows from \eqref{jc2b}--\eqref{jc2d}, \eqref{m3-cB},
and \eqref{m3-Pkm-sector} (see also \eqref{qaB}) that
$\alpha>a/b$ and $\beta<c/d$. This gives
\eqref{m3-open-farey}.

Finally, 
the four displayed equivalences in \eqref{ele1} follow from
\eqref{m3-Pkm-sector} for $(k,m)=(b,a)$, $(v,u)$, $(d,c)$,
using \eqref{jc2v} and \eqref{jc2b}--\eqref{jc2d}.
\end{proof}

\begin{remark}\label{R1/2}
  We have defined Farey cells $\cR_{u/v}$ for $\frac uv\in(0,\frac12)$, and
  shown that they form a partition of $[K,L,M]$ 
(see Figures~\ref{fig:m3-param-00} and \ref{fig:m3-farey-cells}).
It is convenient to define also the Farey cell
$\cR_{1/2}:=S_{1/2}=[J,K,M]$, so that $\set{\cR_{u/v}:\frac uv\in(0,\frac12]}$
is a partition of the full parameter region $\Delta$ in \eqref{Delta}. 
However, note that
this cell is exceptional in several ways; in particular, it is a triangle
and not a rectangle. (It is the intersection of the rectangle
\eqref{m3farey2} for $\frac uv=\frac12$ and $\Delta$.)
\end{remark}

\subsection{The invariant polygon for parameters in a Farey cell}
\label{SSm3construction}
We now construct the invariant polygon 
for $\vp$ in an arbitrary open Farey cell.  In
the ergodic case the invariant compact set is unique by \refT{Tinv}; for other
parameters in the cell the same construction still gives a compact invariant
polygonal domain.

We continue to assume \eqref{m3-parameter-region}.
We further continue with the assumptions and notation of \refT{Tm3-Farey-entry}.
In particular, we fix $(p,q)\in\cR_{u/v}$, where $u/v$ 
is reduced.
We assume also $u/v<\frac12$ as in \refT{Tm3-Farey-entry}
for convenience at a few places in the arguments,
but it should be noted
that the results below also hold for $u/v=1/2$, with $\cR_{1/2}=S_{1/2}$ as
in \refR{R1/2}, with a few minor modifications of the arguments;
the construction in \eqref{cD} below then yields exactly the same
polygon $\cD$ as was found by a different method in \refSS{SSm3S12},
which easily can be verified directly.
 
We also define
\begin{align}\label{m3-lambda-mu}
 \delta_-:=\frac uv-\alpha,\qquad \delta_+:=\beta-\frac uv.
\end{align}
Since $a/b$ and $u/v$, and also $u/v$ and $c/d$, are neighbouring,
\eqref{m3-open-farey} gives
\begin{align}\label{me-lambda-mu-b}
 0<\delta_-<\frac1{bv},\qquad 0<\delta_+<\frac1{dv}.
\end{align}

Define $\kk_j$, for $j=0,\dots,v$, as the integer in $\set{1,\dots,v}$
that satisfies
\begin{align}  \label{qa1}
\kk_j\equiv jb+1\pmod v.
\end{align}
In particular, $\kk_0=\kk_v=1$.
Note that $b$ and $v$ are relatively prime, since $bu-av=1$, 
and thus $\kk_1,\dots,\kk_v$ is a permutation of $1,\dots,v$.
Let
\begin{align}\label{qa2}
 W_0:=X,\qquad W_j:=P_{\kappa_j}\quad(1\le j<v),\qquad W_v:=Z.
\end{align}
Thus, $\set{W_j:1\le j\le v}$ is a permutation of the points 
$\set{P_k:1\le k\le v}$, recalling that we
in \refT{Tm3-Farey-entry} made the rather arbitrary choice
$P_1:=Z=T_3U$ for the multiple-valued $TU$. 
We here also include (as $W_0$) the alternative $X=T_1U$.
(We do not yet consider $T_2U=Y$.)
The following lemma shows that $W_0,\dots,W_v$ are the points
$X,P_1,\dots,P_v$ arranged in decreasing order of $f$, or in increasing
order of $g$ or $h$.
(Since $g(x,y)-f(x,y)=3x$, this is also the order of increasing $x$, i.e.,
from left to right.)
\begin{lemma}\label{Lm3-shadowing}
We have 
\begin{align}
 f(W_0)&>f(W_1)>\dots>f(W_v),\label{qa3f}
\\
 g(W_0)&<g(W_1)<\dots<g(W_v),\label{qa3g}
\\
 h(W_0)&<h(W_1)<\dots<h(W_v).\label{qa3h}
\end{align}
\end{lemma}

\begin{proof}
  Let
\begin{align}\label{qa4}
 t_k:=k\frac uv-\nu(k) \quad (1 \le k \le v)
  \end{align}
and
  \begin{align}\label{qa5}
 \widetilde t_0:=\frac uv-1,
 \qquad \widetilde t_j:=t_{\kappa_j}\quad(1\le j\le v).
\end{align}
We claim that then
\begin{align}\label{qa8}
 \widetilde t_j=\frac uv-1+\frac jv
=\frac{u+j-v}{v},
 \qquad 0\le j\le v.
\end{align}
In fact,
for $1\le j\le v$, we have by 
\eqref{qa5}, \eqref{qa4}, \eqref{qa1}, and 
$bu-av=1$, 
\begin{align}\label{qa6}
  \ttt_j = t_{\kk_j} 
\equiv \kk_j \frac{u}{v}
\equiv(jb+1)\frac{u}{v}
=\frac{u+jbu}v
\equiv\frac{u+j}{v}
\pmod\bbZ.
\end{align}
Furthermore, for $1\le k\le v$, \eqref{qa4} and \eqref{m3-nu-Pk} yield
\begin{align}\label{qa7}
  t_k-\frac{u}{v}=(k-1)\frac{u}{v} - \lrceil{(k-1)\frac{u}{v}}
\in(-1,0],
\end{align}
which together with \eqref{qa6} yields \eqref{qa8} 
for $1\le j\le v$. 
Moreover,
\eqref{qa8} holds also for $j=0$ by
the special definition of $\ttt_0$.

Equations \eqref{m3-Pkm-sector}, \eqref{m3-nu-Pk}, 
\eqref{m3-lambda-mu}, and \eqref{qa4}
yield, for $1\le k\le v$, 
\begin{align}\label{qb1f}
  f(P_k)&=\nu(k)-k\ga = \nu(k)-k\frac uv+k\delta_- =-t_k+k\delta_-,
\\\label{qb1h}
 h(P_k)&=k\gb-\nu(k) = k\frac uv-\nu(k) +k\delta_+ = t_k+k\delta_+.
\end{align}
Hence, recalling the definitions \eqref{qa2} and \eqref{qa5},
for $1\le j\le v$,
\begin{align}\label{qb2f}
 f(W_j)&=-\widetilde t_j+\kappa_j\delta_-,
\\\label{qb2h}
 h(W_j)&= \widetilde t_j+\kappa_j\delta_+.
\end{align}
Furthermore, the special definitions of $W_0=X=P(1,1)$ and 
$\ttt_0=\frac uv-1$ above  
imply, using \eqref{m3-Pkm-sector} and \eqref{m3-lambda-mu},
that \eqref{qb2f}--\eqref{qb2h} hold also for $j=0$.

It follows from the definition \eqref{qa1} and $b+d=v$ that
$\kappa_{j+1}-\kappa_j$ equals either $b$ or $b-v=-d$.
Hence, \eqref{qb2f} and \eqref{qb2h} imply,
using \eqref{qa8} and \eqref{me-lambda-mu-b},
\begin{align}
 f(W_{j+1})-f(W_j)
 &=-\frac1v+(\kappa_{j+1}-\kappa_j)\delta_-<0,
 \label{m3-f-W-diff}\\
 h(W_{j+1})-h(W_j)
 &= \frac1v+(\kappa_{j+1}-\kappa_j)\delta_+>0.
 \label{m3-h-W-diff}
\end{align}
This proves \eqref{qa3f} and \eqref{qa3h}; 
furthermore, \eqref{qa3g} follows from these since $g=h-f$.
\end{proof}

\begin{lemma}\label{LQ4}
  Let $j_0:=v-u\in\set{2,\dots,v-1}$.
Then
\begin{align}\label{qb3}
  \kk_{j_0-1}=d,
\qquad 
\kk_{j_0}=v,
\qquad
  \kk_{j_0+1}=b.
\end{align}

Furthermore, for $0\le j\le v$, let $\jx\in\set{1,\dots,v}$ with
\begin{align}\label{qb4}
  \jx=j+u\pmod v.
\end{align}
Then, $\jox=v$.  
If $0\le j<v$ and $j\neq j_0$, then
\begin{align}
  \label{qbj0}
(j+1)^* = \jx+1.
\end{align}
If $0\le j\le v$ and $j\neq j_0$, so $\kk_j\neq v$, we have 
$\kk_{\jx}=\kk_j+1$, and thus
\begin{align}\label{qb5}
T(W_j)=T(P_{\kk_j})=P_{\kk_{j}+1}=  W_{\jx}.
\end{align}
\end{lemma}

\begin{proof}
Note that $j_0\ge2$, since $u/v<1/2$.
We have, by \eqref{qa1} and $bu-av=1$,
\begin{align}\label{qb6}
  \kk_{j_0}\equiv (v-u)b+1 \equiv -bu+1 \equiv 0 \pmod v,
\end{align}
and thus $\kk_{j_0}=v$.
Similarly,
\begin{align}\label{qb7-}
  \kk_{j_0-1}&\equiv \kk_{j_0}-b = v-b = d \pmod v,
\\\label{qb7+}
  \kk_{j_0+1}&\equiv \kk_{j_0}+b = v+b \equiv b \pmod v,
\end{align}
and thus $\kk_{j_0-1}=d$ and $\kk_{j_0+1}=b$, 
showing \eqref{qb3}.

It is immediate that $\jox=v$.
Furthermore, it follows from \eqref{qb4} that $(j+1)^*\equiv\jx+1\pmod v$,
and thus \eqref{qbj0} follows when $j\neq j_0$ and thus $\jx\neq v$.

Finally, \eqref{qb4} and \eqref{qa1} imply
\begin{align}
  \kk_{\jx}\equiv(j+u)b+1
\equiv \kk_j+ub
=\kk_j+av+1
\equiv\kk_j+1
\pmod v,
\end{align}
and thus $\kk_{\jx}=\kk_j+1$ when $\kk_j\neq v$.
If $1\le j< v$ and $j\neq j_0$, \eqref{qb5} now follows from
\eqref{qa2} and \refT{Tm3-Farey-entry}.

For $j=0$, we have $W_0=X$, whereas $P_{\kappa_0}=P_1=Z$.
Here $0<\alpha<\beta<1$ implies $X\in\cC\inter$ and $Z\in\cA\inter$,
by \eqref{qaA}, \eqref{qaC} and \eqref{m3-Pkm-sector}. 
Moreover, $0^*=u$ and $\kappa_u=2$, since $bu-av=1$.  Thus
\begin{align}\label{cec1}
 T(W_0)=T_3(X)=XZ=T_1(Z)=T(P_1)=P_2=W_{0^*},
\end{align}
which proves \eqref{qb5} when $j=0$.
Similarly, for $j=v$ we have $W_v=Z\in\cA\inter$, $\kk_v=\kk_0=1$, and
$v^*=0^*$, and \eqref{qb5} again follows from \eqref{cec1}.
\end{proof}

\begin{lemma}\label{LQ3}
We have 
\begin{align}\label{qbb5}
  \begin{cases}
    W_j\in\cC\inter, & j=0,\dots,j_0-1,
\\
    W_j\in\cB\inter, & j=j_0,
\\
    W_j\in\cA\inter, & j=j_0+1,\dots,v.
  \end{cases}
\end{align}
\end{lemma}
\begin{proof}
We have, by \eqref{qa8} and $j_0=v-u$, 
\begin{align}\label{qb8}
  \ttt_{j_0}=0
\qquad\text{and}\qquad
\ttt_{j_0\pm1}=\pm\frac{1}v.
\end{align}
Hence, \eqref{qb2f}--\eqref{qb2h} yield,
using \eqref{qb3} and \eqref{me-lambda-mu-b},
\begin{align}\label{qb9f}
  f(W_{j_0})&=0+v\delta_->0,
\\\label{qb9f+}
  f(W_{j_0+1})&=-\frac{1}{v}+b\delta_-<0,
\\\label{qb9h}
  h(W_{j_0})&=0+v\delta_+>0,
\\\label{qb9h-}
  h(W_{j_0-1})&=-\frac{1}v+d\delta_+<0.
\end{align}

\refL{Lm3-shadowing} now shows that 
for $0\le j\le j_0-1$, we have $f(W_j)>0$ and $h(W_j)<0$, and thus 
$g(W_j)=h(W_j)-f(W_j)<0$; hence $W_j\in\cC\inter$ by \eqref{m3-cC}.
Similarly, 
for $j_0+1\le j\le v$, we have $f(W_j)<0$ and $h(W_j)>0$, and thus 
$g(W_j)=h(W_j)-f(W_j)>0$; hence $W_j\in\cA\inter$ by \eqref{m3-cA}.
Finally, using \eqref{qa2} and \eqref{qb3}, $W_{j_0}=P_v\in\cB\inter$
by \refT{Tm3-Farey-entry}.
\end{proof}

Define, for $1\le j\le v$, the point $Q_j$
as the intersection of the $f$-line through $W_{j-1}$ and the $h$-line
through $W_{j}$: thus
\begin{align}\label{m3-Qk}
        f(Q_j)=f(W_{j-1}),\qquad h(Q_j)=h(W_{j}).
\end{align}
Here $Q_j$ denotes a point, whereas $Q_1,Q_2,Q_3$ in
\refSS{SSm3S12} denote auxiliary polygons.
Let also $B$ be the intersection of the $f$-line through $W_{j_0-1}$ and the
$h$-line through $P_0=U$. (The latter line is just  $\gD_h$). Thus
\begin{align}\label{jdb}
  f(B)=f(W_{j_0-1}), \qquad h(B)=0.
\end{align}
Similarly, let $A$ be the intersection of the 
$f$-line through $P_0=U$ (i.e.,  $\gD_f$) 
and the $h$-line through $W_{j_0+1}$. Thus
\begin{align}\label{jda}
  f(A)=0, \qquad h(A)=h(W_{j_0+1}) .
\end{align}
Note that we thus have
$f(W_{j_0-1})=f(B)=f(Q_{j_0})$
and, by \eqref{qb9h}--\eqref{qb9h-},
$h(W_{j_0-1})<h(B)<h(W_{j_0})=h(Q_{j_0})$.
Hence $B$ is an interior point in the line segment $W_{j_0-1}Q_{j_0}$.
Similarly,
$A$ is an interior point in the line segment $Q_{j_0+1}W_{j_0+1}$.

\begin{figure}[htbp]
\centering
\begin{tikzpicture}[x=6.4cm,y=6.4cm,font=\small]
\fill[blue!22] (0.0366667,0.0366667) -- (0.0800000,0.0150000) -- (0.1166667,0.0516667) -- (0.1600000,0.0300000) -- (0.1966667,0.0666667) -- (0.2400000,0.0450000) -- (0.2533333,0.0583333) -- (0.3600000,0.0050000) -- (0.8333333,-0.9416667) -- (0.8000000,-0.9750000) -- (0.7566667,-0.9533333) -- (0.7200000,-0.9900000) -- (0.6133333,-0.9366667) -- (0.6000000,-0.9500000) -- (0.5566667,-0.9283333) -- (0.5200000,-0.9650000) -- (0.4766667,-0.9433333) -- (0.4733333,-0.9466667) -- (0.0000000,0.0000000) -- cycle;
\fill[red!25] (-0.1066667,0.0533333) -- (-0.0833333,0.0766667) -- (-0.0400000,0.0550000) -- (-0.0266667,0.0683333) -- (0.0366667,0.0366667) -- (0.0000000,0.0000000) -- cycle;
\fill[green!35] (-0.6400000,0.0050000) -- (-0.6033333,0.0416667) -- (-0.5600000,0.0200000) -- (-0.5233333,0.0566667) -- (-0.4800000,0.0350000) -- (-0.4433333,0.0716667) -- (-0.4000000,0.0500000) -- (-0.3866667,0.0633333) -- (-0.2800000,0.0100000) -- (-0.2433333,0.0466667) -- (-0.2000000,0.0250000) -- (-0.1633333,0.0616667) -- (-0.1200000,0.0400000) -- (-0.1066667,0.0533333) -- (0.0000000,0.0000000) -- (0.4733333,-0.9466667) -- (0.4400000,-0.9800000) -- (0.3966667,-0.9583333) -- (0.3600000,-0.9950000) -- (0.2533333,-0.9416667) -- (0.2400000,-0.9550000) -- (0.1966667,-0.9333333) -- (0.1600000,-0.9700000) -- (0.1166667,-0.9483333) -- (0.0800000,-0.9850000) -- (-0.0266667,-0.9316667) -- (-0.0400000,-0.9450000) -- (-0.0833333,-0.9233333) -- (-0.1200000,-0.9600000) -- (-0.1633333,-0.9383333) -- (-0.1666667,-0.9416667) -- cycle;
\draw[black,line width=0.6pt] (-0.6400000,0.0050000) -- (-0.6033333,0.0416667) -- (-0.5600000,0.0200000) -- (-0.5233333,0.0566667) -- (-0.4800000,0.0350000) -- (-0.4433333,0.0716667) -- (-0.4000000,0.0500000) -- (-0.3866667,0.0633333) -- (-0.2800000,0.0100000) -- (-0.2433333,0.0466667) -- (-0.2000000,0.0250000) -- (-0.1633333,0.0616667) -- (-0.1200000,0.0400000) -- (-0.0833333,0.0766667) -- (-0.0400000,0.0550000) -- (-0.0266667,0.0683333) -- (0.0800000,0.0150000) -- (0.1166667,0.0516667) -- (0.1600000,0.0300000) -- (0.1966667,0.0666667) -- (0.2400000,0.0450000) -- (0.2533333,0.0583333) -- (0.3600000,0.0050000) -- (0.8333333,-0.9416667) -- (0.8000000,-0.9750000) -- (0.7566667,-0.9533333) -- (0.7200000,-0.9900000) -- (0.6133333,-0.9366667) -- (0.6000000,-0.9500000) -- (0.5566667,-0.9283333) -- (0.5200000,-0.9650000) -- (0.4766667,-0.9433333) -- (0.4400000,-0.9800000) -- (0.3966667,-0.9583333) -- (0.3600000,-0.9950000) -- (0.2533333,-0.9416667) -- (0.2400000,-0.9550000) -- (0.1966667,-0.9333333) -- (0.1600000,-0.9700000) -- (0.1166667,-0.9483333) -- (0.0800000,-0.9850000) -- (-0.0266667,-0.9316667) -- (-0.0400000,-0.9450000) -- (-0.0833333,-0.9233333) -- (-0.1200000,-0.9600000) -- (-0.1633333,-0.9383333) -- (-0.1666667,-0.9416667) -- cycle;
\draw[->,gray] (-0.6900000,0.0000000) -- (0.8833333,0.0000000) node[right] {$x$};
\draw[->,gray] (0.0000000,-1.0450000) -- (0.0000000,0.1766667) node[above] {$y$};
\draw[gray] (-0.5000000,-0.0050000) -- (-0.5000000,0.0050000) node[below=3pt,black,font=\scriptsize] {$-0.5$};
\draw[gray] (-0.2500000,-0.0050000) -- (-0.2500000,0.0050000) node[below=3pt,black,font=\scriptsize] {$-0.25$};
\draw[gray] (0.2500000,-0.0050000) -- (0.2500000,0.0050000) node[below=3pt,black,font=\scriptsize] {$0.25$};
\draw[gray] (0.5000000,-0.0050000) -- (0.5000000,0.0050000) node[below=3pt,black,font=\scriptsize] {$0.5$};
\draw[gray] (-0.0050000,-1.0000000) -- (0.0050000,-1.0000000) node[left=3pt,black,font=\scriptsize] {$-1$};
\draw[gray] (-0.0050000,-0.7500000) -- (0.0050000,-0.7500000) node[left=3pt,black,font=\scriptsize] {$-0.75$};
\draw[gray] (-0.0050000,-0.5000000) -- (0.0050000,-0.5000000) node[left=3pt,black,font=\scriptsize] {$-0.5$};
\draw[gray] (-0.0050000,-0.2500000) -- (0.0050000,-0.2500000) node[left=3pt,black,font=\scriptsize] {$-0.25$};
\begin{scope}
\clip (-0.6900000,-1.0450000) rectangle (0.8833333,0.1766667);
\draw[red!70!black,thin] (0,0) -- (2.0000000,2.0000000);
\draw[red!70!black,thin] (0,0) -- (-2.0000000,1.0000000);
\draw[red!70!black,thin] (0,0) -- (1.0000000,-2.0000000);
\end{scope}
\node at (0.4128333,-0.4605394) {$\cP_1$};
\node at (-0.1876047,-0.4593286) {$\cP_3$};
\node (middlelabel) at (0.0934667,0.1366667) {$\cP_2$};
\draw[thin] (middlelabel) -- (-0.0315333,0.0402333);
\node[below left,font=\scriptsize] at (0,0) {$U$};
\end{tikzpicture}
\caption{The invariant polygon 
for $(p,q)=(9/25,1/200)$
in $\cR_{4/11}$. It has $47$ vertices by \refT{Tm3-4v+3}.
The colours indicate the three pieces.}
\label{fig:m3-S411}
\end{figure}
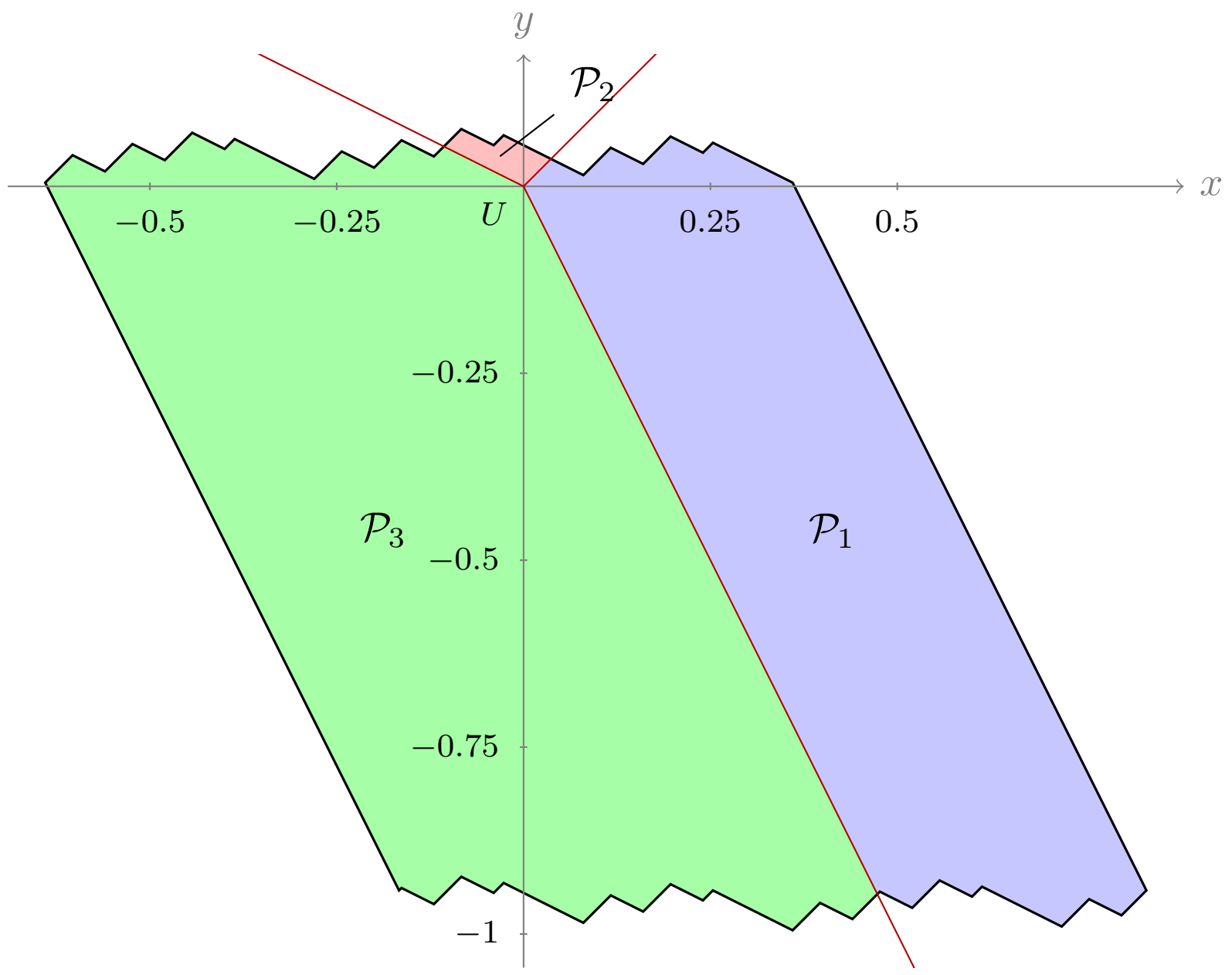

Now define the zig-zag line $\gG^+$
consisting of the segments
$W_{j-1}Q_j$  and 
$Q_jW_{j}$ ($1\le j\le v$), i.e.,
\begin{align}\label{maj34}
\gG^+:=\bigcup_{j=1}^v W_{j-1}Q_jW_{j}.
\end{align}
(This is shown below to be the upper boundary of $\cD$, see
Figures~\ref{fig:m3-S13} and~\ref{fig:m3-S411}.)
Define also
the zig-zag line $\gG^-$:
\begin{align}\label{maj35}
\gG^-:=\bigcup_{j=1}^{j_0-1} W_{j-1}Q_jW_{j}
\cup
W_{j_0-1}BUAW_{j_0+1}
\cup
\bigcup_{j=j_0+2}^{v} W_{j-1}Q_jW_{j},
\end{align}
where we thus have replaced $Q_{j_0}$, $W_{j_0}$, $Q_{j_0+1}$
by $B, U, A$. Hence,
the gap between $\gG^-$ and $\gG^+$ is 
the hexagon
\begin{align}\label{maj36}
  \cP_2:=[B,Q_{j_0},W_{j_0},Q_{j_0+1},A,U].
\end{align}

\begin{lemma}\label{LQ5}
 The zig-zag paths\/ $\gG^+$ and $\gG^-$  are  simple paths from $X$ to $Z$.
For either of the paths, $g$ is strictly increasing as we move
from $X$ to $Z$. 
Furthermore, on each line segment in the paths, either
$f$ is constant and $h$ is increasing (so the segment is on an $f$-line)
or
$h$ is constant and $f$ is decreasing (so the segment is on an $h$-line).

Similarly, in the $xy$-plane,
$x$
is strictly increasing  along $\gG^+$ and $\gG^-$, i.e., these paths are
graphs in the $xy$-plane of functions. 
Furthermore, since all line segments in them lie on $f$-lines and
$h$-lines, they have slope $1$ or $-\frac12$. 
Both $\gG^+$ and $\gG^-$ are contained in the upper halfplane
$\set{(x,y):y\ge0}$.
\end{lemma}
\begin{proof}
  This follows easily from the definitions and results above.
First, \eqref{m3-Qk} and \eqref{qa3h}
imply that on the segment $W_{j-1}Q_j$,
we have $f$ constant while $h$ increases from $h(W_{j-1})$ to $h(W_j)$.
Similarly, \eqref{m3-Qk} and \eqref{qa3f} imply that
on 
$Q_jW_{j}$,
we have $h$ constant while $f$ decreases from $f(W_{j-1})$ to $f(W_j)$.
In both cases, $g=h-f$ is (strictly) increasing.
Hence $g$ is strictly increasing on $\gG^+$, and in particular, $\gG^+$ is
simple (i.e., non-selfintersecting).

For the exceptional segments
in $\gG^-$, we note that \eqref{jdb} and \eqref{jda} imply that
$W_{j_0-1}B$ and $UA$ are on $f$-lines, while
$BU$ and $AW_{j_0+1}$ are on $h$-lines.
Furthermore, $h(W_{j_0-1})<0=h(B)$ by \eqref{qb9h-} and \eqref{jdb}.
The remaining inequalities $f(B)>f(U)=0$, $h(A)>h(U)=0$, 
and $f(W_{j_0+1})<f(A)=0$
follow similarly from \eqref{qb9f}--\eqref{qb9h},
together with
\eqref{qa3f} and \eqref{qa3h}.
Hence, the conclusions in the first paragraph hold for $\gG^-$ too.

Since $3x=h(x,y)-2f(x,y)$, it follows also that $x$ is strictly increasing
along $\gG^+$ and $\gG^-$. It follows from \eqref{m3-f} and \eqref{m3-h}
that the slope of an $f$-line is 1,
and the slope of an $h$-line is $-\frac12$.
Finally, 
it follows from \eqref{m3-nu-Pk} and \eqref{m3-Pkm} that
$\set{W_j:0\le j\le v}=\set X \cup\set{P_k:1\le k\le v}
\subset\set{y>0}$. 
And if we travel along $\gG^+$ from $X$ to $Z$, then the fact
on the slopes just noted shows that 
$y$ increases on each segment $W_{j-1}Q_j$ and then
decreases on $Q_jW_j$; hence $y>0$ everywhere on $\gG^+$.
The same argument holds for $\gG^-$, except that $y=0$ at $U$.
\end{proof}

We note from the proof that in $\gG^+$ and $\gG^-$,
each appearing $W_j$ is a valley, and so is $U\in\gG^-$,
while each $Q_j$ is a peak, and so are $B$ and $A$ in $\gG^-$.

Recall that $U$ is the triple point $\cA\cap\cB\cap\cC$ where all three
choices $T_1,T_2,T_3$ are possible for $T$.

\begin{lemma}\label{LQ6}
  If we exclude the choice $T_2$ for $T$, then
\begin{align}\label{jd3}
        T(\gG^-)=\gG^+ .
\end{align}
\end{lemma}

\begin{proof}
We treat the different terms in \eqref{maj35} separately.

\pfitemnr
Let $1\le j\le j_0-1$.
Then $W_{j-1},W_j\in\cC\inter$ by \eqref{qbb5}.
Furthermore,
it follows from \eqref{m3-Qk},
\refL{Lm3-shadowing}, and
\eqref{qb9f}--\eqref{qb9h-}, 
that
$f(Q_j)>0$ and $h(Q_j)<0$;
hence also $g(Q_j)=h(Q_j)-f(Q_j)<0$ and thus $Q_j\in\cC\inter$ by \eqref{m3-cC}.
Since $\cC\inter$ is convex, the path
$W_{j-1}Q_jW_j\subset\cC\inter$, and thus $T=T_3$ there. 
Furthermore, \eqref{qbj0} implies that $(j-1)^*=\jx-1$, and thus
\eqref{qb5} yields $T(W_{j-1})=W_{\jx-1}$  and $T(W_{j})=W_{\jx}$. 
Since $T_3$ is a
translation, and any translation maps $f$-lines to $f$-lines
and $h$-lines to $h$-lines,
it follows from the definition of $Q_j$ that  $T(Q_j)=T_3(Q_j)=Q_{\jx}$.
Consequently,
\begin{align}\label{qc1}
T(W_{j-1}Q_jW_j)=W_{\jx-1}Q_{\jx}W_\jx\subset\gG^+.
\end{align}

\pfitemnr
Similarly, if $j_0+2\le j\le v$, then
$W_{j-1}Q_jW_j\subset\cA\inter$, and \eqref{qc1} holds again.

\pfitemnr
We have, directly from \eqref{qb4}, $(j_0-1)^*=v-1$ and $(j_0+1)^*=1$.
Hence, $T(W_{j_0-1})=W_{v-1}$ and 
$T(W_{j_0+1})=W_{1}$ by \eqref{qb5}.
Furthermore, $T_3U=Z=W_v$ and $T_1U=X=W_0$.
We now argue as in (i), using $U\in\cA\cap\cC$.
It follows that
$W_{j_0-1}BU\subset \cC$ and $UAW_{j_0+1}\subset\cA$.
Furthermore, since we have excluded $T_2$,
 and with the ``right'' choice of $T$ at $U$,
\begin{align}
T(W_{j_0-1}BU) & = T_3(W_{j_0-1}BU) =W_{v-1}Q_vW_v,\label{qc2a}
\\
T(UAW_{j_0+1})&=T_1(UAW_{j_0+1})=W_0Q_1W_1.\label{qc2b}
\end{align}

Recalling \eqref{maj34}--\eqref{maj35}, 
the images in
\eqref{qc1}--\eqref{qc2b}
cover all the teeth of $\gG^+$, and hence prove \eqref{jd3}.
\end{proof}

To summarize the construction so far, we have constructed zig-zag paths
$\gG^+$ and $\gG^-$ from $X$ to $Z$, with gap between them given by 
\eqref{maj36}.
$\gG^+$ will be the top part of the boundary of $\cD$. 
It remains to add sides and a bottom part.

The sides are just the two
$g$-lines through $X$ and $Z$,
given by
\begin{align}
L_X:2x+y&=  g(x,y)=g(X)=2p+q-2,\label{maj37}
\\
L_Z:2x+y&=  g(x,y)=g(Z)=2p+q.\label{maj38}
\end{align}
Let $\cS$ be the strip between these two lines.
(Thus the horizontal width of $\cS$ is exactly~1.)

To construct the lower part of the boundary of $\cD$,
note that the endpoints $X$ and $Z$ of $\gG^-$ differ by $\ve_1=(1,0)$.
Extend $\gG^-$ to a periodic infinite zig-zag line
$\tgG^-:=\bigcup_{\ell\in\bbZ}(\gG^-+\ell\ve_1)$.
Translate this curve by $Y$ to the periodic infinite zig-zag line
$\tgG_-:=\tgG^-Y$, and let $\gG_-$ be the part of $\tgG_-$ 
that lies in the strip $\cS$.
Do the same with $\gG^+$ and call the resulting zig-zag lines 
$\tgG^+$, $\tgG_+$, and $\gG_+$.  
We can now complete our construction, which will be verified below:
\begin{align}\label{cD}
\cD  
\text{\emph{ is the polygon bounded by $\gG^+$, $\gG_-$, and (parts of) 
the lines $L_X$ and $L_Z$}.}
\end{align}

For examples, see Figures \ref{fig:m3-S13} and \ref{fig:m3-S411}.

\begin{figure}[htbp]
\centering
\begin{tikzpicture}[x=6.4cm,y=6.4cm,font=\small]
\fill[blue!22] (0.1200000,0.1200000) -- (0.2000000,0.0800000) -- (0.6133333,-0.7466667) -- (0.6000000,-0.7600000) -- (0.5200000,-0.7200000) -- (0.4133333,-0.8266667) -- (0.0000000,0.0000000) -- cycle;
\fill[red!25] (-0.5066667,0.2533333) -- (-0.4800000,0.2800000) -- (-0.4000000,0.2400000) -- (-0.3066667,0.3333333) -- (0.1200000,0.1200000) -- (0.0000000,0.0000000) -- cycle;
\fill[green!35] (-0.8000000,0.0800000) -- (-0.6800000,0.2000000) -- (-0.6000000,0.1600000) -- (-0.5066667,0.2533333) -- (0.0000000,0.0000000) -- (0.4133333,-0.8266667) -- (0.4000000,-0.8400000) -- (0.3200000,-0.8000000) -- (0.2000000,-0.9200000) -- (-0.3066667,-0.6666667) -- (-0.3866667,-0.7466667) -- cycle;
\draw[black,line width=0.6pt] (-0.8000000,0.0800000) -- (-0.6800000,0.2000000) -- (-0.6000000,0.1600000) -- (-0.4800000,0.2800000) -- (-0.4000000,0.2400000) -- (-0.3066667,0.3333333) -- (0.2000000,0.0800000) -- (0.6133333,-0.7466667) -- (0.6000000,-0.7600000) -- (0.5200000,-0.7200000) -- (0.4000000,-0.8400000) -- (0.3200000,-0.8000000) -- (0.2000000,-0.9200000) -- (-0.3066667,-0.6666667) -- (-0.3866667,-0.7466667) -- cycle;
\draw[->,gray] (-0.8500000,0.0000000) -- (0.6633333,0.0000000) node[right] {$x$};
\draw[->,gray] (0.0000000,-0.9700000) -- (0.0000000,0.4333333) node[above] {$y$};
\draw[gray] (-0.7500000,-0.0050000) -- (-0.7500000,0.0050000) node[below=3pt,black,font=\scriptsize] {$-0.75$};
\draw[gray] (-0.5000000,-0.0050000) -- (-0.5000000,0.0050000) node[below=3pt,black,font=\scriptsize] {$-0.5$};
\draw[gray] (-0.2500000,-0.0050000) -- (-0.2500000,0.0050000) node[below=3pt,black,font=\scriptsize] {$-0.25$};
\draw[gray] (0.2500000,-0.0050000) -- (0.2500000,0.0050000) node[below=3pt,black,font=\scriptsize] {$0.25$};
\draw[gray] (0.5000000,-0.0050000) -- (0.5000000,0.0050000) node[below=3pt,black,font=\scriptsize] {$0.5$};
\draw[gray] (-0.0050000,-0.7500000) -- (0.0050000,-0.7500000) node[left=3pt,black,font=\scriptsize] {$-0.75$};
\draw[gray] (-0.0050000,-0.5000000) -- (0.0050000,-0.5000000) node[left=3pt,black,font=\scriptsize] {$-0.5$};
\draw[gray] (-0.0050000,-0.2500000) -- (0.0050000,-0.2500000) node[left=3pt,black,font=\scriptsize] {$-0.25$};
\draw[gray] (-0.0050000,0.2500000) -- (0.0050000,0.2500000) node[left=3pt,black,font=\scriptsize] {$0.25$};
\begin{scope}
\clip (-0.8500000,-0.9700000) rectangle (0.6633333,0.4333333);
\draw[red!70!black,thin] (0,0) -- (2.0000000,2.0000000);
\draw[red!70!black,thin] (0,0) -- (-2.0000000,1.0000000);
\draw[red!70!black,thin] (0,0) -- (1.0000000,-2.0000000);
\end{scope}
\node at (0.2938667,-0.3464533) {$\cP_1$};
\node at (-0.3106667,-0.3256000) {$\cP_3$};
\node (middlelabel) at (-0.0403333,0.3933333) {$\cP_2$};
\draw[thin] (middlelabel) -- (-0.1653333,0.1698667);
\node[below left,font=\scriptsize] at (0,0) {$U$};
\end{tikzpicture}
\caption{The invariant polygon 
for $(p,q)=(1/5,2/25)$ in $S_{1/3}\SC\subset\cR_{1/3}$
(see \eqref{m3-S13SC}).
The colours indicate the pieces
$\cP_1$, $\cP_2$, and $\cP_3$.}
\label{fig:m3-S13}
\end{figure}
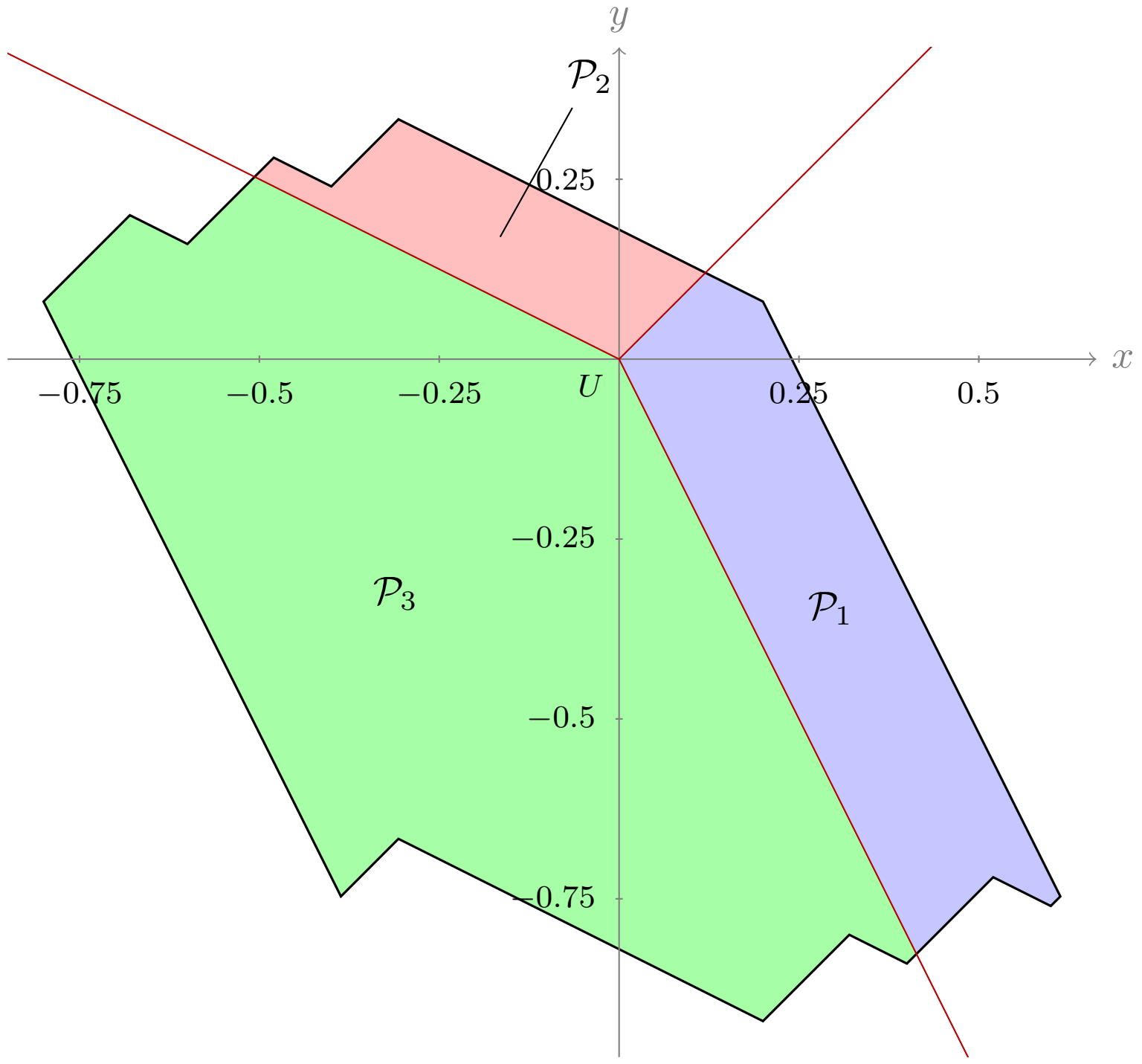

\begin{theorem}\label{TD}
Let $(p,q)$ belong to the open Farey cell $\cR_{u/v}$, and let $\cD$ be the
polygon constructed in \eqref{cD}.  Then $\cD$ is a compact polygonal domain,
and
\begin{align}
 T(\cD)=\cD.
\end{align}
Furthermore,  $\cD\cap\cB=\cP_2$ defined in \eqref{maj36},
and if\/ $\cP_1:=\cD\cap\cA$ and $\cP_3:=\cD\cap\cC$, then the three images
$T_1\cP_1,T_2\cP_2,T_3\cP_3$ have pairwise disjoint interiors and their union
is $\cD$.
If\/ $p,q,1$ are rationally independent, then $\cD$ is the global attractor
$\Kx$ and equals $\cL(\vu)$ for every initial point $\vu\in\Ho$.
\end{theorem}

Before giving the proof, we first give some lemmas.
\begin{lemma}\label{LQ8}
We have
\begin{align}\label{lq8a}
  \tgG_+,\tgG_-\subset\set{(x,y):y<0}.
\end{align}
In particular,
\begin{align}\label{lq8b}
 \tgG_+\cap\cB=\tgG_-\cap\cB=\emptyset. 
\end{align}
It follows that
\begin{align}\label{qc3}
  T(\tgG_-)\subseteq\tgG_+.
\end{align}
\end{lemma}
\begin{proof}
  We have $f+h=3y$ by \eqref{m3-f} and \eqref{m3-h}.
Hence, on both $\gG^+$ and $\gG^-$ we have, using \refL{LQ5},
\begin{align}\label{ql1}
  3y \le f(W_0)+h(W_v)=f(X)+h(Z)=(1-p+q)+(p+2q)=1+3q.
\end{align}
The same holds on the infinite extensions $\tgG^+$ and $\tgG^-$
obtained by horizontal translations. Translating by $Y=(p,q-1)$, we see that
on $\tgG_+$ and $\tgG_-$, we have
\begin{align}\label{ql2}
  y\le \frac{1+3q}3 + q-1 = 2q -\frac{2}{3} <0,
\end{align}
since \eqref{m3-parameter-region} implies $q<\frac13$. 
This shows \eqref{lq8a}.

On the other hand, in $\cB$ we have, by the definition \eqref{m3-cB},
$3y=f+h\ge0$. Hence, \eqref{lq8b} follows.

The translations $T_1$ and $T_3$ differ by the horizontal period $\ve_1$.
Thus 
taking periodic extensions in
\eqref{jd3} gives
\begin{align}\label{qc2}
  T_1(\tgG^-)=T_3(\tgG^-)=\tgG^+.
\end{align}
By \eqref{lq8b}, $T_2$ is not a choice on $\tgG_-$.
Hence we obtain \eqref{qc3} from \eqref{qc2} by translating by $Y$
(since different translations commute).
\end{proof}

\begin{lemma}\label{LQ9}
  $\gG^+$ and $\gG_-$ are disjoint, and together with parts of 
the $g$-lines $L_X$ and  $L_Z$ they form a polygon $\cD$.
(Thus, the definition in \eqref{cD} makes sense.)

Furthermore,
if $L$ is any $g$-line that intersects $\cD$ at all, then
$L\cap\cD$ is the segment $P^+P_-$ where 
$P^+:=L\cap\gG^+$ and
$P_-:=L\cap\tgG_-$. 
\end{lemma}
\begin{proof}
We have $\gG^+\subset\set{(x,y):y\ge0}$ by \refL{LQ5} 
and $\gG_-\subset\tgG_-\subset\set{(x,y):y<0}$ by \refL{LQ8};
consequently, $\gG^+$ and $\gG_-$ are disjoint.
Furthermore, 
$g$ is strictly increasing along both these paths 
by \refL{Lm3-shadowing}.
Hence, 
if we temporarily switch to the coordinates $(g(x,y),y)$, then
$\gG^+$ and $\gG_-$
are disjoint graphs of 
piecewise linear
functions over the same interval $[g(X),g(Z)]$.
Together with the joining segments on $L_X$ and $L_Z$, they therefore
bound a compact polygonal domain. 

The final claim clearly holds when $L=L_X$ or $L=L_Z$.
Let $L$ be any other $g$-line; then $L$ is parallel to the sides $L_X$ and
$L_Z$ and does not intersect them.
Furthermore, $g$ is strictly increasing on $\gG^+$ and $\gG^-$ by \refL{LQ5},
and as a consequence  also on $\tgG_-$. Hence 
the $g$-line $L$ can intersect $\gG^+$ and $\tgG_-$ in at most one point each,
and these are the only possible intersections of $L$ and
the boundary $\partial\cD$.
Since $\cD$ is bounded, it follows that if $L\cap\partial\cD\neq\emptyset$,
then  $L\cap\cD=P^+P_-$ with $P^+$ and $P_-$ as stated.
\end{proof}

\begin{lemma}\label{LQ7}
  We have 
  \begin{align}\label{lq70}
\cD\cap\cB = \cP_2.
  \end{align}
Furthermore, 
  \begin{align}\label{lq7}
T_2(\cP_2)=\cP_2Y\subset\cS\inter.    
  \end{align}
\end{lemma}
\begin{proof}
It follows by \eqref{jdb}, \eqref{jda}, \eqref{qb9f}, \eqref{qb9h}, and 
\refL{Lm3-shadowing} that $f(B)>0$, $h(B)=0$, $f(A)=0$, $h(A)>0$.
Hence, the definition \eqref{m3-cB} implies that
$A,B\in\partial\cB$, and it follows from the monotonicity in \refL{LQ5}
that
\begin{align}\label{lq77}
  \gG^+\cap\cB=B Q_{j_0} W_{j_0} Q_{j_0+1} A
\qquad\text{and}\qquad
  \gG^-\cap\cB=B U A.
\end{align}
Note that $BUA \subset \gD_h\cup\gD_f=\partial\cB$.

Let $L$ be a $g$-line that intersects $\cD$, and let $L\cap\cD=P^+P_-$ with
$P^+:=L\cap\gG^+$ and
$P_-:=L\cap\tgG_-$ as in \refL{LQ9}.
Both $f$ and $h$ are strictly decreasing if we move along $L$
downwards from $P^+$ to $P_-$. 
Hence, \eqref{m3-cB} shows that if $P^+\notin\cB$, then
$P^+P_-\cap\cB=\emptyset$. Furthermore, if $P^+\in\cB$, 
and thus $P^+\in\gG^+\cap\cB$, 
then $P^+P_-\cap\cB=P^+P^-$ where $P^-=L\cap\partial\cB$
and thus $P^-=L\cap BUA$. Hence, 
\begin{align}\label{qd1}
L\cap\cD\cap\cB=P^+P^-=L\cap\cP_2.  
\end{align}
This holds for every $g$-line $L$, so taking the union over all $L$ yields
$\cD\cap \cB=\cP_2$.

For \eqref{lq7}, we note that
the linear function $g$ attains its maximum and minimum on $\cP_2$ at 
some points
on the boundary, which consists of the parts of $\gG^+$ and $\gG^-$ between
$B$ and $A$. Since $g$ is increasing on both $\gG^+$ and $\gG^-$ by
\refL{LQ5}, the minimum is $g(B)$ and the maximum is $g(A)$.
Since $g=h-f$, we have, 
using  \eqref{jda}, \eqref{qb2h}, \eqref{qb8}, \eqref{qb3},
\eqref{m3-lambda-mu},
and \eqref{m3-open-farey},
\begin{align}\label{qc4}
  g(A) &
=h(A) -f(A)
=h(W_{j_0+1})
=\ttt_{j_0+1}+\kk_{j_0+1}\delta_+
=\frac1v+b\delta_+
\notag\\&
=\frac{bu-av}v+b\Bigpar{\gb-\frac uv}
=b\gb-a
<b\frac cd - a
=\frac{bc-ad}{d}=\frac{1}d
\le 1,
\end{align}
and similarly, using \eqref{jdb}, \eqref{qb2f}, \eqref{qb8}, \eqref{qb3},
\begin{align}\label{qc5}
  g(B) 
=d\ga-c
>d\frac{a}{b}-c
=\frac{ad-bc}b=-\frac{1}{b}\ge-1.
\end{align}
This gives
\begin{align}\label{qc6}
g(T_2A)=g(AY)=g(A)+g(Y)<g(Y)+1=g(Z)
\end{align}
and
\begin{align}\label{qc7}
g(T_2B)=g(BY)=g(B)+g(Y)>g(Y)-1=g(X).
\end{align}
Hence, for any $P\in\cP_2$ we have  $g(X)<g(T_2P)<g(Z)$, and thus
$T_2P\in\cS\inter$.
\end{proof}

\begin{proof}[Proof of \refT{TD}]
Since $\cP_2$ is the region between $\gG^+$ and $\gG^-$, 
$T_2(\partial\cP_2)$ is contained in 
$T_2(\tgG^+\cup\tgG^-)=
\tgG_+\cup\tgG_-$.  
By \refL{LQ7}, the
translated hexagon $T_2(\cP_2)=\cP_2Y$ lies strictly inside $\cS$.  
Hence $\cP_2Y$ is
exactly the gap between $\gG_-$ and $\gG_+$.  In particular,
$\gG_-$ and $\gG_+$ have common endpoints on $L_X$ and $L_Z$ and coincide
close to the endpoints.
Together with \refL{LQ5} and \refL{LQ8}, this also shows that the four paths
in \eqref{cD} form the boundary of a compact polygonal domain.

Let $\cD_-$ be the region bounded by $\gG^-$, $L_X$, $L_Z$, and $\gG_-$,
i.e., the closed part of $\cD$ below $\gG^-$. 
Similarly, let $\cD_+$ be the closed subset of $\cD$ 
bounded by $\gG^+$, $L_X$, $L_Z$, and $\gG_+$. 
It follows that
\begin{align}\label{m3-global-two}
 \cD=\cD_-\cup\cP_2=\cD_+\cup\cP_2Y.
\end{align}

Let again $L$ be a $g$-line that intersects $\cD$,
and let 
$P^\pm=L\cap\gG^\pm$ and
$P_\pm=L\cap\gG_\pm$.
Then \refL{LQ9} shows that $L\cap\cD=P^+P_-$.
The same argument shows that 
$L\cap\cD_-=P^-P_-$
and
$L\cap\cD_+=P^+P_+$.

It follows from \eqref{qd1} that $P^-P_-\cap\cB\inter=\emptyset$.  Since
$g$ is constant on $L$, this segment lies in $\cA$ when $g(L)>0$ and in
$\cC$ when $g(L)<0$; for $g(L)=0$ both branches are retained.  Its endpoints
$P^-$ and $P_-$
are mapped by $T$, by \eqref{jd3} and \eqref{qc3}, to the corresponding
endpoints of the intersection of $\cD_+$ and the $g$-line $T(L)$.

We check that all sections of $\cD_+$ and a $g$-line occur.  
We have $g(X)<0<g(Z)$.  
The intersections of $\cD_-$ and a $g$-line $L$
with $0\le g(L)\le g(Z)$  are sent by
$T_1$ to $g$-lines $L':=T(L)$
with  $g(L')\in[g(X),g(X)+g(Z)]$.  
Similarly, the
intersections of $\cD_-$ and a $g$-line $L$
with $g(X)\le g(L)\le0$ are sent by $T_3$ to $g$-lines
$L'$
with  $g(L')\in[g(X)+g(Z),g(Z)]$.  
Together these intervals cover
$[g(X),g(Z)]$.  
Consequently, taking the union of the corresponding intersections gives
\begin{align}\label{qd2}
 T_1(\cD_-\cap\cA)\cup T_3(\cD_-\cap\cC)=\cD_+.
\end{align}
Finally, \eqref{m3-global-two}, \eqref{lq70}, \eqref{lq7}, and \eqref{qd2} yield
\begin{align}
  T(\cD)
&=T_1(\cD_-\cap\cA)\cup T_3(\cD_-\cap\cC)\cup T_2(\cP_2)
\notag\\&
=\cD_+\cup\cP_2Y=\cD.
\end{align}
If $\cP_1:=\cD\cap\cA$ and $\cP_3:=\cD\cap\cC$, then
$\cD=\cP_1\cup\cP_2\cup\cP_3$, and  the interiors of $\cP_1,\cP_2,\cP_3$
are pairwise disjoint.
Since translations preserve area, the equality just proved also shows that
the three images $T_1\cP_1,T_2\cP_2,T_3\cP_3$ have pairwise disjoint interiors
and fill $\cD$.
The last assertion follows from Proposition~\ref{Pm3bridge}.
\end{proof}

\begin{proof}[Proof of the boundary assertion in \refR{Rrand}]
Fix a Farey cell $\cR_{u/v}$  
and a parameter pair $(p,q)$ on the boundary $\partial\cR_{u/v}$.
We will show that the construction above of $\cD$ still works and defines 
an invariant polygon $\cD$, although some vertices in the construction
might coincide.
(This includes $\cR_{1/2}=S_{1/2}$, for which the same construction applies
as noted at the beginning of this subsection.)
Fix also some parameters $(p_n,q_n)$ in the open Farey cell such that
$(p_n,q_n)\to(p,q)$ as \ntoo, and note that the paths $\gG^\pm$ and $\gG_\pm$ 
for $(p,q)$ are limits of the corresponding paths for $(p_n,q_n)$; this is
perhaps simplest seen in the coordinates $(g,y)$ where 
these paths are piecewise linear over $[g(X),g(Z)]$ with slopes $\pm\frac13$
(because $f=(3y-g)/2$ and $h=(3y+g)/2$).
(Some segments may get length zero in the limit, making other segments
coalesce.)
It follows by continuity that
the inequalities in the proof of \refL{LQ5} remain valid in non-strict form;
moreover, $\gG^+\subseteq\{y\ge0\}$ and
$\gG^-\subseteq\cA\cup\cC$, while \eqref{ql2} (in non-strict form) yields
\begin{align}\label{ysup}
 y\le 2q-\tfrac23\le 0 \qquad\text{on }\gG_-.
\end{align}
In fact, we have strict inequality in \eqref{ysup}, since 
equality would require $q=\frac13$ and thus 
$(p,q,r)=(\frac13 ,\frac13,\frac13 )$ by 
\eqref{m3-parameter-region} (in non-strict form); 
in this case, as noted in \refR{REX}, $\cD$ reduces to the hexagon in
\eqref{cH} and $\gG_-=H_3YH_1$ (see \refFig{fig:m3-hexagon}), for which
we have $y<0$ as a consequence of \eqref{hi1} and \eqref{hi3}.
Thus the  upper and lower paths  $\gG^+$ and $\gG_-$
remain separated also for the limit $(p,q)$ and bound 
together with the two sides of $\cS$,
a compact polygonal domain $\cD$.
Let $\cD_n$ be the  polygons constructed for $(p_n,q_n)$; then 
$\cD_n$ converge to $\cD$ in the Hausdorff metric.

Let $\Tn$ denote the operator $T$ for $(p_n,q_n)$, with branches $T_k^{(n)}$.
Then it follows by continuity from the invariance $\Tn(\cD_n)=\cD_n$,
looking at the three cases separately,
\begin{align}\label{tor}
 T_1(\cD\cap\cA\inter)\subseteq\cD,\qquad
 T_2(\cD\cap\cB\inter)\subseteq\cD,\qquad
 T_3(\cD\cap\cC\inter)\subseteq\cD.
\end{align}
Change coordinates again to the $(g,y)$-plane. Then the construction gives
\begin{align}
\cD=\bigl\{(s,y):g(X)\le s\le g(Z),\,
 \gamma_-(s)\le y\le\gamma^+(s)\bigr\},
\end{align}
where $\gamma^\pm$ and $\gamma_\pm$ denote the functions with graphs
$\gG^\pm$ and $\gG_\pm$.
The construction of $\gG^\pm$ in \eqref{maj34}--\eqref{maj35}
gives, see also \eqref{lq77} and note that $\gD_h$ and $\gD_f$ have slopes
$\pm\frac13$ in the $(g,y)$-plane,
\begin{align}\label{vidar}
\gamma^-(s)=\min\{\gamma^+(s),|s|/3\}\ge0,
\qquad g(X)\le s\le g(Z),
\end{align}
and this remains valid at the boundary parameter by continuity.
Since $g(X)<0<g(Z)$ and $\gamma_-<0\le\gamma^-$, the winner-sector
conditions \eqref{m3-cA} and \eqref{m3-cC} give
\begin{align}
\cD\cap\cA
&=\bigl\{(s,y):0\le s\le g(Z),\,
\gamma_-(s)\le y\le\gamma^-(s)\bigr\},\\
\cD\cap\cC
&=\bigl\{(s,y):g(X)\le s\le0,\,
        \gamma_-(s)\le y\le\gamma^-(s)\bigr\}.
\end{align}
Both sets are closures of their interiors. Hence the inclusions \eqref{tor} 
already
proved for $\cD\cap\cA\inter$ and $\cD\cap\cC\inter$, together with
continuity of the translations, yield
\begin{align}
T_1(\cD\cap\cA)\subseteq\cD,
\qquad
T_3(\cD\cap\cC)\subseteq\cD.
\end{align}
For the middle piece, write $\cP_{2,\infty}$ for the Hausdorff limit of
$\cP_{2,n}=\cD_n\cap\cB$.
This limit exists also when it has empty interior: the six designated
vertices converge, and in the $(g,y)$-plane the middle hexagon is the
region between the two graphs on the interval from $g(B)$ to $g(A)$.
Thus
\begin{align}
\cP_{2,\infty}
=\bigl\{(s,y):g(B)\le s\le g(A),\,
         \gamma^-(s)\le y\le\gamma^+(s)\bigr\}.
\end{align}
For every $P\in\cP_{2,\infty}$, choose $P_n\in\cP_{2,n}$ with $P_n\to P$.
Since $T_2^{(n)}P_n\in\cD_n$ and $\cD_n\to\cD$, continuity gives
\begin{align}
T_2(\cP_{2,\infty})\subseteq\cD.
\end{align}
At the boundary parameter, $\cD\cap\cB$ may contain additional segments
on the tie-lines. To check these, recall that $\gG^-$ follows $BUA$ on
the interval $[g(B),g(A)]$, where its graph is $y=|g|/3$, and coincides
with $\gG^+$ outside this interval. Together with
$\cB=\{(s,y):y\ge|s|/3\}$, 
which follows from \eqref{m3-f}--\eqref{m3-h} and \eqref{m3-cB},
this shows that
\begin{align}
\cP_{2,\infty}\subseteq\cD\cap\cB,
\qquad
(\cD\cap\cB)\setminus\cP_{2,\infty}\subseteq\gG^-.
\end{align}
Let $P$ belong to the latter set.
The monotonicity along $\gG^-$, in non-strict form
since we are at a boundary parameter, gives
\begin{align}
f(P)\le f(X)=1-\alpha\le1,
\qquad
h(P)\le h(Z)=\beta\le1,
\end{align}
where $\alpha=p-q\ge0$ and $\beta=p+2q\le1$ on the closed parameter
triangle. Since $P\in\cB$, we also have $f(P),h(P)\ge0$.
Consequently,
\begin{align}
-1\le g(P)=h(P)-f(P)\le1.
\end{align}
Moreover, $g(Y)-1=g(X)$ and $g(Y)+1=g(Z)$, so
\begin{align}
g(X)\le g(PY)=g(P)+g(Y)\le g(Z).
\end{align}
Thus $PY$ belongs to the strip $\cS$.
On the other hand, $P\in\gG^-$ implies $PY\in\tgG_-$ by the definition
of the periodic lower path. Hence
\begin{align}
T_2P=PY\in\tgG_-\cap\cS=\gG_-\subseteq\cD.
\end{align}
This checks $T_2$ on all of $\cD\cap\cB$, including any additional
tied segments. We have therefore proved
\begin{align}\label{tyr}
T_1(\cD\cap\cA)\subseteq\cD,\qquad
T_2(\cD\cap\cB)\subseteq\cD,\qquad
T_3(\cD\cap\cC)\subseteq\cD,
\end{align}
and thus, for the multi-valued map, $T(\cD)\subseteq\cD$ also in the
boundary cases.

Conversely, every point $\vv$ of $\cD$ is a limit of a sequence of points
$\vv_n\in\cD_n$. Since $\cD_n$ is invariant for  $\Tn$
we may choose preimages $\vu_n\in \cD_n$ with $\Tn \vu_n=\vv_n$
for suitable branches.
By passing  to a  subsequence, we may assume that $\vu_n\to \vu$ for some
$\vu$, and also that $T_k^{(n)}(\vu_n)=\vv_n$ for some fixed 
branch index $k\in\set{1,2,3}$. 
Then $\vu\in\cD$, since $\cD_n\to\cD$, and $T_k(\vu)=\vv$.
Furthermore, the winner sectors are closed, 
so $k$ is a (possibly tied) winner at $\vu$ and thus
$\vv=T_k(\vu)\in T(\vu)$.
Thus $\cD\subseteq T(\cD)$, proving invariance of $\cD$ 
on the closed Farey cell.

Note that this includes the vertex $(u/v,0)$ of $\cR_{u/v}$ where $q=0$; 
in this simple case 
(which we usually exclude by our assumption \eqref{m3-parameter-region}),
the upper and lower boundaries of $\cD$ become periodic zig-zag lines, see
also the proof of \refT{T-first-wall},
but
it is easy to see that there are many other invariant polygons.
\end{proof}

\begin{corollary}\label{Careas}
  Let $\cD$ and $\cP_1,\cP_2,\cP_3$ be as in \refT{TD}.
Then these polygons have areas $\gl(\cD)=1$, $\gl(\cP_1)=p$, $\gl(\cP_2)=q$,
and $\gl(\cP_3)=r=1-p-q$.
\end{corollary}
\begin{proof}
  The ergodic case follows by \refC{Cgl}. The upper and lower paths 
$\gG^+$ and $\gG_-$
and their intersections with the strip $\cS$ 
vary continuously throughout each open Farey
cell, also when the cut point on $\tgG_-$  
passes through a vertex. Hence the areas of $\cD$
and its intersections with the fixed winner sectors are continuous there.
Approximating any parameter in the open cell by ergodic parameters proves
the general case, including parameters on internal walls.
\end{proof}

\begin{remark}\label{Rcon}
 We may summarize the construction above as the following  algorithm.  
Given
$(p,q)\in\cR_{u/v}$, compute $\nu(k)$ from \eqref{m3-nu-Pk}, order the points
$P_k$ by the permutation $\kappa_j$ in \eqref{qa1}, and form the vertices
$Q_j,A,B$.  This gives $\gG^+$ and $\gG^-$.  Extend $\gG^-$ periodically,
translate it by $Y$, and cut it by the strip $\cS$.  The upper path, the two
side segments and the resulting lower path, read in cyclic order, give the
complete vertex list of $\cD$.
\end{remark}

\begin{remark}\label{RC}
Define the intersection point 
\begin{align}
  \label{m3-C}
C:=\gG_-\cap\gD_g.
\end{align}
$C$ is thus the point on $\gG_-$ where $g=0$; this point
is unique since $g$ is strictly increasing along $\gG_-$,
see the proof of \refL{LQ9}. 
At this tie point the two values of $T(C)$ are $T_1(C)=CX$ and $T_3(C)=CZ$.
These points lie
on the two boundary lines \eqref{maj37} and \eqref{maj38} of $\cS$,
and both points lie on $\tgG_+$ by \eqref{qc3};
hence the two (common) endpoints of $\gG_-$ and $\gG_+$ are $CX$ and $CZ$.
These points are thus the leftmost and rightmost vertices of $\cD$
on the lower part of the boundary.
\end{remark}

In this subsection, we have fixed $(p,q)\in\cR_{u/v}$.
In the next subsection we compare
the polygons $\cD$ constructed for different $(p,q)$
in the same Farey cell $\cR_{u/v}$.

\subsection{Stability domains and walls}\label{SSSD}

We next show that
the construction above also describes the changes of combinatorial type inside
a Farey cell.
Fix again a reduced fraction $u/v\in(0,\frac12]$,
but now let $(p,q)$ vary in $\cR_{u/v}$.
We continue to use the notation introduced above;
in particular,  the Farey parents of $u/v$ are $a/b$ and $c/d$.
Note that most variables defined above depend on $(p,q)$, but not the
functions $f,g,h$ and the winner domains $\cA,\cB,\cC$, 
see \eqref{m3-f}--\eqref{m3-cC}.

\begin{theorem}\label{TSD}
Fix $u/v\in(0,\frac12]$, and let $(p,q)$ vary in $\cR_{u/v}$.
Then the points $P_k$, $W_j$, $Q_j$, $A$, and $B$ in the construction above
have coordinates given by affine functions of $(p,q)$.
Hence, the upper paths $\gG^\pm$ and the periodic lower paths $\tgG_\pm$
have the same symbolic description for all $(p,q)\in\cR_{u/v}$.
On the other hand, $\tgG_-$ may be cut in different ways by the strip $\cS$
to produce the lower boundary $\gG_-$; this is determined by which line
segment of\/ $\tgG_-$ the line $\gD_g$ intersects, i.e.,
in which line segment of\/ $\tgG_-$ 
the intersection point $C$ in \eqref{m3-C} lies.

Consequently, for each fixed (in a symbolic description) 
line segment $\gU$ of $\tgG_-$, 
the set of all $(p,q)\in\cR_{u/v}$ (if any)
such that $\gD_g$ intersects the
interior $\gU\inter$ forms an open stability domain in $\cR_{u/v}$,
where
the symbolic description of 
the cyclic vertex list of $\cD$, including the winner sector of
each boundary piece, is constant.  

The remaining, exceptional, case is when 
$\gD_g$ goes through a vertex $V$ of\/ $\tgG_-$. In this case, $(p,q)$
is on the boundary between two stability domains, so small changes of
$(p,q)$ may change the combinatorial type.
For each vertex $V\in\tgG_-$, the set $\LL_V$
of all such $(p,q)\in\cR_{u/v}$
is either empty or given by an affine equation, i.e., it is the intersection
of $\cR_{u/v}$ with some line in the plane; 
if $\LL_V$ is non-empty we call it a \emph{wall}.
There is only a finite number of walls in each Farey cell, and
different walls in a Farey cell do not intersect.

Thus, there is a finite number of
stability domains inside a Farey cell $\cR_{u/v}$, and these are 
the components of the complement of the walls in $\cR_{u/v}$.
In particular, each stability domain is an open convex polygonal domain.
\end{theorem}

\begin{remark}\label{Rwallsmeet}
Recall that Farey cells are defined to be open;
it is possible for (the extensions of) two walls to intersect on the
boundary of the Farey cell, 
see 
\refT{ConjFarey+lines} and
Figures \ref{fig:m3-farey-triangles} and \ref{fig:m3-farey-zoom}.
\end{remark}

Label the vertices of $\gG^-$ (in order) as
$E_0=X,\dots,E_{2v}=Z=E_0+\ve_1$, and note that
these are given by \eqref{maj35}.

\begin{remark}\label{Rnowall}
 The special Farey cell $\cR_{1/2}=S_{1/2}$ has no walls.
Indeed, the construction in \refSS{SSm3S12} is valid throughout this Farey cell,
with $C$ distinct from all vertices. 
(We have $C\in\linje{Y,AY}\inter$, see \refFig{fig:m3-S12-main}.)
Alternatively, this follows from the proof below, since it is easily verified
algebraically that \eqref{m3-cut-wall} below cannot be satisfied for 
$(p,q)\in\cR_{1/2}$ and any $E_i$ and $\ell$; note that in this case the
possible $E_i$ are $B,U,A,Z$.

We will see in \refT{T-first-wall} that every other Farey cell has at least
one wall, and thus at least two stability domains.
\end{remark}

\begin{proof}[Proof of \refT{TSD}]
The points $P_k$ and $W_j$ defined in \eqref{m3-nu-Pk} and \eqref{qa2}
are obviously given by affine functions of $(p,q)$ (for a fixed $u/v$).
The same is true for $Q_j$, $B$, and $A$ by \eqref{m3-Qk}--\eqref{jda},
since we may use $(f,h)$ as a fixed coordinate system in the plane.
These points determine $\gG^\pm$ and $\tgG_\pm$.

By \refR{RC}, $\gG_-$ is the part of $\tgG_-$ between
the cut points $CX$ and $CZ=CX+\ve_1$. 
The symbolic description of $\gG_-$, and thus
the invariant polygon $\cD$, therefore changes if we change $(p,q)$ such that
$CX$ (and thus also $CZ$) passes through a vertex of $\tgG_-$.
We have shown in the proof above 
(see the beginning of the proof of \refT{TD})
that $\tgG_+$ and $\tgG_-$ coincide in a neighbourhood of $\partial\cS$, and thus
in a neighbourhood of $CX$. 
Furthermore, $T_1(\tgG_-)=\tgG_+$, see \eqref{qc3} and \eqref{qc2},
and thus it follows that $CX$ is a vertex of $\tgG_-$ if and only if $C$ is a
vertex of $\tgG_-$.
Hence, the symbolic description of $\cD$ depends only on 
which vertex $V$ or open line segment $\gU\inter$ of $\tgG_-$ that $\gD_g$
intersects. 

Fix (in the symbolic description) a  line segment 
$\gU=\linje{U_1,U_2}$ of $\tgG_-$, with $U_1$ to the left of $U_2$
and thus $g(U_1)<g(U_2)$ (since $g$ is strictly increasing along $\tgG_-$).
Then, for any
$(p,q)\in\cR_{u/v}$, 
we have by \eqref{m3-C} $C\in\gU\inter$ if and only if 
\begin{align}\label{wa1}
  g(U_1)<0<g(U_2).
\end{align}
We have shown that the coordinates of the vertices of $\gG^-$ 
(listed in \eqref{maj35})  
are affine functions of $(p,q)\in\cR_{u/v}$.
Thus the same is true for all vertices of $\tgG^-$ and of $\tgG_-$.
Consequently, \eqref{wa1} defines the intersection of the convex Farey cell
$\cR_{u/v}$ with two open affine half-planes. This set is open and convex,
and hence connected if non-empty; it is thus a stability domain.

Similarly, for each vertex $V$ in $\gG_-$, the set of $(p,q)$ such that $C=V$
is given by the affine equation $g(V)=0$. 
We will find these equations explicitly 
in \refT{TW2} below; in particular we will see that the equation
$g(V)=0$
is non-trivial and thus defines a line, which may or may not
intersect $\cR_{u/v}$. If it does, this line thus defines a wall.
We will also see in \refL{LW1}
that only a finite number of these intersections are non-empty, so there is
only a finite number of walls.

Finally, two walls cannot intersect at some $(p,q)\in\cR_{u/v}$,
because that would mean that for this $(p,q)$, the intersection point $C$
would coincide with two different vertices of $\tgG_-$, which clearly is
impossible.
\end{proof}

To study the walls further, and in particular prove the results
completing the proof of \refT{TSD},
we describe the vertices of $\tgG_-$ more explicitly.
The vertices of $\tgG^-$ are 
\set{E_i+\ell\ve_1:1\le i\le 2v,\,\ell\in\bbZ},
and thus the vertices of $\tgG_-=\tgG^-+Y$ are
\set{E_i+Y+\ell\ve_1:1\le i\le 2v,\,\ell\in\bbZ}.
(We may take $i\ge1$ since $E_0=E_{2v}-\ve_1$.)
Furthermore, using \eqref{m3-C} and $g(\ve_1)=2$, 
we have 
that $C$ coincides with the vertex
$E_i+Y+\ell\ve_1$ if and only if
\begin{align}\label{m3-cut-wall}
  g(E_i)+g(Y)+2\ell=0,
\qquad \text{where $i\in\set{1,\dots,2v}$ and $\ell\in\bbZ$}.
\end{align}

\begin{lemma}\label{LW1}
Let $u/v\in(0,\frac12]$. 
Then, for every $i=0,\dots,2v$,
and every $(p,q)\in\cR_{u/v}$, we have
  \begin{align}\label{lw1}
    -3 < g(E_iY)<1,
  \end{align}
and  consequently,
the only possible $\ell$ in \eqref{m3-cut-wall} are $\ell=0$ and $\ell=1$.
\end{lemma}
\begin{proof}
We have $g(E_i)\in [g(X),g(Z)]$ by \refL{Lm3-shadowing}.
Furthermore, 
\begin{align}\label{lw1a}
  0<g(Z)=2p+q < r+p+q=1,
\end{align}
and $g(X)=g(Z)-2$, $g(Y)=g(Z)-1$. Hence
\begin{align}\label{lw1b}
  g(E_iY)\in[g(X)+g(Y),g(Z)+g(Y)]
=[2g(Z)-3,2g(Z)-1]
\subset(-3,1).
\end{align}
\end{proof}

This lemma reduces \eqref{m3-cut-wall} 
to a set of $4v$ lines, of which the ones that
meet $\cR_{u/v}$ give the walls.
These are given more explicitly as follows.
For $0\le k\le v$, let $\kx\in\set{0,\dots,v-1}$ with
\begin{align}\label{kx}
  \kx=k+b\pmod v.
\end{align}
In other words, more explicitly,
\begin{align}\label{kx2}
  \kx=
  \begin{cases}
    k+b&\text{if }0\le k<d,
\\
    k-d&\text{if }d\le k\le v.
  \end{cases}
\end{align}

\begin{theorem}\label{TW2}
Let $u/v\in(0,\frac12]$.
Then the walls in the Farey cell $\cR_{u/v}$ are the intersections of 
$\cR_{u/v}$ with those of the following lines that yield a non-empty
intersection: 
\begin{align}
  2p+q&=\frac{2\nu(k)\pm1}{k+1}, \qquad k=1,\dots, v-1,\label{jw0}
\\
(k+\kx+2)p+(2\kx-k+1)q&
=\nu(k)+\nu(\kx)+\indic{k=1}\pm1,
\qquad k=0,\dots v-1
.\label{jw1}
\end{align}
More generally, any line that intersects $\cR\uv$ and is of one of the types
\begin{align}\label{jw0+}
  2p+q&=\frac{2m+1}{k+1}
\qquad\text{or}
\\\label{jw1+}
(k+\kx+2)p+(2\kx-k+1)q&
=\nu(k)+\nu(\kx)+\indic{k=1}+1+2m,
\end{align}
with $k$ as in \eqref{jw0} or \eqref{jw1} and $m\in\bbZ$,
defines a wall in $\cR\uv$ and is thus of the form \eqref{jw0} or
\eqref{jw1}.
\end{theorem}

We call walls given by \eqref{jw0} \emph{even walls}
and walls given by \eqref{jw1} \emph{odd walls}. 
These names are explained by the proof below, 
see also \refT{T-first-wall}.
The proof shows also that an even wall corresponds to $C$ coinciding with 
a valley in the lower boundary $\gG_-$, while an odd wall corresponds to $C$
coinciding with a peak in $\gG_-$.

\begin{proof}
Assume $u/v< 1/2$.
When $u/v=1/2$ there are no walls as said in \refR{Rnowall} and it can
be checked directly
that none of the lines \eqref{jw0}--\eqref{jw1+} intersects $\cR_{1/2}$;
alternatively the reader may verify that the arguments below hold with minor
modifications also in this exceptional case, yielding the same result
\eqref{jw0}--\eqref{jw1}, with no solutions in $\cR_{1/2}$.

Consider first \eqref{m3-cut-wall} for even $i$.
We call such walls \emph{even walls}; we will see that this agrees with
the definition just before the proof.
From the definition \eqref{maj35} of $\gG^-$, we see that the even walls are
defined by vertices $E_i$ in the set
\begin{align}\label{cw1}
  \set{E_{2j}:1\le j\le v}
=  \set{W_{j}:1\le j\le v,\, j\neq j_0}\cup\set{U}.
\end{align}
We may here ignore $U$, since, 
using \eqref{lw1a},
\begin{align}\label{cw11}
g(UY)=g(Y)=g(Z)-1\in(-1,0) 
\end{align}
in $\cR_{u/v}$, 
and thus \eqref{m3-cut-wall} cannot hold.
Recall that $\set{W_j:1\le j\le v}$ is a permutation of
$\set{P_k:1\le k\le v}$, see \eqref{qa2}, and that $\kk_{j_0}=v$ by
\eqref{qb3} and thus 
$W_{j_0}=P_v$. Hence, the interesting vertices are
\begin{align}\label{cw3}
 \set{W_{j}:1\le j\le v,\, j\neq j_0}
=
 \set{P_{k}:1\le k\le v-1}.
\end{align}
Each such $k$ gives, using \eqref{m3-nu-Pk}, \eqref{m3-Pkm}, and \eqref{m3-XYZ},
\begin{align}\label{cw4}
  P_kY=P(k,\nu(k))+Y
=\bigpar{(k+1)p-\nu(k),(k+1)q-1},
\end{align}
and thus
\begin{align}\label{cw5}
 g(P_kY)=2(k+1)p +(k+1)q-2\nu(k)-1.
\end{align}
Hence, \eqref{m3-cut-wall} for even $i$ yields the equations \eqref{jw0+}
(with $m=\nu(k)-\ell$)
for even walls, and
the two choices $\ell=0$ and $\ell=1$ in \eqref{m3-cut-wall} 
permitted by \refL{LW1} restricts this to the equations in \eqref{jw0}.

We turn to \emph{odd walls}, i.e., walls defined by \eqref{m3-cut-wall}
for an odd $i$.
By \eqref{maj35} again, these are defined by the vertices
\begin{align}\label{ow1}
  \set{E_{2j-1}:j=1,\dots,v}
=\bigset{Q_j:j\in\set{1,\dots,v}\setminus\set{j_0,j_0+1}}
\cup\set{B,A}.
\end{align}
By \eqref{m3-Qk} and \eqref{qa2}, we have, for $2\le j\le v$,
\begin{align}\label{ow2}
  g(Q_j) =h(Q_j)-f(Q_j) = h(W_j) - f(W_{j-1})
=h(P_{\kk_j})-f(P_{\kk_{j-1}}).
\end{align}
For the exceptional case $j=1$ we have instead, since we have defined
$W_0=X$ but $P_1=Z$,
\begin{align}\label{ow3}
  g(Q_1) = h(W_1) - f(W_{0})
=h(P_{\kk_1})-f(X)
=h(P_{\kk_1})-f(Z)-1.
\end{align}
We have also the exceptional cases $j=j_0$ and $j=j_0+1$, for which we
recall from \eqref{qb3} that $\kk_{j_0-1}=d$,  $\kk_{j_0}=v$, and $\kk_{j_0+1}=b$.
Thus,
by \eqref{jdb},
\begin{align}\label{ow4}
  g(B)=h(B)-f(B) = -f(W_{j_0-1})=-f(P_d),
\end{align}
and by \eqref{jda},
\begin{align}\label{ow5}
  g(A)=h(A)-f(A) = h(W_{j_0+1})=h(P_b).
\end{align}

Let $j\in\set{1,\dots,v}$ and let $k:=\kk_{j-1}$.
By the definition \eqref{qa1}, we have
$\kk_{j}\equiv \kk_{j-1}+b=k+b \pmod v$, and thus, recalling \eqref{kx},
$\kk_j=\kx$, except in the case $k=d$ when $\kk_j=v$ and $\kx=0$.
This exception thus occurs for $j=j_0$, which is excluded in \eqref{ow1}.
Hence, for $j\in\set{1,\dots,v}\setminus\set{j_0,j_0+1}$
we have by \eqref{ow2}--\eqref{ow3},
\begin{align}\label{ow6}
  g(Q_j) = h(P_\kx)-f(P_k) -\indic{j=1}
= h(P_\kx)-f(P_k) -\indic{k=1},
\end{align}
where $k\in\set{1,\dots,v}\setminus\set{\kk_{j_0-1},\kk_{j_0}}
=\set{1,\dots,v}\setminus\set{d,v}$.
Furthermore,
since
$P_0=P(0,0)=(0,0)$,
\eqref{ow4}--\eqref{ow5} may be written
\begin{align}\label{ow6B}
  g(B)&=h(P_0)-f(P_d) = h(P_{\kxx{d}})-f(P_d),
\\\label{ow6A}
g(A)&=h(P_b)-f(P_0)=h(P_{\kxx{0}})-f(P_0).
\end{align}
Consequently, since $j=1$ gives $k=\kk_0=1$, the equations
\eqref{m3-cut-wall} for odd $i$ are equivalent to the equations
\begin{align}\label{ow7}
  h(P_\kx)-f(P_k)-\indic{k=1}+g(Y)+2\ell=0,
\qquad k=0,\dots,v-1,\;\ell\in\bbZ.
\end{align}
By \eqref{m3-Pkm-sector}, \eqref{m3-nu-Pk},  and $g(Y)=g(Z)-1=\ga+\gb-1$, 
we may write \eqref{ow7} as
\begin{align}\label{ow8}
(k+1)\ga +(\kx+1)\gb =\nu(k)+\nu(\kx)+\indic{k=1}+1-2\ell,
\end{align}
which by \eqref{m3-Phi} is the same as
\begin{align}\label{ow9}
(k+\kx+2)p +(2\kx-k+1)q =\nu(k)+\nu(\kx)+\indic{k=1}+1-2\ell.
\end{align}
This gives the equations \eqref{jw1+} for odd walls (with $m=-\ell$), and
since necessarily $\ell\in\set{0,1}$ by \refL{LW1}, we may restrict to
\eqref{jw1}.
\end{proof}

We end this section with some further observations on the construction.

\begin{theorem}\label{Tm3-4v+3}
For every stability domain inside a Farey cell $\cR_{u/v}$, the invariant
polygon $\cD$ has $4v+3$ vertices.
For $(p,q)$ on a wall inside $\cR_{u/v}$, the polygon $\cD$ has $4v+2$
vertices. 
\end{theorem}

\begin{proof}
The upper boundary $\gG^+$ has the $v+1$ valleys
$W_0,\dots,W_v$ and the $v$ peaks $Q_1,\dots,Q_v$, and hence $2v+1$
vertices.  

Let us project the construction above to the cylinder $(\bbR/\bbZ)\times\bbR$.
This identifies $X$ and $Z$, and thus $\gG^-$ becomes a closed
polygonal curve $\chgG^-$ with $2v$ vertices on the cylinder;
hence, $\tgG_-$ becomes a translated closed polygonal curve $\chgG_-$ with
$2v$ vertices.
We obtain $\gG_-$ by cutting this curve $\chgG_-$ at the image of $CX$
(which equals the image of $CZ$).
We saw in the proof of \refT{TSD} that when $(p,q)$ lies 
on a wall in $\cR_{u/v}$, then 
$CX$ coincides with a vertex of $\tgG_-$; hence, cutting there yields
a lower boundary $\gG_-$ with $2v+1$ vertices,
and the total number of vertices on $\cD$ is $4v+2$.
On the other hand, for $(p,q)$ in a stability
domain, $CX$ lies on a side, and cutting $\chgG_-$ there yields
a lower boundary $\gG_-$ with $2v+2$ vertices,
and the total number of vertices on $\cD$ is $4v+3$.
\end{proof}

\refT{Tm3-4v+3} is valid only inside the open Farey cells. For parameters
$(p,q)$ on the boundaries between Farey cells, further coincidences of
vertices may occur, see \refR{Rrand}. 

\begin{example}\label{ER13}
Consider the Farey cell $\cR_{1/3}$. 
The Farey parents of $\frac13$ are $\frac ab=\frac01$ and $\frac cd=\frac 12$.
The vertices of $\cR_{1/3}$ are, by
\eqref{m3farey2} and \eqref{m3-fcu}--\eqref{m3-fcm},
\begin{align}\label{R13}
V_{1/3}&=\Bigl(\frac13,0\Bigr),&
L_{1/3}&=\Bigl(\frac19,\frac19\Bigr),\notag\\
U_{1/3}&=\Bigl(\frac16,\frac16\Bigr),&
M_{1/3}&=\Bigl(\frac7{18},\frac1{18}\Bigr).
\end{align}
Put also
\begin{align}
 L'_{1/3}:=\Bigl(\frac4{18},\frac1{18}\Bigr).
\end{align}
Recall also that in the alternative coordinates $(\ga,\gb)$,
 $\cR_{1/3}$ is given by the rectangle $0<\alpha<1/3<\beta<1/2$, 
see \eqref{m3-farey-rect2} and \eqref{m3-open-farey}.
It is easily checked that of the 10 lines in \eqref{jw0}--\eqref{jw1},
only the following 2 intersect $\cR_{1/3}$ and thus define walls:
\begin{align}
p+q=\frac13
 &\quad\Longleftrightarrow\quad 
 \alpha+2\beta=1
 \quad\text{i.e.,  }\linje{V_{1/3},U_{1/3}},\\
2p+q=\frac12
 &\quad\Longleftrightarrow\quad 
 \alpha+\beta=\frac12
 \quad\text{i.e., }\linje{L'_{1/3},U_{1/3}}.
\end{align}
(The first wall is odd and the second is even.
The first line is the one denoted $\LL_{1/3}$ in \refT{T-first-wall}.)

The three stability domains in this cell are thus the triangles
\begin{align}
S_{1/3}\SA&=[V_{1/3},U_{1/3},M_{1/3}],\label{m3-S13SA}\\
S_{1/3}\SB&=[V_{1/3},U_{1/3},L'_{1/3}],\label{m3-S13SB}\\
S_{1/3}\SC&=[L_{1/3},L'_{1/3},U_{1/3}],\label{m3-S13SC}
\end{align}
where, as usual, only the interiors are meant.
The three stability domains give invariant polygons $\cD$ with different
types of lower boundary, 
but by \refT{Tm3-4v+3} all polygons have $4\cdot3+3=15$ vertices.

The construction in \refSS{SSm3construction} yields 
for $\cR_{1/3}$
the following results 
by straightforward calculations.
We have $j_0=2$, and thus 
$B$ is defined by \eqref{jdb} as the intersection of $\gD_h$ and
the $f$-line through $W_1=P_2=XZ$
(for all three stability domains);
equivalently, 
\begin{align}\label{emB}
 B&:=\gD_h\cap T_3T_1\gD_f.
\end{align}
The intersection point $C$ on the lower boundary defined by \eqref{m3-C} 
lies on different sides of the lower boundary of $\cD$ 
(recall that this is what defines the different stability domains in a Farey
cell); for
$S_{1/3}\SA$ , $S_{1/3}\SB$, and $S_{1/3}\SC$ we have $C$ equal to, respectively,
\begin{align}\label{emC1}        
C_1&:=\gD_g\cap T_2\gD_f 
\phantom{T_0}\qquad\text{on }\linje{Y,AY},    
\\\label{emC2}
C_2&:=\gD_g\cap T_2T_3\gD_h
\qquad\text{on }\linje{AY,YZ},
\\\label{emC3}
C_3&:=\gD_g\cap T_2T_3\gD_f
\qquad\text{on }\linje{YZ,AYZ}.
\end{align}
The construction now yields the
explicit formulas below for $\cD$,
where we list the vertices in cyclic
order around the boundary, starting
with the upper boundary from left to right.
In $S_{1/3}\SA$ , $S_{1/3}\SB$, and $S_{1/3}\SC$ we have, respectively,
\begin{align}\label{DS13A}
\cD&=[X,AX,XZ,AXZ,XZ^2,BZ,Z,C_1Z,YZ,AY,Y,BY,XYZ,AXY,C_1X],
\\\label{DS13B}
\cD&=[X,AX,XZ,AXZ,XZ^2,BZ,Z,C_2Z,AYZ,YZ,AY,Y,BY,XYZ,C_2X],
\\\label{DS13C}
\cD&=[X,AX,XZ,AXZ,XZ^2,BZ,Z,C_3Z,YZ^2,AYZ,YZ,AY,Y,BY,C_3X].
\end{align}
This illustrates \refT{TSD}:
the upper boundary is the same Farey zig-zag in all three domains,
but the lower boundary changes.

See \refFig{fig:m3-S13} for one example from $S_{1/3}\SC$.
\end{example}

\begin{remark}\label{Rm3-vertex-orbits}
For a parameter in a stability domain, the vertices can also be grouped into
truncated branch orbits.  The points
$U,A,B,C$ are vertices of the pieces, but not of $\cD$: $U$ belongs to all
three pieces, $A$ to $\cP_1$ and $\cP_2$, $B$ to $\cP_2$ and $\cP_3$, and
$C$ to $\cP_1$ and $\cP_3$; 
these points are thus mapped to (including all branches)
$X,Y,Z,AX,AY,BY,BZ,CX,CZ$, which all are vertices of $\cD$.
Starting with these vertices and following their orbits
as long as they remain vertices of $\cD$, we obtain
(as can be seen from the construction in \refSS{SSm3construction})
the following partial orbits of vertices:
\begin{align}
\begin{array}{c|cccccccc}
\text{initial vertex}&X,Z&Y&AX&AY&BZ&BY&CX&CZ\\ \hline
\text{number of vertices}&v+1&v&d&d&b&b&1&1
\end{array}
\end{align}
Here $X$ and $Z$ have the common next image $XZ$, 
followed by $v-2$ further vertices (all in $\gG^+$);
hence each of them starts a partial orbit with $v$ vertices, but these coincide
from the second vertex so they cover together $v+1$ vertices.
In contrast, the
two images of $C$ have the common next image $CXZ$, which is not a vertex of
$\cD$; hence $CX$ and $CZ$ each yields a partial orbit of size 1.  
Note that this grouping gives again
$(v+1)+v+2d+2b+2=4v+3$ vertices, agreeing with \refT{Tm3-4v+3}.
On a wall this grouping has to be modified,
since $C$ is then itself a vertex; see the proof of \refT{TSD}.
\end{remark}

\begin{remark}\label{Rword}
The $X/Z$ word along the branch
$U,P_1,\dots,P_v$ is the standard rational rotation word of slope $u/v$:
the prefix ending at $P_k$ contains $\nu(k)$ letters $X$ and
$k-\nu(k)$ letters $Z$,
see \eqref{m3-nu-Pk} and \eqref{m3-Pkm}.  
Thus $u/v$ determines this $X/Z$ word, 
and \eqref{m3-nu-Pk} and \eqref{qa1}
give the orbit points and their order for the construction above.
\end{remark}

Finally, we make some observations on the  hexagonal piece $\cP_2$,
including some relations between different Farey cells.
To keep track of the data attached to the Farey cell, write
$A_{u/v}:=A$ and $B_{u/v}:=B$.

\begin{proposition}\label{Pm3-cell-data}
The symbolic form of $A_{u/v}$, $B_{u/v}$ and $\cP_2$ is fixed by the Farey
cell $\cR_{u/v}$, and is thus the same for all stability domains inside the
cell.  More precisely,
\begin{align}\label{m3-AB-cell-data}
 \bigl(f(A_{u/v}),h(A_{u/v})\bigr)&=(0,b\beta-a),
 &
 \bigl(f(B_{u/v}),h(B_{u/v})\bigr)&=(c-d\alpha,0),
\end{align}
the two 
peaks of the middle hexagon satisfy
\begin{align}\label{m3-middle-peaks}
 Q_{j_0}&=A_{u/v}X^cZ^{d-c},
 &
 Q_{j_0+1}&=B_{u/v}X^aZ^{b-a},
\end{align}
and consequently
\begin{align}\label{m3-middle-explicit}
 \cP_2=
 [U,B_{u/v},A_{u/v}X^cZ^{d-c},X^uZ^{v-u},
 B_{u/v}X^aZ^{b-a},A_{u/v}].
\end{align}
The tie vertices have the following inheritance property along the Farey tree:
\begin{align}\label{m3-AB-inheritance}
 b>d&\quad\Longrightarrow\quad B_{u/v}=B_{a/b},
 &
 d>b&\quad\Longrightarrow\quad A_{u/v}=A_{c/d}.
\end{align}
\end{proposition}

\begin{proof}
For $u/v=1/2$, the assertions follow directly from
\eqref{pda}--\eqref{pdb} and \eqref{pd2}, with Farey parents $0/1$ and $1/1$;
\eqref{m3-AB-inheritance} is vacuous since $b=d=1$.

Assume henceforth that $u/v<1/2$.
By \eqref{qa2} and \eqref{qb3}, the three central valleys 
in the upper boundary $\gG^+$
are
$W_{j_0-1}=P_d$, $W_{j_0}=P_v$ and $W_{j_0+1}=P_b$.  
Equations \eqref{jdb}, \eqref{jda}, \eqref{jc2b}--\eqref{jc2d},
and \eqref{m3-Pkm-sector} therefore give
\eqref{m3-AB-cell-data}.  Moreover, \eqref{jc2b}--\eqref{jc2d} and
\eqref{jc2v} yield $P_b+P_d=P_v$, since $b+d=v$ and $a+c=u$.  
Since $f$ and $h$ are
linear, 
and furthermore may be regarded as coordinates in $\bbR^2$,
\eqref{jdb}, \eqref{jda} and \eqref{m3-Qk} now yield, after simple calculations,
\begin{align}
 Q_{j_0}=A_{u/v}+P_d,
 \qquad
 Q_{j_0+1}=B_{u/v}+P_b.
\end{align}
This proves \eqref{m3-middle-peaks} and \eqref{m3-middle-explicit},
recalling \eqref{jc2b}--\eqref{jc2d} and \eqref{maj36}.
In particular, the symbolic expressions for all six vertices 
of $\cP_2$ are fixed by the Farey cell.

Suppose that $b>d$.  Then
\begin{align}
 \frac ab=\frac{(a-c)+c}{(b-d)+d},
\end{align}
and $bc-ad=1$ shows that $(a-c)/(b-d)$ is non-negative and neighbouring to
both $a/b$ and $c/d$.  Hence $(a-c)/(b-d)$ and $c/d$ are the Farey parents
of $a/b$.  Thus the right
parent used to define $B_{a/b}$ is again $c/d$, and
\eqref{m3-AB-cell-data} implies $B_{a/b}=B_{u/v}$.  If $d>b$, then similarly
\begin{align}
 \frac cd=\frac{a+(c-a)}{b+(d-b)},
\end{align}
where $(c-a)/(d-b)$ is positive and neighbouring to both $a/b$ and
$c/d$.  Thus the left parent used to define $A_{c/d}$ is $a/b$ and
$A_{c/d}=A_{u/v}$.  This proves \eqref{m3-AB-inheritance}.
\end{proof}

Recall that, in contrast, the symbolic form of $C$ depends on the stability
domain. 

\subsection{The system of walls}\label{SSCW}

\begin{figure}[htbp]
\centering
\includegraphics[width=0.9\textwidth]{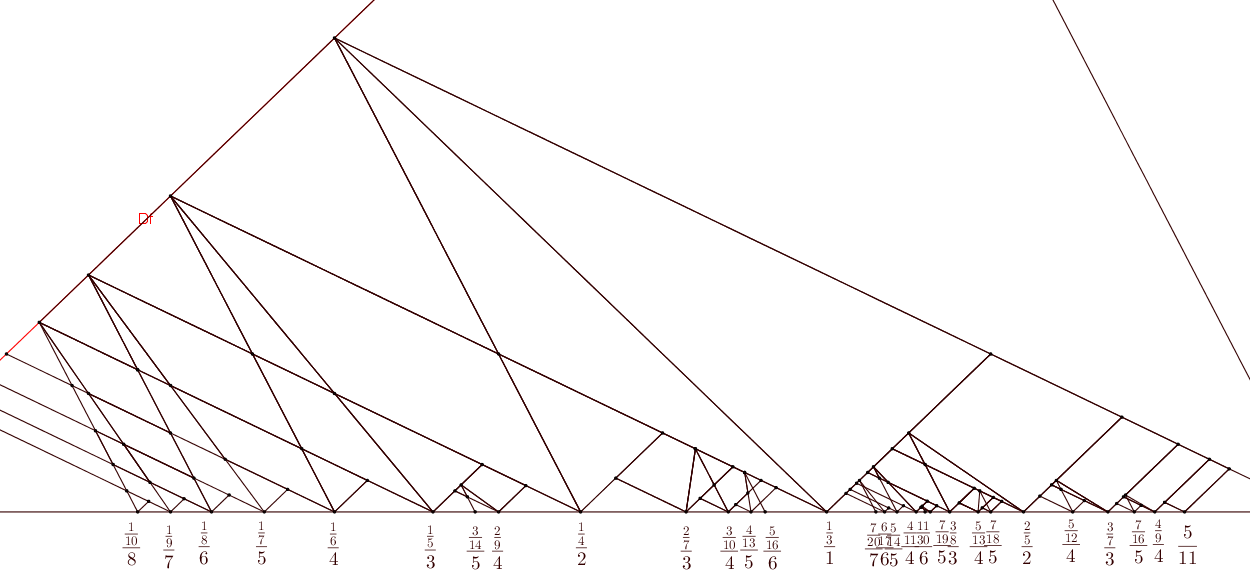}
\caption{Farey cells together with the wall segments $\LG_{u/v}$.
The number below a fraction gives its generation in the Farey tree.}
\label{fig:m3-farey-triangles}
\end{figure}

\begin{figure}[htbp]
\centering
\includegraphics[width=0.9\textwidth]{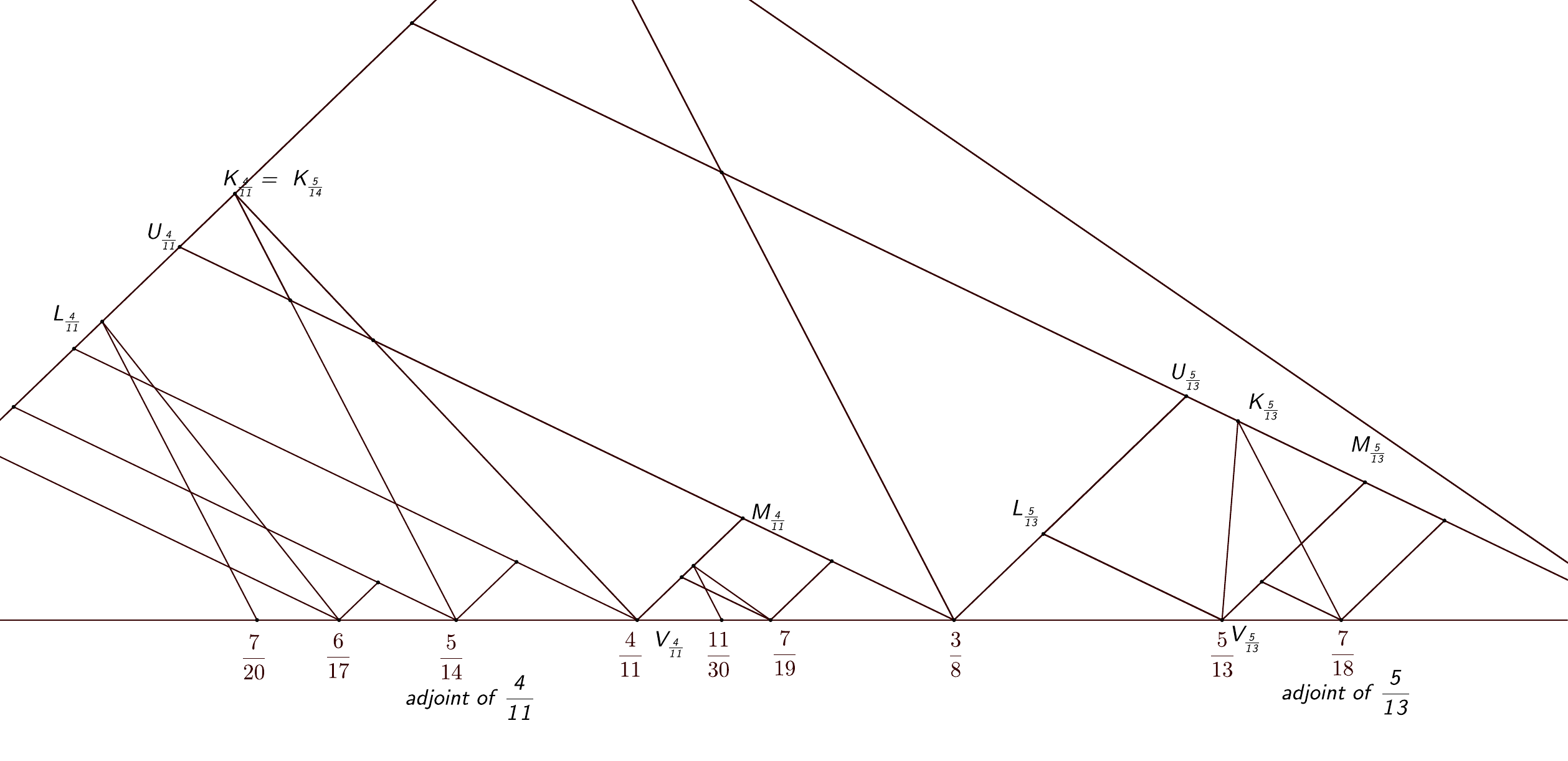}
\caption{A magnified part of the wall pattern in
\refFig{fig:m3-farey-triangles}.}
\label{fig:m3-farey-zoom}
\end{figure}

We have in the preceding subsection studied the walls inside a fixed
Farey cell.
A glance at \refFigs{fig:m3-farey-triangles} and \ref{fig:m3-farey-zoom}
shows that the walls typically extend in a straight line across several
Farey cells. In this section we study this pattern of walls.
Recall that there is no wall in $\cR_{1/2}$, so we will only consider
Farey cells $\cR_{u/v}$ with $u/v<\frac12$.

We first show that in every such Farey cell $\cR\uv$, 
there is exactly one wall that ends at the vertex
$V\uv=(\frac uv,0)$ on the axis $q=0$.
(Equivalently, the line extending the wall goes through $V\uv$.)

For a reduced fraction $u/v\in(0,\frac12)$, let
$a/b<c/d$ be the Farey parents of $u/v$, and define the line
\begin{align}\label{m3-H}
 \HG_{u/v}:\qquad v\left(p-\frac uv\right)+\tau_{u/v}q=0,
\end{align}
where
\begin{align}\label{m3-tau}
\tau_{u/v}=
\begin{cases}
\displaystyle \frac v2,
        & v \text{ even},\\[5pt]
\displaystyle 2b+\frac d2,
        & v \text{ odd},\ b \text{ odd},\ d \text{ even},\\[5pt]
\displaystyle \frac b2-d,
        & v \text{ odd},\ b \text{ even},\ d \text{ odd}.
\end{cases}
\end{align}

\begin{theorem}\label{T-first-wall}
  Let $u/v$ be a reduced fraction in $(0,\frac12)$.
Then there is exactly one wall in $\cR\uv$ with endpoint
$V\uv=(\frac{u}v,0)$.
Consequently, there are exactly two stability domains in $\cR_{u/v}$ 
that have $V_{u/v}$ on their boundaries. (Or, equivalently, contain $(p,q)$
with $q$ arbitrarily small.)

Furthermore, the wall
with endpoint $V\uv=(\frac{u}v,0)$
is the intersection of $\cR\uv$ and the line $\HG\uv$ in \eqref{m3-H}.
The wall is even if and only if  $v$ is.
\end{theorem}

\begin{proof}
By \refT{TW2} and its proof, the walls in $\cR\uv$ are the intersections
of $\cR\uv$ and the lines defined by \eqref{m3-cut-wall}.
Now let $(p,q)\in\cR\uv$ converge to the boundary point $V\uv$.
Then, by continuity (see \refR{Rrand} and its proof) the point $E_i$
converges to a limit $E_i\oo$, where we as for other points below
use the superscript ${}\oo$ to denote the limiting case
$(p,q)=(u/v,0)$.

We consider the construction in \refSS{SSm3construction} for this limiting
case. 
Then
$\ga=\gb=u/v$ and
$\delta_-=\delta_+=0$, see \eqref{m3-lambda-mu},
and thus \eqref{qb2f}--\eqref{qb2h} simplify to, recalling \eqref{qa8}, 
\begin{align}\label{em4}
  -f(W_j\oo)=h(W_j\oo)=\ttt_j = \frac{u+j-v}{v},
\qquad j=0,\dots,v.
\end{align}
We have $\ttt_{j_0}=0$ (since $j_0=v-u$), and thus
$f(W_{j_0}\oo)=h(W_{j_0}\oo)=0$; hence it follows from
\eqref{m3-Qk}--\eqref{jda} that 
$B\oo=Q_{j_0}\oo$ and $A\oo=Q_{j_0+1}\oo$. 
Hence, $\gG^+=\gG^-$ in this limiting case, and
the points $E_i\oo$ in \eqref{m3-cut-wall} are simply
$\set{W_j\oo,Q\oo_j:j=1,\dots,v}$.
We have by \eqref{em4} and \eqref{m3-Qk}, since now $g(Y\oo)=g(u/v,-1)=2u/v-1$,
\begin{align}
   \label{em5w}
g(W_j\oo Y\oo)&=-f(W_j\oo)+h(W_j\oo)+g(Y\oo)
=2\ttt_j+\frac{2u-v}{v}=\frac{4u-3v+2j}{v}.
\\\label{em5q}
g(Q_j\oo Y\oo)&=-f(W_{j-1}\oo)+h(W_j\oo)+g(Y\oo)
=\frac{4u-3v+2j-1}{v}.
\end{align}
(At this limiting parameter, the upper and lower boundaries become
periodic zig-zag lines.)
The numerators in \eqref{em5w} and \eqref{em5q} are exactly the integers in 
$[4u-3v+1,4u-3v+2v]$, and
it follows that there is exactly one vertex 
$V\oo_i$ such that 
this numerator is divisible by $2v$, and thus $g(V\oo_iY\oo)$ is an even
integer $-2\ell$. This means that \eqref{m3-cut-wall} holds 
for $(p,q)=(u/v,0)$, and thus the proof of \refT{TW2} shows that 
$(u/v,0)$ lies on exactly one of the lines \eqref{jw0+}--\eqref{jw1+}.

It remains to show that this line, $\LL$ say, intersects the open Farey
cell $\cR\uv$, and thus defines a wall by \refT{TW2}.
To see this, we 
move along the line $\LL$ with $q$ increasing from $0$, and use
the coordinates $\ga$ and $\gb$ in \eqref{m3-Phi}. 
In the case of a line \eqref{jw0+}, we write
\eqref{jw0+} as 
\begin{align}\label{em6}
\ga+\gb=\frac{2m+1}{k+1},
\end{align}
and
in the case of a line \eqref{jw1+}, we recall from the proof of \refT{TW2}
that it can be written as \eqref{ow8}.
In both cases, the coefficients of $\ga$ and $\gb$ are both positive.
Recall also that  $\gb-\ga=3q$ by \eqref{Phi1}.
It follows that as $q$ increases, $\ga$ and $\gb$ move in opposite
directions, so $\ga$ decreases and $\gb$ increases from their common initial
value $\frac{u}v$ when $q=0$. 
Consequently, for small positive $q$, the
inequalities \eqref{m3-open-farey} are satisfied, and thus these
points on $\LL$ lie in $\cR\uv$. (See also \eqref{m3-farey-rect2}.)
This shows that the line $\LL$ intersects $\cR\uv$, and thus it defines a
wall. 
This proves the first statement in the theorem.

To identify the line $\LL$ explicitly, we consider the cases when it 
gives an even or odd wall separately.

\pfcase1{even walls.}
 First, suppose that $\LL$ is  given by  \eqref{jw0+},
and thus defines an even wall.
Then, since $\LL$ goes through $(u/v,0)$,
\begin{align}\label{bay1}
    2\frac{u}v&=\frac{2m+1}{k+1},
\end{align}
which is the same as
\begin{align}\label{bay2}
  2(k+1)u=(2m+1)v.
\end{align}
In particular, $v$ has to be even. 
Furthermore, the line \eqref{jw0+} has slope $-2$ and can thus
be written as \eqref{m3-H}
with $\tau\uv=v/2$ as given in \eqref{m3-tau}, 
so the line $\LL$ equals  $\HG\uv$.

\pfcase2{odd walls.}
Suppose that $\LL$ is  given by  \eqref{jw1+},
and thus defines an odd wall.
Denoting the \rhs{} of \eqref{jw1+} by $z$, we see that 
$(u/v,0)\in\LL$ means that
\begin{align}\label{bay3}
  (k+\kx+2)\frac{u}v = z
.\end{align}
Since $u$ and $v$ are coprime, this requires that $k+\kx+2$ is a multiple of
$v$.
Furthermore, $0<k+\kx+2\le (v-1)+(v-1)+2=2v$, 
with equality in the second inequality only
if $k=\kx=v-1$, which is impossible by \eqref{kx}. 
Consequently, \eqref{bay3} implies 
\begin{align}\label{bay5}
  k+\kx+2=v.
\end{align}
Moreover, \eqref{kx} implies that either $\kx=k+b$ or $\kx=k+b-v=k-d$.

\pfcase{2a}{$\kx=k+b$.}
In this case, \eqref{bay5} is $2k+b+2=v=b+d$, which has the solution
\begin{align}\label{bay6d}
  k=\tfrac12d-1
\qquad\text{and thus}\qquad
\kx=b+\tfrac12d-1.
\end{align}
Thus, in this case $d$ has to be even.
Furthermore, \eqref{bay6d} yields
\begin{align}\label{bay6d+}
  2\kx-k+1=k+2b+1=2b+\tfrac12 d.
\end{align}

\pfcase{2b}{$\kx=k-d$.}
In this case, \eqref{bay5} is $2k-d+2=v=b+d$, which has the solution
\begin{align}\label{bay6b}
  k=\tfrac12b+d-1
\qquad\text{and thus}\qquad
\kx=\tfrac12b-1.
\end{align}
Thus, in this case $b$ has to be even.
Furthermore, \eqref{bay6b} yields
\begin{align}\label{bay6b+}
  2\kx-k+1=k-2d+1=\tfrac12b- d.
\end{align}

Consequently, in both cases one of $b$ and $d$ is even, and then the other
is necessarily odd, and thus $v$ is odd.
Furthermore, \eqref{bay6d+} and \eqref{bay6b+} show that in both cases, 
$2\kx-k+1=\tau\uv$ defined in \eqref{m3-tau}. This and \eqref{bay5}
show that the line $\LL$ equals $\HG\uv$ in \eqref{m3-H} in the case of an odd
wall too.

We have also seen that $v$ is necessarily even if the wall is even,
and odd if the wall is odd.
\end{proof}

To describe the walls  away from the axis $q=0$,
we state first a simple criterion for finding the Farey cell that a given
$(p,q)$ belongs to.

\begin{lemma}\label{LUV}
  Let $(p,q)\in\gD$, and let as usual $(\ga,\gb):=(p-q,p+2q)$;
 thus $0<\ga<\gb<1$ by \eqref{Phi2}. 
Let $u/v$ be the reduced fraction in $[\ga,\gb]$ with smallest denominator. 
If $u/v\in(\ga,\gb)$, then $(p,q)\in\cR_{u/v}$.
Otherwise, i.e., if $u/v=\ga$ or $u/v=\gb$, then $(p,q)$ lies on the boundary
$\partial\cR_{u/v}$. 
\end{lemma}
\begin{proof}
  Note first that there is always a unique $u/v \in[\ga,\gb]$ with smallest
  denominator. 
To see this, 
let $v$ be the smallest denominator of any reduced fraction in $[\ga,\gb]$,
and consider the set of all reduced fractions in $[0,1]$ with denominator at
most $v$. 
(This is called a \emph{Farey sequence} (or \emph{Farey series}); 
Farey sequences are closely connected to the Farey tree, but consider the
reduced fractions in a different order.)
If $[\ga,\gb]$ would contain two or more reduced fractions with denominator
$v$, then there cannot be any fraction with smaller denominator between them, so
they would be neighbours, but neighbours in a Farey sequence are Farey
neighbours in the sense above, see \eqref{sb1}, which is impossible if they
have the same denominator $>1$.
(For this and other properties of Farey sequences
see e.g.\ \cite[p.~118--119]{CM}.)

Thus, let $u/v$ be the unique reduced fraction in $[\ga,\gb]$ with
smallest denominator, and let its Farey parents be $a/b$ and $c/d$.
Since $b,d<v$, it follows that $a/b$ and $c/d$ cannot belong to $[\ga,\gb]$,
and thus
\begin{align}\label{jz0}
  \frac ab <\ga \le \frac uv \le\gb <\frac cd.
\end{align}
The result now follows from \eqref{m3-farey-rect2}, see also
\eqref{m3-open-farey}. 
\end{proof}

We can now state our main result on walls, describing all walls in
all Farey cells explicitly. We first need some more notation.

Let $u/v$ be a reduced fraction in $(0,\frac12)$,
and let its Farey parents be $a/b$ and $c/d$.
Recall the line $\LL\uv$ through $V\uv=(\frac uv,0)$ defined in
\eqref{m3-H}. 
Define $K_{u/v}$ as the point on the line $\LL\uv$ whose height $q$ is
\begin{align}\label{qK}
 q(K_{u/v}):=
  \begin{cases}
\displaystyle\frac{2}{3v \min(b,d)},
        & v \text{ even},\\[5pt]
\displaystyle\frac{2}{3b(b+v)},
        & v \text{ odd},\ b \text{ odd},\ d \text{ even},\\[5pt]
\displaystyle\frac{2}{3d(d+v)},
        & v \text{ odd},\ b \text{ even},\ d \text{ odd}.
\end{cases}
\end{align}
and put
\begin{align}\label{guv}
 \LG_{u/v}:=\linje{V_{u/v},K_{u/v}}.
\end{align}
which thus is a segment of the line $\LL\uv$.

\begin{theorem}\label{Twalls}
  If $u/v\in(0,\frac12)$ is a reduced fraction, then
the line segment $\LG\uv$ intersects a finite number of Farey cells,
and the intersection of $\LG\uv$ with each of them is a wall in
that cell. 
This wall is even (odd) if and only if the
denominator $v$ is.

Conversely, any wall in an open Farey cell
is the intersection of that Farey cell and some line segment $\LG\uv$ 
with $u/v\in(0,\frac12)$.
\end{theorem}

For clarity, we prove \refT{Twalls} for odd and even walls separately, although 
the two proofs are rather similar, and some arguments are the same.

\begin{proof}[Proof of \refT{Twalls} for even walls.]
  First, by \refT{TW2}, any even wall in a Farey cell is
  the intersection of the Farey cell and a line $\LL$ given by \eqref{jw0+},
where,
since $0<2p+q<1$ in $\gD$, we have $0<2m+1<k+1$.
The line $\LL$ 
intersects the $p$-axis at $(u/v,0)$ with $u/v:=(2m+1)/(2(k+1))\in(0,\frac12)$,
and since it has slope $-2$, it equals $\LL\uv$.

The even walls in a Farey cell $\cR\uvi$  are thus given by the
intersections with  a line $\LL\uv$ such that this intersection is nonempty
and furthermore $\LL\uv$  can be written as \eqref{jw0+} for the Farey cell
$\cR\uvi$, which simply means that $k\le v'-1$.

We switch to the coordinates $(\ga,\gb)$. Then $\LL\uv$ is given by
\begin{align}\label{jz1}
  \Phi(\LL\uv)=\bigset{(\tfrac uv-s,\tfrac uv+s):s\in\bbR},
\end{align}
where we note that we have, by \eqref{Phi1},
\begin{align}\label{jz2}
  q=\frac{\gb-\ga}3=\frac{2s}3.
\end{align}
Thus $s=0$ gives the point $V\uv$. 
Consider a point $P\in\gD$ on $\LL\uv$ with $s>0$, and suppose for simplicity that
$s$ is irrational. 
Then, by \refL{LUV}, $P\in \cR\uvi$ where $u'/v'$ is the reduced fraction in 
$(\tfrac uv -s,\tfrac uv+s)$ with smallest denominator.
For small $s>0$ we thus have $u'/v'=u/v$, but if $s$ grows, eventually
$v'$ decreases. 
In particular, $v'\le v$, which implies that there is only a finite number
of different Farey cells $\cR\uvi$ that are hit by $\LL\uv$.

The equation \eqref{m3-H} for $\LL\uv$ can be written as 
\begin{align}\label{jz33}
2p+q=\frac{2u}v = \frac{u}{v/2},  
\end{align}
which (since $u$ is odd because $v$ is even)
is of the type \eqref{jw0+}
with $k=v/2-1$.
Hence, $P$ is on an even wall in $\cR\uvi$ if and only if 
$v/2-1\le v'-1$.
(Note that the fraction \eqref{jz33} is irreducible.)
This is equivalent to $v'\ge v/2$.
Consequently:

\emph{$P$ lies on an even wall if and only if 
there is no fraction $u'/v'\in(\tfrac uv -s,\tfrac uv+s)$ with denominator
$v'<v/2$}.

Equivalently, we have:

\emph{$P$ lies on an even wall if and only if $0<s<s\uv$, where
$s\uv$ is the distance from $u/v$ to the closest fraction $u'/v'$ with
denominator $v'<v/2$.}

It remains to find the closest $u'/v'$ with $v'<v/2$.
Let the Farey parents of $u/v$ be $a/b$ and $c/d$.
Assume first that $b<d$, and thus $b<v/2$. 
Since $a/b$ and $c/d$ are Farey neighbours, it
follows that the Farey parents of $c/d$ are $a/b$ and $(c-a)/(d-b)$.
If $d-b>b$, we can continue, and by induction we obtain an increasing
sequence $r_i:=(u-ia)/(v-ib)$, $i=0,\dots,i_1:=\ceil{v/b}-1$, such that
$r_0=u/v$, $r_1=c/d$, and 
the Farey parents of $r_i$ are $a/b$ and $r_{i+1}$ for $0\le i<i_1$.
Let $i_2:=\floor{v/(2b)}+1\le i_1$. Then the denominator of $r_i$ is
$v-ib\ge v/2$ for $i<i_2$, but $v-i_2b<v/2$.
It follows from \refL{Lm3-no-small} applied to
$(\frac{a}{b},\frac{u}{v})$ and all $(r_i,r_{i+1})$ for $0\le i< i_2$
that every reduced fraction in
$(a/b,r_{i_2})$ has denominator that is at least equal to one of the
denominators of $r_i$, $i< i_2$, and thus at least $v/2$.
Consequently, the reduced fractions with denominator $<v/2$ closest to $u/v$
to the left and right are $a/b$ and $r_{i_2}$.
Furthermore, since by definition $i_2>v/(2b)$ and thus $i_2b>v/2$, 
\begin{align}\label{vc1}
  r_{i_2}-\frac{u}v 
= \frac{u-i_2a}{v-i_2b}-\frac{u}v
=\frac{i_2(bu-av)}{v(v-i_2b)}
=\frac{i_2}{v(v-i_2b)}
>\frac{v/(2b)}{v\cdot v/2}
=\frac1{bv}
=\frac{u}v-\frac{a}b.
\end{align}
Hence, if $b<d$, then $a/b$ is the reduced fraction with denominator $<v/2$
closest to $u/v$, and the distance is $1/vb$.

If $b>d$, then 
a similar argument with $r_i:=(u-ic)/(v-id)$ shows that
$c/d$ is the reduced fraction with denominator
$<v/2$ 
closest to $u/v$, and the distance is $1/vd$.

Consequently, we have in both cases
\begin{align}\label{vc2}
  s\uv=\frac{1}{v\min(b,d)}.
\end{align}
The corresponding value of $q$ is by \eqref{jz2} and \eqref{qK}
\begin{align}\label{vc3}
q= q\uv :=\frac{2s\uv}{3}
=\frac{2}{3v\min(b,d)}
=q(K\uv).
\end{align}
The corresponding value of $p$ follows from \eqref{jz33}:
\begin{align}\label{vc4}
p\uv :=\frac{u}v -\frac12q\uv
=\frac{u}v - \frac{1}{3v\min(b,d)}
.\end{align}

At $K\uv$, the parent of $u/v$ with the smaller denominator is an endpoint of
$[\ga,\gb]$. The closest-fraction comparison above shows that it is the
unique fraction of smallest denominator in this interval.
Hence, by \refL{LUV}, 
$K\uv$ belongs to no open Farey cell.
The point $V\uv$ also belongs to no open Farey cell, 
since $q=0$ and $\ga=\gb$ there.

Consequently, if
$\LG\uv$ intersects a Farey cell $\cR\uvi$, then (by considering a point on
the intersection with irrational $q$), we see that the intersection is a
wall in $\cR\uvi$. (And also that $v'\ge v/2$.) 

Conversely, if 
the intersection of $\LL\uv$ and $\cR\uvi$ is a wall
in the latter, then $q\le q(K\uv)$ for any point $(p,q)$ on the intersection
with irrational $q$, and thus for every point on the intersection.
Hence the wall is contained in $\LG\uv$, and equals $\LG\uv\cap\cR\uvi$.
\end{proof}

\begin{proof}[Proof of \refT{Twalls} for odd walls.]
We split the proof into two parts.

\medskip
\noindent\emph{Part I.}
Consider an odd wall $\LLL$ in a Farey cell $\cR\uvx1$.
We will often use the coordinates $(\ga,\gb)$ defined by \eqref{m3-Phi};
recall also \eqref{Phi1}. 
(It will always be clear from the context whether we use $(p,q)$ or
$(\ga,\gb)$ as coordinates.) 
The line $\LL$ through the wall $\LLL$ 
is given by an equation of the form \eqref{jw1+}, which as shown in the
proof of \refT{TW2} is equivalent to  \eqref{ow8}. We write this as
\begin{align}\label{va1}
(k+1)\ga +(\kx+1)\gb = z
\end{align}
for some $z\in\bbZ$. 
Note that both coefficients $k+1$ and $\kx+1$ are strictly positive.
The $p$-axis is given by $q=0$ or, equivalently, see \eqref{Phi1}, by $\ga=\gb$.
Hence, $\LL$ intersects the $p$-axis at the unique rational point
\begin{align}\label{va2}
  \frac{u}v=\ga=\gb=\frac{z}{k+\kx+2},
\end{align}
where we as usual assume that $u/v$ is written in reduced form.
We can then write $\LL$ in parametric form as
\begin{align}\label{va3}
(\ga,\gb)=\bigpar{\ga(s),\gb(s)}=\Bigpar{\frac{u}v -(\kx+1)s,\frac{u}v+(k+1)s},
\qquad s\in\bbR,  
\end{align}
where we are interested in $s>0$ which gives $\ga<\gb$ and thus $q>0$.

By assumption, $\LL$ intersects $\cR\uvx1$, and thus there exists $(\ga,\gb)$
of the form \eqref{va3} with $s>0$ such that $(\ga,\gb)\in\Phi(\cR\uvx1)$,
which implies that, see \eqref{m3-farey-rect2} and \eqref{m3-open-farey},
\begin{align}\label{va4}
  0 < \ga < \frac{u_1}{v_1} <\gb <\frac12.
\end{align}
Hence, \eqref{va1} implies $0<z<(k+\kx+2)\frac12$, and thus \eqref{va2} shows
$\frac{u}v\in(0,\frac12)$. We denote as usual the Farey parents of
$\frac{u}v$ by $\frac{a}b$ and $\frac{c}d$.

Let us follow the line $\LL$ from a starting point in $\cR\uvx1$ downwards
to $u/v$. In the parametric representation \eqref{va3}, then $s$ decreases
from some initial value $s_1$ to 0, and this means that $\ga(s)$ 
increases from $\ga(s_1)$
to $u/v$ and $\gb(s)$ decreases from $\gb(s_1)$ 
to $u/v$. \refL{LUV} shows that at the start, 
$u_1/v_1\in(\ga(s_1),\gb(s_1))$, and that as $s$ decreases and the interval 
$\bigpar{\ga(s),\gb(s)}$ shrinks, the line passes through a sequence of
Farey cells $\cR\uvx{j}$, $j=1,\dots,N$, with (strictly)
increasing denominators $v_j$,
and ending at $u_N/v_N=u/v$ for some finite $N\ge1$.
Indeed, the fraction $u/v$ remains in every interval, so $v_j\le v$.
There are only finitely many fractions with denominator at most $v$,
and for sufficiently small $s$ the only such fraction in the interval is $u/v$.
It follows also from \refL{LUV} that $\LL$ crosses the boundaries of these
Farey cells transversally; the intersections are single points and the 
line $\LL$
immediately proceeds into the next open Farey cell.

Suppose that $N>1$, so that $\LL$ crosses from $\cR\uvx1$ into $\cR\uvx2$ at
some point $D_2=(p_2,q_2)$ given by $\bigpar{\ga(s_2),\gb(s_2)}$.

Recall that $P_0=U=(0,0)$.
By \refT{Tm3-Farey-entry}, for $(p,q)\in\cR\uvx1$, $P_1,\dots,P_{v_1-1}$ 
is the beginning of a branch orbit of $T$, where none of the points lies in
$\cB$. By continuity, $P_1,\dots,P_{v_1-1}$ is a branch orbit also
at the boundary point $D_2$, with none of the points in
$\cB\inter$. 
Furthermore, all these points lie in $\set{y>0}$, and therefore they are not
on the tie line $\gD_g$. Hence, if we allow only $T_1$ and $T_3$, there are
no ties at these points, so $P_1,\dots,P_{v_1-1},P_{v_1}$ is the unique
(partial) orbit started at $P_1$.
The same applies if we start inside $\cR\uvx2$ and approach $D_2$ from the
other side, now considering the orbit $P_1,\dots,P_{v_2-1}$.
Since $v_2>v_1$, it follows that at parameters $D_2$,
the points $P_1,\dots,P_{v_1}$ have the same limit from both sides.
Consequently, the algebraic form of
$P_1,\dots,P_{v_1}$  is the same in $\cR\uvx1$ and $\cR\uvx2$.
(In other words, see \eqref{m3-nu-Pk}, $\nu(k)$ is the same in
$\cR\uvx1$ and $\cR\uvx2$ for $1\le k \le v_1$.
The value $\nu(0)=0$ is fixed as well.)
This comparison uses only the shorter of the two partial orbits,
so it applies in either direction at a cell crossing.

Write $a_j/b_j<c_j/d_j$ for the Farey parents of $u_j/v_j$.
We next show that the left Farey parent stays the same for all $j$ when
$k<\kx$, and the right Farey parent stays the same when $k>\kx$.
Put
\begin{align}\label{va5}
 H:=k+\kx+2,\qquad
 D_L:=b_1z-a_1H,\qquad D_R:=c_1H-d_1z.
\end{align}
The Farey relation $b_1c_1-a_1d_1=1$ gives
\begin{align}\label{va6}
 d_1D_L+b_1D_R=H.
\end{align}
Moreover, since \eqref{va1} intersects the open Farey rectangle
$\cR\uvx1$,
we have by \eqref{m3-open-farey}
\begin{align}\label{va7}
 \frac{(k+1)a_1}{b_1}+\frac{(\kx+1)u_1}{v_1}
 <z<
 \frac{(k+1)u_1}{v_1}+\frac{(\kx+1)c_1}{d_1}.
\end{align}
Using $b_1u_1-a_1v_1=c_1v_1-d_1u_1=1$, this yields
\begin{align}\label{va8}
 v_1D_L>\kx+1,\qquad v_1D_R>k+1.
\end{align}
In particular, $D_L$ and $D_R$ are positive integers.

Let $\nu_1$ denote the function $\nu$ for $u_1/v_1$.
The definition \eqref{m3-nu-Pk} 
gives
\begin{align}
 \nu_1(j+b_1)&=\nu_1(j)+a_1+\indic{j=1},
 &&0\le j<d_1,\label{va9a}\\
 \nu_1(j+d_1)&=\nu_1(j)+c_1,
 &&0\le j<b_1.\label{va9b}
\end{align}
Indeed, $b_1u_1/v_1=a_1+1/v_1$ and
$d_1u_1/v_1=c_1-1/v_1$. 
Adding $1/v_1$ to
$(j-1)u_1/v_1$ changes its ceiling precisely when $j=1$ in
\eqref{va9a}. Subtracting $1/v_1$ in \eqref{va9b} cannot change the
ceiling: a fractional part $1/v_1$ would require
$j-1\equiv b_1\pmod{v_1}$, which is outside the stated range.
The odd-wall equation \eqref{ow8} also gives
\begin{align}\label{va10}
 z\equiv \nu_1(k)+\nu_1(\kx)+\indic{k=1}+1\pmod2.
\end{align}

\pfcase1{$k<\kx$.}
In this case,
$\kx=k+b_1$ and $k<d_1$, and thus
\begin{align}\label{va71}
H=2(k+1)+b_1\le2d_1+b_1
\qquad\text{with}\qquad 
H\equiv b_1 \pmod2.
\end{align}
Equations \eqref{va9a} and \eqref{va10} yield $z\equiv a_1+1\pmod2$,
and hence, using \eqref{va71},
\begin{align}\label{va72}
 D_L
&\equiv b_1(a_1+1) -a_1b_1
= b_1\pmod2,
\\ \label{va73}
D_R&\equiv c_1b_1-d_1(a_1+1)
=1-d_1\pmod2.
\end{align}
If $b_1$ were even, then $d_1$ would be odd and both determinants
$D_L$ and $D_R$
would be positive even integers. This would give
$H\ge2d_1+2b_1$ in \eqref{va6}, in contradiction to \eqref{va71}.
Thus $b_1$ is odd. 
Hence $H$ and $D_L$ are odd by \eqref{va71} and \eqref{va72}. 
Furthermore, by \eqref{va6} and \eqref{va71},
$d_1D_L\le H-b_1=2(k+1)\le2d_1$.
Since $D_L$ is positive and odd, we obtain
$D_L=1$ and thus by \eqref{va5}
\begin{align}\label{va11a}
 b_1z-a_1H=1.
\end{align}
Consequently $z$ and $H$ are coprime, 
and thus the fraction $z/H$ in \eqref{va2} is reduced; 
hence $u=z$ and $v=H$. In particular, $v$ is odd.

Furthermore,
at the starting point, $\ga(s_1)>a_1/b_1$; 
hence \eqref{va1}, \eqref{va11a}, and \eqref{va5} imply 
\begin{align}\label{va73a}
  b_1(\kx+1)\gb(s_1)
&=b_1z-b_1(k+1)\ga(s_1)
<b_1z-(k+1)a_1
\notag\\&
=a_1H+1-(k+1)a_1=(\kx+1)a_1+1,
\end{align}
and thus
\begin{align}\label{va12a}
 \frac{a_1}{b_1}<\ga(s_1)<\gb(s_1)
 <\frac{a_1}{b_1}+\frac{1}{b_1(\kx+1)}.
\end{align}
These bounds remain valid as $s$ decreases.
For any open Farey cell $\cR\uvx j$ met below the starting point,
\refL{LUV} gives 
$u_j/v_j\in\bigpar{\ga(s),\gb(s)}$ and further
$v_j\le v=H$, since $u/v$ remains in
$(\ga(s),\gb(s))$.
Put $\delta_j:=b_1u_j-a_1v_j$.
This is a positive integer, and \eqref{va12a} gives
\begin{align}
 \frac{\delta_j}{b_1v_j}
 =\frac{u_j}{v_j}-\frac{a_1}{b_1}
 <\frac{1}{b_1(\kx+1)}.
\end{align}
Thus $v_j>\delta_j(\kx+1)$.
But $v_j\le H=(k+1)+(\kx+1)<2(\kx+1)$, so $\delta_j\ge2$ is impossible.
Hence $\delta_j=1$ and 
thus $a_1/b_1$ and $u_j/v_j$ are Farey neighbours;
furthermore
$v_j>\kx+1>b_1$.
It follows that the left Farey parent of $u_j/v_j$ is $a_1/b_1$:
the other parent is $(u_j-a_1)/(v_j-b_1)$, and their mediant is $u_j/v_j$.
In particular, $b_j=b_1$ for every cell met on the way down.

\pfcase2{$k>\kx$.} 
In this case, $\kx=k-d_1$ and $\kx< b_1$, and thus
\begin{align}\label{va81}
H=2(\kx+1)+d_1\le2b_1+d_1
\qquad\text{with}\qquad 
H\equiv d_1 \pmod2.
\end{align}
Here $k\ge d_1\ge2$, so \eqref{va9b} and \eqref{va10} give
$z\equiv c_1+1\pmod2$. Thus
\begin{align}\label{va82}
 D_L&\equiv b_1(c_1+1) -a_1d_1
= 1+b_1\pmod2,
\\ \label{va83}
D_R&\equiv c_1d_1-d_1(c_1+1)
\equiv d_1\pmod2.
\end{align}
If $d_1$ were even, both determinants would be positive even integers,
contrary to \eqref{va81}
and \eqref{va6}.
Hence $d_1$ is odd, and thus $H$ and $D_R$ are odd by \eqref{va81} and
\eqref{va83}. 
Furthermore, \eqref{va6} and \eqref{va81} yield
$b_1D_R\le H-d_1=2(\kx+1)\le2b_1$,
and since $D_R$ is an odd positive integer, we have $D_R=1$ and thus by
\eqref{va5} 
\begin{align}\label{va11b}
 c_1H-d_1z=1.
\end{align}
Again,
$z$ and $H$ are coprime,  and \eqref{va2} implies that $u=z$ and $v=H$;
in particular, $v$ is odd.

Furthermore, in this case
\eqref{va1}, \eqref{va11b}, and \eqref{va5} imply
\begin{align}\label{va73b}
  d_1(k+1)\ga(s_1)&
=d_1z-d_1(\kx+1)\gb(s_1)
>d_1z-(\kx+1)c_1
\notag\\&
=c_1H-1-(\kx+1)c_1
=(k+1)c_1-1
\end{align}
and thus
\begin{align}\label{va12b}
 \frac{c_1}{d_1}-\frac{1}{d_1(k+1)}
 <\ga(s_1)<\gb(s_1)<\frac{c_1}{d_1}.
\end{align}
Now $\delta_j:=c_1v_j-d_1u_j$ is a positive integer,
$v_j>\delta_j(k+1)$, and $v_j\le v=H<2(k+1)$.
Consequently $\delta_j=1$, $v_j>k+1>d_1$, and the right Farey parent
of every $u_j/v_j$ is $c_1/d_1$.
Thus $d_j=d_1$ throughout.

In both cases, 
it follows for any $j=1,\dots,N$,
since $v_j\ge v_1$, that 
the original indices $k,\kx$ lie in
$\set{0,\dots,v_j-1}$ and retain the relation \eqref{kx2}.
The orbit comparison above shows, successively at each cell crossing,
that $\nu_j(k)$ and $\nu_j(\kx)$ retain their values.
Hence \eqref{va10} remains valid (with $\nu_j$)
in every cell met below the starting point, 
and thus the same line $\LL$ is of the form \eqref{ow8}, 
and thus \eqref{jw1+}, for that cell too. 
By \refT{TW2}, $\LL$ is therefore an odd wall in each of these cells.
In the last cell, $\cR\uv$, the line goes through $V\uv$, so
\refT{T-first-wall} identifies it as $\LL\uv$.
We have also proved that its coefficients in \eqref{va1} sum to the
reduced denominator $v$, a fact used below.

So far, we have proved that every odd wall in a Farey cell is the
intersection with some line $\LL\uv$, with $v$ odd; 
furthermore, for any such line $\LL\uv$, 
the union of the closures of the walls on $\LL\uv$ form a line segment on
$\LL\uv$, starting at $V\uv$. It remains to see how far this line segment
stretches.

\bigskip
\medskip
\noindent\emph{Part II.}
Fix $u/v\in(0,\frac12)$ with Farey parents $a/b<c/d$.
Assume $v$ is odd. Thus either $b$ is odd and $d$ is even, or $b$ is even and
$d$ is odd.


\pfcase1{$b$ is odd and $d$ even.}
As in the proof above for even walls,
we construct an increasing sequence of reduced fractions
$r_i:=u_i/v_i:=(u-ia)/(v-ib)$, with $i=0,\dots,i_1:=\ceil{v/b}-1$, such that
$r_0=u/v$, $r_1=c/d$, and
the Farey parents of $r_i$ are $a/b$ and $r_{i+1}$ for $0\le i<i_1$.
We write $\LL\uv$ as \eqref{va1} and parametrize it by
\eqref{va3}, where in this case
we have \eqref{bay6d} and thus
\begin{align}\label{vb1}
 k=\tfrac12d-1,\qquad
 \kx=k+b=b+\tfrac12d-1=v-\tfrac12d-1.
\end{align}
The right-hand side of \eqref{va1} is now $z=u$.
Note also that the coefficients $k+1=\frac12d$ and
$\kx+1=v-\frac12d$ are relatively prime,
since any common divisor would divide both $d$ and their sum $v$,
whereas $\gcd(d,v)=1$.

We move along $\LL\uv$ starting at $s=0$, where $\ga=\gb=u/v$, and
increase $s$. Then $\ga$ decreases and $\gb$ increases.
The coordinate $\ga$ reaches $a/b$ when
\begin{align}\label{vb2}
 s=\sx:=\frac{1}{\kx+1}\Bigpar{\frac{u}v-\frac{a}b}
 =\frac{1}{v-\frac12d}\cdot\frac{1}{vb}
 =\frac{1}{b(v-\frac12d)v}.
\end{align}
Extend $r_t=(u-ta)/(v-tb)$ to real $t\in[0,v/b)$.
Then $\gb$ reaches $r_t$ when
\begin{align}\label{vb3}
 s=\frac{1}{k+1}\Bigpar{\frac{u-ta}{v-tb}-\frac{u}{v}}
 =\frac{t(bu-av)}{(d/2)v(v-tb)}
 =\frac{t}{(d/2)v(v-tb)}.
\end{align}
Equations \eqref{vb3} and \eqref{vb2} agree when $bt=d/2$.
Consequently, before $\ga$ reaches $a/b$, $\gb$ hits
$r_1,\dots,r_I$, where
\begin{align}\label{vbT}
 I:=\ceil{d/(2b)}-1
< i_1
.\end{align}
If $d/(2b)$ is an integer, then  $r_{I+1}$ is reached only at $s=\sx$.
For $0<s<\sx$, the open cells met are exactly
$\cR\uvx i$, $0\le i\le I$, in this order.
Indeed, between the successive crossings we have
$a/b<\ga<r_i<\gb<r_{i+1}$, which is the open Farey rectangle criterion
\eqref{m3-open-farey} for $\cR\uvx i$.
(In particular, if $I=0$, the segment meets only the original cell
$\cR\uv$.)
For these indices,
$i<d/(2b)$ and thus
\begin{align}\label{vb4}
 v_i=v-ib>v-\tfrac12d=\kx+1.
\end{align}
The left Farey parent of every $u_i/v_i$ is $a/b$, so $b_i=b$.
Hence $\kx=k+b_i$ in each of these cells, and both indices are in
$\set{0,\dots,v_i-1}$.
The orbit comparison in Part I, applied at each crossing, shows that
the values of $\nu(k)$ and $\nu(\kx)$ do not change.
Thus 
$\LL\uv$
is of the form \eqref{ow8} and \eqref{jw1+} in every cell met.
By \refT{TW2}, each non-empty intersection is an odd wall.

At $s=\sx$, we have
equality in \eqref{va73a} (with $a$ and $b$ instead of $a_1$ and $b_1$)
and thus
\begin{align}\label{vb-end-ab}
 \ga=\frac ab,\qquad
 \gb=\frac ab+\frac{1}{b(\kx+1)}.
\end{align}
This point belongs to the closure of the lower parameter triangle
\eqref{Phi4}:
$\ga\ge0$ and $\gb\le1/2$, the latter because
$(\kx+1)(b-2a)\ge2$.
Let $m/n\in[\ga,\gb]$ be a reduced fraction different from $a/b$.
Since $\ga=a/b$, we have
\begin{align*}
 \frac{1}{bn}
 \le \frac{m}{n}-\frac{a}{b}
 \le \gb-\ga
 =\frac{1}{b(\kx+1)}.
\end{align*}
Hence $n\ge\kx+1=b+d/2>b$.
Thus the unique fraction of smallest denominator in $[\ga,\gb]$ is $a/b$
and occurs at an endpoint.
If $a>0$, \refL{LUV} shows that this point belongs to no open Farey cell.
If $a=0$, the same conclusion follows immediately from $\ga=0$.
 
For $s>\sx$, as long as the parameters remain in the lower parameter
triangle, $a/b$ belongs to $(\ga,\gb)$.
Therefore any open cell met has denominator $v'\le b$.
Suppose that the line $\LL$ is an odd wall there.
By Part I, any representation of this line as \eqref{ow8} has
right-hand side $u$ and coefficients of $\ga$ and $\gb$ summing to $v$.
Thus its coefficients are those in \eqref{vb1}.
But $\kx+1>b\ge v'$, contrary to the index range in \eqref{ow8}.
Hence the line $\LL$ cannot be an odd wall there.
It also cannot be an even wall, since the two coefficients are unequal.

Thus the wall intersections occur exactly for $0<s<\sx$, apart from the
isolated crossings of cell boundaries. Their closures form the segment
$0\le s\le\sx$.
By \eqref{va3}, \eqref{Phi1}, \eqref{vb1}, and \eqref{vb2},
the endpoint height is
\begin{align}\label{vb5}
 q\uv:=\frac{k+\kx+2}{3}\sx
 =\frac{v}{3}\sx
 =\frac{1}{3b(v-\frac12d)}
 =\frac{2}{3b(b+v)}.
\end{align}
This is $q(K\uv)$ in \eqref{qK}.
The corresponding $p$-value is, by \eqref{m3-H},
\begin{align}\label{vb6}
 p(K\uv)=\frac{u}v-\frac{2b+d/2}{v}q\uv
 =\frac{u}v-\frac{4b+d}{3bv(b+v)}.
\end{align}

\pfcase2{$b$ is even and $d$ is odd.}
This is a reflected version of Case 1.
Now 
\eqref{bay6b} gives
\begin{align}\label{vb7}
 k=b/2+d-1,\qquad \kx=b/2-1=k-d.
\end{align}
In the parametrization \eqref{va3}, $\gb$ reaches $c/d$ when
\begin{align}\label{vb8}
 s=\sx:=\frac{1}{k+1}\Bigpar{\frac cd-\frac uv}
 =\frac{1}{d(b/2+d)v}.
\end{align}
To find the cells met for $0<s<\sx$, we now use the decreasing
sequence
\begin{align}\label{vb12}
 r_i:=u_i/v_i:=\frac{u-ic}{v-id},
 \qquad 0\le i\le i_1:=\ceil{v/d}-1.
\end{align}
Then $r_0=u/v$, $r_1=a/b$, and the Farey parents of $r_i$ are
$r_{i+1}$ and $c/d$ for $i<i_1$.
For the extension to real $t\in[0,v/d)$, $\ga$ reaches $r_t$ when
\begin{align}\label{vb13}
 s=\frac{1}{\kx+1}\Bigpar{\frac uv-\frac{u-tc}{v-td}}
 =\frac{t}{(b/2)v(v-td)}.
\end{align}
This agrees with \eqref{vb8} when $dt=b/2$.
Thus the open cells met before the endpoint are exactly
$\cR\uvx i$, $0\le i\le I$, where
\begin{align}\label{vb14a}
 I:=\ceil{b/(2d)}-1 < i_1,
\end{align}
and for these indices
\begin{align}\label{vb14b}
 v_i=v-id>v-b/2=k+1.  
\end{align}
The strict rectangle inequalities for $\cR_{\uvx i}$ are now
$r_{i+1}<\ga<r_i<\gb<c/d$.
All these cells have right Farey parent $c/d$, so $d_i=d$ and
$\kx=k-d_i\equiv k+(v_i-d_i)\pmod{v_i}$.
The orbit comparison in Part I shows that $\nu(k)$ and $\nu(\kx)$
are  preserved
at the crossings. Hence the line satisfies \eqref{ow8} in every cell met,
and each intersection is an odd wall by \refT{TW2}.

At the endpoint $\sx$, we have equality in \eqref{va73b} (with $c$ and $d$)
and thus
\begin{align}\label{vb9}
 \ga=\frac cd-\frac{1}{d(k+1)},\qquad \gb=\frac cd.
\end{align}
Here $\ga>0$ and $\gb<1/2$, because $c(k+1)>1$ and $d$ is odd.
Let $m/n\in[\ga,\gb]$ be a reduced fraction different from $c/d$.
Since $\gb=c/d$, we have
\begin{align*}
 \frac{1}{dn}
 \le \frac{c}{d}-\frac{m}{n}
 \le \gb-\ga
 =\frac{1}{d(k+1)}.
\end{align*}
Hence $n\ge k+1=b/2+d>d$.
Thus the unique fraction of smallest denominator in $[\ga,\gb]$ is $c/d$
and occurs at an endpoint. By \refL{LUV}, this point belongs to no open
Farey cell.

For larger $s$, the interval $[\ga,\gb]$ contains $c/d$, and every open cell met
has denominator $v'\le d$. 
By Part I, any representation as \eqref{ow8} must have the
coefficients in \eqref{vb7}, which is impossible since $k+1>d\ge v'$.
The unequal coefficients also exclude an even wall.

Hence, again the wall intersection occur exactly for $0<s<\sx$, apart from a
finite number of boundary crossings.
The endpoint coordinates are
\begin{align}
 q(K\uv)&=\frac{v}{3}\sx
 =\frac{1}{3d(b/2+d)}
 =\frac{2}{3d(d+v)},\label{vb10}\\
 p(K\uv)&=\frac uv-\frac{b/2-d}{v}q(K\uv)
 =\frac uv+\frac{2d-b}{3dv(d+v)}.\label{vb11}
\end{align}
Thus \eqref{qK} holds in this case as well.

To summarize: 
In Part II we have shown that, in both cases, 
for any odd reduced fraction $u/v\in(0,\frac12)$, 
the line $\LL\uv$ forms odd walls in a finite number of Farey cells, and the
union of these walls (+ boundary points) is exactly the segment  $\LG\uv$.
Moreover,
in Part I, we have shown that any odd wall in a Farey cell is 
obtained in this way as the the intersection of the Farey cell and some line
$\LL\uv$ with $v$ odd. 
This completes the proof.
%
\end{proof}

\begin{corollary}\label{Cwalls}
The walls in any open Farey cell are precisely the non-empty
intersections of that cell with the segments $\LG_{u/v}$,
$u/v\in(0,1/2)$.
\end{corollary}

\begin{proof}
  This follows immediately from \refT{Twalls}.
\end{proof}

\subsection{The global wall pattern}\label{SSm3-wall-pattern}

Call a reduced fraction even or odd according to the parity of its
denominator. An odd fraction has a unique even child in the Farey tree;
call the two fractions \emph{adjoint}.
Conversely, every even fraction in $(0,1/2)$ is the child in the tree
of a unique odd fraction in $(0,1/2)$.
Hence, every reduced fraction in $(0,1/2)$ has a unique adjoint reduced
fraction.

The following theorem shows that the walls form a pattern of triangles,
as is clearly visible in 
\refFigs{fig:m3-farey-triangles} and \ref{fig:m3-farey-zoom}.

\begin{theorem}\label{ConjFarey+lines}
For adjoint fractions $u/v$ and $u'/v'$, 
the endpoints of the segments $\LG\uv$ and $\LG\uvi$ coincide:
$K_{u/v}=K_{u'/v'}$.
This common endpoint is the intersection of first-wall lines $\LL\uv$ and
$\LL\uvi$, 
and the two segments thus form the sides of a triangle with base on $q=0$.
\end{theorem}

\begin{proof}
Let $v$ be odd and let $a/b<c/d$ be the Farey parents of $u/v$.

If $b$ is odd, the even child is $u'/v'=(u+a)/(v+b)$,
whose Farey parents are $a/b$ and $u/v$.
By \eqref{qK}, \eqref{vc4}, and \eqref{vb6}, its endpoint has coordinates
\begin{align}
 q(K_{u'/v'})&=\frac{2}{3b(v+b)}=q(K_{u/v}),\\
 p(K_{u'/v'})&=\frac{u+a}{v+b}-\frac{1}{3b(v+b)}
 =\frac uv-\frac{4b+d}{3bv(v+b)}
 =p(K_{u/v}),
\end{align}
where the second equality uses $bu-av=1$ and $v=b+d$.

If $d$ is odd, the even child is $u'/v'=(u+c)/(v+d)$.
Now its parents are $u/v$ and $c/d$, so
\eqref{qK}, \eqref{vc4}, and \eqref{vb11}  yield
\begin{align}
 q(K_{u'/v'})&=\frac{2}{3d(v+d)}=q(K_{u/v}),\\
 p(K_{u'/v'})&=\frac{u+c}{v+d}-\frac{1}{3d(v+d)}
 =\frac uv+\frac{2d-b}{3dv(v+d)}
 =p(K_{u/v}),
\end{align}
using $cv-du=1$. 

The two first-wall lines  $\LL\uv$ and $\LL\uvi$
thus meet at $K\uv$, and together with the $p$-axis $q=0$ they  form a
triangle with vertices $V\uv,V\uvi,K\uv$, and sides $\LG\uv$ and $\LG\uvi$.
\end{proof}

\begin{remark}
  Recall that each wall in a Farey cell signifies that some vertex on
  the lower boundary of $\cD$ lies on the tie-line $\gD_g$.
Thus, when two walls meet at a point on the boundary of the Farey cell 
(recall that they cannot meet in the interior), this means that as $(p,q)$
approaches this boundary point, the two corresponding vertices of $\cD$ merge;
these vertices are necessarily a valley and a peak on the lower boundary,
and when they meet, they annihilate each other, which explains why there is
no continuation of these walls as walls in the Farey cell on the
other side of the boundary.
\end{remark}

\section{Dynamical systems for other election  methods}
\label{Sdiv}

We end by noting that the same dynamical point of view also applies to the
standard divisor methods.
Consider, as in \refR{Rprop}, an election where several representatives are
chosen from party lists.  Two of the most common rules used in practice
are \dH's method and \SL's method.  
They are examples of linear divisor methods, which can be defined as
follows.  We write $p_j$ for
the vote proportion of party $j$, and assume here that all $p_j$ are positive.
For other descriptions of divisor methods, see \cite{BY}, \cite{Pukelsheim},
and \cite[Appendix A]{SJ262}.

\def\MethodName{Linear divisor methods}
\begin{method}
Seats are awarded one at a time.  Let $r_j(t)$ be the number of seats already
given to party $j$ after $t$ rounds.  In round $t$, the winning party is
\begin{align}
        E(t):=\argmax_{j\in[m]} \frac{p_j}{r_j(t-1)+\beta}.
\end{align}
Here $\beta>0$ is fixed.  The choice $\beta=1$ gives \dH{}'s method, and
$\beta=1/2$ gives \SL{}'s method.
\end{method}

Other values of $\beta$ have also been used or discussed.  For example,
$\beta=1/3$ gives the Danish method and $\beta=2$ gives Imperiali's method.
Although the definition above assumes $\beta>0$, the limiting value
$\beta=0$ corresponds to Adams's method, with the usual extra convention for
the first seats.

The rule can be put in almost the same form as Frege's rule.  Since all
$p_j>0$, we may invert the fractions and write
\begin{align}
E(t)
&=\argmin_{j\in[m]} \frac{r_j(t-1)+\beta}{p_j}
\notag\\&
=\argmax_{j\in[m]} \frac{-r_j(t-1)-\beta}{p_j}
\notag\\&
=\argmax_{j\in[m]} \frac{t p_j-r_j(t-1)-\beta}{p_j}.
\end{align}
In the last line we have only added the same constant $t$ to all the
fractions.  If $z_j(t)$ is the discrepancy, defined as for Frege's method 
in \refS{SFrege} so that
$r_j(t)=t p_j-z_j(t)$, this becomes
\begin{align}
        E(t)=\argmax_{j\in[m]} \frac{p_j+z_j(t-1)-\beta}{p_j}.
\end{align}
If we make a uniform shift of the variables,
replacing $z_j(t)-\gb$ by $z_j(t)$, this becomes 
the following, which is very similar to the
Frege recursion studied in the present paper.

\begin{theorem}\label{Tdivisor}
A linear divisor method with parameter $\beta$ can be written as the recursion
\begin{align}\label{div1}
        s_j(t)=z_j(t-1)+p_j,
        \qquad
        z_j(t)=s_j(t)-{\bf 1}_{\{E(t)=j\}},
\end{align}
with initial condition
\begin{align}\label{div2}
        z_j(0)=-\beta,
        \qquad 1\le j\le m,
\end{align}
and with winner rule
\begin{align}\label{div3}
        E(t):=\argmax_{j\in[m]} \frac{s_j(t)}{p_j}.
\end{align}
Ties may be treated by any fixed tie-breaking rule, or by the multi-valued
convention used earlier in the paper.
\end{theorem}

\begin{proof}
The calculation preceding the theorem proves the displayed recursion and the
winner rule after the uniform shift of variables.  The statement about ties is
just the corresponding choice of convention for the displayed argmax.
\end{proof}

Equivalently, define
\begin{align}\label{thetap}
        \theta_{\vp}(\vx):=\argmax_{1\le i\le m} \frac{x_i}{p_i}.
\end{align}
Then the corresponding piecewise translation is
\begin{align}\label{div4}
        T_{\vp}^{\rm div}(\vx)=\vx+\vp-\ve_{\theta_{\vp}(\vx+\vp)},
        \qquad
        \vz(0)=-\beta\1 .
\end{align}
Thus all linear divisor methods have the same map $T_{\vp}^{\rm div}$ for a
fixed vote vector $\vp$;  the parameter $\beta$ changes only the starting point.
Note also that the recursion \eqref{div1} has the same form as the Frege
recursion \eqref{sj}--\eqref{zj};
however, the winner rule \eqref{div3} is different and depends on $\vp$;
the winner is now chosen by the weighted quantities $x_i/p_i$, rather than
by the coordinates $x_i$ themselves.
(One can also make a linear change of variables, replacing $x_i/p_i$ by $x_i$.
Then the winner is determined by \eqref{theta} as for the Frege recursion,
but the translations \eqref{div4} will be different and depend explicitly on
$\vp$.)

\begin{example}
For the case $m=3$, we may as in \refS{Sm=3}
regard the recursion as a dynamical system in $\bbR^2$.
Use again the notation $\vp=(p,q,r)$. The invariant hyperplane for
\eqref{div1}--\eqref{div2} is $z_1+z_2+z_3=-3\beta$. Put
\begin{align}
\vw:=\vz+3\beta\vp,
\end{align}
so that $w_1+w_2+w_3=0$, and identify
$\vw=(x,y,-x-y)$ with $(x,y)\in\bbR^2$. Since
\begin{align}
\frac{z_i+p_i}{p_i}=\frac{w_i}{p_i}+1-3\beta,
\end{align}
the winner \eqref{div3} equals $\argmax_{i\in[m]} (w_i/p_i)$. 
Moreover, the update of $\vw$ is
$\vw\mapsto\vw+\vp-\ve_i$. Thus the planar map is the piecewise
translation \eqref{m3-T}, but with the winner domains
$\cA,\cB,\cC$ given, as in \eqref{m3-cA}--\eqref{m3-cC}, by the functions
\begin{align}
 f^{\rm div}(x,y)&:=-qx+py,\\
 g^{\rm div}(x,y)&:=(1-p-q)x-p(-x-y)=(1-q)x+py,\\
 h^{\rm div}(x,y)&:=(1-p-q)y-q(-x-y)=qx+(1-p)y.
\end{align}
We expect that there is an invariant polygon for this system, as in
\refS{Sm=3},
but we do not know whether there are several stability domains or not.
However, the simple result for the bias in \cite{SJ262}, 
compared to the more complicated result found in
\refS{Sbias},
suggests that the situation is simpler than for Frege's method.
\end{example}

The map \eqref{div4} is a piecewise translation, and thus another example of
the maps studied by \citet{gun}. However, it seems to
fall outside the class of maps studied by Adler et al.\
\cite{adler0,adler10,adler17,adler18}.
It would be interesting to study this system 
further, and compare it to the results above for Frege's system,
but we leave this for future research.

\section*{Acknowledgement}

The third author wishes to thank Anders Johansson for many good discussions.

\newcommand\AAP{\emph{Adv. Appl. Probab.} }
\newcommand\JAP{\emph{J. Appl. Probab.} }
\newcommand\JAMS{\emph{J. \AMS} }
\newcommand\MAMS{\emph{Memoirs \AMS} }
\newcommand\PAMS{\emph{Proc. \AMS} }
\newcommand\TAMS{\emph{Trans. \AMS} }
\newcommand\AnnMS{\emph{Ann. Math. Statist.} }
\newcommand\AnnPr{\emph{Ann. Probab.} }
\newcommand\CPC{\emph{Combin. Probab. Comput.} }
\newcommand\JMAA{\emph{J. Math. Anal. Appl.} }
\newcommand\RSA{\emph{Random Structures Algorithms} }
\newcommand\DMTCS{\jour{Discr. Math. Theor. Comput. Sci.} }

\newcommand\AMS{Amer. Math. Soc.}
\newcommand\Springer{Springer-Verlag}
\newcommand\Wiley{Wiley}

\newcommand\vol{\textbf}
\newcommand\jour{\emph}
\newcommand\book{\emph}
\newcommand\inbook{\emph}
\def\no#1#2,{\unskip#2, no. #1,} 
\newcommand\toappear{\unskip, to appear}

\newcommand\arxiv[1]{\texttt{arXiv}:#1}
\newcommand\arXiv{\arxiv}

\end{document}